\documentclass[11pt]{article}

\usepackage[numbers]{natbib}
\usepackage{enumitem}
\usepackage{amsmath}
\usepackage{amsfonts}
\usepackage{amssymb}
\usepackage{epsfig}
\usepackage{dsfont}
\usepackage{mathrsfs}
\usepackage{bbm}
\usepackage{bm}
\usepackage{multirow}
\usepackage{makecell}
\usepackage{threeparttable}
\usepackage{authblk}
\usepackage{mathabx} 
\usepackage{amsxtra}
\usepackage{graphicx}
\usepackage{animate}
\usepackage{tikz}
\usepackage{algorithmic}
\usepackage{float}  
\usepackage{lipsum}
\usepackage{caption}
\usepackage{soul}
\usepackage{url}
\usepackage{algorithm}
\usepackage[title]{appendix}
\usepackage{boldline}
\usepackage{diagbox}
\usepackage{xcolor}
\usepackage{pifont}
\usepackage{harmony}
\usepackage{scalerel}
\usepackage[colorlinks=true]{hyperref}
\hypersetup{urlcolor=blue, citecolor=blue,linkcolor=blue}
\usepackage[capitalise,nameinlink]{cleveref}

\RequirePackage{xcolor}[2.11]
\colorlet{inlinkcolor}{green!60!black}
\colorlet{exlinkcolor}{red!50!black}
\colorlet{reviewcolor}{black!50}
\crefname{equation}{}{}
\crefname{figure}{Figure}{Figures}
\crefname{assumption}{Assumption}{Assumptions}
\crefname{subappendix}{Appendix}{Appendices}

\allowdisplaybreaks[4]
\usepackage{subfigure}
\counterwithin{algorithm}{section}
\newtheorem{theorem}{Theorem}[section]
\newtheorem{lemma}{Lemma}[section]
\newtheorem{proposition}{Proposition}[section]
\newtheorem{example}{Example}[section]
\newtheorem{remark}{Remark}[section]
\newtheorem{definition}{Definition}[section]
\newtheorem{assumption}{Assumption}[section]
\newtheorem{corollary}{Corollary}[section]

\numberwithin{equation}{section}
\numberwithin{definition}{section}

\newcommand{\dd}{\mathsf {d\kern -0.07em l}} 
\newcommand{\bgeqn}{\begin{eqnarray}}
\newcommand{\edeqn}{\end{eqnarray}}
\newcommand{\bgeq}{\begin{eqnarray*}}
\newcommand{\edeq}{\end{eqnarray*}}
\newcommand{\bec}{\begin{center}}
\newcommand{\enc}{\end{center}}
\newcommand{\R}{{\rm I\!R}}

\newcommand{\inmat}[1]{\mbox{\rm {#1}}}

\renewcommand{\Box}{\framebox{\rule{0.3em}{0.0em}}}

\newcommand{\be}{\begin{equation}}
\newcommand{\ee}{\end{equation}}

\newcommand{\vt}{{\vartheta}}

\newcommand{\proof}{\noindent{\bf Proof. }}

\def\bbe{{\mathbb{E}}}

\providecommand{\BIBand}{and}
\newcommand{\F}{{\cal F}}
\usepackage{adjustbox}
\usepackage{booktabs} 
\def\supp{{\rm supp}}
\usepackage{makecell}

\title{Adaptive Preference Elicitation in Preference Robust CPT-Based Shortfall\footnote{This project is 
 supported by Hong Kong RGC grant 14204624 and The Chinese University of Hong Kong start-up grant.}
}

\author[1]{Sainan Zhang}
 \author[2]{Erick Delage}
 \author[1]{Huifu Xu}
 \affil[1]{Department of Systems Engineering and Engineering Management,
           The Chinese University of Hong Kong, Shatin, N.T., Hong Kong.}
\affil[2]{GERAD and Department of Decision Sciences, HEC Montr\'eal, Montreal, Quebec H3T 2A7, Canada.}         
\affil[]{\textsf{snzhang.m@gmail.com,
                erick.delage@hec.ca,
                hfxu@se.cuhk.edu.hk}}

\begin{document}
  \maketitle
\vspace{-3em}
\begin{abstract} 
This paper investigates the cumulative prospect theory (CPT)-based shortfall risk measure where the value (utility) function and the distortion (probability weighting) functions are ambiguous. 
We propose a worst-case CPT-based shortfall risk measure where
the capital requirement is based on the worst-case value function and distortion functions from an ambiguity set
constructed with user-specific information (e.g.,~preference between risky and certain lotteries) and well-known general information such as the $S$-shape of the value functions.
When the underlying random loss is discretely distributed and parameters in pairs of lotteries remain undetermined, 
we demonstrate that the worst-case value function and the distortion functions have a piecewise linear structure.
Based on this observation, we develop an adaptive preference elicitation approach extended from the well-known polyhedral method to simultaneously elicit the distortion functions and the positive and negative parts of the value function by solving a sequence of mixed-integer linear programs.
Based on the elicited distortions, 
we design new queries to construct a feasible set of the holistic value functions.
The proposed framework iteratively shrinks the ambiguity set through designed queries, yielding a computationally implementable robust CPT-based shortfall risk measure.
To examine the in-sample and out-of-sample performance of the elicited CPT-based preference functions and the resulting preference robust shortfall risk measure,
we conduct numerical experiments under an $S$-shaped value function and inverse $S$-shaped distortion functions and report the corresponding~results.
\end{abstract}

\vspace{-0.15cm}

\section{Introduction}
\setcounter{equation}{0}
 
Let $\xi: (\Omega,{\cal F},\mathbb{P})\to\R$ be a random variable representing a financial loss, and let $l:\R\to \R$ be a convex, increasing, and non-constant {\em loss} function.
The utility-based shortfall risk measure (SRM) of $\xi$ is defined~as
\begin{equation}
\label{eq:defi-SR}
\inf_{x\in \R} \{x:
\bbe_P[l(\xi-x)]\leq \lambda\},
\end{equation}
where the expectation is taken with respect to (w.r.t.) the probability distribution of $\xi$, 
i.e., $P=\mathbb{P}\circ \xi^{-1}$, and $\lambda$ is an interior point of the range of $l$.
The SRM in \eqref{eq:defi-SR} is the minimum cash injection into the financial position ($-\xi$) such that the expected loss of the resulting position does not exceed the specified level $\lambda$.
Since its introduction by \cite{FoS02}, 
SRM has received wide attention in risk management and finance not only because it covers some well-known risk measures,
such as value-at-risk, conditional value-at-risk, 
and the entropic risk measure, but also because it better captures tail losses (through a proper choice of loss function $l$, see \citep{FoS04}). 
Moreover,  it satisfies invariance under randomization and is more suitable for dynamic measurement of risks over time (see \citep{GSW08,We06}).

The shortfall risk model \eqref{eq:defi-SR} is essentially based on the expected utility theory (EUT, see \citep{VNM47}) if we regard the loss function as a disutility function ($l(y):=-u(-y)$).  
The EUT asserts that if a decision maker's (DM's) preference over random prospects satisfies some axiomatic properties (including completeness, transitivity, continuity, and independence), 
then there exists a unique utility function up to a positive affine transformation such that the DM's preference over the random prospects can be described by the expected utility.
While EUT is recognized as a normative benchmark in decision theory, 
it is also found that the theory contradicts some experimental studies.
Classic studies such as \citet{Al53} and \citet{El1961} reveal behavior inconsistent
with the independence axiom and expected utility, 
which prompts various remedies/modifications of the EUT.
Rank-dependent expected utility (RDEU) theory (see \citep{Qui82}, \citep[Chapter 5]{Qui12}, \citep{Sch89} and \citep{Abd02}) and cumulative prospect theory (CPT, see \citep{TvK92}, \citep[p.~158]{WaT93} and \citep[Definition 1]{KSW11}) are best known among them.
In the RDEU theory, a DM's preference is not only described by a utility function 
but also by a distortion function  (probability weighting function) which modifies the 
cumulative probabilities of the underlying random prospects.
Consequently, a DM's preference for a random prospect $-\xi$ is represented by distorted expected utility $\int_{-\infty}^{\infty}u(x)dg(\mathbb{P}(-\xi\leq x))$, 
as opposed to the undistorted expected utility of $-\xi$ under the EUT, where the distortion function $g:[0,1]\to [0,1]$ is increasing and continuous with $g(0)=0$ and $g(1)=1$. 
The CPT takes a step forward by allowing the DM's preference to be dependent on 
the relative wealth as opposed to the overall wealth in the RDEU theory and subsequently the utility function (called the value function) is defined on the deviations from a reference point (representing relative wealth), 
i.e., on the additional gains and losses.
The distortion function is also modified with so-called sign dependence which consists
of two separate distortion functions,
one transforming probabilities of gains and 
the other transforming probabilities of losses.
Mathematically, the DM's preference for $-\xi$ under CPT can be written as
\begin{equation}
\label{eq:defi-bbe-wwv}
\begin{aligned}
 \bbe_{g_1 g_2}[v(-\xi)]
 :=&\int_{0}^{\infty}g_2(\mathbb P(v(-\xi)\geq t))dt
 -\int_{-\infty}^{0} g_1(\mathbb P(v(-\xi)\leq t))dt.
 \end{aligned}
\end{equation}
From the formulation, we can see that $g_2$ modifies the decumulative distribution function (also known as the survival function or the complementary cumulative distribution function (complementary CDF)) of $v(-\xi)$ when $v(-\xi)$ is positive whereas $g_1$ modifies the CDF of $v(-\xi)$ when $v(-\xi)$ is negative. 
Note that when $g_1(t)=1-g_2(1-t)$, for $t\in [0,1]$, which excludes sign dependence, the CPT reduces to the RDEU theory \cite[p.~19]{TuF02}.
In general, $g_1(t)\neq 1-g_2(1-t)$ because the DM's view over losses and gains is not necessarily ``symmetric''.
The definition of the SRM inherits the drawbacks of the EUT and this prompts the introduction of a generalized shortfall risk (GSR) measure by \cite{MaC18}:
\begin{equation}
\label{eq:PR-IPC-CPT}  
\inf_{x\in \R}\left\{ x: \mathbb
E_{g_1g_2}[v(\xi-x)]\leq \lambda\right\},
\end{equation}
where $v$ is the ``disvalue'' function (disutility function)\footnote{If $u(\cdot)$ is a value function (utility function) in the CPT, then $-u(-\cdot)$ is a ``disvalue'' function (disutility function).}.
Throughout the paper, we continue to call $v$ a value function although we actually mean the ``disvalue'' function to ease the exposition. 
In this paper, we take on this stream of research but with a very different focus. 
We consider a situation where the functional $\bbe_{g_1g_2}[v(\cdot)]$ in \eqref{eq:defi-bbe-wwv} is ambiguous due to incomplete information on the value and
distortion functions.
For example, there could be a case in which the preference of a market regulator may be described by CPT with a tuple of the value and distortion functions but there is a lack of complete information (particularly from the modeller's perspective) on these functions.
If we know that the true functions lie in a finite set of plausible value and distortion functions, 
then we may use model \eqref{eq:PR-IPC-CPT} to compute the GSR one by one using all possible tuples of the value and distortion functions and subsequently identify the worst tuple which gives rise to the largest shortfall risk measure.
However, this approach is not feasible in practice when the set is infinite.

A potential way to tackle the issue is to use  partially available information to construct a set of plausible value and distortion functions, 
called an ambiguity set, 
and obtain the worst GSR measure.
This extends the preference robust version of the utility-based shortfall risk measure \eqref{eq:defi-SR} (see \citep{DGX22}) to the CPT-based shortfall.
By assuming that the DM's risk preference can be characterized by shortfall or, 
more generally, by convex/coherent risk measures (when $l$ is convex, the shortfall risk measure is convex),
\citet{DGX22} use the DM's risk preference information elicited via certainty equivalents and the DM's risk preference over large tail losses to construct  an ambiguity set of loss functions.
This paper takes a step further to consider a preference robust version of model \eqref{eq:PR-IPC-CPT}.

A key element of the research is how to elicit the DM's preferences and use them to construct the ambiguity set and reduce the size of the set effectively. 
Research on preference elicitation is well documented primarily in the marketing literature (see \citep{TSHD03,THS04,THG07,SaV19} and references therein).
For instance, \citet{THS04} propose a polyhedral method for eliciting customers' preferences characterized by multivariate linear utility functions in discrete choice models.
\citet{SaV19} develop an ellipsoidal approach which extends the effectiveness of the polyhedral approach to the high response-error setting.
\citet{VYMDR20} propose an objective-oriented method for eliciting a multivariate linear utility function.  
Unfortunately, none of these approaches can be readily applied to elicit the value and distortion functions in CPT because they rely heavily on the linearity of the function to be elicited. 
\citet{ArD15} apply the so-called random relative utility split (RRUS) scheme to elicit a nonlinear univariate utility function under the EUT framework.
However, this approach cannot be applied to simultaneously elicit the value and distortion functions in CPT.
\citet{Abd00} develops a two-step method to elicit the positive and negative parts of nonlinear univariate utility functions and the distortion functions under CPT by constructing a standard sequence of outcomes and a standard sequence of probabilities,
where a tuple of the utility and distortion functions is estimated.
\citet{ABP07} develop a four-step method to estimate the utility for gains and losses under prospect theory.
Here, we follow this stream of research by proposing modified polyhedral methods that enable the use of pairwise comparisons to effectively elicit the DM's preference and subsequently construct the joint ambiguity set of the value and distortion functions in~CPT.

\subsection{Contributions}
The main contributions of the paper can be summarized as follows.
 \begin{itemize}
 \item 
\textbf{Modeling}.
We consider the robust GSR based on CPT.
Unlike preference robust utility-based shortfall (see \citep{DGX22}) which corresponds to identity distortion functions in the GSR framework, 
the robust GSR considers general inverse $S$-shaped distortion functions. 
Moreover, unlike the preference robust distortion risk measure (DRM, \citep{WaX23}) corresponding to an identity value function in the GSR framework, the robust GSR considers an $S$-shaped value function.
The new model focuses on the case where the DM's preferences can be described by CPT with non-identity value and distortion functions.

\item
\textbf{Methodology}.
We propose a novel polyhedral method, modified from the well-known polyhedral method 
originally developed by \cite{THS04} for eliciting linear multivariate utility functions, 
to elicit the nonlinear univariate value function and the distortion functions in CPT. 
Specifically, we consider piecewise linear value and distortion functions, 
represent each of them as a linear function of a vector of increments corresponding to linear pieces between consecutive breakpoints, 
and implement the polyhedral cuts in the space of the increment vectors. 
Each cut enables us to  purposely (rather than randomly, as in the RRUS method of \cite{ArD15}) generate a query for eliciting the DM's preferences and effectively reduce the ambiguity of the value and distortion functions. 
Mathematically, each cut requires solving a coupled, specially structured maximization problem.
To handle the bilinear terms in the constraints of the CPT-based robust GSR model, 
we develop four types of well-structured queries such that the resulting constraints in the ambiguity set become linear.

\item \textbf{Numerical tests}.
We demonstrate how the modified polyhedral method can be applied to  elicit the pair of the positive (resp. negative) part of the value function and the corresponding distortion function,
and how to apply the modified polyhedral method to elicit the holistic value function. 
Using a sequence of synthetic queries,
we conduct tests of the out-of-sample performance of the resulting CPT-value (based on the worst-case value function and distortion functions) in terms of hit rates in comparison with the true CPT-value.
\end{itemize}

It is important to note that this paper focuses on ambiguity in the value function and the distortion functions in the CPT models, 
how to design appropriate preference elicitation procedures to estimate these functions and to derive the computationally feasible evaluation of the preference robust shortfall, rather than exploring sufficient conditions for the convexity or coherence of GSR.
Our restriction that the value function is $S$-shaped and the distortion functions are inverse $S$-shaped is in alignment with the findings in the literature (see  \citep{TvK92,Pre98,WuG96,Abd00,TvF95,BlP00,BBEJ22}).

The rest of the paper is organized as follows.
Section~\ref{Sec:robust-GSR-CPT} introduces the worst-case GSR measure based on CPT and discusses the equivalent definition, 
a shortfall risk measure with a worst-case constraint, 
referred to as the preference robust GSR based on CPT (PRGSR-CPT),
and its main properties.
Section~\ref{sec:tractability} presents a computationally feasible reformulation of the PRGSR-CPT.
Section~\ref{sec:Estimating-dis} develops two types of queries to estimate a pair consisting of the positive (resp. negative) part of the value function and the corresponding distortion function by a modified polyhedral method with theoretical guarantees.
Based on these queries, 
it further discusses the elicitation of the holistic value functions using a modified polyhedral method tailored for the value function.
The resulting bilinear optimization problem is approximated by applying a piecewise McCormick relaxation in Section~\ref{sec:computation}.
Section~\ref{sec5:case_study} reports numerical test results including the efficiency of preference elicitation.
Finally, Section~\ref{sec:conclusion} concludes.
All proofs are relegated to the Appendix.

\subsection{Notation}
For any integer $N\geq 1$, we write $[N]$ for $\{1,\ldots,N\}$.
We use $|S|$ to denote the cardinality of a set $S$.
We use a boldface lowercase letter to denote a column vector and a boldface uppercase (capital) letter to denote a matrix. 
To ease the exposition, we write $(x_1,\ldots,x_n)$ for an $n$-dimensional column vector ${\bm x}$ and $({\bm a}, {\bm b})$ for $({\bm a}^\top, {\bm b}^\top)^\top$ where $\top$ denotes the transpose.
Given any matrix ${\bm A} = (a_{ij})_{i\in[m],j\in[n]} \in \R^{m\times n}$, we let ${\bm a}_i^\top$ be its $i$-th row vector.
Let ${\bm e}$ be the vector with all the entries being $1$.
Finally, we use ${\cal L}^\infty(\R)$ to represent the space of essentially bounded random variables $\xi: (\Omega,\F,\mathbb{P})  \rightarrow \R$, equipped with the essential supremum norm $\| \xi\|_\infty :=\inf\{a\in \R :\mathbb{P}(|\xi|>a)=0\}$.

\section{Preference Robust Generalized Shortfall Risk Under CPT}
\label{Sec:robust-GSR-CPT}

In this section, we formally introduce the worst-case CPT-based shortfall,
the PRGSR-CPT, and establish their equivalence. This equivalence paves the way for calculating the worst-case CPT-based shortfall by evaluating PRGSR-CPT.
In Section~\ref{sec:properties}, we analyze the convexity and coherence of the PRGSR-CPT based on results established in the non-robust cases \citep{MaC18} and explain the relation of these properties to the notions of so-called strong risk aversion or strong risk taking behaviors in the behavioral economics literature \citep{RoS70}.

\subsection{PRGSR-CPT}
\label{sec:CPT-shortfall}
We begin with a formal definition of the worst-case CPT-based shortfall.
\begin{definition}
[Worst-Case CPT-Based Shortfall]
\label{def:Worst-CPT-shortfall}
Let $\underline{x}, \bar{x}$ be such that $\underline{x}<0<\bar{x}$, 
and let $\mathscr V$ denote the set of all value functions $v: [\underline{x},\bar{x}]\rightarrow \R$ which are strictly increasing with $v(0)=0$. 
Let $\mathscr{G}$ denote the set of all distortion functions $g:[0,1]\rightarrow [0,1]$ which are strictly increasing with $g(0)=0$ and $g(1)=1$.
For each tuple $(v,g_1,{g_2})\in \mathscr{V}\times \mathscr{G} \times \mathscr{G}$,
the CPT-based shortfall of a loss $\xi\in {\cal L}^\infty(\R)$ is defined as follows:
\begin{equation}
\label{eq:rho-CPT}
\rho_{v,g_1,{g_2}}(\xi):=\inf\limits_{x\in \R} \left\{ x: \bbe_{g_1{g_2}}[v(\xi-x)]\leq v(0) \right\}.
\end{equation}
Let $\mathbb{ V}\mathbb{G}_1\mathbb{G}_2\subset {\mathscr V}\times {\mathscr G}\times {\mathscr G}$ be a joint ambiguity set of the value functions and the distortion functions.
The worst-case (largest) CPT-based shortfall of $\xi$ associated with $\mathbb{V}\mathbb{G}_1\mathbb{G}_2$ is defined as
\begin{equation}
\label{eq:PRGRS-CPT-equiv}
\Pi_{\mathbb{ V}\mathbb{G}_1\mathbb{G}_2}(\xi) 
:= \sup\limits_{(v,g_1,{g_2})\in \mathbb{ V}\mathbb{G}_1\mathbb{G}_2} \rho_{v,g_1,{g_2}}(\xi).
\end{equation}
\end{definition}

In Definition~\ref{def:Worst-CPT-shortfall}, 
we may set $\bar{x}=\max\{\bar{x},-\underline{x}\}$ so that the value function is defined over a symmetric interval $[-\bar{x},\bar{x}]$.
Note also that the interval $[\underline{x},\bar{x}]$ can be chosen as the convex hull of the supports of $\xi-x$
for all $x$ between the essential lower and upper bounds of $\xi$.
There is an alternative way to define the robust version of CPT-based shortfall.
\begin{definition}
[PRGSR-CPT]
\label{defi:worst-CPT-constraint}
The preference robust generalized shortfall risk based on CPT of a loss $\xi\in {\cal L}^{\infty}(\R)$ associated with $\mathbb{ V}\mathbb{G}_1\mathbb{G}_2$ is defined as
\begin{equation}
\label{eq:PRGRS-CPT}
\begin{aligned}
\rho_{\mathbb{V}\mathbb{G}_1\mathbb{G}_2}(\xi)
:=\inf_{x\in\R}\left\{x: 
\sup_{(v,g_1,g_2)\in \mathbb{V}\mathbb{G}_1\mathbb{G}_2}
\mathbb E_{g_1g_2}[v(\xi-x)]
\leq 0\right\}.
\end{aligned}
\end{equation}
\end{definition}

In Definition~\ref{defi:worst-CPT-constraint}, the worst-case distorted expected value of $v(\xi-x)$ is used in the constraint.
Note that we set the reference point in CPT to be zero in \eqref{eq:rho-CPT} and \eqref{eq:PRGRS-CPT} since the preference functional $\bbe_{g_1 g_2}[v(\cdot)]$  in the shortfall \eqref{eq:rho-CPT} and the robust shortfall \eqref{eq:PRGRS-CPT} is applied to $\xi-x$ instead of $\xi$.
Here we make a blanket assumption that the robust constraint function in \eqref{eq:PRGRS-CPT} is well-defined. 
A sufficient condition is that $x$ is bounded and $\xi\in {\cal L}^{\infty}(\R)$ is essentially bounded and $(v,g_1,g_2)\in \mathbb{V}\mathbb{G}_1\mathbb{G}_2$ are continuous and strictly increasing.
Note that not all ambiguity sets $\mathbb{ V}\mathbb{G}_1\mathbb{G}_2$  considered in this paper are Cartesian product sets.

In the rest of the paper, CPT-based shortfall risk always refers to problem \eqref{eq:rho-CPT}.
The next proposition states that the two definitions in \eqref{eq:PRGRS-CPT-equiv} and \eqref{eq:PRGRS-CPT} are equivalent.
\begin{proposition}
\label{prop:PGSR-robust-equal}
Let $\Pi_{\mathbb{ V}\mathbb{G}_1\mathbb{G}_2}(\xi)$ and $\rho_{\mathbb{ V}\mathbb{G}_1\mathbb{G}_2}(\xi)$ be defined as in \eqref{eq:PRGRS-CPT-equiv} and \eqref{eq:PRGRS-CPT}, respectively. 
Then 
\begin{equation}
\label{eq:PRGSR-equivalent}
\Pi_{\mathbb{ V}\mathbb{G}_1\mathbb{G}_2}(\xi)
=\rho_{\mathbb{ V}\mathbb{G}_1\mathbb{G}_2}(\xi).
\end{equation}
\end{proposition}

In Section~\ref{exa:finite}, we give a simple example to illustrate the relationship between the worst-case CPT-based shortfall and PRGSR-CPT and highlight the advantages of PRGSR-CPT.
We further discuss the case that the reference point $B$ is nonzero.
Empirical investigations of prospect theory largely take the status quo as the reference point; 
for the case of multiple reference points, see \citep{markle2018goals}.

The PRGSR-CPT framework covers two important preference robust risk measures as special cases.
\begin{itemize}
\item[1.]{\bf Preference Robust Utility-Based Shortfall Risk Measure.}
If the ambiguity set of the distortion functions only contains a singleton identity
function, i.e., $g_i(p)=p$, for $p\in [0,1]$, $i=1,2$, 
then $\Pi_{{\cal V} }(\xi):= \sup_{v\in {\cal V}} \rho_{v,g_1,{g_2}}(\xi)$, 
where ${\cal V}$  denotes the coordinate projection of $\mathbb{V}{\mathbb{G}_1\mathbb{G}_2}$ onto $\mathscr{V}$ with identity distortion.
It follows from Proposition~\ref{prop:PGSR-robust-equal} that $\Pi_{{\cal V}}(\xi)$
is equivalent to the preference robust utility-based shortfall $\rho_{{\cal V}}(\xi):=\inf_{x\in \R} \left\{x:\sup_{v\in {\cal V}} \bbe[v(\xi-x)] \leq 0\right\}$,
which has been fully investigated by \cite{DGX22}.

\item[2.]{\bf  Preference Robust Distortion Risk Measure.}
If the ambiguity set of the value functions consists of a singleton identity function, 
i.e., $v(t)=t$, for $ t\in \R$, and $g_1(p)=1-{g_2}(1-p)$, then $\bbe_{g_1{g_2}}[v(\xi)]=\bbe_{g_1}[\xi]$.
The explicit integral representation of $\bbe_{g_1 g_2}[v(\xi)]$ follows from \eqref{eq:defi-bbe-wwv} and integration by parts.
Let ${\cal G}_1$ be the ambiguity set of $g_1$.
Define $\rho_{{\cal G}_1}(\xi):=\inf_{x\in \R} \left\{ x:  \sup_{g_1\in {\cal G}_1} 
\bbe_{g_1}[\xi-x] \leq 0 \right\}$.
Hence, $\rho_{{\cal G}_1}(\xi)$ may be interpreted as a shortfall-type risk measure based on Yaari's dual theory of choice.
It is well-known that $\bbe_{g_1}[\xi]$ defines a coherent risk measure when $g_1(p)$ is convex.
In this case, the translation invariance property holds, which implies that $\bbe_{g_1}[\xi-x]=\bbe_{g_1}[\xi]-x$, 
and thus $\rho_{g_1}(\xi)=\inf_{x\in \R}\{x:\bbe_{g_1}[\xi-x]\leq 0\}=\bbe_{g_1}[\xi]$, 
which coincides with the classical DRM.
Finally, letting $\Pi_{{\cal G}_1}(\xi):= \sup_{g_1\in {\cal G}_1}\rho_{g_1}(\xi)$,
we obtain that $\Pi_{{\cal G}_1}$ collapses to the preference robust DRM $\sup_{g_1\in {\cal G}_1}\int_{\R} tdg_1(F_{\xi}(t))$ (see \citep{WaX23}).
Moreover, by Proposition~\ref{prop:PGSR-robust-equal}, $\Pi_{{\cal G}_1}(\xi)=\rho_{{\cal G}_1}(\xi)$, 
and consequently $\rho_{{\cal G}_1}(\xi)$ also collapses to the preference robust DRM.
\end{itemize}

\subsection{Properties of PRGSR-CPT}
\label{sec:properties}
Next, we discuss the properties of the newly introduced PRGSR-CPT.
First, it is known that the CPT-based shortfall $\rho_{v,g_1,g_2}(\cdot)$ is a law-invariant monetary risk measure (satisfying monotonicity and translation invariance, see \citep[Proposition 2.6]{MaC18}). 
Since these properties are preserved under the supremum operation, $\Pi_{\mathbb{ V}\mathbb{G}_1\mathbb{G}_2}(\cdot)$ is also a law-invariant monetary risk measure. 
We can therefore conclude, via Proposition~\ref{prop:PGSR-robust-equal}, 
that the PRGSR-CPT $\rho_{\mathbb{ V}\mathbb{G}_1\mathbb{G}_2}(\cdot)$ is a law-invariant monetary risk measure. 
Second, it might be interesting to discuss whether PRGSR-CPT is a convex (or coherent) risk measure.
Recall that \cite[Corollary 2.8, Theorem 2.10]{MaC18} derive conditions under which the CPT-based shortfall is a convex (or coherent) risk measure. 
The following two propositions state that if we use the set of all the value and distortion functions satisfying these conditions to define the joint ambiguity set
$\mathbb{V}\mathbb{G}_1\mathbb{G}_2$, 
then the resulting PRGSR-CPT is also a convex (or coherent) risk measure.

Define ${\cal R}_{\rm con}:=\{\rho_{v,g_1,{g_2}}:(v,g_1,g_2)\in \mathscr{V}\times \mathscr{G}\times \mathscr{G}, \; \rho_{v,g_1,{g_2}} \inmat{ is a convex risk measure}\}$ and  $\rho_{{\cal R}_{\rm con}}(\xi):=\sup_{\rho_{v,g_1,{g_2}}\in {\cal R}_{\rm con}}\rho_{v,g_1,{g_2}}(\xi)$.
The following proposition establishes a characterization of convex CPT-based shortfall risk measures.
\begin{proposition}
\label{Prop:Convex}
Let ${\cal V}_{1}\subset {\mathscr V}$ be the set of all value functions $v$ that are convex on $\R_-$ and convex on $\R_+$, 
let ${\cal G}_{\rm con}\subset {\mathscr G}$ be the set of all convex distortion functions, 
and let ${\cal G}_{\rm cav}\subset {\mathscr G}$ be the set of all concave distortion functions.
Define
\begin{equation}
\label{eq:constraint-vww}
\begin{aligned}
\mathbb{V}\mathbb{G}_1\mathbb{G}_2^{\rm con}
:= \left\{(v,g_1,{g_2})\in {\cal V}_1\times {\cal G}_{\rm con}\times {\cal G}_{\rm cav}: \frac{v'_+(0)}{v'_-(0)} \geq   \sup_{p\in(0,1)}   \frac{(g_1)'_-(p)}{(1-{g_2}(1-\cdot))'_+(p)}  \right\},
\end{aligned} 
\end{equation}
where $v'_-(0)$ and $v'_+(0)$ are the left and right derivatives of $v$ at $0$, and
$(g_1)'_-(p)$ and $(1-{g_2}(1-\cdot))'_+(p)$ are the left and right derivatives of $g_1$ and $1-{g_2}(1-\cdot)$ at point $p$, respectively.
Then $(v,g_1,g_2)\in\mathbb V\mathbb G_1\mathbb G_2^{\rm con} \Leftrightarrow \rho_{v,g_1,g_2}(\cdot)$ is a convex risk measure.
Consequently, $\Pi_{\mathbb{V}\mathbb{G}_1\mathbb{G}_2^{\rm con}}(\xi)=\rho_{{\cal R}_{\rm con}}(\xi)$.
\end{proposition}

Define ${\cal R}_{\rm coh}:=\{\rho_{v,g_1,{g_2}}: (v,g_1,g_2)\in \mathscr{V}\times \mathscr{G}\times \mathscr{G},\; \rho_{v,g_1,{g_2}}\;\inmat{is a coherent risk measure}\}$
and $\rho_{{\cal R}_{\rm coh}}(\xi)\hspace{-0.1cm}:=\sup_{\rho_{v,g_1,{g_2}}\in {\cal R}_{\rm coh}}\rho_{v,g_1,{g_2}}(\xi)$.
We now turn to the coherent case and provide a characterization.
\begin{proposition}
\label{prop:robust-RM-Pi}
Define ${\cal V}_{\rm coh}:=\left\{v\in {\mathscr V} :\, v(x)=a_1x \inmat{ for } x\geq 0, v(x)=a_2x \inmat{ for } x<0, a_1,a_2>0 \right \}$.
Let ${\cal G}_{\rm con}$ and ${\cal G}_{\rm cav}$ be defined as in Proposition~\ref{Prop:Convex}, 
and ${\mathbb{V}\mathbb{G}_1\mathbb{G}_2^{\rm coh}}:=\{(v,g_1,g_2):v\in{\cal V}_{\rm coh},
(g_1,g_2)\in{\cal G}_{\rm pair}(v)\}$,
where ${\cal G}_{\rm pair}(v)$ consists of all pairs $(g_1,g_2) \in {\cal G}_{\rm con}\times {\cal G}_{\rm cav}$ such that $\Upsilon(p):=a_1(1-g_2(1-p))-a_2g_1(p)$ is increasing in $p\in [0,1]$.
Then $(v,g_1,g_2)\in \mathbb V\mathbb G_1\mathbb G_2^{\rm coh} \Leftrightarrow \rho_{v,g_1,g_2}(\cdot)$ is a coherent risk measure.
Consequently, $\Pi_{\mathbb{V}\mathbb{G}_1\mathbb{G}_2^{\rm coh}}(\xi)=\rho_{{\cal R}_{\rm coh}}(\xi)$.
\end{proposition}

It might be helpful to make some comments on the condition about $\Upsilon(p)$. Consider the case that $1-g_2(1-p)=g_1(p)$. 
Then $\Upsilon(p)=(a_1-a_2)g_1(p)$. 
Consequently, the increasing property is equivalent to $a_1-a_2>0$. 
Moreover, this condition is related to the inequality in \eqref{eq:constraint-vww}.
In fact, the inequality implies that $v'_+(0)(1-{g_2}(1-\cdot))'_+(p)-v'_-(0)(g_1)'_-(p) \geq 0$ for all $p\in (0,1)$.
Since $g_2$ is concave over $[0,1]$, it follows that $1-g_2(1-p)$ is convex over $[0,1]$.
Together with the convexity of $g_1$, we assert that $v'_+(0)(1-{g_2}(1-p))-v'_-(0)g_1(p)$ is continuous and differentiable almost everywhere over $[0,1]$.
For every \(p\in[0,1]\), we have 
$v'_+(0)(1-{g_2}(1-p))-v'_-(0)g_1(p) 
 = \int_{0}^{p} (v'_+(0)(1-{g_2}
(1-\cdot))'_+(t)-v'_-(0)(g_1)'_-(t))dt$.
This implies that the function $v'_+(0)(1-{g_2}(1-p))-v'_-(0)g_1(p)$ is increasing in $p\in [0,1]$.

\begin{remark}
\label{rem:risk-taking}
The results in Propositions~\ref{Prop:Convex} and \ref{prop:robust-RM-Pi} are related to the notions of strong risk aversion and strong risk taking.
A prospect is defined as a random variable with a finite probability distribution, represented by $A:=(p_1,x_1;\cdots;p_n,x_n)$, 
where probability $p_j$ is assigned to outcome $x_j$ for $j\in [n]$.
The outcomes are ranked in decreasing order, i.e., $x_1\geq \cdots\geq x_n$.
\begin{itemize}
\item[(i)]
We begin by recalling the notion of ``strong risk aversion (SRA)'' introduced by \cite{RoS70}. 
A DM's preference $\succeq$ is said to exhibit strong risk aversion if, 
for any prospect $A=(p_1,x_1;\cdots;p_n,x_n)$ and $\delta>0$, $A\succeq (p_1,x_1;\cdots;p_i,x_i+\frac{\delta}{p_i};\cdots;p_j,x_j-\frac{\delta}{p_j};\cdots;p_n,x_n)$, 
for some $j>i$ with $p_i,p_j>0$, where the parameter $\delta$ is chosen such that the rank ordering of outcomes is maintained.
There are some necessary conditions on $g_1,g_2$ and $v$ when the DM's preference is characterized by CPT with an SRA risk attitude. 
Specifically, it follows from \cite[Theorem 1]{ScZ08} that if SRA holds, 
then the gain-weighting function $g_2$ must be concave whereas the loss-weighting function $g_1$ must be convex.
The utility function is concave for losses and also concave for gains, 
and the following relationship is satisfied $\frac{v'_-(0)}{v'_+(0)}\geq \sup_{p\in (0,1)} \frac{(1-g_2(1-\cdot))'_+(p)}{(g_1)'_-(p)}$.
Next, we discuss how inequality in \eqref{eq:constraint-vww} can be derived when the DM's preference is ``strong risk taking (SRT)''.

Analogous to SRA, we can define SRT as follows.
A DM's preference $\succeq$ is said to exhibit SRT if, for all prospects $A=(p_1,x_1;\cdots;p_n,x_n)$ and $\delta>0$, it holds that 
\begin{equation}
\label{eq:ineq5}
\begin{aligned}
A \preceq \left(
p_1,x_1;\cdots;p_i,x_i+\frac{\delta}{p_i};
\cdots;p_j,x_j-\frac{\delta}{p_j};
\cdots;p_n,x_n\right),
\end{aligned}
\end{equation}
where $p_i,p_j>0$.
The parameter $\delta$ is chosen such that the rank ordering of outcomes is maintained.
We assume that the gain-weighting function $g_2$ is concave and the loss-weighting function $g_1$ is convex. 
The utility function is convex for losses and also convex for gains, and the relationship \eqref{eq:ineq5} is satisfied.
In that case, \eqref{eq:ineq5} is equivalent to $\frac{v'_+(0)}{v'_-(0)}\geq \sup_{s\in(0,1)} \frac{(g_1)'_-(s)}{(1-{g_2}(1-\cdot))'_+(s)}$, coinciding with the inequality in \eqref{eq:constraint-vww} (see Section~\ref{SC:risk-seeking} for the proof of the equivalence).

\item[(ii)]
The coherence of monetary risk measures can be deduced from their positive homogeneity and convexity.
Loosely speaking, if you double your portfolio, then you double your risk. 
In financial risk management, positive homogeneity implies that the risk of a position is proportional to its size.
Specifically, for prospects $A=(p_1,x_1;\cdots;p_n,x_n)$ and any $r>0$, ${\rm CPT}(rA) =r{\rm CPT}(A)$, where ${\rm CPT}(\cdot)$ is defined as in Section~\ref{SC:risk-seeking}.
Thus, the positive homogeneity of the value function holds (see Section~\ref{SC:risk-seeking} for details).
\end{itemize}
\end{remark}

\subsection{Classical Value Function and Distortion Functions in CPT}
\label{sec:Literature_CPT}
In Proposition~\ref{Prop:Convex}, the value function is required to be convex over $\R_-$ and convex over $\R_+$, 
and the distortion function $g_1$ is required to be convex over $[0,1]$ whereas $g_2$ is required to be concave over $[0,1]$.
These requirements are usually not fulfilled in practice.
A variety of methods have been proposed to determine the shapes of the value function
and the distortion function in the CPT literature.
\citet{KaT79} assume that the value function is concave for gains and convex for losses. 
The assumption is supported by the experimental results in \cite{Abd00}.
\citet{TvK92} claim that individuals exhibit risk seeking for gains and risk aversion for losses of low probability, and risk aversion for gains and risk seeking for losses of high probability.
They infer that this pattern of risk attitude suggests an inverse $S$-shaped nonlinear transformation of the probability scale as well as an $S$-shaped utility function.
Preference conditions concerning the shape of the probability weighting function were formulated by \cite{TvW95,WuG96,WuG98,Pre98}.

One of the approaches to elicit these functions is to consider prespecified parametric forms for them (see \citep{TvK92,CaH94,TvF95}).
The experimental values of parameters $\alpha$, $\beta$, $\lambda$ of the value function and $\gamma_-$, $\gamma_+$, $\delta_-$, $\delta_+$ of the distortion functions are shown in Table~\ref{tab:Distortion} based on \cite{WaX23,RiW06,Abd00,BlP00,Sto06,BBEJ22}.
In the table, Type (I) is Tversky and Kahneman's single-parameter form, 
and Type (II) is the Lattimore et al.~function \citep{LBW92}, 
which takes into account curvature and elevation through parameters $\gamma$ and $\delta$, respectively.
The parameter $\gamma$ determines whether the weighting function exhibits an inflection point and the direction of the inflection (e.g., concave-convex).
The estimates of $\delta$ in \cite{Abd00} show that the probability weighting function for median data exhibits more elevation for losses ($\delta=0.84$) than for gains ($\delta=0.65$).
Type (III) is proposed by~\cite{Pre98}.

\begin{table}[hbt!]
\centering
\caption{Some typical parametric forms of the value function $v$ and the distortion function $g$, ``$-$" means the reference does not report the parameter}
\setlength\extrarowheight{1.6pt}
\begin{adjustbox}{width=1\textwidth}
\small{
\begin{threeparttable}
\begin{tabular}{l c c l l}
\toprule
References  &    
Forms of $v(t)$ &  Parameters $\alpha$, $\beta$, $\lambda$ & Forms of $g(t)$ & Parameters $\gamma_-$, $\gamma_+$, $\delta_-$, $\delta_+$ \\
\midrule
\multirow{3}{*}{\citet{TvK92}} & \multirow{3}{*}{$\left\{ 
\begin{array}{ll}
x^{\alpha} &~~\inmat{if}~~ x\geq 0\\
-\lambda (-x)^{\beta} &~~\inmat{if}~~ x<0
\end{array}
\right.$} & $\lambda=2.25$
& \multirow{3}{*}{(I)~~$\dfrac{t^{\gamma}}{\left(t^{\gamma}+(1-t)^{\gamma}\right)^{\frac{1}{\gamma}}}$}
&  \\
&& losses: $\beta=0.88$ &&  losses: $\gamma_-=0.69$\\
&&  gains: $\alpha=0.88$ && gains: $\gamma_+=0.61$ \\
\midrule
\multirow{5}{*}{\citet{Abd00}} & \multirow{5}{*}{$\left\{ 
\begin{array}{ll}
x^{\alpha} &~~\inmat{if}~~ x\geq 0\\
-\lambda (-x)^{\beta}&~~\inmat{if}~~ x<0
\end{array}
\right.$} & 
&  \multirow{3}{*}{(I)~~$\dfrac{t^{\gamma}}{\left(t^{\gamma}+(1-t)^{\gamma}\right)^{\frac{1}{\gamma}}}$}
&  \\
&& 
&&  losses: $\gamma_-=0.70$\\
&&  losses: $\beta=0.92$ && gains: $\gamma_+=0.60$ \\
&& gains: $\alpha=0.89$ 
& \multirow{2}{*}{(II)~~$\dfrac{\delta t^{\gamma}}{\delta t^{\gamma}+(1-t)^{\gamma}}$} & 
losses: $\gamma_-=0.65$, $\delta_-=0.84$\\
&&&& gains: $\gamma_+=0.60$, $\delta_+=0.65$\\
\midrule
\multirow{7}{*}{\citet{WuG96}} & \multirow{7}{*}{$x^{\alpha}$} & \multirow{3}{*}{gains: $\alpha=0.52$}
& \multirow{3}{*}{(I)~~$\dfrac{t^{\gamma}}{\left(t^{\gamma}+(1-t)^{\gamma}\right)^{\frac{1}{\gamma}}}$} & \multirow{3}{*}{gains: $\gamma_+=0.71$}\\
&&&&\\
&&&&\\
&& $-$ & \multirow{2}{*}{(II)~~$\dfrac{\delta t^{\gamma}}{\delta t^{\gamma}+(1-t)^{\gamma}}$} & \multirow{2}{*}{gains: $\gamma_+=0.68$, $\delta_+=0.84$} 
\\
&&&&\\
&& \multirow{2}{*}{gains: $\alpha=0.48$}&\multirow{2}{*}{(III)~~$\exp(-\delta (-\ln{t})^{\gamma})$}& \multirow{2}{*}{gains: $\gamma_+=0.74$, $\delta_+=1$}\\
&&&&\\
\midrule    
\multirow{3}{*}{\citet{CaH94}} & \multirow{3}{*}{$x^{\alpha}$} & \multirow{3}{*}{gains: $\alpha=0.37$ \footnote{\cite{CaH94} do not report $\alpha$, but \cite{WuG96} repeat their estimation procedure using data generously provided by Teck-Hua Ho and found $\alpha=0.37$, see \citep[p.~1686]{WuG96}.}}
& \multirow{3}{*}{(I)~~$\dfrac{t^{\gamma}}{\left(t^{\gamma}+(1-t)^{\gamma}\right)^{\frac{1}{\gamma}}}$}
&  \multirow{3}{*}{gains: $\gamma_+=0.56$}\\
&&&&\\
&&&&\\
\midrule
\multirow{3}{*}{\citet{TvF95}} & \multirow{3}{*}{$x^{\alpha}$} & \multirow{3}{*}{gains: $\alpha=0.88$} 
& \multirow{3}{*}{(II)~~$\dfrac{\delta t^{\gamma}}{\delta t^{\gamma}+(1-t)^{\gamma}}$}
& \multirow{3}{*}{gains: $\gamma_+=0.69$, $\delta_+=0.77$} \\
&&&&\\
&&&&\\
\midrule
\multirow{3}{*}{\citet{GoW99}} & \multirow{3}{*}{$x^{\alpha}$} & \multirow{3}{*}{gains: $\alpha=0.49$} 
& \multirow{3}{*}{(II)~~$\dfrac{\delta t^{\gamma}}{\delta t^{\gamma}+(1-t)^{\gamma}}$}
& \multirow{3}{*}{gains: $\gamma_+=0.44$, $\delta_+=0.77$} \\
&&&&\\
&&&&\\
\midrule
\multirow{7}{*}{\citet{BlP00}} & \multirow{7}{*}{$x^{\alpha}$} & \multirow{7}{*}{ $\alpha=0.779$}
& \multirow{3}{*}{(I)~~$\dfrac{t^{\gamma}}{\left(t^{\gamma}+(1-t)^{\gamma}\right)^{\frac{1}{\gamma}}}$} & \multirow{3}{*}{$\gamma =0.674$}\\
&&&&\\
&&&&\\
&&  & \multirow{2}{*}{(II)~~$\dfrac{\delta t^{\gamma}}{\delta t^{\gamma}+(1-t)^{\gamma}}$} & \multirow{2}{*}{$\gamma=0.550$, 
$\delta=0.816$} \\
 &&&&\\
&& &\multirow{2}{*}{(III)~~$\exp(-\delta(-\ln{t})^{\gamma})$}& \multirow{1}{*}{$\gamma =0.533$, $\delta=1$}\\
&& && \multirow{1}{*}{$\gamma =0.534$, $\delta=1.083$}\\
\midrule
\end{tabular}
\end{threeparttable}
}
\end{adjustbox}
\label{tab:Distortion}
\end{table}

This means that our robust CPT-based shortfall risk measure is generally non-convex.
Indeed, this is a departure from \cite{DGX22}, 
which is based on the expected utility theory with convex disutility functions.
To facilitate computation of the robust shortfall, 
we will propose a bisection-based algorithm with piecewise McCormick relaxations to compute the robust shortfall risk measure.

\section{Finite-Dimensional Reformulation of the Inner Problem in PRGSR-CPT}
\label{sec:tractability}
In this section, we discuss computational schemes for the proposed PRGSR-CPT with the ambiguity set $\mathbb{V}\mathbb{G}_1\mathbb{G}_2$.
We begin by discussing how to compute the CPT-based shortfall $\rho_{v,g_1,g_2}(\xi)$.
Throughout this section, we restrict our discussion to the case that $\xi$ is discretely distributed. 
Let $\Xi:=\{\xi_1,\ldots,\xi_K\}$. 
\begin{assumption}
\label{ass:discrete}
$P$ is a discrete distribution with $P(\{\xi_k\})=p_k>0$, for $k=1,\cdots,K$.
\end{assumption}

In the case when $\xi$ is continuously distributed, we may view the discrete distribution as a sample average approximation (SAA). 
Convergence of the SAA in terms of the optimal value and optimal solutions can be established under some moderate conditions analogous to those in the literature such as \cite{DGX22,ZXW22,ZhX21}. 
We omit the details of the theoretical analysis as our focus is on preference elicitation in this paper and leave the nonconvex case for future research.
We will illustrate the performance of SAA in Section~\ref{sec:sensitivity} from a numerical perspective.

\begin{assumption}
\label{ass:discrete_Lip}
Each $v\in \mathscr{V}$ is Lipschitz continuous over $[\underline{x}, \bar{x}]$ with Lipschitz modulus being bounded by $L$.
\end{assumption}

The constant $L$ gives an upper bound on the absolute value of the marginal value $v'(t)$ over $[\underline{x},\bar{x}]$ (particularly on $v'_+(0)$ and $v'_-(0)$).
Assumption~\ref{ass:discrete_Lip} is needed when we measure the difference between 
the piecewise linear value function $v_J^*$ and the true value function $v^*$ 
but the Lipschitz modulus $L$ is not necessarily known in preference elicitation and numerical tests.
The assumption is satisfied by many value functions of practical interest, 
including,
e.g., the exponential form $v(x)= \lambda(1-e^{-\beta x})/\beta$ for $x<0$, 
$v(x)= (1-e^{-\alpha x})/\alpha$ for $x\geq 0$, where $\alpha>0$, $\beta>0$, and $\lambda\geq 1$.

\subsection{Computation of CPT-Based Shortfall}
\label{Sec:computation_gsr}
We begin by ordering the elements $\{\xi_{k}\}_{k=1}^{K}$ in increasing order  ($\xi_{k_1}\leq \xi_{k_2}$, $k_1<k_2$), and obtaining $\xi_{k_1}-x\leq \xi_{k_2}-x$ with $k_1<k_2$ for a fixed $x$.
By convention, for the outcomes satisfying $\xi_{k}-x<0$, 
we label them with indices $k\in \{1,\cdots,K_0(x)-1\}$ and for those $\xi_k-x\geq 0$, we label them with $k\in \{K_0(x),\cdots,K\}$.
Here, $K_0(x)=K+1$ if all outcomes are negative,
and empty sums are understood to be zero.
Consequently, we can rewrite the integral
$\bbe_{g_1g_2}[v(\xi-x)]
=\int_{-\infty}^{0} v(t) d{g_1}(F_{\xi-x}(t))
-\int_{0}^{\infty} v(t) d({g_2}(1-F_{\xi-x}(t))) $
as
\begin{subequations}
\begin{align}
 \bbe_{{g_1}{g_2}}[v(\xi-x)] 
& = \sum_{k=1}^{K_0(x)-1} v(\xi_k-x) ({g_1}(q_k)-{g_1}(q_{k-1}))\nonumber \\
& \quad + \! \! \sum_{k={K_0(x)}}^{K} \! \!
v(\xi_k-x)({g_2}(1-q_{k-1})-{g_2}(1-q_k)) \label{eq:saa_ww} \\
& =: h_{v,g_1,g_2}^1(x)+h_{v,g_1,g_2}^2(x), 
\label{eq:sum_h-expectation}
\end{align}
\end{subequations}
where $F_{\xi-x}(t)=\sum_{k=1}^Kp_k\mathds{1}_{t\geq \xi_k-x}$, 
$\mathds{1}_{t\geq \xi_k-x}$ is the indicator function of $[\xi_k-x,\infty)$, 
$q_{0} :=0,$ and $q_k :=\sum_{j=1}^{k}p_j$ for $k\in [K]$.
We refer readers to a concrete example of  calculating $\bbe_{g_1g_2}[v(\xi-x)]$ in Section~\ref{SC:exa-bbe}.
Since $v$, $g_1$, and $g_2$ are monotonically increasing, $\bbe_{{g_1}{g_2}}[v(\xi-x)]$ is also monotonically decreasing w.r.t.~$x$.
This monotonicity property prompts us to use a bisection method (see \citep[Section 2.1]{BFB15}) to compute the CPT-based shortfall.
We refer readers to Section~\ref{sec:Bisection-norobust} for the bisection method used to compute $\rho_{v,g_1,g_2}(\xi)$.
If there exists a computationally manageable reformulation of the constraint function $\sup_{(v,g_1,{g_2})\in\mathbb{V}\mathbb{G}_1\mathbb{G}_2} \mathbb{E}_{g_1{g_2}}[v(\xi-x)]$ in \eqref{eq:PRGRS-CPT}, 
then the PRGSR-CPT $\rho_{\mathbb{V}\mathbb{G}_1\mathbb{G}_2}(\xi)$ in \eqref{eq:PRGRS-CPT} can also be solved by the bisection method as described in Algorithm~\ref{algorithm:GS}.

\subsubsection{Comparison with Alternative Shortfalls}
To illustrate the differences among the utility-based shortfall, 
the DRM-based shortfall,
the RDEU-based shortfall,
and the CPT-based shortfall, 
we consider a simple case where there is no ambiguity in the value and distortion functions. 
\begin{example}[Calculation of Four Shortfalls]
\label{exa:comparison_three}
We consider the case that the true $S$-shaped value function is given by $v^*(t)=(2(1-e^{-8t}))/8$ for $t\in [0,0.5]$ and $v^*(t)=(e^{3t}-1)/3$ for $t\in [-0.5,0]$.
The true distortion functions are $g_i^*(p)=p^{\gamma_i}/(p^{\gamma_i}+(1-p)^{\gamma_i})$,
$i=1,2$, where $\gamma_1=0.89$ and $\gamma_2=0.6$.
Let $\xi=(0.01,0;0.66,0.4;0.33,0.5)$.
We consider the following four models:
$\rho_{v^*}(\xi):=\inf_{x\in \R} \{ x:\bbe[v^*(\xi-x)] \leq 0\}$,
$\rho_{g_1^*}(\xi):=\inf_{x\in \R} \{ x:\bbe_{g_1^*}[\xi-x] \leq 0\}$,  
$\rho_{v^*,g^*_1}(\xi):=\inf_{x\in \R} \{ x:  \bbe_{g_1^*}[v^*(\xi-x)] \leq 0\}$,
and 
$\rho_{v^*,g_1^*,g_2^*}(\xi):=\inf_{x\in \R} \{ x: \bbe_{g_1^*g_2^*}[v^*(\xi-x)] \leq 0\}$.

Let 
$H_{v^*}(x):= p_1v^*(\xi_1-x)+p_2v^*(\xi_2-x) +p_3v^*(\xi_3-x)$,
$H_{g_1^*}(x):=g_1^*(q_1)(\xi_1-x)
  +(g_1^*(q_2)-g_1^*(q_1))(\xi_2-x)
+(1-g_1^*(q_2))(\xi_3-x)$,
$H_{v^*,g_1^*}(x):= g_1^*(q_1)v^*(\xi_1-x)
  +(g_1^*(q_2)-g_1^*(q_1))v^*(\xi_2-x)
  +(1-g_1^*(q_2))v^*(\xi_3-x)$,
and $H_{v^*,g_1^*,g_2^*}(x) = \sum_{k=1}^{K_0(x)-1}v^*(\xi_k-x)s_k^-
  +\sum_{k=K_0(x)}^3v^*(\xi_k-x)s_k^+$,
where $q_0=0$, $q_1=p_1$, $q_2=p_1+p_2$, $q_3=1$, $s_k^-={g^*_1}(q_k)-{g^*_1}(q_{k-1})$ and $s_k^+= {g^*_2}(1-q_{k-1})-{g^*_2}(1-q_k)$.

The resulting shortfalls are
$\rho_{v^*}(\xi)=\inf_{x\in \R} \{ x: H_{v^*}(x)\leq 0\}=0.4441$,
$\rho_{g_1^*}(\xi)=\inf_{x\in \R} \{ x: H_{g_1^*}(x)\leq 0\}=0.4282$,
$\rho_{v^*,g_1^*}(\xi)=\inf_{x\in \R} \{ x: H_{v^*,g_1^*}(x)\leq 0\}=0.4451$,
and
$\rho_{v^*,g_1^*,g_2^*}(\xi)=\inf_{x\in \R} \{ x: H_{v^*,g_1^*,g_2^*}(x)\leq 0\}=0.4491$.
We can see that $\rho_{v^*,g_1^*,g_2^*}(\xi)$ is larger than $\rho_{v^*}(\xi)$ due to the effects of the distortion functions $g_1$ and $g_2$: $s_1^-=0.0165>0.01=p_1$, $s_2^-=
0.6361<0.66=p_2$, $s_3^- = 1-0.6526 = 0.3474>0.33=p_3$, and $s_1^+ = 0.0597>0.01=p_1$, $s_2^+ =0.5450<0.66=p_2$, $s_3^+=0.3953>0.33=p_3$.
Therefore, it is inappropriate to replace $\rho_{v,g_1,g_2}(\xi)$ with $\rho_{v}(\xi)$, $\rho_{g_1}(\xi)$ or $\rho_{v,g_1}(\xi)$ if the DM's preference can only be described by CPT-based shortfall instead of EUT-based, DRM-based, or RDEU-based shortfalls.
\end{example}

\subsection{A Joint Ambiguity Set for Value and Distortion Functions}
\label{sec3:ambiguity_VW}
A key ingredient of the PRGSR-CPT in \eqref{eq:PRGRS-CPT} is the construction of 
the ambiguity set $\mathbb{V}\mathbb{G}_1\mathbb{G}_2$ of the value and distortion functions since it concerns the collection of a DM's preference information and 
the computation of PRGSR-CPT.
There are two types of information: 
(a) general information, such as risk taking over losses and risk aversion over gains considered in the literature of behavioral economics ($S$-shaped value functions), 
and overweighting of events from impossible to possible or possible to certainty and underweighting of those from possible to more possible by many empirical tests (inverse $S$-shaped distortion functions);
(b) user-specific information, such as pairwise comparisons between a risky lottery and a certain lottery under CPT.
Here, we assume that the inflection point (concave-convex) of the distortion function is unknown, 
and thus in the following preference elicitation, we do not impose the inverse $S$-shape constraint of the distortion functions.
Here we describe how to use the pairwise comparisons of lotteries to elicit $v$, $g_1$, and $g_2$.
Let $Z_1^m$ be a risky lottery and $Z_2^m$ be a certain lottery for $m\in [M]$:
\begin{equation}
\label{eq:questions}
Z_1^m=\left\{\begin{aligned}
r_1^m \; & \mbox{w.p. } 1-p^m\\
r_3^m \; & \mbox{w.p. } p^m
\end{aligned}
\right.
\quad \mbox{and}\quad Z_2^m=r_2^m,
\end{equation}
where parameters $r_1^m, r_2^m, r_3^m$ satisfy $r_1^m<r_3^m$, 
and $r_1^m,r_2^m,r_3^m\in [\underline{x},\bar{x}]$, and $p^m\in (0,1)$.
These parameters remain undetermined throughout Section~\ref{sec:tractability},
and we will discuss in Section~\ref{sec:Estimating-dis} how to identify them through preference elicitation with the DM.
We use a superscript $m$ to denote the $m$-th query. 
By replacing $\xi-x$ with the lottery $Z_1^m$ in \eqref{eq:saa_ww}--\eqref{eq:sum_h-expectation} and setting $q_0 :=0$, $q_1:=1-p^m$, $q_2 :=1$, $K:=2$, we obtain the distorted expected value $\bbe_{{g_1g_2}}[v(Z_1^m)]$ as
$v(r_1^m)g_1(1-p^m)+v(r_3^m)(1-g_1(1-p^m))$
if $r_1^m,r_3^m<0$;
$v(r_1^m)g_1(1-p^m)+v(r_3^m)g_2(p^m)$
if $r_1^m<0\leq r_3^m$; and
$v(r_1^m)(1-g_2(p^m))+v(r_3^m)g_2(p^m)$
if $r_1^m,r_3^m\geq0$.
Combining these three cases, we obtain
\begin{equation}
\label{eq:distorted_exp-2}
\begin{aligned}
\bbe_{g_1g_2}[v(Z_1^m)]
={}&v(r_1^m)\Big(
g_1(1-p^m)\mathds{1}_{r_1^m<0} +(1-g_2(p^m))\mathds{1}_{r_1^m\geq0}\Big)\\
&+v(r_3^m)\Big(
(1-g_1(1-p^m))\mathds{1}_{r_3^m<0} +g_2(p^m)\mathds{1}_{r_3^m\geq0}\Big),
\end{aligned}
\end{equation}
and $\bbe_{{g_{1}g_{2}}}[v(Z_2^m)]=v(r_2^m)$.
Let $\mathscr{V}_{\rm S}\subset \mathscr{V}$ be the set of all value functions with inflection point $0$, 
which are concave over $\R_+$ and convex over $\R_-$.
Let $\mathscr{G}_{\rm invS}^i\subset \mathscr{G}$ be the set of the distortion functions $g_i$, for $i=1,2$.
For each $m\in[M]$ and $(v, g_{1},\;g_{2})\in \mathscr{V}_{\rm S}\times
\mathscr{G}_{\rm invS}^1 \times\mathscr{G}_{\rm invS}^2$, 
define $H_m({v},g_{1},g_{2}):= h_m\big(\bbe_{g_{1}g_{2}}[v(Z_1^m)] - \bbe_{g_{1}g_{2}}[v(Z_2^m)]\big)$.
The joint ambiguity set of the value functions and the distortion functions is defined as follows:
\begin{equation}
\label{eq:ambiguity set}
\begin{aligned}
\mathbb{V}\mathbb{G}_1\mathbb{G}_2 = \Big\{(v,g_1,g_2)\in \mathscr{V}_{\rm S}\times
\mathscr{G}_{\rm invS}^1 \times\mathscr{G}_{\rm invS}^2: 
 H_m({v},g_{1},g_{2}) \leq 0,\; m\in[M]\Big\},
\end{aligned}
\end{equation}
where $h_m=1$ if $Z_2^m$ is preferred to $Z_1^m$ by the DM whose preference is described by CPT, and $h_m=-1$ otherwise.
In practice, we need to design questionnaires to elicit the function $v$ and distortion functions $g_1$ and $g_2$ efficiently and determine the parameters $r_1^m,r_2^m,r_3^m,p^m$ to construct a specific ambiguity set (denoted by $(\mathbb{V}\mathbb{G}_1\mathbb{G}_2)_M$ in Section~\ref{sec:Estimating-dis}).
Three criteria prevail in constructing the ambiguity set $\mathbb{V}\mathbb{G}_1\mathbb{G}_2$ of the value and distortion functions.

\begin{remark}
[Three Criteria for Constructing $\mathbb{V}\mathbb{G}_1\mathbb{G}_2$]
First, $\mathbb{V}\mathbb{G}_1\mathbb{G}_2$ should contain the unknown true value and distortion functions or functions sufficiently close to them.
Second, the size of $\mathbb{V}\mathbb{G}_1\mathbb{G}_2$ should decrease as more queries are asked.
Third, $\mathbb{V}\mathbb{G}_1\mathbb{G}_2$ should admit a finite-dimensional reformulation of the inner problem associated with PRGSR-CPT.
\end{remark}

\begin{remark}[Preference Inconsistency]
To ensure $\mathbb{V}\mathbb{G}_1\mathbb{G}_2\neq \emptyset$, we implicitly assume 
that the DM's preferences are consistent. 
In practice, however, 
inconsistency may occur particularly in the process of preference elicitation either because the DM's preferences are intrinsically inconsistent 
or there are some elicitation errors (see \citep{BeO13,HZXZ22,ChL26}).
Here we give a comprehensive treatment of the issue. 
We begin, in Section~\ref{Sec:inconsistency},  
by considering three types of inconsistencies related respectively to the evaluation of CPT value of lotteries,
the perception of lottery outcomes and the DM's erroneous responses to queries, 
and discuss in detail how we may tackle these inconsistencies.
Next, we move on to show, in Example~\ref{exa:PRE_INC} of Section~\ref{Sec:inconsistency},
the consequence of using an incorrect model, 
such as EUT, DRM, or RDEU, to construct an ambiguity set of the value and distortion functions when the DM's preferences can only be described by CPT.
Finally, in Section~\ref{sec:illustration_Inconsis}, 
we  investigate the consequence of mistakenly using EUT-based polyhedral cuts to elicit the preferences of a CPT DM (where the DM’s preference can be characterized by CPT). 
\end{remark}

\subsection{Piecewise Linear Value Functions and Distortion Functions}
\label{sec:PLFs0}
In the forthcoming Theorem~\ref{thm:reformulation-constraint-function-VW}, 
we show that constructing  the ambiguity set using piecewise linear value and distortion functions yields the same constraint function values.
Accordingly, we elicit piecewise linear value and distortion functions in Section~\ref{sec:Estimating-dis}. 
Next, we introduce these piecewise linear value and distortion functions, 
which paves the way for the finite-dimensional reformulation of the inner problem associated with PRGSR-CPT in Section~\ref{Sec:tactability}.

It is well known that convex functions can be recovered by taking the supremum of their supporting affine functions (see \citep{BoV04}).
Analogous to \cite{ArD15} and \cite{WaX23}, 
we can use the support function approach to reformulate the robust CPT-based preference functional used in the definition of PRGSR-CPT as a bilinear programming problem when the value function is $S$-shaped (convex over $\R_-$ and concave over $\R_+$), 
and the distortion function is inverse $S$-shaped (concave over $[0,t^*]$ and convex over $[t^*,1]$ with an inflection point $t^*$).
However, we will take a different approach because $t^*$ is unknown and our model is much more complex (since the constraint function contains both the value and distortion functions) and using support function-based methods would require us to go through a number of combinations where the value function and the distortion functions take particular shapes. 
Specifically, we will demonstrate in Theorem~\ref{thm:reformulation-constraint-function-VW} that the worst-case value function and the worst-case distortion functions are both piecewise linear, 
and that their values at the breakpoints (where the functions switch from one piece to another) are determined by a set of points related to the scenarios of $\xi$.
This observation prompts us to use piecewise linear value and distortion functions to develop formulations.

Let ${\cal X}:=\{x_1,\cdots,x_{J}:\underline{x}=x_1<x_2<\cdots<x_{J}=\overline{x}\}$ be the set of breakpoints of a piecewise linear continuous value function $v_J$.
Let ${\bm a}:=(a_1,\cdots,a_{J-1})$, where  $a_j:=\frac{v(x_{j+1})-v(x_j)}{x_{j+1}-x_j}$ represents the slope of the value function over the interval $[x_j,x_{j+1}]$ for $j\in [J-1]$.
Denote ${\bm b}:=(b_1,\cdots,b_{J-1})$ as the vector of intercepts.
The piecewise linear value function $v_{J}:[\underline{x},\bar{x}]\to \R$ is defined as 
\begin{equation}
\label{eq:PLA_utility}
v_{J}(x)=\sum_{j=1}^{J-1} (a_j x +b_j) \phi_j(x)=(x{\bm a}+{\bm b})^\top{\bm \phi}(x),
\end{equation}
where ${\bm \phi}(x):=(\phi_1(x),\cdots,\phi_{J-1}(x))$, $\phi_{1}(x) := \mathds{1}_{[x_{1},x_{2}]}(x)$, and $\phi_{j}(x):=\mathds{1}_{(x_{j},x_{j+1}]}(x),
$ $j\in [J-1]\backslash{\{1\}}$ (see \citep{GXZ23}).
The coefficients ${\bm a}$ and ${\bm b}$ represent the slopes and intercepts of the linear pieces, respectively.
To develop a new polyhedral method to elicit the value and distortion functions, 
we begin by introducing the piecewise linear distortion functions based on vectors of increments,
and by normalizing the value functions over  $[\underline{x},0]$ and $[0,\bar{x}]$ respectively for the convenience of elicitation separately.

\subsubsection{Piecewise Linear Distortion Functions}
Let ${\cal S}:=\{t_1,\cdots,t_{N}:0=t_1<\cdots <t_{N}=1\}$ be the set of breakpoints
of the distortion function.
Define
$
{\bm h}_i^{{\cal S}}(t):=
 \big(\underbrace{1,\ldots,1}_{i-1},\frac{t-t_i}{t_{i+1}-t_{i}},0,\ldots,0\big) \in \R^{N-1}$, 
 for all $t\in (t_{i},t_{i+1}]$. 
for $i\in [N-1]$ and define ${\bm h}^{{\cal S}}(t) :=\sum_{i=1}^{N-1} {\bm h}_i^{{\cal S}}(t) \mathds{1}_{(t_i,t_{i+1}]}(t)$. 
Let ${\bm h}^{\cal S}(0):={\bm 0}$.
The piecewise linear distortion functions $g_{1N}:[0,1]\to [0,1]$ and $g_{2N}:[0,1]\to [0,1]$ are defined as
\begin{equation}
\label{eq:PLA_distortion}
g_{1N}(t):={\bm \sigma}_1^\top{\bm h}^{\cal S}(t)
\quad \inmat{and} \quad 
g_{2N}(t):={\bm \sigma}_2^\top{\bm h}^{\cal S}(t),
\end{equation}
where ${\bm \sigma}_i:=(\sigma_{i,1},\cdots,\sigma_{i,N-1})$ with $\sigma_{i,s}:=g_i(t_{s+1})-g_i(t_s)$ representing the increment of the distortion function $g_i$ over the interval $[t_s,t_{s+1}]$ for $i=1,2$ and $s\in [N-1]$.
Unlike \cite{HZXZ22}, where this kind of function is used to approximate continuous utility functions, 
we use them to approximate the distortion functions. 
To facilitate computation in later discussions, we reformulate ${\bm h}^{{\cal S}}(t) \in \R^{N-1}$ as 
\begin{equation}
\label{eq:g(t)-PRO}
\begin{aligned}
{\bm h}^{{\cal S}}(t) =  \Bigg(  l_{1}(t) + \sum_{s\in [N-1] \setminus \{1\}} \psi_{s}(t), 
l_{2}(t)  +  \sum_{s \in [N-1] \setminus \{1, 2\}}\psi_{s}(t),\cdots, l_{N-1}(t)\Bigg),
\end{aligned}
\end{equation} 
where $\psi_{s}(t):= \mathds{1}_{(t_{s},t_{s+1}]}(t)$ and $l_{s}(t):=\left(\frac{t - t_{s}}{t_{s+1} - t_{s}}\right)\psi_{s}(t)$, for $s\in [N-1]$.
Note that ${\bm h}^{{\cal S}}(t)$ can be formulated as a linear mapping subject to mixed-integer linear constraints (see \citep[p.~535]{HZXZ22}).

\subsubsection{Normalization, Piecewise Linear Approximation, and Scalarization of Value Functions}
\label{sec:PLFs}
In this section, we introduce four forms of the value function in terms of normalization, 
piecewise linear approximation (PLA) and scalarization, 
in preparation for the forthcoming discussions in Section~\ref{sec:Estimating-dis}.
Normalization of the  positive and negative parts of the value function is needed for
the modified polyhedral method to be presented in Section~\ref{sec:PM_modified}; 
PLA of the normalized value functions is required in the design of queries to elicit these functions in Section~\ref{sec:PM_modified};
and finally scalarization is applied to them to construct a holistic value function over $[\underline{x},\bar{x}]$ provided that the function values $v(\underline{x})$ and $v(\bar{x})$ are known. 
In the case that the values are unknown, 
we adopt specifically scaled value functions to be estimated in Section~\ref{sec:ratior}.
\begin{itemize}
\item[(i)] {\bf Normalized Value Functions $v_-$ and $v_+$}.
Let $v_-:[\underline{x},0]\to [-1,0]$ and $v_+:[0,\bar{x}]\to [0,1]$ be two normalized value functions defined as follows:
\begin{equation}
\label{eq:vminus_vplus}
\begin{aligned}
v_-(x)&:=-\frac{v(x)}{v(\underline{x})}, \inmat{ for }x\in [\underline{x},0] 
\quad \inmat{and} \quad 
v_+(x):=\frac{v(x)}{v(\bar{x})}, \inmat{ for } x\in [0,\bar{x}].
\end{aligned}
\end{equation}
Then $v_-(\underline{x})=-1$, $v_-(0)=v_+(0)=0$ (since $v(0)=0$), and $v_+(\bar{x})=1$.
Throughout the paper, 
we use subscripts ${-}$ and ${+}$ to signify the negative part and the positive part of the value function, respectively.
One of the challenges in eliciting the value function is that both $v(\bar{x})$ and $v(\underline{x})$ are unknown.

\item[(ii)] {\bf Piecewise Linear Value Functions $v_{J_-}:[\underline{x},0]\to [-1,0]$, $v_{J_+}:[0,\bar{x}]\to [0,1]$ and $v_J:[\underline{x},\bar{x}]\to [0,1]$}.
Let ${\cal X}_+:=\{x^+_{1},\cdots,x^+_{J_+}: 0=x^+_{1}<\cdots<x^+_{J_+}=\bar{x}\}$ be the set of breakpoints of $v_+$, and ${\bm h}^{{\cal X}_+}(x):= \Big(\underbrace{1,\ldots,1}_{j-1},\frac{x-x_j^+}{x_{j+1}^+-x_{j}^+},0,\ldots,0\Big) \in \R^{J_+-1}$, 
for all $x\in (x_{j}^+,x_{j+1}^+]$ and $j\in[J_+-1]$, and let ${\bm h}^{{\cal X}_+}(0):={\bm 0}$.
Define
\begin{equation}
\label{eq:PLA_utility_2}
v_{J_+}(x):={\bm v}_+^\top{\bm h}^{{\cal X}_+}(x),
\end{equation}
where ${\bm v}_+:=(v_{+,1},\cdots,v_{+,{J_+}-1})\in \R_{+}^{{J_+}-1}$ (the $j$-th element of a vector ${\bm v}_+$ is denoted by $v_{+,j}$),
$v_{+,j}:=v_+(x^+_{j+1})-v_+(x^+_{j})$ represents the increment of $v_+$ over the interval $[x^+_{j},x^+_{j+1}]$ for $j\in [{J_+}-1]$.
Note that $v_{J_+}$, as defined, maps from $[0,\bar{x}]$ to $[0,1]$.
Likewise, let ${\cal X}_-:=\{x_1^-,\cdots,x_{J_-}^-: \underline{x}=x_1^-<\cdots<x_{J_-}^-=0\}$ be the set of breakpoints of $v_-$.
Define ${\bm v}_-\in \R_{+}^{{J_-}-1}$, ${\bm h}^{{\cal X}_-}(x)\in \R^{{J_-}-1}$ and the piecewise linear function (PLF) of $v_-:[\underline{x},0]\rightarrow [-1,0]$ by
\begin{equation}
\label{eq:PLA_utility_minus}
v_{J_-}(x):={\bm v}_-^\top{\bm h}^{{\cal X}_-}(x)-1.
\end{equation}
Let ${\bm h}^{{\cal X}_-}(\underline{x}):={\bm 0}$.
We can also define ${\bm v}\in \R^{J-1}$, ${\bm h}^{\cal X}(x)\in \R^{J-1}$ and the PLF of $v:[\underline{x},\bar{x}]\to [0,1]$ by ${v}_J(x)={\bm v}^\top{\bm h}^{\cal X}(x)$.
Note that $v_{J_-}(x)$, $v_{J_+}(x)$, and $v_J(x)$ will be used in the design of Question Types~1--2 and 4 in Section~\ref{sec:Estimating-dis}.

\item[(iii)] {\bf Scaling of the Holistic Value Function $v$}.
As for the whole value function $v$ over $[\underline{x},\bar{x}]$, scaling $v$ by $\alpha>0$ will not affect the feasibility of $x$ in PRGSR-CPT \eqref{eq:PRGRS-CPT},
that is, $\bbe_{g_1g_2}[v(\xi-x)]\leq 0$ is equivalent to  $\bbe_{g_1g_2}[\alpha v(\xi-x)]\leq 0$ for any $\alpha>0$.
Note that $v$ is strictly increasing and $v(0)=0$ with $\underline{x}<0<\bar{x}$.
Hence, $v(\underline{x})<0$.
Note that the value function $v(x)$ can be recovered from $v_-$ and $v_+$ as follows:
\begin{equation}
\label{eq:v(x)_scale}
v(x)=\left\{\begin{aligned}
-v(\underline{x}) v_-(x) \; & \inmat{for}~ x\in [\underline{x},0],\\
v(\bar{x})v_+(x) \; & \inmat{for}~x\in [0,\bar{x}].
\end{aligned}
\right.
\end{equation}
However, $v(\underline{x}) $ and $v(\bar{x})$ are unknown.
Thus, we cannot recover the holistic value function $v$ over $[\underline{x},\bar{x}]$ 
if only the normalized functions $v_-(x)$ over $[\underline{x},0]$ and $v_+(x)$ over $[0,\bar{x}]$ are obtained.
If $r:=-\frac{v(\bar{x})}{v(\underline{x})}$ is obtained, 
we can construct the following scaled value function.

\item[(iv)] {\bf A Special Scaled Value Function $\tilde{v}:[\underline{x},\bar{x}]\to [-1,r]$}.
Take a specific scaling parameter $\bar{\alpha}:=-\frac{1}{v(\underline{x})}>0$, 
and scale the value function $v$ by $\bar{\alpha}$, 
i.e., $\tilde{v}(x):={\bar{\alpha}}{v}(x)$
for $x\in [\underline{x},\bar{x}]$ and thus $\tilde{v}(\underline{x})=-1$, 
$\tilde{v}(0)=0$, and $\tilde{v}(\bar{x})=r$.
By \eqref{eq:v(x)_scale}, we have 
\begin{equation}
\label{eq:value_tilde}
\tilde{v}(x):=\left\{\begin{array}{cl}
v_-(x) &\inmat{ for }x<0,\\
0 &\inmat{ for }x=0,\\
r  v_+(x) &\inmat{ for }x>0.
\end{array}
\right.
\end{equation}
Note that the proportion $r$ is related to the level of loss aversion when $-\underline{x}=\bar{x}$.
In CPT, loss aversion is interpreted as value (utility) being steeper for losses than for corresponding gains, 
that is, $v(-z)-v(-z-\delta)\geq v(z+\delta)-v(z)$ for all $z\geq 0$ and $\delta>0$, 
which means $v(0)-v(\underline{x})\geq v(\bar{x})-v(0)$ if $\underline{x}=-\bar{x}$ and thus $r\leq 1$ (see \citep{BeT95}).
\end{itemize}

\begin{figure}[!ht]
\begin{minipage}[t]{0.49\linewidth}
\centering
\includegraphics[width=3.2in]{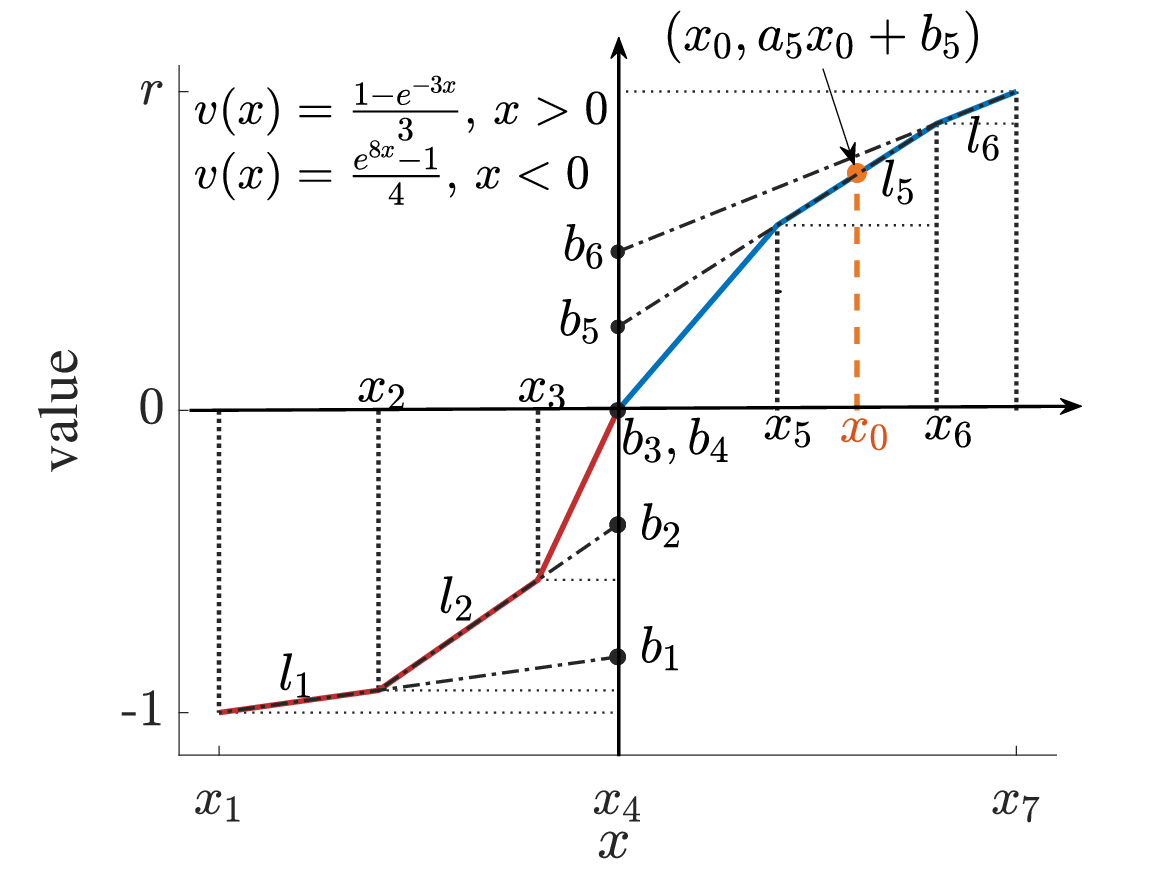}
\text{\footnotesize{ \hspace{0.7cm}\textbf{(a)}}} 
\end{minipage}
\begin{minipage}[t]{0.49\linewidth}
\centering
\includegraphics[width=3.2in]{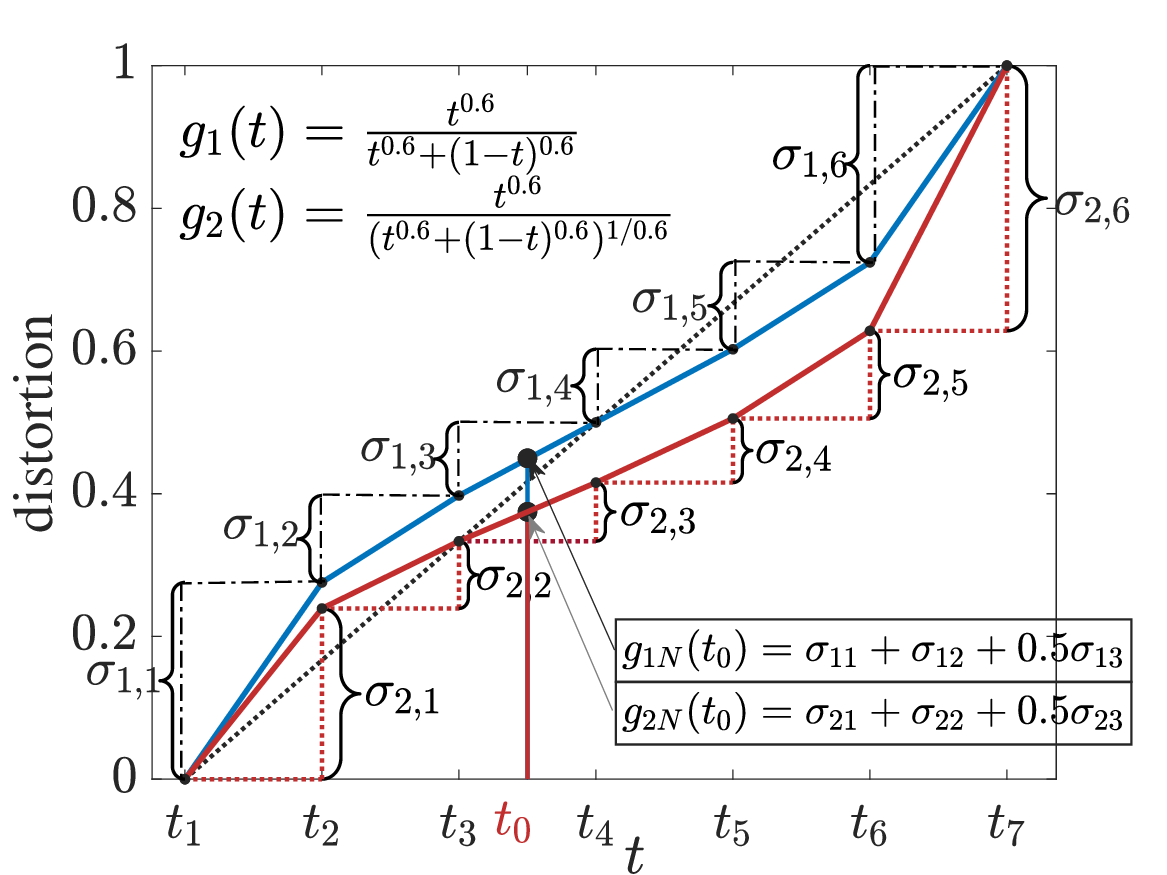}
\text{\hspace{0.7cm}\footnotesize{\textbf{(b)}}}
\end{minipage}
\caption{\textbf{(a)} Piecewise linear value function $v_J(x)$. $l_i(x)=a_ix+b_i$ for $i=1,\cdots,6$.
Set of breakpoints ${\cal X}=\{x_1,\cdots,x_7\}$.
 Consider $x_0:=0.5({x_5+x_6})$.
Since $x_0\in (x_5,x_6]$, then
$\phi_5(x_0)=\mathds{1}_{(x_5,x_6]}(x_0)=1$,
$\phi_1(x_0)=\mathds{1}_{[x_1,x_2]}(x_0)=0$,
and $\phi_i(x_0)=\mathds{1}_{(x_{i},x_{i+1}]}(x_0)=0$ for $i=2,3,4,6$.
In this case, ${\bm \phi}(x_0)=(\phi_1(x_0),\phi_2(x_0),\cdots,\phi_6(x_0))=(0,0,0,0,1,0)$,
and thus $v_J(x_0)=(x_0{\bm a}+{\bm b})^T{\bm \phi}(x_0)=x_0a_5+b_5$.
\textbf{(b)} Piecewise linear distortion function $g_{1N}(t)$ and $g_{2N}(t)$.
Consider $t_0:=0.5(t_3+t_4)$.
Since $t_0
\in (t_3,t_4]$,
then $\psi_3(t_0)=\mathds{1}_{(t_3,t_4]}(t_0)=1$, $l_3(t_0)=\left(\frac{t_0-t_3}{t_4-t_3}\right)\psi_3(t_0)=0.5$, and $\psi_j(t_0)=\mathds{1}_{(t_j,t_{j+1}]}(t_0)=0$, $l_j(t_0)=\left(\frac{t_0-t_j}{t_{j+1}-t_j}\right)\psi_j(t_0)=0$ for $j=1,2,4,5,6$.
By (\ref{eq:g(t)-PRO}),
${\bm h}^{\cal S}(t_0)=(l_1(t_0)+\sum_{j=2}^6\psi_j(t_0),
l_2(t_0)
+\sum_{j=3}^{6}\psi_j(t_0),l_3(t_0)+\sum_{j=4}^6\psi_j(t_0),l_4(t_0)+\sum_{j=5}^6\psi_j(t_0),l_5(t_0)+\psi_6(t_0),l_{6}(t_0))=(\psi_3(t_0),\psi_3(t_0),l_3(t_0),0,0,0)=(1,1,0.5,0,0,0)$.
$g_{iN}(t_0)={\bm \sigma}_i^T{\bm h}^{\cal S}(t_0)=\sigma_{i1}+\sigma_{i2}+0.5\sigma_{i3}$ for $i=1,2$.}
\label{fig:exa_value_distortion_2}
\end{figure}

We confine our discussions to the value function defined in \eqref{eq:value_tilde}.
To ease the exposition, 
we use $v(\cdot)$ to denote $\tilde{v}(\cdot)$ in the forthcoming discussions with a slight abuse of notation.
See an illustration of the piecewise linear value function $v_J(x)$ in Figure~\ref{fig:exa_value_distortion_2}~(a) and the piecewise linear distortion functions $g_{1N}$ and $g_{2N}$ in Figure~\ref{fig:exa_value_distortion_2}~(b). 
Unless specified otherwise, 
we use the piecewise function $v_J$  of the form \eqref{eq:PLA_utility} to approximate the value function in terms of theoretical results and use $v_{J_-}$,  $v_{J_+}$, 
and $v_-$, $v_+$ defined as in Parts (i) and (ii) during preference elicitation.

\subsection{Finite-Dimensional Representation of the Inner Problem}
\label{Sec:tactability}
In this section, we derive a reformulation of the proposed PRGSR-CPT  under Assumptions~\ref{ass:discrete}--\ref{ass:discrete_Lip}.
To this end, we first introduce the joint ambiguity set of the piecewise linear value 
function and distortion functions and show that any three-tuple from   ${\mathbb{V}\mathbb{G}_1\mathbb{G}_2}$ has a corresponding PLA from the ambiguity set with PLFs.
Based on this result,
we develop a finite-dimensional reformulation of the constraint functions.
Recall that the constraint function in problem \eqref{eq:PRGRS-CPT} is
\begin{equation}
\label{eq:constraint_fun}
h_{\mathbb{V}\mathbb{G}_1\mathbb{G}_2} (\xi-x):=\sup_{(v,{g_1},{g_2})\in 
\mathbb{V}\mathbb{G}_1\mathbb{G}_2} \bbe_{{g_1}{g_2}}[v(\xi-x)],
\end{equation}
where $\bbe_{{g_1}{g_2}}[v(\xi-x)]$ is defined as in \eqref{eq:saa_ww}. 
Let $\bar{\cal J} := \bigcup_{m=1}^M( \supp(Z_1^m)\cup \supp(Z_2^m))\bigcup\{\underline{x},\bar{x},0\}$, 
$\bar{\cal S}_1:= \{1-p^m\}_{m=1}^M \bigcup \{0,1,t_1^*\}$, 
$\bar{\cal S}_2:=\{p^m\}_{m=1}^M \bigcup \{0,1,t_2^*\}$ with $\bar{J}:=|\bar{\cal J}|=3M+3$, 
$\bar{N}_i:=|\bar{\cal S}_i|$, where $t_i^*$, $i=1,2$, is the inflection point of $g_i$ (changes from being concave to convex) and it is unknown in advance.
Let
\begin{equation}
\label{eq:J}
{\cal J}(x) := \bar{\cal J}\bigcup \{\xi_k-x\}_{k=1}^K,
\end{equation}
and denote $J:=|{\cal J}(x)|=\bar{J}+K$ with $x_j$ denoting the $j$-th entry of the ranked support set ${\cal J}(x)$.
To ease the exposition, we write ${\cal J}$ for ${\cal J}(x)$.
We observe from \eqref{eq:sum_h-expectation} and \eqref{eq:distorted_exp-2} that 
the value of $h_{\mathbb{V}\mathbb{G}_1\mathbb{G}_2}(\xi-x)$ in \eqref{eq:constraint_fun} depends on the function values of $g_1$ at $q_k$ for $k\in [K_0(x)-1]$, 
and the function values of $g_2$ at $1-q_k$ for $k=K_0(x)-1,\cdots,K$, 
where $q_{0}=0,$ $q_k=\sum_{j=1}^{k}p_j$ for $k\in [K]$.
Consequently, we consider the two sets of breakpoints of the piecewise linear distortion functions: 
\begin{equation}
\label{eq:S}
\begin{aligned}
{\cal S}_1:=\bar{\cal S}_1\bigcup \{q_k\}_{k=1}^{K_0(x)-1}
\quad \inmat{and} \quad 
{\cal S}_2:=  \bar{\cal S}_2\bigcup\{1-q_k\}_{k=K_0(x)-1}^{K}.
\end{aligned}
\end{equation}
For clarity of the exposition, scenarios in $\Xi$ will be indexed by $k$, elements in ${\cal J}$ will be indexed by $j$, 
and  those in ${\cal S}_i$ will be indexed by $s$.
Thus, the size of the problem will be determined by the number of pairwise comparisons $M$, 
the number of scenarios $K$, 
and the total numbers of breakpoints $J=|{\cal J}|$ and $N_i=|{\cal S}_i|$ for $i=1,2$.
Next, we will discuss that the ambiguity set with PLFs using breakpoints ${\cal J}$, 
${\cal S}_1$, and ${\cal S}_2$ is equivalent to that using $\bar{\cal J}$, $\bar{\cal S}_1$, and $\bar{\cal S}_2$.

\subsubsection{Joint Ambiguity Set of PLFs}
\begin{figure}[!ht]
\begin{minipage}[t]{0.49\linewidth}
\centering
\includegraphics[width=3.2in]{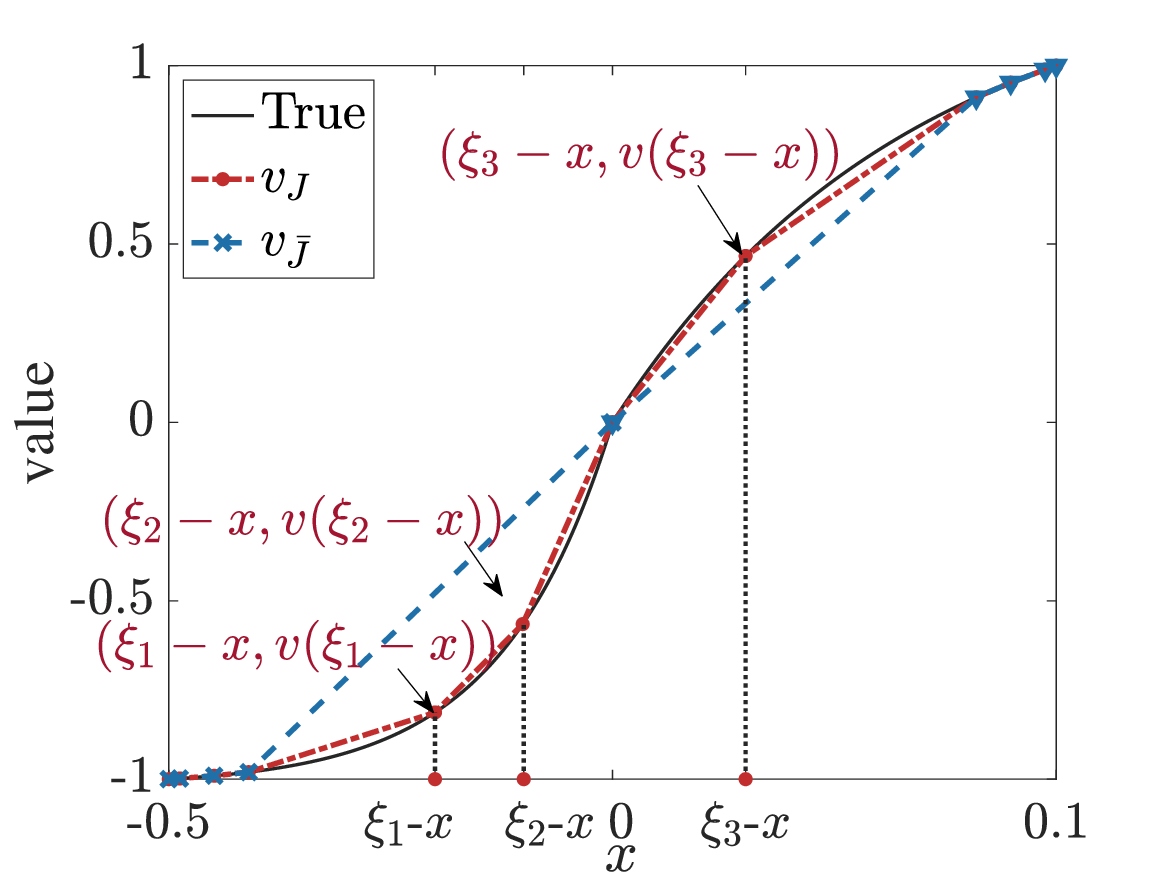}
\text{\footnotesize{\hspace{0.7cm}\textbf{(a)}}}
\end{minipage}
\begin{minipage}[t]{0.49\linewidth}
\centering
\includegraphics[width=3.2in]{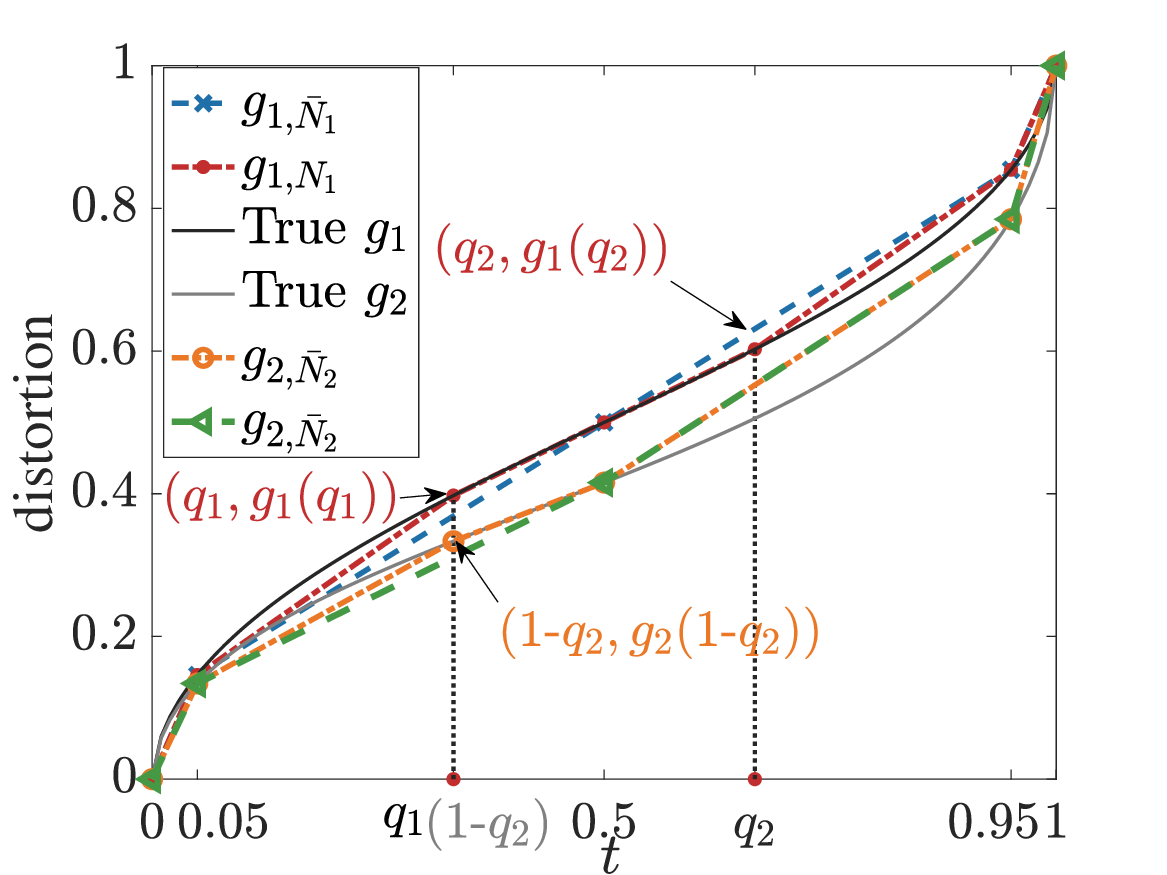}
\text{\footnotesize{\hspace{0.7cm}\textbf{(b)}}}
\end{minipage}
\caption{\footnotesize{\textbf{(a)} Piecewise linear value functions $v_{\bar{J}}$ and $v_{J}$ with $\bar{\cal J}=\{-0.5,\allowbreak
-0.487,\allowbreak
-0.448,\allowbreak
-0.409,\allowbreak
0,\allowbreak
0.409,\allowbreak
0.448,\allowbreak
0.487,\allowbreak
0.5\}$ and
${\cal J}=\bar{\cal J}\cup\{(\xi_k-x)\}_{k=1}^3$,
where
$\{(\xi_k-x)\}_{k=1}^3=\{-0.2,-0.1,0.15\}$.
We find that $v(x)=v_J(x)=v_{\bar{J}}(x)$ for $x\in \bar{\cal J}$,
and that $v_{J}$ has three additional breakpoints 
$-0.2,-0.1,0.15$
compared to $v_{\bar{J}}$.
\textbf{(b)} Piecewise linear distortion functions $g_{i\bar{N}_i}$ and $g_{iN_i}$  with $\bar{\cal S}_i=\{0,0.05,0.5,0.95,1\}$, for $i=1,2$,
and ${\cal S}_1=\bar{\cal S}_1\cup \{q_1,q_2\}$,
${\cal S}_2=\bar{\cal S}_2\cup \{1-q_2\}$ with $q_k=k/3$ for $k=1,2,3$.
Likewise, $g_{i}(t)=g_{iN_i}(t)=g_{i\bar{N}_i}(t)$ at $t\in \bar{\cal S}_i$ for $i=1,2$.
$g_{1N_1}$ (red dotted line) has two additional breakpoints $q_1=1/3$, $q_2=2/3$ compared to $g_{1\bar{N}_1}$ (blue dotted line),
whereas $g_{2N_2}$ (orange dotted line) has one additional breakpoint $1-q_2=1/3$ compared to 
$g_{2\bar{N}_2}$ (green dotted line).}}
\label{fig:value_J_barJ}
\end{figure}

Let $\mathscr{V}_{\rm S}$ and $\mathscr{G}_{\rm invS}^{i}$ be defined as in \eqref{eq:ambiguity set}, 
$\mathscr{V}_{\rm S}^{J} \subset \mathscr{V}_{\rm S}$ be the subset of piecewise linear value functions with breakpoints in ${\cal J}$, 
and  $\mathscr{G}_{\rm invS}^{iN_i}\subset \mathscr{G}_{\rm invS}^{i}$ be the subsets of piecewise linear distortion functions with breakpoints in the sets ${\cal S}_i$, for $i=1,2$.
Let $\mathscr{U}:=\mathscr{V}_{\rm S}^{J}\times \mathscr{G}_{\rm invS}^{1N_1}\times \mathscr{G}_{\rm invS}^{2N_2}$.
Define
\begin{equation}
\label{eq:Question_VJ}
\mathbb{V}_J\mathbb{G}_{N_1}\mathbb{G}_{N_2} := \left\{ (v_{J},g_{1N_1},g_{2N_2})\in\mathscr{U}: 
\begin{aligned}
&\;  \exists\; (v_{\bar{J}}, g_{1\bar{N}_1},\;g_{2\bar{N}_2})\in \mathscr{U} \;\inmat{such that}  \\
&\; v_{\bar{J}}(y)= v_{J}(y), \; y\in  \bar{\cal J}\\
&\;  g_{i\bar{N}_i}(t)=g_{i{N}_i}(t), \; t\in \bar{\cal S}_i, \forall i=1,2  \\
&\;  H_m({v}_{\bar{J}},g_{1\bar{N}_1},g_{2\bar{N}_2})\leq 0, m\in [M]
\end{aligned}
\right\}, 
\end{equation} 
where $H_m$ is defined before \eqref{eq:ambiguity set},  $v_{\bar{J}}$ and $v_{J}$ are piecewise linear value functions with breakpoints in the sets $\bar{\cal J}$ and ${\cal J}$, respectively.
Since ${\cal J}=\bar{\cal J}\cup \{\xi_k-x\}_{k=1}^K$, it follows from \eqref{eq:Question_VJ} that $v_{\bar{J}}(x)=v_{J}(x)$ for $x\in \bar{\cal J}$.
Moreover, since $\bar{\cal J}$ is only a subset of ${\cal J}$, we conclude that $v_{J}$ is a finer PLA than $v_{\bar{J}}$ (see Figure~\ref{fig:value_J_barJ}~(a)).
Likewise, $g_{i\bar{N}_i}$ and $g_{iN_i}$ are piecewise linear distortion functions with breakpoints in $\bar{\cal S}_i$ and ${\cal S}_i$ respectively.
Since ${\cal S}_1=\bar{\cal S}_1\cup\{q_k\}_{k=1}^{K_0(x)-1}$ (resp. ${\cal S}_2=\bar{\cal S}_2\cup \{1-q_k\}_{k=K_0(x)-1}^{K}$), 
we conclude from \eqref{eq:Question_VJ} that $g_{i\bar{N}_i}(t)=g_{iN_i}(t)$ at $t\in \bar{\cal S}_i$, 
and that $g_{iN_i}$ is a finer piecewise linear distortion function than $g_{i\bar{N}_i}$ for $i=1,2$ (see Figure~\ref{fig:value_J_barJ}~(b)).
Since both ${\cal J}$ and $\bar{{\cal J}}$ encompass all of the outcomes of $Z_1^m$ and $Z_2^m$, 
$\bbe_{g_{1\bar{N}_1}g_{2\bar{N}_2}}[v_{\bar{J}}(Z_1^m)] \leq \bbe_{g_{1\bar{N}_1}g_{2\bar{N}_2}}[v_{\bar{J}}(Z_2^m)]$
is equivalent to $\bbe_{g_{1{N}_1}g_{2{N}_2}}[v_{{J}}(Z_1^m)] \leq \bbe_{g_{1{N}_1}g_{2{N}_2}}[v_{{J}}(Z_2^m)]$,
and thus $\mathbb{V}_J\mathbb{G}_{N_1}\mathbb{G}_{N_2}$ is equal to the following set, with $H_m$ defined before \eqref{eq:ambiguity set}:
\begin{equation}
\label{eq:Question_N}
\begin{aligned}
\mathbb{V}_J\mathbb{G}_{N_1}\mathbb{G}_{N_2} = \{(v_{J},g_{1N_1},g_{2N_2})\in \mathscr{U}: H_m({v}_{J},g_{1{N}_1},g_{2{N}_2})\leq 0,\; m\in [M]\big\}.
\end{aligned}
\end{equation}

Since $\mathscr{U}=\mathscr{V}_{\rm S}^{J}\times \mathscr{G}_{\rm invS}^{1N_1}\times \mathscr{G}_{\rm invS}^{2N_2} \subset \mathscr{V}_{\rm S}\times \mathscr{G}_{\rm invS}^{1}\times \mathscr{G}_{\rm invS}^{2}$ (see the definitions of the sets in \eqref{eq:ambiguity set}) and all tuples $(v_J,g_{1N_1},g_{2N_2})$ from $\mathbb{V}_J\mathbb{G}_{N_1}\mathbb{G}_{N_2}$ satisfy the constraint in $\mathbb{VG}_1\mathbb{G}_2$ (see \eqref{eq:ambiguity set}), 
it follows that all these tuples lie within $\mathbb{V}\mathbb{G}_1\mathbb{G}_2$, 
that is, $\mathbb{V}_J\mathbb{G}_{N_1}\mathbb{G}_{N_2}\subset \mathbb{V}\mathbb{G}_1\mathbb{G}_2$.
Conversely, for any $(v,g_1,g_2)\in \mathbb{V}\mathbb{G}_1\mathbb{G}_2$, 
we can construct a three-tuple of PLFs $(v_{J},g_{1N_1},g_{2N_2})\in \mathbb{V}_J\mathbb{G}_{N_1}\mathbb{G}_{N_2}$
by connecting the function values of $v$ and $g_i$, $i=1,2$, at two consecutive breakpoints.
The following proposition states this.

\begin{proposition}[$\mathbb{V}_J\mathbb{G}_{N_1}\mathbb{G}_{N_2}$ vs. $\mathbb{V}\mathbb{G}_1\mathbb{G}_2$]
\label{prop:ambiguity_N_opt}
Let ${\cal J}$ and ${\cal S}_i$, $i=1,2$, be defined as in \eqref{eq:J} and \eqref{eq:S} respectively.
For any $(v,g_1,g_2)\in \mathbb{V}\mathbb{G}_1\mathbb{G}_2$, 
there exists a three-tuple $(v_{J},g_{1N_1},g_{2N_2})\in  \mathscr{V}_{\rm S}^{J}\times \mathscr{G}_{\rm invS}^{1N_1}\times \mathscr{G}_{\rm invS}^{2N_2}$ such that $v_J(x_j)=v(x_j)$ for all $x_j\in {\cal J}$ and $g_{iN_i}(t_s)=g_i(t_s)$ for all $t_s\in {\cal S}_i$, $i=1,2$, 
and $(v_{J},g_{1N_1},g_{2N_2})\in \mathbb{V}_J\mathbb{G}_{N_1}\mathbb{G}_{N_2}$. 
Moreover, $v_J$ and $g_{iN_i}$ are the piecewise linear interpolations of $v$ and $g_i$ on ${\cal J}$ and ${\cal S}_i$, respectively.
Specifically, $(v_{J},g_{1N_1},g_{2N_2})$ can be constructed by
\begin{equation}
\label{eq:PLA_v_g1g2}
\begin{aligned}
v_{J}(x) & :=v(x_{j-1})+\frac{v(x_j)-v(x_{j-1})}{x_j-x_{j-1}}(x-x_{j-1}), \;\;  \mbox{for } x\in [x_{j-1},x_j], \; j\in[J]\backslash{\{1\}}, \\
g_{iN_i}(t) & :=g_{i}(t_{s-1})+\frac{g_i(t_s)-g_i(t_{s-1})}{t_s-t_{s-1}}(t-t_{s-1}), \;\;  \mbox{for } t\in [t_{s-1},t_s],\; s\in[N_i]\backslash{\{1\}},i=1,2.
\end{aligned}
\end{equation}
\end{proposition}

Proposition~\ref{prop:ambiguity_N_opt} means that any three-tuple from  ${\mathbb{V}\mathbb{G}_1\mathbb{G}_2}$ has a PLA in $\mathbb{V}_J\mathbb{G}_{N_1}\mathbb{G}_{N_2}$.
We are now ready to state the main result of this section.

\subsubsection{Finite-Dimensional Reformulation of the Constraint Function}
Under the assumption that the underlying loss $\xi$ is discretely distributed and the value function is Lipschitz continuous, 
we establish the finite-dimensional formulation of $h_{\mathbb{V}\mathbb{G}_1\mathbb{G}_2}(\xi-x)$.
The next theorem states~this.

\begin{theorem}[Finite-Dimensional Reformulation of the Inner Problem of PRGSR-CPT]
\label{thm:reformulation-constraint-function-VW}
Let ${\cal J}(x)$ be defined in \eqref{eq:J}, 
and ${\cal S}_1$ and ${\cal S}_2$ be defined in \eqref{eq:S}. 
Set $J = |{\cal J}(x)|$, $N_1 = |{\cal S}_1|$, and $N_2 = |{\cal S}_2|$.
Let $\mathbb{V}_J\mathbb{G}_{N_1}\mathbb{G}_{N_2}$ be defined in \eqref{eq:Question_N}, 
and let $\mathbb{V}\mathbb{G}_1\mathbb{G}_2$ be defined in \eqref{eq:ambiguity set}. 
Define $h_{\mathbb{V}\mathbb{G}_1\mathbb{G}_2}(\xi-x)$ as in \eqref{eq:constraint_fun}.
Suppose that Assumptions~\ref{ass:discrete}--\ref{ass:discrete_Lip} hold.
For any fixed $x$, the following assertions hold.

(i) $h_{\mathbb{V}_J\mathbb{G}_{N_1}\mathbb{G}_{N_2}}(\xi-x)=h_{\mathbb{V}\mathbb{G}_1\mathbb{G}_2}(\xi-x)$.

(ii) The constraint function  $h_{\mathbb{V}\mathbb{G}_1\mathbb{G}_2}(\xi-x)$ in \eqref{eq:constraint_fun} is the optimal value of the following  program:
\begin{subequations}
\label{robust_v_o}
\begin{align}
\sup_{{\bm a},{\bm b},{\bm \sigma}_1,{\bm \sigma}_2}  & \Bigg\{\sum_{k=1}^{K_0(x)-1} 
((\xi_k-x) {\bm a}+{\bm b})^\top{\bm \phi}(\xi_k-x)  {\bm \sigma}_1^\top ({\bm h}^{{\cal S}_1}(q_k)-{\bm h}^{{\cal S}_1}(q_{k-1}))  \nonumber \\
&\;\;  + \sum_{k=K_0(x)}^{K} ((\xi_k-x) {\bm a}+{\bm b})^\top{\bm \phi}(\xi_k-x)  {\bm \sigma}_2^\top ({\bm h}^{{\cal S}_2}(1-q_{k-1})-{\bm h}^{{\cal S}_2}(1-q_{k})) \Bigg\} \label{robust_Wr0} \\
{\rm s.t.} \;\;\; & \frac{\sigma_{1,s+1}}{t_{s+2}-t_{s+1}}\leq \frac{\sigma_{1,s}}{t_{s+1}-t_s},\; \forall s=1,\cdots,s^1-2, \label{robust_Wrb-1}  \\
& \frac{\sigma_{2,s+1}}{t_{s+2}-t_{s+1}}\leq \frac{\sigma_{2,s}}{t_{s+1}-t_s},\; \forall s=1,\cdots,s^2-2, \label{robust_Wrb-2}  \\
& \frac{\sigma_{1,s}}{t_{s+1}-t_{s}}\leq \frac{\sigma_{1,s+1}}{t_{s+2}-t_{s+1}},\; \forall  s=s^1,\cdots,N_1-2,\label{robust_Wrd-1} \\
& \frac{\sigma_{2,s}}{t_{s+1}-t_{s}}\leq \frac{\sigma_{2,s+1}}{t_{s+2}-t_{s+1}},\; \forall s=s^2,\cdots,N_2-2,\label{robust_Wrd-2} \\
& {\bm e}^\top{\bm \sigma}_i=1,\; \forall i=1,2,\label{robust_Wre} \\
& a_{j-1} x_{j}+b_{j-1}=a_{j} x_{j}+b_{j}, \; \forall j\in[J-1]\backslash{\{1\}}, \label{robust_Wrf}\\
& a_{1}\underline{x}+b_{1}=-1, \;
b_{J_0}=0,\; a_{J-1}\bar{x}+b_{J-1}=r, \label{robust_Wrg}\\
& a_j-a_{j+1} \leq 0,\; \forall j=1,\cdots,J_0-2, \label{robust_Wri-1}\\
& a_{j+1}-a_j \leq 0, \; \forall j=J_0,\cdots,J-2, \label{robust_Wri-2}\\
& h_m(R^m({\bm a},{\bm b},{\bm \sigma}_1,{\bm \sigma}_2)  -(r_2^m{\bm a}+{\bm b})^\top{\bm \phi}(r_2^m))\leq 0,\; \forall m\in [M], \label{robust_Wrj}\\
&{\bm a}\geq 0,\; {\bm a}\leq L{\bm e},\;{\bm \sigma}_1\geq 0, \;{\bm \sigma}_2\geq 0, \label{robust_Wrk}
\end{align}
\end{subequations}
where ${\bm a},{\bm b}\in \R^{J-1}$, ${\bm e}$ is a vector with components $1$, 
${\bm \sigma}_i\in \R^{N_i-1}$ for $i=1,2$, $s^i$ is the index of the point $t_i^*$, 
where $g_{i}$ changes from being concave to being convex for $i=1,2$, $J_{0}$ is the index of the point $0$, where the value function $v$ changes from being convex to being concave, 
$q_{0}=0,$ $q_k=\sum_{j=1}^{k}p_j$ for $k\in [K]$.
Let $w_1^m:={\bm \sigma}_1^\top{\bm h}^{{\cal S}_1}(1-p^m)$,
$w_2^m:={\bm \sigma}_2^\top{\bm h}^{{\cal S}_2}(p^m)$
and define $R^m({\bm a},{\bm b},{\bm \sigma}_1,{\bm \sigma}_2):=R_1^m({\bm a},{\bm b},{\bm \sigma}_1,{\bm \sigma}_2) +R_2^m({\bm a},{\bm b},{\bm \sigma}_1,{\bm \sigma}_2)$ with
$R_1^m({\bm a},{\bm b},{\bm \sigma}_1,{\bm \sigma}_2)
:=(r_1^m{\bm a}+{\bm b})^\top{\bm \phi}(r_1^m) (w_1^m\mathds{1}_{r_1^m<0}
 +(1-w_2^m)\mathds{1}_{r_1^m\geq0})$,
 and 
$R_2^m({\bm a},{\bm b},{\bm \sigma}_1,{\bm \sigma}_2)
:=(r_3^m{\bm a}+{\bm b})^\top{\bm \phi}(r_3^m)
((1-w_1^m)\mathds{1}_{r_3^m<0}
+w_2^m\mathds{1}_{r_3^m\geq0})$.
\end{theorem}

The proof is deferred to Section~\ref{EC:Theorem1}.
It is helpful to explain the constraints of program \eqref{robust_v_o}.
Constraints \eqref{robust_Wrb-1}--\eqref{robust_Wrb-2} characterize the concavity of $g_{1N_1}$ and $g_{2N_2}$ over $[0,t_{s^1}]$ and $[0,t_{s^2}]$, respectively, 
and constraints \eqref{robust_Wrd-1}--\eqref{robust_Wrd-2} characterize the convexity of $g_{1N_1}$ and $g_{2N_2}$ over $[t_{s^1},1]$ and $[t_{s^2},1]$, respectively.
Constraint \eqref{robust_Wre} stipulates that $g_{iN_i}(1)=1$ and $g_{iN_i}(0)=0$ for $i=1,2$.
Constraint~\eqref{robust_Wrf} ensures the continuity of the PLF $v_J$.
Constraint~\eqref{robust_Wrg} requires $v_J(\bar{x})=v(\bar{x})=r$, $v_J(0)=v(0)=0$,
and $v_J(\underline{x})=v(\underline{x})=-1$,
by the definition of $v_J$ and $v$ (see \eqref{eq:value_tilde}).
Constraints \eqref{robust_Wri-1}--\eqref{robust_Wri-2} reflect the convexity and concavity of $v_J$ over $[\underline{x},0]$ and $[0,\bar{x}]$, respectively, 
via the slopes of two consecutive linear pieces.
Constraint~\eqref{robust_Wrk} enforces monotonicity of $v_J$, $g_{1N_1}$ and $g_{2N_2}$ and the Lipschitz continuity of $v_J$.
We will explain the objective function and constraint \eqref{robust_Wrj} of the program in the proof.
Note that the objective function in \eqref{robust_Wr0} and the constraint function in \eqref{robust_Wrj} are bilinear in $({\bm a},{\bm \sigma}_i)$ and $({\bm b},{\bm \sigma}_i)$ for $i=1,2$, 
which makes problem \eqref{robust_v_o} computationally challenging.

From the theorem, we can see that the worst-case value function and the worst-case distortion functions lie in  $\mathbb{V}_J\mathbb{G}_{N_1}\mathbb{G}_{N_2}$
(see \eqref{eq:Question_N} for the definition of the latter), and the robust constraint function can be computed by solving a finite-dimensional bilinear program, where the value function is characterized by slopes ${\bm a}$ and intercepts  ${\bm b}$, and the distortion functions by increments ${\bm \sigma}_1,{\bm \sigma}_2$.
The resulting bilinear program is nonconvex and therefore not polynomial-time solvable in general.
When the bilinear constraints~\eqref{robust_Wrj} are simplified into linear constraints according to Question Types~1--4 in Section~\ref{sec:Estimating-dis}, 
the finite-dimensional problem~\eqref{robust_v_o} still retains a bilinear objective function.
For instance, we solve the resulting problem for a fixed value $x = 0.3$, $\xi_k\in [0,0.5]$ for $k\in [K]$ with $K = 10$ scenarios, $N_1 = N_2 = 25$, $J= 25$, and $M = 50$
(with $M_1 = M_2 = 20$ for Question Types~1--2 and $M_4 = 10$ for Question Type~4).
The problem is solved within approximately $1946$ seconds using Gurobi.
To tackle the issue, we will use the well-known piecewise McCormick relaxation to approximate the bilinear term in Section~\ref{sec:McCormick}. 
In Section~\ref{sec:sensitivity}, 
we then show that instances under identical experimental settings, with up to $10$ scenarios, 
can be solved within reasonable computational time using Gurobi. 
Specifically, the solution time increases from approximately $10$ seconds with one McCormick segment, 
to about $30$ seconds with $10$ segments, and to roughly $120$ seconds with $20$ segments.

If the robust constraint function can be calculated, 
then we can calculate PRGSR-CPT using the bisection method.
If we can identify a lower bound $t^0_l$ and an upper bound  $t^0_u$ for the PRGSR-CPT measure $t^*=\rho_{\mathbb{VG}_1\mathbb{G}_2}(\xi)=\inf\{t\in \R: \sup_{(v,g_1,g_2)\in \mathbb{VG}_1\mathbb{G}_2}\bbe_{g_1g_2}[v(\xi-t)]\leq 0\}$, i.e., $\rho_{\mathbb{V}\mathbb{G}_1\mathbb{G}_2}(\xi)=t^*$ and $t^0_l\leq t^* \leq t^0_u$,
then we can carry out a bisection search to figure out $t^*$ (see the algorithmic details in~Section~\ref{sec:bisection-robust}).

\subsubsection{Motivation for the Modified Polyhedral Method for Preference Elicitation}
The parameters $p^m, r_1^m, r_2^m, r_3^m$ in the constraint \eqref{robust_Wrj} are undetermined.
Program \eqref{robust_v_o} is difficult to solve because of the bilinear terms appearing in the objective and the constraint. 
In practice, one needs to design the questionnaire in a smart way. 
In the next section, we will discuss how to properly design lotteries $Z_1^m$, $Z_2^m$  so that the parameters $p^m, r_1^m, r_2^m, r_3^m$ in the constraint \eqref{robust_Wrj} 
can be determined efficiently.
One of the key issues to be resolved is that the standard RRUS scheme which is widely 
used for preference elicitation in the literature on expected utility models is not applicable to CPT.
The design of the questionnaire in RRUS is essentially determining the probability of a risky lottery as the midpoint between the minimum and maximum possible values of the normalized utility function at a certain outcome point based on the answers to the previous queries.
The scheme cannot be used here for eliciting the value function unless the distortion functions are given.
The other challenge is that normalization of the value function, 
a key step in RRUS,
does not work under the CPT framework.
Specifically, let $\hat{v}(x):=\frac{v(x)-v(\underline{x})}{v(\bar{x})-v(\underline{x})}$.
$Z_1^m$ and $Z_2^m$ are defined as in~\eqref{eq:questions}, with $r_1^m < 0 < r_3^m$.
Then
\begin{equation}
\label{eq:normalization_extra}
\begin{aligned}
 \bbe_{g_1g_2}[v(Z_1^m)]-\bbe_{g_1g_2}[v(Z_2^m)]\leq 0 
& \Leftrightarrow  v(r_1^m)g_{1}(1-p^m) +v(r_3^m)g_{2}(p^m) -v(r_2^m)\leq 0 \\
&\Leftrightarrow  \hat{v}(r_1^m)g_{1}(1-p^m)  + \hat{v}(r_3^m)g_{2}(p^m)-\hat{v}(r_2^m) \\
& \quad \;\; +\frac{v(\underline{x})}{v(\bar{x})-v(\underline{x})}(g_1(1-p^m)+g_2(p^m)-1)\leq 0 \\
& \not\Leftrightarrow  \hat{v}(r_1^m)g_{1}(1-p^m) +\hat{v}(r_3^m)g_{2}(p^m)
-\hat{v}(r_2^m)\leq 0.
\end{aligned}
\end{equation}
This is essentially because there are two distortion functions $g_1$ and $g_2$ in CPT such that the distorted expectation of a risky lottery with outcomes $r_1^m<0$ and $r_3^m>0$ yields a term $g_1(1-p^m)+g_2(p^m)-1\neq 0$.

Note that \citet{HFD16} consider stochastic programs with ambiguity in both the utility function and the underlying probability distribution, 
whose computation is based on a bilinear programming problem.
They provide an exact evaluation of the resulting robust optimization when the utility function is concave and the probability distributions are represented by vectors in a
polyhedron.
The problem can be approximately computed by using the reformulation-linearization technique and the semidefinite programming relaxation for the general ambiguity set of probability distributions.
All these results are based on the assumption that the utility function is concave and can be represented by vectors using the support function-based method.
Different from this, we consider the $S$-shaped utility function and the worst-case piecewise linear utility function represented by ${\bm a},{\bm b}$ in \eqref{robust_v_o}
that is parameterized by $x$, even though the underlying distortion functions can be characterized by vectors in a polyhedron.
The issue of bilinear terms does not pose an immediate difficulty to the solvability of \eqref{robust_Wr0} when the problem size is relatively small. 
For example, we can use Gurobi to solve the problem. 
However, the problem size of \eqref{robust_v_o} depends heavily on the number of queries used for eliciting the DM's risk preference. 
This brings us to another query as to how to design questionnaires optimally so that the joint ambiguity set is reduced efficiently (see \citep{VYMDR20}).
In practice, we may add queries pair by pair. 
For instance, let $(\mathbb{V}\mathbb{G}_1\mathbb{G}_2)_{m-1}$ be a joint ambiguity set constructed with the DM's answers to the first $m-1$ queries, 
let ${\cal M}_{m}:=\{(v,g_{1},g_{2}) \in \mathscr{U}:
h_m(\bbe_{g_{1\bar{N}_1}g_{2\bar{N}_2}}[v_{\bar{J}}(Z_1^m)]
-\bbe_{g_{1\bar{N}_1}g_{2\bar{N}_2}}[v_{\bar{J}}(Z_2^m)])\leq 0\}$
denote the set obtained from the $m$-th query.
Then 
$
(\mathbb{V}\mathbb{G}_1\mathbb{G}_2)_{m}:=(\mathbb{V}\mathbb{G}_1\mathbb{G}_2)_{m-1}\cap {\cal M}_{m}.
$
Here we may regard ${\cal M}_{m}$ as a ``hypersurface'' cutting off part of the ambiguity set $(\mathbb{V}\mathbb{G}_1\mathbb{G}_2)_{m-1}$.

\begin{remark}[Framework of the Preference Elicitation]
To address the first issue, we design Question Types 1--2 in Section~\ref{sec:Estimating-dis} to obtain the ambiguity set of the vectors of increments of $g_1$, ${g}_2$, $v_-$ and $v_+$, 
by exploiting the ideas of the well-known polyhedral method in the literature on linear multivariate utility models \citep{THS04} in Section~\ref{sec:modified_PM}. 
To address the second issue, we utilize the breakpoints of  $g_1$ and $g_2$ in Question Types 1--2 to obtain the ambiguity set of the vectors of the increments of  $g_1$, $g_2$ and the ratio of the value functions $r=-\frac{v(\bar{x})}{v(\underline{x})}$ via Question Type $3$ in Section~\ref{Sec:Estimation_r}.
The flowchart in Figure~\ref{fig:flow_chart_Q} illustrates the procedures for estimating $g_1,g_2$, $v_-$, $v_+$, and $r$ through asking the DM Question Types 1--4.

\begin{figure}[!ht]
 \begin{minipage}[t]{\linewidth}
\centering
\includegraphics[width=6.2in]{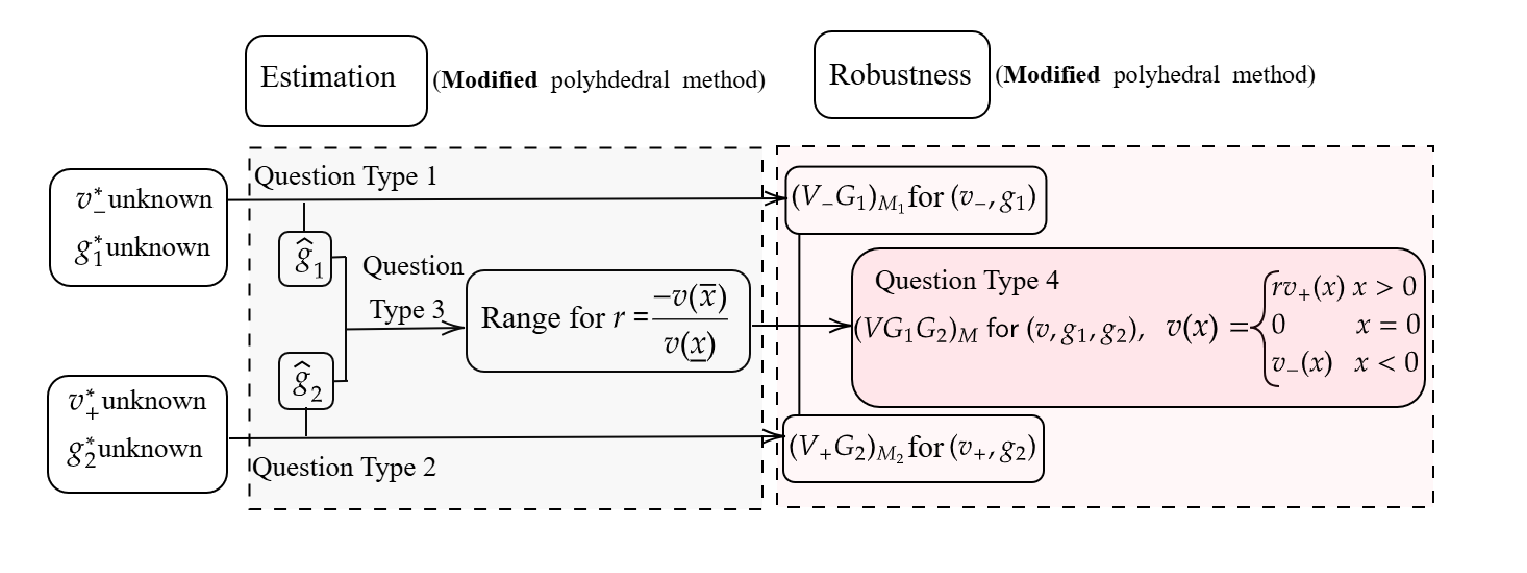}
\end{minipage}
\caption{\footnotesize{The procedure of constructing the ambiguity set $(\mathbb{V}\mathbb{G}_1\mathbb{G}_2)_M$ of the value function $v$ over $[\underline{x},\bar{x}]$ and the distortion functions $g_1$ and $g_2$.
\underline{Estimation portion}:
$(\mathbb{V}_-\mathbb{G}_1)_{M_1}$ is the set of the value function $v_-$ over $[-\underline{x},0]$ and the distortion function $g_1$ with byproduct,  
estimate $\widehat{g}_1$, 
through Question Type 1 when their true counterparts $v_-^*,g_1^*$ are unknown. 
Likewise,
$(\mathbb{V}_+\mathbb{G}_2)_{M_2}$ is the set of the value function $v_+$ over $[0,\bar{x}]$ and the distortion function $g_2$ with byproduct,  
estimate $\widehat{g}_2$, 
through Question Type 2 when their true counterparts $v_+^*,g_2^*$ are unknown. 
The sets $(\mathbb{V}_-\mathbb{G}_1)_{M_1}$ and $(\mathbb{V}_+\mathbb{G}_2)_{M_2}$  are used in the follow-up robustness procedure.
Next, $\widehat{g}_1,\widehat{g}_2$ are used to estimate the ratio $-\frac{v(\bar{x})}{v(\underline{x})}$ through Question Type 3.
\underline{Robustness portion}: 
The modified
polyhedral method will be used in the estimation part
whereas 
the modified polyhedral method for the ratio of the value function, combined with the ambiguity sets $(\mathbb{V}_-\mathbb{G}_1)_{M_1}$ and $(\mathbb{V}_+\mathbb{G}_2)_{M_2}$ will be used in the part of ambiguity set construction.
}}
\label{fig:flow_chart_Q}
\end{figure}
\end{remark}

\section{Preference Elicitation for the Value Function and the Distortion Functions}
\label{sec:Estimating-dis}
In this section, we discuss the design of questionnaires for constructing and reducing
the ambiguity set $\mathbb{V}\mathbb{G}_1\mathbb{G}_2$ defined as in \eqref{eq:ambiguity set} in Section~\ref{sec:tractability}. 
We do so in several steps.
In Section~\ref{sec:bilinear}, 
we design Question Types $1$ and $2$ to tackle the bilinear terms in constraint \eqref{robust_Wrj}, converting them into linear terms; 
in Section~\ref{sec:PM_modified}, 
we develop a modified polyhedral method to estimate the negative (resp. positive) part of the value function and the distortion functions;
in Section~\ref{Sec:Estimation_r}, 
we design Question Type $3$ for the estimation of a ratio $r$ related to the value functions, and in Section~\ref{sec:ratior}, 
we develop a modified polyhedral method to estimate the holistic value function $v$.
Based on the estimation of the distortion functions and the value functions, 
we construct the specific ambiguity set $(\mathbb{V}\mathbb{G}_1\mathbb{G}_2)_M$ in Section~\ref{sec:AS}.

\subsection{Question Types 1 and 2 for the Estimation of \texorpdfstring{$(v_-,g_1)$}{} and \texorpdfstring{$(v_+,g_2)$}{}}
\label{sec:bilinear}
We begin by developing a strategy to tackle the bilinear terms as noted in Theorem~\ref{thm:reformulation-constraint-function-VW}. 
To this end, we design two types of queries which will transform the DM's answers to the pairs of lotteries $(Z_1^m, Z_2^m)$, for $m\in [M_1]$ and $m\in [M_2]$, into a linear system of inequalities w.r.t. $(v_-,g_1)$ and $(v_+,g_2)$, respectively.

\vspace{0.1cm}

\noindent{\bf Question Type 1.}
This type of question comprises a risky lottery $Z_1^m$ and a certain lottery $Z_2^m$:
\begin{equation*}
Z_1^m=\left\{\begin{aligned}
\underline{x} \; & \mbox{ w.p. } 1-p^m\\
0 \; & \mbox{ w.p. } p^m
\end{aligned}
\right.
\quad \mbox{and} \quad 
Z_2^m=r_{2-}^m,
\end{equation*}
where $r_{2-}^m\in [\underline{x},0]$ and $p^m\in (0,1)$ are two parameters to be determined. 
The index $m$ denotes the $m$-th query, and there are $M_1$ such queries in total to be designed.

The key point to note is that the outcomes of $Z_1^m$ are fixed and they are both in the region $[\underline{x},0]$. 
This will effectively enable us to formulate the DM's answer to the query as a linear inequality. 
To see this, by \eqref{eq:distorted_exp-2}, $\bbe_{g_1g_2}[v(Z_1^m)]=v(\underline{x})g_1(1-p^m)$ and $\bbe_{g_1g_2}[v(Z_2^m)]=v(r_{2-}^m)$. 
Let $v$ and $v_-$ be defined as in \eqref{eq:value_tilde}. 
Then $v(\underline{x})=v_-(\underline{x})=-1$ and $v(r_{2-}^m)=v_-(r_{2-}^m)$.
Based on the DM's preference over the two lotteries, we will be able to derive a linear inequality in terms of $g_1(1-p^m)$ and $v_-(r_{2-}^m)$,
$h_m(-g_1(1-p^m)-v_-(r_{2-}^m))\leq 0$, for $m\in [M_1]$,
where $h_m=1$ if $Z_2^m\succeq Z_1^m$ and $h_m=-1$ otherwise.

\vspace{0.1cm}

\noindent{\bf Question Type 2.}
The second type of question is designed with the same consideration as in Question Type $1$ except that one of the outcomes of $Z_1^m$ is changed to the right endpoint $\bar{x}$:
\begin{equation*}
Z_1^m=\left\{\begin{aligned}
0 \;& \mbox{ w.p. } 1-p^m\\
\bar{x} \;& \mbox{ w.p. } p^m
\end{aligned}
\right.
\quad \mbox{and} \quad
Z_2^m=r_{2+}^m,
\end{equation*}
where $r_{2+}^m\in [0,\bar{x}]$ and $p^m\in (0,1)$.
In this case, $\bbe_{g_1g_2}[v(Z_1^m)]=v(\bar{x})g_2(p^m)$ and $\bbe_{g_1g_2}[v(Z_2^m)]=v(r_{2+}^m)$.
Let $v$ and $v_+$ be defined as in \eqref{eq:value_tilde}. 
Then $v(\bar{x})=rv_+(\bar{x})=r$, $v(r_{2+}^m)=rv_+(r_{2+}^m)$, 
and the DM's answer to each of the queries will enable us to derive the following inequalities:
$
h_m(rg_{2}(p^m)-rv_+(r_{2+}^m))\leq 0,\;\inmat{for }m\in [M_2] 
\Leftrightarrow  
h_m(g_{2}(p^m)-v_+(r_{2+}^m))\leq 0,\;\inmat{for}\; m\in [M_2]$,
where $h_m=1$ if $Z_2^m\succeq Z_1^m$ and $h_m=-1$ otherwise.

In the rest of this section, we propose to use the basic idea of the well-known polyhedral method (see \citep{THS04}) to estimate $(v_-,g_1)$ and $(v_+,g_2)$. 
There is a gap though: the polyhedral method works only for linear utility functions with a finite number of parameters whereas here the four functions are all nonlinear and infinite dimensional.
To fill the gap, we will focus on PLF representations of $(v_-,g_1)$ and $(v_+,g_2)$. 
As we discussed in Section~\ref{sec:PLFs0}, 
a normalized PLF is determined by the vector of increments over a specified set of intervals.
As such, we can use  the linear system of inequalities established through questionnaires to specify a feasible set for the vector of increments and then develop an appropriate strategy 
to reduce the set by implementing a sequence of cuts.

\subsection{Modified Polyhedral Method for the Estimation of \texorpdfstring{ $(v_-,g_1)$}{} and \texorpdfstring{$(v_+,g_2)$}{}}
\label{sec:PM_modified}
In this section, we discuss how to generate the $m$-th query for Question Types $1$ and $2$ based on the previous $m-1$ queries by extending the polyhedral method (see \citep{THS04})
to the elicitation of the value function and distortion function simultaneously.
In Section~\ref{sec:modified_PM}, we construct the initial polyhedron of the increments of the positive part of the value function $v_+$ and the distortion function $g_2$.
In Section~\ref{sec:generation_new}, we develop the steps, 
including the models for the analytic center and the models for generating parameters in the $m$-th query.
In Section~\ref{sec:thm-gua}, we establish some theoretical results which guarantee the efficiency of the elicitation method.

\subsubsection{Description of the Method}
\label{sec:modified_PM}
We concentrate on developing a modified polyhedral method to design questionnaires for eliciting $(v_+,g_2)$ since the same strategy can be applied to the elicitation of $(v_-,g_1)$.
The method relies on piecewise linear approximations of $v_+$ and $g_2$: 
\begin{equation}
\label{eq:g2v+}
v_{J_+}(x)={\bm v}_+^\top {\bm h}^{{\cal X}_+}(x)
\quad \inmat{and} \quad
g_{2N_2}(t)={\bm \sigma}_2^\top{\bm h}^{{\cal S}_2}(t),
\end{equation}
where ${\cal X}_+=\{x_1^+,\cdots,x_{J_+}^+:0=x_1^+<x_2^+<\cdots<x_{J_+}^+=\bar{x}\}$, ${\cal S}_2:=\{t_1,\cdots,t_{N_2}:0=t_1<t_2<\cdots<t_{N_2}=1\}$, $J_+=|{\cal X}_+|$,
$N_2=|{\cal S}_2|$, ${\bm v}_+\in \R^{J_+-1}$ and ${\bm \sigma}_2\in \R^{N_2-1}$,
${\bm h}^{{\cal X}_+}(x)$ is defined as in \eqref{eq:PLA_utility_2}, and
${\bm h}^{{\cal S}_2}{(t)}$ is defined as in \eqref{eq:g(t)-PRO} with ${\cal S}$ being replaced with ${\cal S}_2$.
In the setup of \eqref{eq:g2v+}, the piecewise linear functions are uniquely determined by the vectors of increments ${\bm v}_+$ and ${\bm \sigma}_2$. 
Our aim is to estimate ${\bm v}_+$ and ${\bm \sigma}_2$ based on elicited information of the DM's preferences.
We begin with an initial polyhedron of the ranges of the two vectors and then try to reduce the ranges effectively based on the fundamental ideas of the polyhedral method (see \citep{THS04})
by strategically designing pairwise comparison lotteries $Z_1^m$, $Z_2^m$ described in Section~\ref{sec:bilinear} and asking the DM to make a choice.

Let ${\bm I}\in \R^{r\times r}$ be the identity matrix of appropriate dimension $r$,
and let $\widetilde{\bm v}_{+}:=\left[\begin{matrix}{\bm I} & 0\end{matrix}\right]{\bm v}_+\in \R^{J_+-2}$, $\widetilde{\bm \sigma}_{2}:=\left[\begin{matrix}{\bm I}&0\end{matrix}\right]{\bm \sigma}_2\in \R^{N_2-2}$,
where the last components of ${\bm v}_+$ and ${\bm \sigma}_2$ are dropped.
Conversely, we can use $\widetilde{\bm v}_+$ and $\widetilde{\bm \sigma}_2$ to represent ${\bm v}_+$ and ${\bm \sigma}_2$. 
Let ${\bm e}\in\R^l$ be the vector with all components equal to $1$, and let the mapping ${\bm R}:\R^{l}\rightarrow \R^{l+1}$ be defined as  ${\bm R}{\bm r}=\left[\begin{matrix}
{\bm r}\\
1-{\bm e}^\top{\bm r}
\end{matrix}\right]$. 
Then ${\bm v}_+={\bm R}\widetilde{\bm v}_+$, and ${\bm \sigma}_2={\bm R}\widetilde{\bm \sigma}_2$.
Let $d_j^+:=\frac{1}{x_{j+1}^+-x_j^+}, \quad j=1,\ldots,J_+-1$ and define ${\bm D}_+\in \R^{(J_+-2)\times(J_+-2)}$ as follows:
\begin{equation}
\label{eq:D_+}
{\bm D}_+:= \left[ \begin{array}{ccccc}
-d_1^+ & d_2^+ & & &\\
&-d_2^+ & d_3^+ & &\\
& &\ddots &\ddots &\\
& & &-d_{J_+-3}^+ &d_{J_+-2}^+\\
-d_{J_+-1}^+ &-d_{J_+-1}^+
&\cdots &-d_{J_+-1}^+
& \displaystyle \sum_{j=J_+-2}^{J_+-1}-d_{j}^+
\end{array}
\right],
\end{equation}
and ${\bm d}_+:=\left(0,\cdots,0,\frac{-1}{x_{J_+}^+-x_{J_+-1}^+}\right)\in \R^{J_+-2}$.
Then ${\bm D}_+\widetilde{\bm v}_+\leq {\bm d}_+$ represents that the piecewise linear value function $v_{J_+}$ defined as in \eqref{eq:g2v+} is concave over $[0,\bar{x}]$.
The initial polyhedron of $(\widetilde{\bm v}_+,\widetilde{\bm \sigma}_2)$ is
\begin{equation}
\label{eq:defi_Delta_0}
\begin{aligned}
\Delta_+^0 := \big \{ \left( \widetilde{\bm v}_+, \widetilde{\bm \sigma}_2 \right)\in \R^{J_+-2}\times \R^{N_2-2}: -\widetilde{\bm \sigma}_2\leq 0,  -\widetilde{\bm v}_+\leq 0, {\bm e}^\top\widetilde{\bm \sigma}_2\leq 1, {\bm e}^\top\widetilde{\bm v}_+\leq 1, {\bm D}_+\widetilde{\bm v}_+\leq {\bm d}_+\big\}.
\end{aligned}
\end{equation}
Now suppose that we have established $\Delta_+^{m-1}$ and update it by asking the DM a question of Type 2 and then constructing $h_m(g_{2N_2}(p^m) - v_{J_+}(r_{2+}^m)) \leq 0$ based on the DM's answer. 
The inequality can be rewritten as $h_m({\bm \sigma}_2^\top {\bm h}^{{\cal S}_2}(p^m) - {\bm v}_+^\top {\bm h}^{{\cal X}_+}(r_{2+}^m)) \leq 0$.
By incorporating ${\bm v}_+ = {\bm R}\widetilde{{\bm v}}_+$ and ${\bm \sigma}_2 = {\bm R}\widetilde{{\bm \sigma}}_2$, we obtain
\begin{equation}
\label{eq:choice-xkyk}
\begin{aligned}
\Delta_+^{m}:=\Delta_+^{m-1} \bigcap  \left \{ \left( \widetilde{\bm v}_+, \widetilde{\bm \sigma}_2 \right)\in \R^{(J_+-2)+{(N_2-2)}}:  h_m\left[\begin{matrix} 
-{\bm h}^{{\cal X}_+}(r_{2+}^m)\\
{\bm h}^{{\cal S}_2}(p^m)
\end{matrix}\right]^\top
\left[\begin{matrix}
{\bm R}\widetilde{\bm v}_+\\
{\bm R}\widetilde{\bm  \sigma}_2
\end{matrix}
\right]\leq 0 \right\}. 
\end{aligned}
\end{equation}
It remains to be explained how the lotteries of Question Type 2 should be designed so that  $\Delta_+^{m-1}$ is effectively cut down to $\Delta_+^{m}$. 
This requires us to strategically set parameters $p^{m}$ and $r_{2+}^{m}$ in the design of the query. 
We propose a modified polyhedral method which mimics the fundamental idea of the polyhedral method in \cite{THS04} to identify $(p^{m}, r_{2+}^{m})$.
To this end, we introduce the notion of an analytic center.
\begin{definition}[Analytic Center]
\label{def:AC}
Consider a polyhedron being defined by a set of linear inequalities $\tilde{\bm a}_i^\top {\bm x} \leq \tilde{b}_i$,  for $i\in [I]$.
The analytic center of the polyhedron is defined as a global minimizer of the function $-\sum_{i=1}^I \log(\tilde{b}_i-\tilde{\bm a}_i^\top {\bm x})$ over $\{{\bm x}:\tilde{\bm a}_i^\top{\bm x}<\tilde b_i,\ i\in[I]\}$.
\end{definition}

\subsubsection{Design of a New Query}
\label{sec:generation_new}
The following steps~1--3 can be used to optimally design the $m$-th query.

{\bf Step 1}: Find the analytic center ${\bm c}^m=({\bm c}^m_1,{\bm c}^m_2)$
with ${\bm c}_1^m\in \R^{J_+-2}$, ${\bm c}_2^m\in \R^{N_2-2}$,
that is, the optimal solution to the following problem w.r.t.~$\widetilde{\bm v}_+$ and
$\widetilde{\bm \sigma}_2$:
\begin{subequations}
\begin{align*}
\max\limits_{{\bm s},\widetilde{\bm \sigma}_2,{\widetilde{\bm v}_+}}~~& \sum_{i=1}^{2J_{+}+N_2-4+m-1} \log(s_i) \\
{\rm s.t.}\quad & {\bm e}^\top \widetilde{\bm \sigma}_2+{s}_1=1,\; {\bm e}^\top \widetilde{\bm v}_++{s}_2=1, \\
& -\widetilde{\sigma}_{2,j}+{s}_{2+j}=0,\;  j\in[N_2-2], \\ 
& -\widetilde{v}_{+,i}+{s}_{N_2+i}=0,\; \forall i\in[J_+-2], \\
&  {\bm D}_+\widetilde{\bm v}_{+}+{\bm s}_{[J_++N_2-1,2J_++N_2-4]}={\bm d}_+,\\
& h_l\left(\left[\begin{matrix}
-{\bm h}^{{\cal X}_+}(r_{2+}^l)\\
{\bm h}^{{\cal S}_2}(p^l)
\end{matrix}
\right]^\top\left[\begin{matrix}
{\bm R}\widetilde{\bm v}_+\\
{\bm R}\widetilde{\bm  \sigma}_2
\end{matrix}
\right]\right)  +s_{2J_{+}+N_2-4+l}=0, \; \forall l\in [m-1], \\
& {\bm s}\in \R_{+}^{2J_{+}+N_2-4+m-1}, \; \widetilde{\bm \sigma}_2\in \R^{N_2-2}, \; \widetilde{\bm v}_+\in \R_+^{J_+-2},
\end{align*}
\end{subequations} 
where $\widetilde{v}_{+,i}$ denotes the $i$-th element of vector $\widetilde{\bm v}_+$ and $\widetilde{\sigma}_{2,j}$ denotes the $j$-th element of vector $\widetilde{\bm \sigma}_2$, 
${\bm s}_{[J_++N_2-1,2J_++N_2-4]}:=(s_{J_++N_2-1},\cdots, s_{2J_++N_2-4})\in \R^{J_+-2}$.

 \vspace{0.2cm}
{\bf Step 2}: Find the two intersection points between the longest axis of the polyhedron $\Delta_+^{m-1}$  (longest axis of the associated Sonnevend’s ellipsoid) and the boundary of the polyhedron $\Delta_+^{m-1}$, 
denoted by $\left(\widetilde{\bm v}^{1m}_+,\widetilde{\bm \sigma}_2^{1m} \right)$ and  $\left(\widetilde{\bm v}^{2m}_+, \widetilde{\bm \sigma}_2^{2m}\right)$ (see \citep{THS04} and \citep{SaV19}).
 
\vspace{0.2cm}
{\bf Step 3} (Improved):
Let ${\cal D}_+:=\left\{\sum_{j=1}^i c_{2,j}^m : i\in [N_2-2]\right\}$.
Find  $Z_{1}^m$, $Z_{2}^m$ by identifying $p^m$ and $r_{2+}^m$:
\begin{subequations}
\label{eq:gene_m_Q-1}
\begin{align}
& \begin{aligned}
p^{m,i}\in \arg\max_{p\in (0,1)} 
\left\{ \left[\begin{matrix}
{\bm R}\widetilde{\bm v}^{1m}_+\\
{\bm R}\widetilde{\bm \sigma}_2^{1m}
\end{matrix}\right]^\top \left[\begin{matrix}
0\\
{\bm h}^{{\cal S}_2}(p)
\end{matrix}\right] :  \left[\begin{matrix}
{\bm R}{\bm c}_1^{m}\\
{\bm R}{\bm c}_2^{m}
\end{matrix}
\right]^\top
\left[\begin{matrix}
0\\
{\bm h}^{{\cal S}_2}(p)
\end{matrix}\right]\leq D_i
\right\},
\end{aligned}\label{eq:gene_m_Q-1-a}\\
& \begin{aligned}
r_{2+}^{m,i}\in \arg\max_{r_2\in (0,\bar{x})}  
\left\{\left[\begin{matrix}
{\bm R}\widetilde{\bm v}^{2m}_+\\
{\bm R}\widetilde{\bm \sigma}_2^{2m}
\end{matrix}\right]^\top 
\left[\begin{matrix}
{\bm h}^{{\cal X}_+}(r_{2})\\
0
\end{matrix}\right] :
\left[\begin{matrix}
{\bm R}{\bm c}_1^{m}\\
{\bm R}{\bm c}_2^{m}
\end{matrix}\right]^\top
\left[\begin{matrix}
{\bm h}^{{\cal X}_+}(r_{2})\\
0
\end{matrix}\right]  \leq D_i 
\right\},
\end{aligned}
\label{eq:gene_m_Q-1-b}
\end{align}
\end{subequations}
where $D_i\in {\cal D}_+$, for $i\in [N_2-2]$.
We will choose a pair of $(p^{m,i_0},r_{2+}^{m,i_0})$ from the set of optimal solutions $(p^{m,i},r_{2+}^{m,i})$,  $i\in [N_2-2]$ obtained from solving problems \eqref{eq:gene_m_Q-1}.
To ease the exposition, denote the selected pair by $(p^m,r_{2+}^m)$.
This pair leads to the following constraint involving $g_2$ and $v_+$:
$ h_m(rg_{2}(p^m)-rv_+(r_{2+}^m))\leq 0
\Leftrightarrow  
h_m(g_{2}(p^m)-v_+(r_{2+}^m))\leq 0$,
where $h_m=1$ if $\bbe_{g_1g_2}[v(Z_1^m)]\leq \bbe_{g_1g_2}[v(Z_2^m)]$, and $h_m=-1$ otherwise.
Note that the left-hand side (lhs) of the constraint function in \eqref{eq:gene_m_Q-1-a} simplifies to $({\bm R}{\bm c}_2^{m})^\top{\bm h}^{{\cal S}_2}(p)$, 
which is increasing in $p$, whereas the lhs of the constraint function in \eqref{eq:gene_m_Q-1-b} simplifies to $({\bm R}{\bm c}_1^{m})^\top{\bm h}^{{\cal X}_+}(r_{2})$, 
which is increasing in $r_2$. 
Then the optimal solutions $p^{m,i}$ and $r_{2+}^{m,i}$ lie on the boundary of the feasible region of problem \eqref{eq:gene_m_Q-1}.
These properties enable us to construct the cut hyperplane that passes through the analytic center $({\bm c}_1^{m}, {\bm c}_2^{m})$ of the polyhedron $\Delta_+^{m-1}$ (see \eqref{eq:cut_postive} and Proposition~\ref{Prop:CB} below). 
Before introducing the cut hyperplane, we first show that $p^m$ lies within the set of breakpoints ${\cal S}_2$.

\begin{proposition}
[${p}^m\in {\cal S}_2$]
\label{prop:pS_2}
Assume that $(\widetilde{\bm \sigma}_2^{1m},1-{\bm e}^\top\widetilde{\bm \sigma}_2^{1m})$ has no zero component.
Then ${p}^m\in {\cal S}_2$.
\end{proposition}

{\bf Update of Set ${\cal X}_+^0$}.
Let ${\cal X}_+^0:=\{x^+_{1},\cdots,x^+_{J_+}: 0=x^+_{1}<\cdots<x^+_{J_+}=\bar{x}\}$ denote the initial set of breakpoints for $v_+$.
Let $ {\cal X}_+^{m-1}$ be the set of breakpoints of $v_{+}$ in the $(m-1)$-th iteration with $J_+^{m-1}:=|{\cal X}_+^{m-1}|$.
Since $r_{2+}^m\notin {\cal X}_+^{m-1}$ may occur, we can update 
${\cal X}_+^m:={\cal X}_+^{m-1}
\cup \{r_{2+}^{m}\}$
and $J_+^m:=|{\cal X}_+^m|$.
Let $d_j^+=1/(x_{j+1}^+-x_j^+)$ for $j=1,\ldots,J_+^m-1$. Then, update the matrix ${\bm D}_+$ in \eqref{eq:D_+} as follows:
\begin{equation*}
{\bm D}_{{\cal X}_+^m}\hspace{-0.01in}:=\hspace{-0.01in}
\left[\begin{array}{ccccc}
-d_1^+ & d_2^+ & & &\\
&-d_2^+ & d_3^+ & &\\
& &\ddots &\ddots &\\
& & &-d_{J_+^m-3}^+ &d_{J_+^m-2}^+\\
-d_{J_+^m-1}^+ &-d_{J_+^m-1}^+
&\cdots &-d_{J_+^m-1}^+
&\displaystyle \sum_{j=J_+^m-2}^{J_+^m-1}-d_{j}^+
\end{array}
\right]
\end{equation*}
and ${\bm d}_{{\cal X}_+^m}:=\left(0,\cdots,0,\frac{-1}{x_{J_+^m}^+-x_{J_+^m-1}^+}\right)\in \R^{J_+^m-2}$.
Then the initial polyhedron $\Delta_+^0$ of $(\widetilde{\bm v}_+,\widetilde{\bm \sigma}_2)$ is updated to 
$\Delta_{{{\cal X}_+^m}}^0:=\{\left(\widetilde{\bm v}_+,\widetilde{\bm \sigma}_2\right)\in \R^{J_+^m-2}\times \R^{N_2-2}:-\widetilde{\bm \sigma}_2\leq 0, -\widetilde{\bm v}_+\leq 0,
{\bm e}^\top\widetilde{\bm \sigma}_2\leq 1, {\bm e}^\top\widetilde{\bm v}_+\leq 1, {\bm D}_{{\cal X}_+^m}\widetilde{\bm v}_+\leq {\bm d}_{{\cal X}_+^m}\}$.
Combining this with Proposition~\ref{prop:pS_2}, 
we can ensure the existence of 
$p^m\in {\cal S}_2$ and $r_{2+}^m\in {\cal X}_+^m$ such that 
\begin{equation}
\label{eq:equal_star}
g^*_{2}(p^m)-v_+^*(r_{2+}^m)=g^*_{2N_2}(p^m)-v^*_{J_+^m}(r_{2+}^m),
\end{equation}
where $v_+^*$ and $g^*_2$ are the true value and distortion functions, 
and $v_{J_+^m}^*$, $g_{2N_2}^*$ are the corresponding piecewise linear counterparts.
Based on the response of the DM, that is,  
$h_m(g^*_{2}(p^m)-v_+^*(r_{2+}^m))\leq 0$,
it follows from \eqref{eq:equal_star} that $h_m(g^*_{2N_2}(p^m)-v_{J_+^m}^*(r_{2+}^m))\leq 0$.
This means there is no difference between $(v_+^*, g^*_2)$ and $(v_{J_+^m}^*,g_{2N_2}^*)$ 
when we use them to characterize the DM's preference to lotteries whose certain outcome $r_{2+}^m$ lie in  ${\cal X}_+^m$.
However, inconsistency may occur when $p^m\notin{\cal S}_2$ and/or
$r_{2+}^m\notin{\cal X}_+^m$ in which case we say there is a direction error (see \citep[Section 5.2]{ZGSX25}).
The consistency prompts us to construct  a ``cut'' in the space of the piecewise linear functions $h_m(g_{2N_2}(p^m)-v_{J_+^m}(r_{2+}^m))\leq 0$.
The resulting polyhedron in the space of increment vectors can be written as
\begin{equation}
\label{eq:poly_Delta+}
\begin{aligned}
\Delta_+^{m}:=\Delta_{{\cal X}_+^m}^0\cap \Big\{\left(\widetilde{\bm v}_+,\widetilde{\bm \sigma}_2\right)\in \R^{J_+^m-2}\times \R^{N_2-2}: l_{\rm cut,+}^{l,m}(\widetilde{\bm v}_+,\widetilde{\bm \sigma}_2)\leq 0,\; l\in [m]\Big\},
\end{aligned}
\end{equation}
where $\Delta_{{\cal X}_+^m}^0$ is updated according to the change of the set of breakpoints ${\cal X}_+^m$ and the dimension of the space of vectors of increments of the PLFs, where
\begin{equation}
\label{eq:cut_postive}
\begin{aligned}
l_{\rm cut,+}^{l,m}(\widetilde{\bm v}_+,\widetilde{\bm \sigma}_2):= h_l\big(({\bm R}\widetilde{\bm \sigma}_2)^\top{\bm h}^{{\cal S}_2}(p^l) - ({\bm R}{\widetilde{\bm v}_+})^\top{\bm h}^{{\cal X}_+^m}(r_{2+}^l)\big)
\end{aligned}
\end{equation}
defines the $l$-th cut hyperplane through
$l_{\rm cut,+}^{l,m}(\widetilde{\bm v}_+,
\widetilde{\bm \sigma}_2)=0$ constructed using the set of breakpoints ${\cal X}_+^m$. 
We will prove that the cut hyperplane $l_{\rm cut,+}^{m,m-1}(\widetilde{\bm v}_+, \widetilde{\bm \sigma}_2)$ passes through the analytic center $({\bm c}_1^m,{\bm c}_2^m)\in \R^{J_+^{m-1}-2}\times \R^{N_2-2}$ of the polyhedron $\Delta_+^{m-1}$, where 
$\Delta_+^{m-1}= {\Delta_{{\cal X}_+^{m-1}}^0}
\cap\Big\{ \left(\widetilde{\bm v}_+,\widetilde{\bm \sigma}_2\right)\in \R^{J_+^{m-1}-2}\times \R^{N_2-2}:
l_{\rm cut,+}^{l,{m-1}}(\widetilde{\bm v}_+,\widetilde{\bm \sigma}_2)\leq 0,\; l\in [m-1]\Big\}$.
Before proceeding to theoretical analysis, we make some comments on the modified polyhedral method.

\begin{remark} 
The proposed modified polyhedral method differs from the original polyhedral method in at least two respects:
(a) the polyhedra are defined in the space of the vectors of the increments 
of the piecewise linear value function and the distortion function, 
(b) the maximization problems in \eqref{eq:gene_m_Q-1} mimic the original polyhedral method but they differ significantly in the definition of the analytic center, the identification of intersection points between the longest axis of the associated Sonnevend’s ellipsoid and the boundaries of the polyhedron in Step 2 and the structures of the objective function and constraint functions.
\end{remark}

\begin{remark}
[Polyhedral Method vs. RRUS]
When there is no information about the shape and the Lipschitz modulus of the utility function in EUT, the polyhedral method exhibits better performance than RRUS, as Example~\ref{exa:comparison} shows.
\end{remark}

\subsubsection{Theoretical Guarantees}
\label{sec:thm-gua}
The next two propositions state that the cut hyperplane passes through the analytic center of the polyhedron, and the size of the polyhedron decreases as the number of queries increases.

\begin{proposition}
[Choice Balance]
\label{Prop:CB}
If the vectors $\widetilde{\bm \sigma}_2^{1m}$ and $\widetilde{\bm v}_+^{2m}$ have no zero components, then the cut hyperplane $l_{\rm cut,+}^{m,{m-1}}(\cdot)$
passes through the analytic center ${\bm c}^m$ of the polyhedron $\Delta_+^{m-1}$ defined as in \eqref{eq:poly_Delta+}.
\end{proposition}

Step $3$ is the key step in the modified polyhedral method.
It determines a hyperplane $l_{\rm cut,+}^{m,{m-1}}(\cdot)$ which cuts the polyhedron ${\Delta}_+^{m-1}$ into two parts through the analytic center. 
The next proposition states that the size of the polyhedron decreases as
$m$ increases.

\begin{proposition}
[Reducing the Size of the Polyhedron $\Delta_+^m$]
\label{prop:decrease_AC}
Assume that ${\cal X}_+^m={\cal X}_+^0$ for all $m\geq0$. 
Let ${\bm c}^m$ denote the analytic center of $\Delta_+^{m-1}$. 
Then the feasible region $\Delta_+^{m-1}$ shrinks as $m$ increases.
\end{proposition}

We refer readers to a concrete example in Section~\ref{SC:exa-polyhedral}
which shows how the algorithm works when ${\cal X}_+^m={\cal X}_+^0$, even though a direction error due to $r_{2+}^m\notin {\cal X}_+^{m-1}$ may occur.

\begin{remark}
\label{rem:improved_PM}
Before discussing Question Type 3, we make the following comments.
\begin{itemize}
\item[(a)] 
The conditions of Proposition~\ref{prop:decrease_AC} may not be satisfied in practice, since $r_{2+}^m\notin {\cal X}_+^{m-1}$ may occur.
However, the actual reduction in the distance from ${\bm c}^m$ to the boundaries of the polyhedron $\Delta_+^{m-1}$ is illustrated in Figure~\ref{fig:distance} even though the condition that ${\cal X}_+^m={\cal X}_+^0$ may not be satisfied.

\item[(b)] 
Like $\Delta_+^m$, we can apply the modified polyhedral method to shrink the polyhedron $\Delta_-^{m}$ of $(\widetilde{\bm v}_-,\widetilde{\bm \sigma}_1)$,
with the initial polyhedron 
$\Delta_-^0 := \big\{ (\widetilde{\bm v}_-, \widetilde{\bm \sigma}_1)\in \R^{J_{-}-2}\times \R^{N_1-2}:
-\widetilde{\bm \sigma}_1\leq 0,
-\widetilde{\bm v}_-\leq 0,
{\bm e}^\top\widetilde{\bm \sigma}_1\leq 1, {\bm e}^\top\widetilde{\bm v}_-\leq 1, {\bm D}_-\widetilde{\bm v}_-\leq {\bm d}_- \big\}$,
where ${\bm D}_-$ and ${\bm d}_-$ can be defined according to ${\bm D}_+$ and ${\bm d}_+$ as in \eqref{eq:D_+}.
After $m-1$ queries, the resulting polyhedron is given by
\begin{equation}
\label{eq:Delta_minus}
\begin{aligned} 
\Delta_-^{m-1}=\Delta^0_{{\cal X}_-^{m-1}} \cap
\Big\{\left(\widetilde{\bm v}_-,\widetilde{\bm \sigma}_1\right) \in \mathbb{R}^{J_-^{m-1}-2}\times \mathbb{R}^{N_1-2}: l_{\rm cut,-}^{l,{m-1}}(\widetilde{\bm v}_-,\widetilde{\bm \sigma}_1)\leq 0,\; l\in [m-1]\Big\}, 
\end{aligned} 
\end{equation}
where $l_{\rm cut,-}^{l,{m-1}}(\widetilde{\bm v}_-,
\widetilde{\bm \sigma}_1):=h_l \big(-({\bm R}\widetilde{\bm \sigma}_1)^\top{\bm h}^{{\cal S}_1}(1-p^l)-({\bm R}\widetilde{\bm v}_-)^\top({\bm h}^{{\cal X}_-^{m-1}}(r_{2-}^l)-{\bm e})\big)$,
$\Delta^0_{{\cal X}_-^{m-1}}$ is the updated version of $\Delta_-^0$ corresponding to the current set ${\cal X}_-^{m-1}$, and $J_-^{m-1} = |{\cal X}_-^{m-1}|$.
Denote the numbers of queries for Question Types~1 and~2 by $M_1$ and $M_2$, respectively, 
and let the corresponding polyhedra be $\Delta_-^{M_1}$ and $\Delta_+^{M_2}$.

\item[(c)] Proposition~\ref{prop:decrease_AC} and the numerical tests pave the way for us to use an estimate of the distortion functions by the analytic center of the polyhedra $\Delta_-^{M_1}$ and $\Delta_{+}^{M_2}$ in the forthcoming design of Question Type~4, 
where only the holistic value function needs to be determined.
To avoid propagation of the errors occurring in the observation of the inflection point of the underlying distortion functions, 
we do not incorporate the inverse $S$-shape of the distortion function when eliciting the polyhedra $\Delta_-^{M_1}$ and $\Delta_{+}^{M_2}$, 
and in the forthcoming estimate of the distortion functions.

\item[(d)] 
It is important to note that the polyhedral method is based on the setup that the DM's value function and the distortion functions have a specified piecewise linear structure. 
The modified polyhedral method provides an efficient (optimal) way to reduce the ambiguity set of piecewise linear value and distortion functions. 
The estimated functions $\widehat{g}_1,\widehat{g}_2$, constructed by the analytic centers of the resulting feasible sets $\Delta_-^{M_1}$ and $\Delta_{+}^{M_2}$ of the vectors of the increments,
are the ``best estimates'' that we can obtain after $M_1$ and $M_2$ cuts from the specified class of piecewise linear functions.

Specifically, after obtaining the DM’s preference information from $M_1$ and $M_2$ queries, we define $J_-^{M_1} := |{\cal X}_-^{M_1}|$ and $J_+^{M_2} := |{\cal X}_+^{M_2}|$ as the cardinalities of the updated breakpoint sets.
We denote by $\Delta_-^{M_1} \subset \mathbb{R}^{J_-^{M_1} + N_1 - 4}$ 
and $\Delta_+^{M_2} \subset \mathbb{R}^{J_+^{M_2} + N_2 - 4}$ the resulting polyhedra corresponding to the variables $(\widetilde{\bm v}_-, \widetilde{\bm \sigma}_1)$ and $(\widetilde{\bm v}_+, \widetilde{\bm \sigma}_2)$, respectively.
Let ${\bm c}_-^{M_1} =({\bm c}_{-,1}^{M_1},{\bm c}_{-,2}^{M_1})\in \Delta_-^{M_1}$
and ${\bm c}_+^{M_2} =({\bm c}_{+,1}^{M_2},{\bm c}_{+,2}^{M_2}) \in \Delta_+^{M_2}$
be the corresponding analytic centers.
Let ${\cal S}_2$ be defined as in \eqref{eq:g2v+} with $|{\cal S}_2|=N_2$.
Then the estimates are given~by
\begin{equation}
\label{eq:hat_g1g2}
\begin{aligned}
\widehat{g}_1(t):=\left({\bm R}{\bm c}_{-,2}^{M_1}\right)^\top{\bm h}^{{\cal S}_1}(t)
\quad \inmat{and} \quad 
\widehat{g}_2(t):=\left({\bm R}{\bm c}_{+,2}^{M_2}\right)^\top{\bm h}^{{\cal S}_2}(t).
\end{aligned}
\end{equation}
\end{itemize}
\end{remark}

In practice, we cannot calculate $v(\bar{x})$ and $v(\underline{x})$ or even estimate them. On the other hand, as discussed in Section~\ref{sec:PLFs},
scaling the value function does not affect the optimal solution of the robust shortfall \eqref{eq:PRGRS-CPT}.
This prompts us to identify the ratio $r=-{v(\bar{x})}/{v(\underline{x})}$ first
and then use the normalized negative part $v_-(x)$ (with range $[-1,0]$) and the normalized positive part $v_+(x)$ (with range $[0,1]$) to construct an ambiguity set of the value functions satisfying ${v}(x)=rv_+(x)$ when $x\geq 0$, and ${v}(x)=v_-(x)$ when $x<0$.
To better understand the design of Question Type~$3$, let us explain its structure in the next subsection.

\subsection{Question Type  \texorpdfstring{$3$}{} for the Estimation of the Range of the Ratio \texorpdfstring{$r$}{}}
\label{Sec:Estimation_r}
We begin by developing a strategy to estimate the ratio~$r$. To this end, we design a type of question that enables us to formulate a linear system of inequalities in the ratio~$r$ based on the DM’s responses to the queries.

\noindent{\bf Question Type 3.}
Generate lotteries 
\begin{equation*}
Z_1^m=\left\{\begin{aligned}
\underline{x} \; & \mbox{w.p. } 1-p^m\\
\bar{x} \; & \mbox{w.p. } p^m
\end{aligned}
\right.
\quad \mbox{and} \quad 
Z_2^m=0,\;\inmat{for } m\in [M_3].
\end{equation*}
Then $\bbe_{\widehat{g}_1\widehat{g}_2}[v(Z_1^m)]=v(\underline{x})\widehat{g}_1(1-p^m)+ v(\bar{x})\widehat{g}_2(p^m)$ and $\bbe_{\widehat{g}_1\widehat{g}_2}[v(Z_2^m)]=0$.
The DM's preference to lotteries can be formulated as $h_m(\bbe_{\widehat{g}_1\widehat{g}_2}[v(Z_1^m)]-\bbe_{\widehat{g}_1\widehat{g}_2}[v(Z_2^m)])\leq 0$,
where  $h_m=1$ if $Z_2^m$ is preferred, and $h_m=-1$ otherwise.
Let $r=-\frac{v(\bar{x})}{v(\underline{x})}$. 
Then the inequality can be recast as $h_m(\widehat{g}_1(1-p^m)-r \widehat{g}_2(p^m))\geq0$.
Let $[\underline{r}^m,\bar{r}^m]:=\{r\in \R_+:h_l(\widehat{g}_1(1-p^l)-r \widehat{g}_2(p^l))\geq0, l\in [m]\}$.
Next, we will reduce the range of $r$ by reducing the range of $p^m\in [0,1]$.
Observe that when $p^m=0$, we have $Z_1^m=\underline{x}$ and $Z_1^m \preceq Z_2^m$, and when $p^m=1$, we have $Z_1^m=\bar{x}$ and $Z_1^m \succeq Z_2^m$.
Moreover, let $\phi(p^m):=
\bbe_{\widehat{g}_1\widehat{g}_2}[v(Z_1^m)]$.
Then
$\phi(p^m):= v(\underline{x})\widehat{g}_1(1-p^m)+ v(\bar{x})\widehat{g}_2(p^m)$.

\begin{figure}[!ht]
\centering
\begin{minipage}[t]{0.48\textwidth}
\vspace{0pt}
\centering
\captionof{figure}{Changes in the Range $[\underline{r}^m,\bar{r}^m]$
of the Ratio $r=-\frac{v(\bar{x})}{v(\underline{x})}$
as the Number of Type 3 Queries $M_3$ Increases from $2$ to $9$.}
\label{fig:range_R}
\includegraphics[width=3.0in]{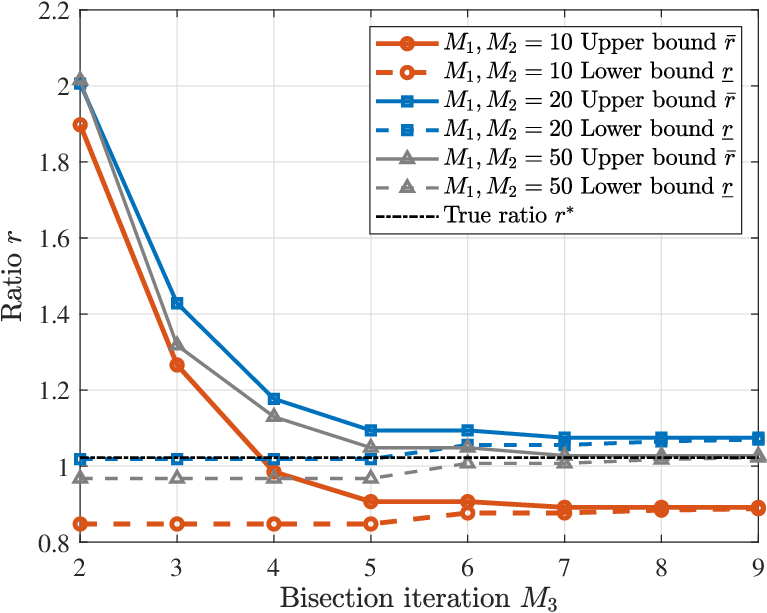}
\end{minipage}
\hfill
\begin{minipage}[t]{0.49\textwidth}
\vspace{0pt}
\centering
\captionof{table}{Bisection Method for Solving $\phi(p^m)=v(\underline{x})\widehat{g}_1(1-p^m)+ v(\bar{x})\widehat{g}_2(p^m)=0$. Note that $r_2^m=0$, $[p^{0,L},p^{0,U}]=[0,1]$.}
\label{tab:bisection_r2_final_23}
{\small
\setlength{\tabcolsep}{3pt}
\begin{tabular}{c|c|c|c|c}
\hline
\multicolumn{5}{c}{\textbf{Bisection Method}} \\
\hline
${m}$ & ${p^{m-1,L}}$ & ${p^{m-1,U}}$ & $p^m$ & $h_m\in \{-1,1\}$\\
\hline
1 & 0.000000 & 1.000000 & 0.500000 & -1 \\
2 & 0.000000 & 0.500000 & 0.250000 & 1 \\
3 & 0.250000 & 0.500000 & 0.375000 & 1 \\
4 & 0.375000 & 0.500000 & 0.437500 & 1 \\
5 & 0.437500 & 0.500000 & 0.468750 & 1 \\
6 & 0.468750 & 0.500000 & 0.484375 & -1 \\
7 & 0.468750 & 0.484375 & 0.476562 & 1 \\
8 & 0.476562 & 0.484375 & 0.480469 & -1 \\
9 & 0.476562 & 0.480469 & 0.478516 & -1\\
10 & 0.476562 & 0.478516 & 0.477539 & $\simeq$ \\
\hline
\end{tabular}
\par}
{\footnotesize
\raggedright
\quad For $m=10$, ``$\simeq$'' means $Z_1^{10}\simeq Z_2^{10}$.
\par}
\end{minipage}
\end{figure}

As $p^m$ increases from $0$ to $1$, $\phi(p^m)$ strictly increases from $v(\underline{x})<0$ to $v(\bar{x})>0$.
There must therefore exist some $p^* \in (0,1)$ such that 
$\phi(p^*)=v(\underline{x})\widehat{g}_1(1-p^*)+ v(\bar{x})\widehat{g}_2(p^*)=0$.
Our strategy to reduce the ambiguity of (the range of) $r$ is through identifying $p^*$ by employing the well-known bisection method. 
This is illustrated in the process in Table~\ref{tab:bisection_r2_final_23}. 
Based on the ten queries in Table~\ref{tab:bisection_r2_final_23} and the conclusion that $\phi(p^{10})= 0$, we can conclude that $\widehat{r}=\frac{\widehat{g}_1(1-p^{10})}{\widehat{g}_2(p^{10})}= 1.01850$,
which is close to the true ratio $r^*=1.02257$, with only $M_3=10$.
Given that the distortion functions $\widehat{g}_1$ and $\widehat{g}_2$ are estimated, we obtain better estimates of $r$ as the number of Question Types~1 and~2 increases (see Figure~\ref{fig:range_R}).

After asking Question Types 1-3 with numbers of queries $M_1,M_2$ and $M_3$,
we can construct the polyhedron $\Delta_-^{M_1}$ in \eqref{eq:Delta_minus} for $(\widetilde{\bm v}_-,\widetilde{\bm \sigma}_1)$ and the polyhedron $\Delta_+^{M_2}$ in \eqref{eq:poly_Delta+} for $(\widetilde{\bm v}_+,\widetilde{\bm \sigma}_2)$ and then estimate the ratio $\widehat{r}$. 
By the definition of the value function $\widetilde{v}$ in \eqref{eq:value_tilde}, 
based on the ranges of $\widetilde{\bm v}_-$, 
$\widetilde{\bm v}_+$ and $\widehat{r}$, 
we can construct the value function over $[\underline{x},\bar{x}]$ via $v(x)=v_-(x)$ for $x<0$ and $v(x)=rv_+(x)$ for $x>0$.
Then the ambiguity set can take the form: 
\begin{equation}
\label{eq:Question_1_3}
(\mathbb{V}{\mathbb{G}_1\mathbb{G}_2})_M^0=\left\{
(v,g_1,g_2)\in \mathscr{V}\times \mathscr{G}\times \mathscr{G}: 
\begin{aligned}
&\; g_{iN_i}(t)=({\bm R}\widetilde{\bm \sigma}_i)^\top{\bm h}^{{\cal S}_i}(t),\; i=1,2\\
&\;  \left({\widetilde{\bm v}_-},\widetilde{\bm \sigma}_1\right)\in \Delta_-^{M_1},\left({\widetilde{\bm v}_+},\widetilde{\bm \sigma}_2\right)\in \Delta_+^{M_2}\\
&\;  v(x)=({\bm R}\widetilde{\bm v}_-)^\top{\bm h}^{{\cal X}_-^{M_1}}(x)  +\widehat{r}({\bm R}\widetilde{\bm v}_+)^\top{\bm h}^{{\cal X}_+^{M_2}}(x)
\end{aligned}
\right\}.
\end{equation}
However, the ambiguity set $(\mathbb{V}{\mathbb{G}_1\mathbb{G}_2})_M^0$ is still 
large and the resulting robust shortfall \eqref{eq:PRGRS-CPT-equiv} is too conservative.
To strengthen the modified polyhedral method, we will propose a new approach in the next subsection to reduce the range of holistic value functions over the whole interval $[\underline{x},\bar{x}]$ by designing Question Type 4,
where the outcomes of $Z_1^m$ and $Z_2^m$ can take any point in $(\underline{x},\bar{x})$.
To apply the modified polyhedral method, we need to normalize the value function over the interval $[\underline{x},\bar{x}]$.  
However, as illustrated in \eqref{eq:normalization_extra},
normalization of the value function does not generally preserve the preference inequalities under CPT, as we explain in the next remark.

\begin{remark}
\label{rem:normalization}
Consider the lottery $Z_1^m$ with outcomes $r_1^m<0$ and $r_3^m>0$ and the certain lottery $Z_2^m$ with $r_2^m$ such that $r_1^m\leq r_2^m \leq r_3^m$.
In this case, we can use the estimated distortion functions $\widehat{g}_1,\widehat{g}_2$  to establish $\bbe_{\widehat{g}_1\widehat{g}_2}[v(Z_1^m)]\leq \bbe_{\widehat{g}_1\widehat{g}_2}[v(Z_2^m)]$ if the DM prefers $Z_2^m$. 
This inequality can be equivalently written as 
\begin{equation}
\label{eq:elicit-given-g-1}
v(r_1^m)\widehat{g}_1(1-p^m)+v(r_3^m)\widehat{g}_2(p^m)\leq v(r_2^m).
\end{equation} 
Let $\hat{v}(x):=\frac{v(x)-v(\underline{x})}{v(\bar{x})-v(\underline{x})}$. 
Since $\frac{v(\underline{x})}{v(\bar{x})-v(\underline{x})}=\frac{-1}{{r}+1}$,
then by normalization of the value function $v$, we have
\begin{equation}
\label{eq:normal-term} 
\begin{aligned}
v(r_1^m)\widehat{g}_{1}(1-p^m) +v(r_3^m)\widehat{g}_{2}(p^m) 
=\hat{v}(r_1^m)\widehat{g}_1(1-p^m)+\hat{v}(r_3^m)\widehat{g}_2(p^m)  - \frac{(\widehat{g}_1(1-p^m) + \widehat{g}_2(p^m))}{{r}+1}
 \end{aligned}
\end{equation}
and $v(r_2^m)=\frac{v(r_2^m)-v(\underline{x})}{v(\bar{x})-v(\underline{x})}+\frac{v(\underline{x})}{v(\bar{x})-v(\underline{x})} =\hat{v}(r_2^m)+\frac{-1}{{r}+1}$.
Since the term $1-\widehat{g}_1(1-p^m)-\widehat{g}_2(p^m)\neq 0$ in general, 
the term on the right-hand side (rhs) of \eqref{eq:normal-term} does not disappear,
which means that \eqref{eq:elicit-given-g-1} is not equivalent to $\hat{v}(r_1^m)\widehat{g}_1(1-p^m)+\hat{v}(r_3^m)\widehat{g}_2(p^m) \leq \hat{v}(r_2^m)$.
Based on the observations, we will introduce the fourth question in the next subsection.
\end{remark}

\subsection{Modified Polyhedral Method for Refining the Holistic Value Function \texorpdfstring{$v$}{v}}
\label{sec:ratior}
In this section, we discuss how to elicit the holistic value function based on the estimates $\widehat{g}_1$ and $\widehat{g}_2$ (see \eqref{eq:hat_g1g2}).
In Sections~\ref{sec:pre_utility}--\ref{sec:Thm-Gue2}, 
we clarify the initial polyhedron of the increments of the value function ${v}$, 
the models for generating the $m$-th query and the theoretical guarantees, respectively.\\
{\bf Question Type 4}.  Generate lotteries: 
\begin{equation*}
Z_1^m=\left\{\begin{aligned}
r_1^m \; & \inmat{w.p.}~1-p^m\\
r_3^m \; & \inmat{w.p.}~p^m
\end{aligned}
\right.
\quad \mbox{and} \quad 
Z_2^m=r_2^m,\;\inmat{for } m\in [M_4],
\end{equation*}
where $r_1^m<0<r_3^m$, $r_1^m\leq r_2^m \leq r_3^m$, $p^m\in (0,1)$ for $m\in [M_4]$.
Given $\widehat{g}_1$ and $\widehat{g}_2$ as described in \eqref{eq:hat_g1g2},
if we choose a value for $p^m$, then $\widehat{g}_1(1-p^m)$ and $ \widehat{g}_2(p^m)$ are known, 
and the only function to be identified is the value function $v$.
Since $\bbe_{\widehat{g}_1\widehat{g}_2}[v(Z_1^m)] = v(r_1^m)\widehat{g}_1(1-p^m)
+v(r_3^m)\widehat{g}_2(p^m)$ and  $\bbe_{\widehat{g}_1\widehat{g}_2}[v(Z_2^m)]=v(r_2^m)$,
the DM's preference to the lotteries can be formulated as 
$h_m\left(\bbe_{\widehat{g}_1\widehat{g}_2}[v(Z_1^m)]-\bbe_{\widehat{g}_1\widehat{g}_2}[v(Z_2^m)]\right)\leq 0$, where  $h_m=1$ if the DM prefers $Z_2^m$ and $h_m=-1$ otherwise.
Based on the observations in Remark~\ref{rem:normalization}, 
we will introduce some preference functionals to design queries in the next subsection.

\subsubsection{Set-up}
\label{sec:pre_utility}
Let ${\cal X}_0={\cal X}_-^{M_1}\cup {\cal X}_+^{M_2}$ be the initial set of breakpoints of the value function over $[\underline{x},\bar{x}]$ and define $J_0:=|{\cal X}_0|=J_-^{M_1}+J_+^{M_2}-1$, 
where ${\cal X}_-^{M_1}\cap {\cal X}_+^{M_2}=\{0\}$.
Based on $\bbe_{\widehat{g}_1\widehat{g}_2}[v(Z_1^m)]$ and  $\bbe_{\widehat{g}_1\widehat{g}_2}[v(Z_2^m)]$, 
the preference functional based on the piecewise linear approximation is defined as follows:
\begin{subequations}
\label{eq:pre_func}
\begin{align}
&\varrho_{1}((r_1^m,r_3^m,p^m); {\bm v}) :=  ({\bm R}{\bm v})^\top {\bm h}^{{\cal X}_0}(r_1^m) \widehat{g}_1(1-p^m)
+({\bm R}{\bm v})^\top{\bm h}^{{\cal X}_0}(r_3^m)\widehat{g}_2(p^m)  -\frac{(\widehat{g}_1(1-p^m)+\widehat{g}_2(p^m))}{\widehat{r}+1},\\
&\varrho_{2}(r_2^m;{\bm v}):=({\bm R}{\bm v})^\top {\bm h}^{{\cal X}_0}(r_2^m)+\frac{-1}{\widehat{r}+1}.
\end{align}
\end{subequations}
Note that we normalize the value function over $[\underline{x},\bar{x}]$
and denote it by $v$ temporarily in this subsection.
Accordingly, the vector ${\bm v}\in \R^{J_0-2}$ in \eqref{eq:pre_func} satisfies ${\bm e}^\top{\bm R}{\bm v}=1$.
Moreover, the vectors $\widetilde{\bm v}_-\in \R^{J_-^{M_1}-2}$, 
$\widetilde{\bm v}_+\in \R^{J_+^{M_2}-2}$ and the ratio $\widehat{r}$ in Question Types 1--3 satisfy ${\bm R}{\bm v}=\left(\frac{1}{1+\widehat{r}}{\bm R}\widetilde{\bm v}_-,\frac{\widehat{r}}{1+\widehat{r}}{\bm R}\widetilde{\bm v}_+\right)$.
Let ${\Delta}_0$ be the initial polyhedron of the increment vector of the value function over $[\underline{x},\bar{x}]$. 
To incorporate the information on the negative and positive parts $v_-$ and $v_+$ of the value function, 
as well as the ratio $\widehat{r}$ obtained from Question Types 1--3 (see \eqref{eq:Question_1_3}), 
and to ensure feasibility of the resulting joint ambiguity set $\mathbb{V}\mathbb{G}_1\mathbb{G}_2$ under Question Types 1--4, 
we define ${\Delta}_0$ as
$
{\Delta}_0=\Big\{  
{\bm v}\in \R^{J_0-2}: {\bm R}{\bm v}=\Big(\frac{1}{1+\widehat{r}}{\bm R}\widetilde{\bm v}_-,\frac{\widehat{r}}{1+\widehat{r}}{\bm R}\widetilde{\bm v}_+\Big),
 (\widetilde{\bm v}_-,
{\bm c}_{-,2}^{M_1}) \in \Delta_-^{M_1},\;
(\widetilde{\bm v}_+,
{\bm c}_{+,2}^{M_2})\in \Delta_+^{M_2} \Big\},
$
where ${\bm c}_{-,2}^{M_1}$ and ${\bm c}_{+,2}^{M_2}$ denote the $\sigma$-components of the analytic centers ${\bm c}_-^{M_1}=({\bm c}_{-,1}^{M_1},{\bm c}_{-,2}^{M_1})\in \R^{J_-^{M_1}+N_1-4}$ and ${\bm c}_+^{M_2}=({\bm c}_{+,1}^{M_2},{\bm c}_{+,2}^{M_2})\in \R^{J_+^{M_2}+N_2-4}$ of the polyhedra $\Delta_-^{M_1}$ and $\Delta_{+}^{M_2}$, respectively (see \eqref{eq:hat_g1g2}). 
Next, we will develop a modified polyhedral method to elicit the value function,
which will circumvent the issue in Remark~\ref{rem:normalization}.

\subsubsection{Design of a New Query} 
\label{sec:gene_guarantee} 
In {\bf Steps 1--2}, let ${\bm c}^m$ be the analytic center of $\Delta_{m-1}$ and ${\bm v}^m_1$ and ${\bm v}^m_2$ be the two intersection points between the longest axis of Sonnevend's ellipsoid associated with the polyhedron $\Delta_{m-1}$ and the boundary of the polyhedron.

\vspace{0.2cm}

{\bf Step 3}.
The two problems generating the $m$-th pair of queries are defined as
\begin{subequations}
\label{eq:problem_Q3}
\begin{align}
r_2^{m,i}
&\in \arg\max\big\{\varrho_2(x;{\bm v}_1^m): x\in[\underline{x},\bar{x}],
\;\varrho_2(x;{\bm c}^m)\leq D_i\big\}, \label{eq:problem_Q3-a}\\
(r_1^{m,i},r_3^{m,i},p^{m,i})
&\in \arg\max\left\{\varrho_1((x,y,p);{\bm v}_2^m):
\begin{aligned}
&x\in[\underline{x},0),\;y\in(0,\bar{x}], p\in(0,1)\\
&\quad \varrho_1((x,y,p);{\bm c}^m) \leq D_i
\end{aligned}
\right\},
\label{eq:problem_Q3-b}
\end{align}
\end{subequations} 
where $\varrho_1$ and $\varrho_2$ are defined as in \eqref{eq:pre_func}, and $D_i\in {\cal D}$ is a constant. 
To guarantee the feasibility of problems~\eqref{eq:problem_Q3}, $D_i$ satisfies $D_i\in \big(\frac{-1}{\widehat{r}+1},\frac{\widehat{r}}{\widehat{r}+1}\big)$.
We select some tuple $(r_1^{m,i_0}, r_2^{m,i_0}, r_3^{m,i_0}, p^{m,i_0})$ 
from the set of optimal solutions $\{(r_1^{m,i}, r_2^{m,i}, r_3^{m,i}, p^{m,i}):i = 1, \dots, |{\cal D}|\}$,
and denote it by $(r_1^{m}, r_2^{m}, r_3^{m}, p^{m})$.

\subsubsection{Theoretical Guarantees}
\label{sec:Thm-Gue2}
The next proposition states that the cut hyperplane passes through the analytic center of the polyhedron.

\begin{proposition}
\label{prop:Q3}
Assume that the intersection points ${\bm v}^m_1$ and ${\bm v}^m_2$ in problem \eqref{eq:problem_Q3} have no zero components and $r_1^{m,i}, r_2^{m,i}, r_3^{m,i} \in {\cal X}_{m-1}$.
Then the optimal solutions $r_2^{m,i}$ and $(r_1^{m,i},r_3^{m,i},p^{m,i})$ satisfy $ \varrho_2(r_2^{m,i};{\bm c}^m)=D_i$ and $\varrho_1((r_1^{m,i},r_3^{m,i},p^{m,i});{\bm c}^m)=D_i$, 
and the resulting cut hyperplane passes through the analytic center of the polyhedron $\Delta_{m-1}$.
\end{proposition}
Proposition~\ref{prop:Q3} can be proved analogously to Proposition~\ref{Prop:CB}.
Proposition~\ref{prop:Q3} assumes that the generated outcomes of $Z_1^m$ and $Z_2^m$ lie within the current set of the breakpoints ${\cal X}_{m-1}$, 
and the feasibility of the scheme in Question Type~4.
In practice, the conditions may not be satisfied, 
and we can update the set ${\cal X}_m:={\cal X}_{m-1} \cup \{r_1^m,r_2^m,r_3^m\}$ with $J_m:=|{\cal X}_m|$, and
\begin{equation}
\label{eq:update_M1M2}
\begin{aligned}
&{\cal X}_-^{M_1} \leftarrow  {\cal X}_-^{M_1}\cup  (\{r_1^m,r_2^m,r_3^m\} \cap [\underline{x},0]),\\
&{\cal X}_+^{M_2} \leftarrow {\cal X}_+^{M_2}\cup (\{r_1^m,r_2^m,r_3^m\}\cap [0,\bar{x}]).
\end{aligned}
\end{equation}
Subsequently, the polyhedral feasible set of ${\bm v}$ is defined as follows:
\begin{equation}
\label{eq:whole_uti}
\Delta_{m}:=\Delta_0^{{\cal X}_{m}} \cap  \left\{{\bm v}\in \R^{J_{m}-2}: 
\phi_{\rm cut}^{l,m}\left({\bm v}\right) \leq 0, l\in [m] \right\},
\end{equation}
where the updated initial polyhedron $\Delta_0^{{\cal X}_m}$ is defined as
\begin{equation}
\label{eq:update_Delta0}
\begin{aligned}
\Delta_0^{{\cal X}_m}
=\Big\{ {\bm v}\in \R^{J_m-2}:
{\bm R}{\bm v}=\left(\frac{1}{1+\widehat{r}}{\bm R}\widetilde{\bm v}_-,\frac{\widehat{r}}{1+\widehat{r}}{\bm R}\widetilde{\bm v}_+\right),   (\widetilde{\bm v}_-,
{\bm c}_{-,2}^{M_1}) \in \Delta_-^{M_1},\;\;
(\widetilde{\bm v}_+,
{\bm c}_{+,2}^{M_2})\in \Delta_+^{M_2}
\Big\}
\end{aligned}
\end{equation}
and for each $l\in [m]$,
\begin{equation}
\label{eq:phi_cutm}
\begin{aligned}
\phi_{\rm cut}^{l,m}({\bm v}):= h_l\Big(({\bm R}{\bm v})^\top({\bm h}^{{\cal X}_m}(r_1^l) \widehat{g}_1(1-p^l) 
+{\bm h}^{{\cal X}_m}(r_3^l)\widehat{g}_2(p^l) -{\bm h}^{{\cal X}_m}(r_2^l)) 
-\frac{\widehat{g}_1(1-p^l)+\widehat{g}_2(p^l)-1}{\widehat{r}+1}\Big).
\end{aligned}
\end{equation}

\begin{remark}
We make some comments on the roles played by the $S$-shaped value function and the inverse $S$-shaped distortion functions in the elicitation process.
In Question Type~1, the convexity of the value function over $[\underline{x},0]$ is used.
But the inverse $S$-shape of the distortion function $g_1$ is not used, 
since we do not know the inflection point (i.e., concavity-convexity).
In Question Type~2, 
the concavity of the value function over $[0,\bar{x}]$ is used.
But the inverse $S$-shape of the distortion function $g_2$ is not used, 
since we do not know the inflection point.
In Question Type~4, the $S$-shape of the value function over the whole interval $[\underline{x},\bar{x}]$ is imposed.
When estimated distortion functions $\widehat{g}_1$ and $\widehat{g}_2$ and numerical approximation are involved, 
relaxing the $S$-shape constraint may improve the feasibility of the elicitation procedure.
Note that we obtain good performance on the holistic value function elicitation in Figures~\ref{fig:value} and~\ref{fig:CPT_PLA-4}.
\end{remark}

Based on Question Types 1--4, 
we can construct the joint ambiguity set $\mathbb{V}{\mathbb{G}_1\mathbb{G}_2}$ of value and distortion functions.
Before doing so, let us clarify the relation between $\Delta_-^{M_1}$, $\Delta_+^{M_2}$ and $\Delta_{M_4}$ in Remark~\ref{rem:improved_PM} (b), \eqref{eq:poly_Delta+}, and \eqref{eq:whole_uti}.
To this end, we map the vector ${\bm v}$ of increments of the holistic value function defined in \eqref{eq:pre_func} into vectors ${\bm v}_-$ and ${\bm v}_+$ for the negative and the positive parts of the value function, defined in \eqref{eq:PLA_utility_minus} and \eqref{eq:PLA_utility_2}.

\begin{remark}
Note that the DM's preference to the certainty lottery $Z_2^l$ in Question Type $4$ is formulated as follows:
\begin{subequations}
\begin{align}
 \bbe_{\widehat{g}_1\widehat{g}_2}[v(Z_1^l)]-\bbe_{\widehat{g}_1\widehat{g}_2}[v(Z_2^l)] \leq 0 
&\quad \Leftrightarrow  v(r_1^l) \widehat{g}_{1}(1-p^l) +v(r_3^l)\widehat{g}_{2}(p^l)
-v(r_2^l)\leq 0 \label{eq:expand_preference}\\
& \quad  \Leftrightarrow   \frac{v(r_1^l)}{-v(\underline{x})}\widehat{g}_{1}(1-p^l)
- \frac{v(r_3^l)}{v(\underline{x})}\widehat{g}_{2}(p^l)
+\frac{v(r_2^l)}{v(\underline{x})} \leq  0 \label{eq:scale_preference}\\
& \quad  \Leftrightarrow \! v_-(r_1^l)\widehat{g}_{1}(1-p^l)+rv_+(r_3^l)\widehat{g}_{2}(p^l) \nonumber \\
& \quad \quad  -v_-(r_2^l)\mathds{1}_{r_2^l<0}-rv_+(r_2^l)\mathds{1}_{r_2^l\geq 0}\leq 0,\label{eq:normalization_extra-Ref} 
\end{align}
\end{subequations}
where \eqref{eq:scale_preference} is obtained by dividing both sides of the inequality by $-v(\underline{x})>0$, 
and \eqref{eq:normalization_extra-Ref} follows from the definitions of $v_-(x)=-\frac{v(x)}{v(\underline{x})}$ and $v_+(x) =\frac{v(x)}{v(\bar{x})}$ in \eqref{eq:vminus_vplus}, 
together with the definition of $r=\frac{-v(\bar{x})}{v(\underline{x})}$.
In the $m$-th iteration, 
the term on the rhs of \eqref{eq:expand_preference},
multiplied by $h_l$, 
is approximated by $\phi_{\rm cut}^{l,m}({\bm v})\leq 0$ in \eqref{eq:phi_cutm}, 
expressed in terms of the increment vector ${\bm v}\in \R^{J_m-2}$ of the holistic value function $v$.
In contrast, the inequality on the rhs of \eqref{eq:normalization_extra-Ref} ensures that the preference relation can be characterized in terms of the increment vectors ${\bm v}_-\in \R^{J_-^{M_1}-1}$ of the negative value function $v_-(x)$ and ${\bm v}_+\in \R^{J_+^{M_2}-1}$ of the positive value function $v_+(x)$.

Specifically, we approximate $v_+(x)$ and $v_-(x)$ in \eqref{eq:normalization_extra-Ref} by the piecewise linear functions
$v_{J_+^{M_2}}(x)={\bm v}_+^\top{\bm h}^{{\cal X}_+^{M_2}}(x)$
and $v_{J_-^{M_1}}(x)={\bm v}_-^\top{\bm h}^{{\cal X}_-^{M_1}}(x)-1$, 
defined in \eqref{eq:PLA_utility_2} and \eqref{eq:PLA_utility_minus}, respectively, where ${\cal X}_+$ and ${\cal X}_-$ are replaced with ${\cal X}_+^{M_2}$ and ${\cal X}_-^{M_1}$.
Substituting these approximations into \eqref{eq:normalization_extra-Ref} yields
$
\phi_l^{v-v+}({\bm v}_-,{\bm v}_+):= h_l\big( ({\bm v}_-^\top{\bm h}^{{\cal X}_-^{M_1}}(r_1^l)-1)\widehat{g}_{1}(1-p^l) 
 + \widehat{r}{\bm v}_+^\top{\bm h}^{{\cal X}_+^{M_2}}(r_3^l)\widehat{g}_2(p^l)
 -({\bm v}_-^\top{\bm h}^{{\cal X}_-^{M_1}}(r_2^l)-1)\mathds{1}_{r_2^l<0} 
-\widehat{r}
{\bm v}_+^\top{\bm h}^{{\cal X}_+^{M_2}}(r_2^l)\mathds{1}_{r_2^l\geq 0}\big),
$
which is linear in the increment vectors $({\bm v}_-,{\bm v}_+)$.
The two representations are equivalent: $\phi_{\rm cut}^{l,M_4}({\bm v})\leq 0
\Leftrightarrow \phi_l^{v-v+}({\bm v}_-,{\bm v}_+)\leq 0$, where $\phi_{\rm cut}^{l,M_4}(\cdot)$ is defined in \eqref{eq:phi_cutm}. 
Consequently, to ensure consistency with the polyhedra $\Delta_-^{M_1}$ and $\Delta_+^{M_2}$ derived from Question Types~1--2, we map $\Delta_{M_4}$ defined in \eqref{eq:whole_uti} to the space of $({\bm v}_-,{\bm v}_+)$,
retain only the information from Question Type~4, and remove the constraints corresponding to ${\Delta}_0^{{\cal X}_{M_4}}$ (see \eqref{eq:update_Delta0} with $m=M_4$):
\begin{equation}
\label{eq:Delta_v-v+}
\begin{aligned}
\Delta_{M_4}^{v-v+}
:= \Big\{({\bm v}_-,{\bm v}_+)\in \R^{J_-^{M_1}-1}\times \R^{J_+^{M_2}-1}:  \phi_l^{v-v+}({\bm v}_-,{\bm v}_+) \leq 0, l\in [M_4] \Big\}.
\end{aligned} 
\end{equation}
\end{remark}
In the next subsection, we construct the ambiguity set $\mathbb{V}{\mathbb{G}_1\mathbb{G}_2}$ of the value and distortion functions.

\subsection{Construction of the Joint Ambiguity Set}
\label{sec:AS}
In this section, we can combine the inequalities obtained in Question Types $1$, $2$, and $4$ to construct the joint ambiguity set $\mathbb{V}{\mathbb{G}_1\mathbb{G}_2}$.
The polyhedra $\Delta_-^{M_1}$ of vectors $(\widetilde{\bm v}_-,\widetilde{\bm \sigma}_1)$, $\Delta_+^{M_2}$ of vectors $(\widetilde{\bm v}_+,\widetilde{\bm \sigma}_2)$, and $\Delta_{M_4}^{v-v+}$ of vectors $({\bm v}_-,{\bm v}_+)$
are established in Remark~\ref{rem:improved_PM} (b), \eqref{eq:poly_Delta+}, and \eqref{eq:Delta_v-v+}, respectively.
Then we can combine the inequalities obtained in $\Delta_-^{M_1}$, $\Delta_+^{M_2}$, and $\Delta_{M_4}^{v-v+}$ to define the set of vectors $\widetilde{\bm v}_-$, 
$\widetilde{\bm v}_+$, $\widetilde{\bm \sigma}_1$ and $\widetilde{\bm \sigma}_2$:
\begin{equation}
\label{eq:Delta1234}
\begin{aligned}
\widehat{\Delta}_{M_1+M_2+M_3+M_4}:= 
\Bigg\{({\bm R}\widetilde{\bm v}_-,{\bm R}\widetilde{\bm v}_+,{\bm R}\widetilde{\bm \sigma}_1,{\bm R}\widetilde{\bm \sigma}_2):
\begin{aligned}
&(\widetilde{\bm v}_-,\widetilde{\bm \sigma}_1)\in \Delta_-^{M_1},\left(\widetilde{\bm v}_+,\widetilde{\bm \sigma}_2\right)\in \Delta_+^{M_2}\\
& ({\bm R}\widetilde{\bm v}_-,{\bm R}\widetilde{\bm v}_+)\in \Delta_{M_4}^{v-v+}
\end{aligned}
\Bigg\}.
\end{aligned}
\end{equation}
To ease the exposition, let $M = M_1 + M_2 + M_3 + M_4$ denote the total number of queries, 
and let $\widehat{\Delta}_M^{v}$ denote the projection of $\widehat{\Delta}_M$ onto the space of $({\bm R}\widetilde{\bm v}_-, {\bm R}\widetilde{\bm v}_+)$.
Then we define the ambiguity set of value and distortion functions as
\begin{equation}
\label{eq:ambiguity_set_definition}
(\mathbb{V}{\mathbb{G}_1\mathbb{G}_2})_{M} := \left\{ (v,g_1,g_2)\in \mathscr{V}\times \mathscr{G}\times \mathscr{G}:
\begin{aligned}
&  g_{iN_i}(t)= g_{i}(t), \; t\in {\cal S}_i,i=1,2 \\
&  v(x)=v_J(x), \; x\in {\cal X}_-^{M_1}\cup {\cal X}_+^{M_2}\\
&  g_{iN_i}(t)=({\bm R}\widetilde{\bm \sigma}_i)^\top{\bm h}^{{\cal S}_i}(t), \; i=1,2\\
&  v_J(x)=({\bm R}\widetilde{\bm v}_-)^\top{\bm h}^{{\cal X}_-^{M_1}}(x)-1, x<0\\
&  v_J(x)=\widehat{r}({\bm R}\widetilde{\bm v}_+)^\top{\bm h}^{{\cal X}_+^{M_2}}(x), \; x\geq 0 \\
&  \left({\widetilde{\bm v}_-},\; \widetilde{\bm \sigma}_1\right)\in \Delta_-^{M_1},\left({\widetilde{\bm v}_+},\widetilde{\bm \sigma}_2\right)\in \Delta_+^{M_2} \\
&\;   ({\bm R}\widetilde{\bm v}_-,{\bm R}\widetilde{\bm v}_+)
\in \Delta_{M_4}^{v-v+}
\end{aligned}
\right\}, 
\end{equation}
where $\Delta_{M_4}^{v-v+}$ is defined as in \eqref{eq:Delta_v-v+} and $v_J(x)$ is a PLA of the value function $\tilde{v}$ in \eqref{eq:value_tilde}, 
which is obtained by scaling the value function $v(x)$ by parameter $-\frac{1}{v(\underline{x})}$.
In Section~\ref{sec3:ambiguity_VW}, 
we present a generic form of the ambiguity set $\mathbb{V}{\mathbb{G}_1\mathbb{G}_2}$ defined as in \eqref{eq:ambiguity set},
and establish Theorem~\ref{thm:reformulation-constraint-function-VW} accordingly.
To effectively reduce the size of the ambiguity set as the number of queries increases and to circumvent the bilinearity in the ambiguity set with respect to $v$ and $g_1,g_2$, 
we design Question Types 1--4 and obtain the polyhedra $\Delta_-^{M_1}$, $\Delta_+^{M_2}$ and $\Delta_{M_4}^{v-v+}$ and construct the specific ambiguity set $(\mathbb{V}{\mathbb{G}_1\mathbb{G}_2})_M$.
In the next subsection, we approximately reformulate the constraint function in \eqref{eq:PRGRS-CPT} based on the specific structure of $(\mathbb{V}{\mathbb{G}_1\mathbb{G}_2})_M$.

\section{Approximate Reformulation of the Constraint Function}
\label{sec:computation}
In this section, we update the reformulation of the constraint function  \eqref{eq:PRGRS-CPT}, 
based on the ambiguity set $(\mathbb{V}{\mathbb{G}_1\mathbb{G}_2})_M$ constructed in Section~\ref{sec:AS} and the finite-dimensional reformulation established in Theorem~\ref{thm:reformulation-constraint-function-VW}.

\subsection{Reformulation of Problem~\texorpdfstring{\eqref{robust_v_o}}{} Using the Constructed Ambiguity Set}
\label{sec:REFORMULATION}
After tackling the issues of bilinearity in the ambiguity set ${\mathbb V}{\mathbb G}_1{\mathbb G}_2$ and questionnaire design listed at the end of Section~\ref{Sec:tactability}, 
we are now in a position to obtain a reformulation of problem~\eqref{robust_v_o}, 
which is down to solving a problem with a bilinear objective function.

\begin{corollary} 
\label{coro:inner_trac}
Let the sets of breakpoints ${\cal X}_-^{M_1}$ and ${\cal X}_+^{M_2}$ be defined in \eqref{eq:update_M1M2}, 
where $J_-^{M_1}$ and $J_+^{M_2}$ denote the updated numbers of breakpoints.
Let ${\cal S}_1$ and ${\cal S}_2$ be the sets of breakpoints of $g_1$ and $g_2$, 
as defined in \eqref{eq:Delta_minus} and \eqref{eq:g2v+}, respectively. 
Analogously to the sets ${\cal J}(x)$ in \eqref{eq:J} and ${\cal S}_1, {\cal S}_2$ in \eqref{eq:S}, we further update the breakpoint sets as follows:
${\cal X}_-^{M_1} \leftarrow {\cal X}_-^{M_1}\cup \{\xi_k-x\}_{k=1}^{K_0(x)-1}$,
${\cal X}_+^{M_2} \leftarrow {\cal X}_+^{M_2}\cup \{\xi_k-x\}_{k=K_0(x)}^{K}$,
${\cal S}_1 \leftarrow {\cal S}_1 \cup \{q_k\}_{k=1}^{K_0(x)-1}$, 
and ${\cal S}_2 \leftarrow {\cal S}_2 \cup \{1-q_k\}_{k=K_0(x)-1}^{K}$.
The corresponding numbers of breakpoints are updated accordingly.
Let the ambiguity set $(\mathbb{V}\mathbb{G}_1\mathbb{G}_2)_M$ be defined in \eqref{eq:ambiguity_set_definition}, 
and let $\widehat{g}_1(\cdot)$ and $\widehat{g}_2(\cdot)$ be defined in \eqref{eq:hat_g1g2}.
Under Assumptions~\ref{ass:discrete}--\ref{ass:discrete_Lip}, 
the constraint function  $h_{\mathbb{V}\mathbb{G}_1\mathbb{G}_2}(\xi-x)$  in \eqref{eq:PRGRS-CPT}, 
with $\mathbb{V}\mathbb{G}_1\mathbb{G}_2$ replaced by $(\mathbb{V}\mathbb{G}_1\mathbb{G}_2)_M$, 
equals the optimal value of the following bilinear program:
\begin{subequations}
\label{eq:bilinear}
\begin{align}
\sup_{\substack{{\bm v}_-,{\bm v}_+\\
{\bm \sigma}_1,{\bm \sigma}_2}}  \Bigg\{ \! & \sum_{k=1}^{K_0(x)-1} 
{\bm v}_-^\top({\bm h}^{{\cal X}_-^{M_1}}(\xi_k-x)-{\bm e}) {\bm \sigma}_1^\top ({\bm h}^{{\cal S}_1}(q_k)-{\bm h}^{{\cal S}_1}(q_{k-1})) \nonumber \\
&\;  + \widehat{r} \sum_{k=K_0(x)}^{K} {\bm v}_+^\top{\bm h}^{{\cal X}_+^{M_2}}(\xi_k-x)   {\bm \sigma}_2^\top ({\bm h}^{{\cal S}_2}(1-q_{k-1})-{\bm h}^{{\cal S}_2}(1-q_{k})) \Bigg\} \label{eq:objective-2}\\
 {\rm s.t.} \;\;
&{\bm e}^\top{\bm \sigma}_i=1,\; {\bm \sigma}_i\geq 0,\; \forall i=1,2, \label{eq:bilinear-b}\\
&{\bm v}_+\geq 0,\;  {\bm v}_-\geq 0,\; 
{\bm e}^\top{\bm v}_-=1, \; 
{\bm e}^\top{\bm v}_+=1,\label{eq:bilinear-c}\\
&\frac{v_{-,j}}{x_{j+1}-x_j}-\frac{v_{-,j+1}}{x_{j+2}-x_{j+1}} \leq 0,\; \forall j\in[J_-^{M_1} -2],\label{eq:bilinear-d}\\
&\frac{v_{+,j+1}}{x_{j+2}-x_{j+1}}-\frac{v_{+,j}}{x_{j+1}-x_{j}} \leq 0,\; \forall j\in [J_+^{M_2}-2],\label{eq:bilinear-e}\\
& h_m ({\bm \sigma}_2^\top{\bm h}^{{\cal S}_2}(p^m)   - {\bm v}_+^\top{\bm h}^{{\cal X}_+^{M_2}}(r_{2+}^m))\leq 0,\; \forall m\in [M_2],\label{eq:bilinear-f}\\
&
h_m\big(-{\bm \sigma}_1^\top{\bm h}^{{\cal S}_1}(1-p^m)  - {\bm v}_-^\top({\bm h}^{{\cal X}_-^{M_1}}(r_{2-}^m)-{\bm e})\big)\leq 0,\; \forall m\in [M_1],\label{eq:bilinear-g}\\
& h_m(({\bm v}_-^\top{\bm h}^{{\cal X}_-^{M_1}}(r_1^m)-1)\widehat{g}_1(1-p^m)  + \widehat{r}{\bm v}_+^\top{\bm h}^{{\cal X}_+^{M_2}}(r_3^m)\widehat{g}_2(p^m)\nonumber \\
&\quad  \;\; -({\bm v}_-^\top{\bm h}^{{\cal X}_-^{M_1}}(r_2^m)-1)\mathds{1}_{r_2^m<0}  - \widehat{r}{\bm v}_+^\top{\bm h}^{{\cal X}_+^{M_2}}(r_2^m)\mathds{1}_{r_2^m\geq 0})\leq 0,\; \forall m\in [M_4],\label{eq:bilinear-h}
\end{align}
\end{subequations}
where $v_{-,j}$ and $v_{+,j}$ denote the $j$-th component of ${\bm v}_-$ and ${\bm v}_+$, respectively.
\end{corollary}

Problem \eqref{eq:bilinear} is a bilinear program with linear constraints
where the objective function consists of a sum of bilinear terms $v_{-,i}\sigma_{1,j}$ and $v_{+,i}\sigma_{2,j}$.
The problem has the same structure as that of the joint-constrained bilinear program considered by \cite{ShA92}. 
To tackle the bilinear terms, 
we adopt the standard McCormick relaxation (see \citep{McC76}) procedure, 
where each bilinear term is replaced by a new variable coupled by four linear inequality constraints.
Let $v_{-,i}^L$ and $v_{-,i}^U$ denote the lower and upper bounds of $v_{-,i}$ over the feasible region defined by \eqref{eq:bilinear-b}--\eqref{eq:bilinear-h},
and define $v_{+,i}^L$, $v_{+,i}^U$, $\sigma_{1,j}^L$, $\sigma_{1,j}^U$, $\sigma_{2,j}^L$, and $\sigma_{2,j}^U$ analogously.
To linearize the bilinear term in the objective function \eqref{eq:objective-2},
we apply piecewise McCormick relaxation to obtain an approximation of problem  \eqref{eq:bilinear} (see problem \eqref{eq:trilinear_relax_2}).
For the bilinear term $v_{+,i}\sigma_{2,j}$, we partition the range of $\sigma_{2,j}$ into $N$ disjoint subintervals and apply McCormick relaxation over each subinterval. 
Likewise, the bilinear term $v_{-,i}\sigma_{1,j}$ can be linearized.
The implementation details of the piecewise McCormick relaxation for approximating problem~\eqref{eq:bilinear} are provided in Section~\ref{sec:Piecewise_Mc}.
Applying this technique yields an MILP approximation of problem~\eqref{eq:bilinear}.

\subsection{Approximation of the Bilinear Program \texorpdfstring{\eqref{eq:bilinear}}{} via Piecewise McCormick Relaxation}
\label{sec:McCormick}
By introducing new variables $w_{ij}^+=v_{+,i}\sigma_{2,j}$ for $i\in [J_+^{M_2}-1]$, $j\in [N_2-1]$ into problem \eqref{eq:bilinear} alongside \eqref{eq:wMcC-2}, 
and introducing new variables $w_{ij}^-=v_{-,i}\sigma_{1,j}$ for $i\in [J_-^{M_1}-1]$, $j\in [N_1-1]$ into problem \eqref{eq:bilinear} along with system \eqref{eq:wMcC-2-minus}, 
we can obtain a relaxed version of problem \eqref{eq:bilinear} in the next proposition.
\begin{proposition}
Under the setting of Corollary~\ref{coro:inner_trac}, 
the bilinear problem~\eqref{eq:bilinear} can be approximated via piecewise McCormick relaxation as follows:
\begin{subequations}
\label{eq:trilinear_relax_2}
\begin{align}
 \sup_{\substack{{\bm w}^-,{\bm w}^+,{\bm v}_-,{\bm v}_+,
 {\bm \sigma}_1\\
 {\bm \sigma}_2,
 {\bm y}^{g_2},
 {\bm y}^{g_1},
 \widehat{\bm V}^+,\widehat{\bm V}^-}}
  \Bigg\{ & \sum_{k=1}^{K_0(x)-1} \sum_{i=1}^{J_-^{M_1}-1} \sum_{j=1}^{N_1-1} ({\bm h}^{{\cal X}_-^{M_1}}(\xi_k-x)-{\bm e})_i  ({\bm h}^{{\cal S}_1}(q_k)-{\bm h}^{{\cal S}_1}(q_{k-1}))_j w_{ij}^- \nonumber \\
& \; + \widehat r \sum_{k=K_0(x)}^{K} \sum_{i=1}^{J_+^{M_2}-1}\sum_{j=1}^{N_2-1} ({\bm h}^{{\cal X}_+^{M_2}}(\xi_k-x))_i  ({\bm h}^{{\cal S}_2}(1-q_{k-1})-{\bm h}^{{\cal S}_2}(1-q_{k}))_j w_{ij}^+\Bigg\} \label{eq:objective-3}\\
{\rm s.t.} \;\;\;  & \eqref{eq:bilinear-b}-\eqref{eq:bilinear-h},\\
&\eqref{eq:wMcC-2},\; \forall i\in [J_+^{M_2}-1], j\in [N_2-1],\\
&\eqref{eq:wMcC-2-minus},\; \forall i\in [J_-^{M_1}-1], j\in [N_1-1].
\end{align}
\end{subequations}
\end{proposition}

Given that problem \eqref{eq:trilinear_relax_2} can be evaluated, 
we can compute the PRGSR-CPT $\rho_{({\mathbb{V}\mathbb{G}_1\mathbb{G}_2})_{M}}(\xi)$ using the bisection method, analogous to Algorithm~\ref{algorithm:thm4-1}.
We note that additional preference robust models could be further identified based on $(\mathbb{V}\mathbb{G}_1\mathbb{G}_2)_M$. 
For example,  one may consider a preference robust multivariate decision-making problem
$\rho_{({\mathbb{V}\mathbb{G}_1\mathbb{G}_2})_{M}} (f({\bm z},{\bm \zeta}))$ with a financial loss function $f({\bm z},{\bm \zeta})$.
We refer readers to Section~\ref{SC:extention_models} for details as to how to derive these models and develop corresponding solution methods.

Now all three issues listed at the end of Section~\ref{Sec:tactability} have been tackled. 
In the next section, 
we will report the efficiency of Question~Types~1--4 for constructing the ambiguity set $(\mathbb{V}\mathbb{G}_1\mathbb{G}_2)_{M}$ of the value and distortion functions, and for solving the preference robust generalized shortfall $\rho_{(\mathbb{V}\mathbb{G}_1\mathbb{G}_2)_{M}}(\xi)$
by problem~\eqref{eq:trilinear_relax_2}.

\section{Numerical Tests}
\label{sec5:case_study}
To examine the performance of the PRGSR-CPT problem and the computational schemes, 
we conducted numerical experiments based on an academic example,
where the true value and distortion functions are known.
In this section, we report the results in detail.

We begin by examining the performance of the estimated distortion functions,
the corresponding worst-case value functions and the PRGSR-CPT measure
$\rho_{(\mathbb{V}\mathbb{G}_1\mathbb{G}_2)_{M}}(\xi)$, 
computed via problem~\eqref{eq:trilinear_relax_2} using a small-scale simulation study, and then move on to analyze the sensitivity of the worst-case value functions and the PRGSR-CPT measures.
Finally, we examine the out-of-sample performance of the CPT-based preference functional $\bbe_{{g}_1{g}_2}[v(\cdot)]$ in terms of hit rate.
Apart from these small-scale simulation studies, 
we also present the convergence of the SAA PRGSR-CPT measures.
All of the numerical experiments are carried out with MATLAB R2024b installed on a PC (16GB RAM, CPU 2.3 GHz) with Intel Core i7 processor.  
We will report the details of the test and the outcomes.

\subsection{Setup} 
\label{sec:setup}
Let $\xi$ be a random variable uniformly distributed over $\Xi=\{\xi_1,\cdots,\xi_K\}$ with $K=10$ and $P(\{\xi_k\})=1/K$ for $k\in [K]$, 
where $\xi_1,\cdots,\xi_K$ are sorted in increasing order. 
The true value function $ v^*(t)$ is defined over $[\underline{x},\bar{x}]=[-0.5,0.5]$, 
where $v^*(t)=(2(e^{8 t}-1))/8$ for $t\in [-0.5,0]$ and
$v^*(t)=(1-e^{-3 t})/3$ for $t\in [0,0.5]$.
The true distortion functions are given by $g_i^*(p)=p^{\gamma_i}/(p^{\gamma_i}+(1-p)^{\gamma_i})$, where $\gamma_i=0.6$ for $i=1,2$.
In this setup, $g_1^*(p)= g_2^*(p)$. 
We consider the special case for the simplification of computation although 
the schemes developed in Section~\ref{sec:Estimating-dis} allow $g_1^*$ and $g_2^*$ to differ.

\subsection{Estimating the Distortion Function \texorpdfstring{$g_2$}{g2}}
\label{sec:estimation_hatg}
\underline{Distortion function $g_2$}. 
Consider a set of breakpoints of $g_2$ denoted by ${\cal S}_2=\{0, 0.025, \dots,0.1\} \cup \{0.15, 0.2, \dots,0.85\} \cup \{0.9, 0.925, \dots,1\}$, 
where $0$ and $1$ are endpoints. 
Let $N_2=|{\cal S}_2|=25$ and let ${\bm \sigma}_2=(\sigma_{2,1},\cdots,{\sigma}_{2,24})$ denote the vector of increments of $g_2$.
Then we define $\widetilde{\bm \sigma}_2=(\sigma_{2,1},\cdots,\sigma_{2,23})\in \R^{23}$.

\noindent
\underline{Value function $v_+$}.
Let ${\cal X}_+^0=\{0,0.05,\ldots,0.5\}$ be a set of breakpoints of the value function $v_+$ and $J_+^0=|{\cal X}_+^0|=11$.
Let ${\bm v}_+=(v_{+,1},{v}_{+,2},\ldots,{v}_{+,10})$ denote the vector of increments of $v_+$.
Then $\widetilde{\bm v}_+=(v_{+,1},v_{+,2},\ldots,{v}_{+,9})\in \R^9$.

\begin{figure}[!ht]
\begin{minipage}[t]{0.49\linewidth}
\centering
\includegraphics[width=3.0in]{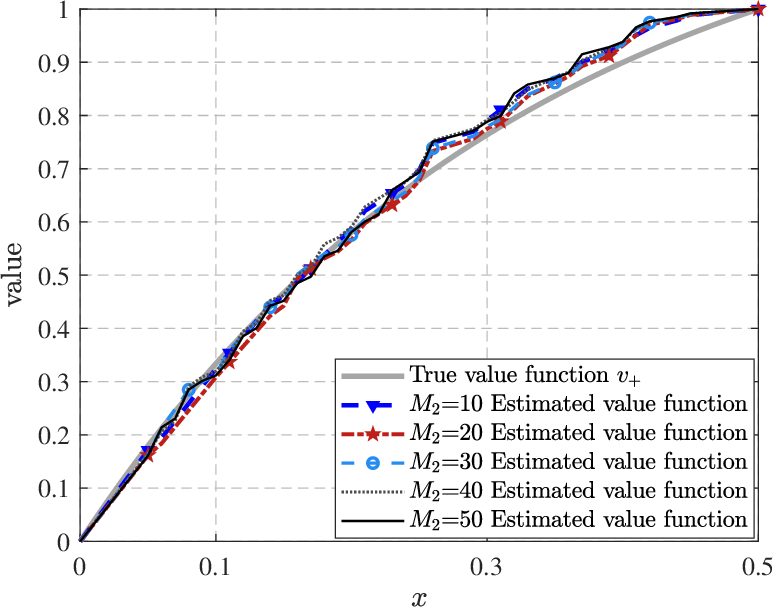}
\text{\footnotesize{\hspace{0.7cm}\textbf{(a)}}}
\end{minipage}
\begin{minipage}[t]{0.49\linewidth}
\centering
\includegraphics[width=3.0in]{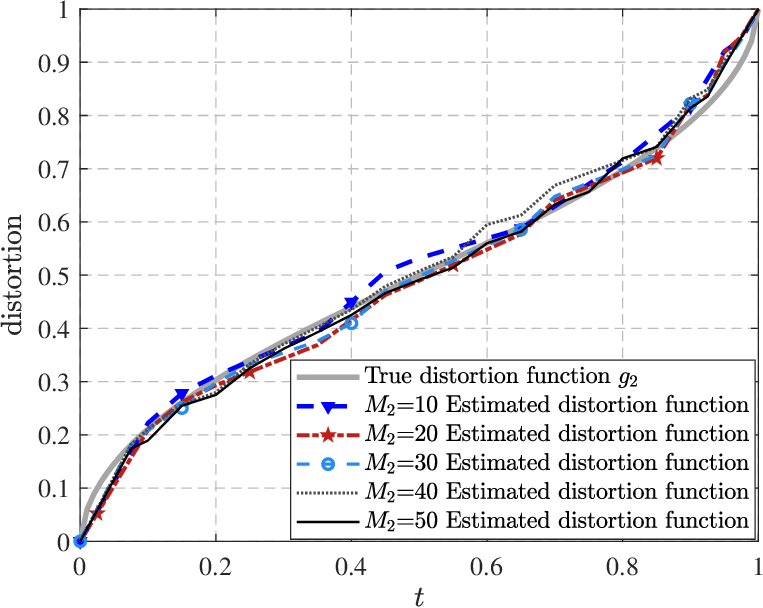}
 \text{\footnotesize{\hspace{0.7cm}\textbf{(b)}}}
\end{minipage}
\caption{\footnotesize{\textbf{(a)} The positive part of the value function $v_+$.
\textbf{(b)} 
The distortion function $g_2$.
Both are constructed using the analytic center
$({\bm v}_+^{\rm AC},{\bm \sigma}_2^{\rm AC})$.}}
\label{fig:analytic}
\end{figure}

The initial feasible set is 
$\Delta_+^0=\{\left(\widetilde{\bm v}_{+},
\widetilde{\bm \sigma}_2\right)\in \R^{9+23}:
-\widetilde{\bm \sigma}_2\leq 0, -\widetilde{\bm v}_{+}\leq 0,
{\bm e}^\top\widetilde{\bm \sigma}_2\leq 1, {\bm e}^\top\widetilde{\bm v}_{+}\leq 1, {\bm D}_+\widetilde{\bm v}_+\leq {\bm d}_+\}
$ (see \eqref{eq:defi_Delta_0}).
Following Steps 1--3 in Section~\ref{sec:generation_new}, we identify the parameters $p^1,r_{2+}^1$ of $Z_1^1,Z_2^1$ via \eqref{eq:gene_m_Q-1} and subsequently construct $\Delta_+^1$ (see \eqref{eq:choice-xkyk} with $m=1$).
Repeating the process $M_2$ times, we obtain the polyhedron $\Delta_+^{M_2}$ (see \eqref{eq:choice-xkyk} for $m=M_2$).
Let $({\bm v}_+^{\rm AC},{\bm \sigma}_2^{\rm AC})$ denote the analytic center of the polyhedron $\Delta_+^{M_2}$. 
The corresponding value and distortion functions induced by $({\bm v}_+^{\rm AC},{\bm \sigma}_2^{\rm AC})$ are given by $v_+^{\rm AC}(x)=({\bm R}{\bm v}_+^{\rm AC})^\top{\bm h}^{{\cal X}_+^{M_2}}(x)$ and $g_2^{\rm AC}(t)=({\bm R}{\bm \sigma}_2^{\rm AC})^\top{\bm h}^{{\cal S}_2}(t)$.
Figure~\ref{fig:analytic} depicts the change of the value function over $[0,0.5]$ induced by the analytic center, 
and the behavior of the distortion function $g_2$ over $[0,1]$.
To ensure feasibility of the resulting polyhedron $\Delta_+^m$, 
we impose the concavity of the value function $v_+$ over ${\cal X}_+^0$ rather than ${\cal X}_+^m$.

\subsection{Estimating the Distortion Function \texorpdfstring{$g_1$}{g1}}
\label{sec:estimation_hatg1}
\begin{figure}[!ht]
\begin{minipage}[t]{0.49\linewidth}
\centering
\includegraphics[width=3.0in]{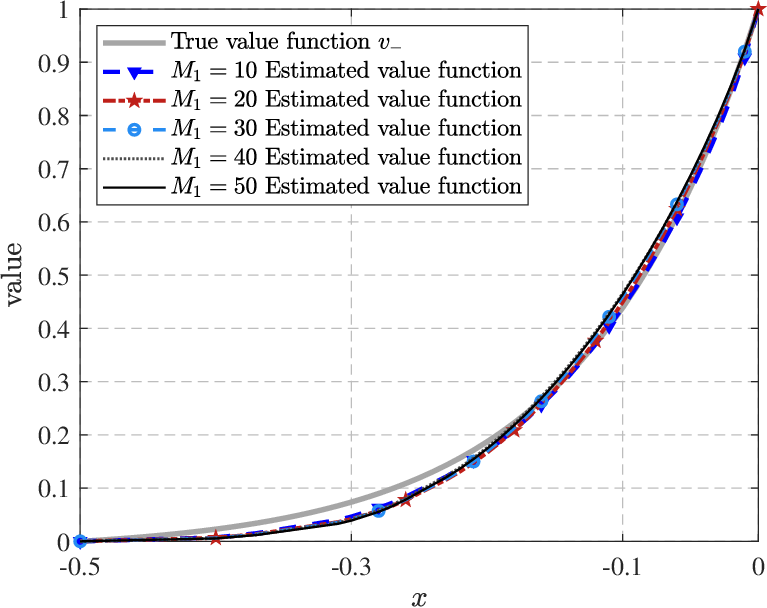}
\text{\footnotesize{\hspace{0.7cm}\textbf{(a)}}}
\end{minipage}
\begin{minipage}[t]{0.49\linewidth}
\centering
\includegraphics[width=3.0in]{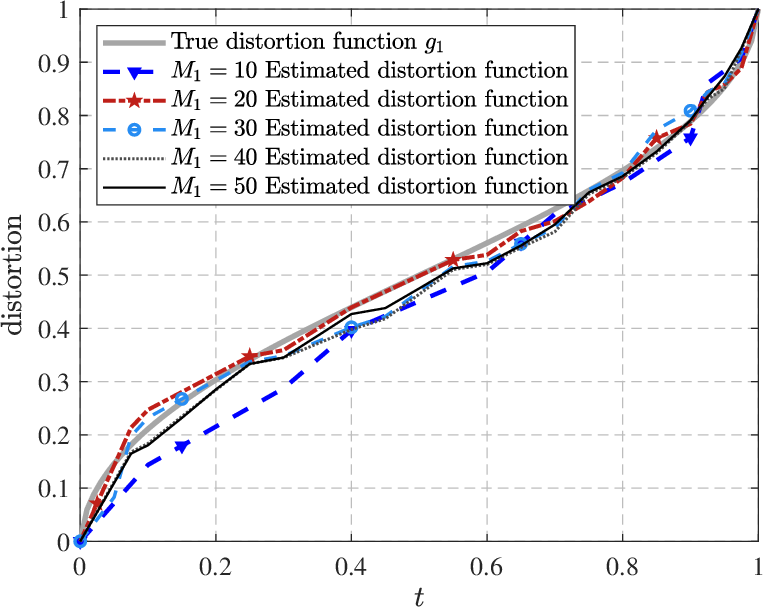}
 \text{\footnotesize{\hspace{0.7cm}\textbf{(b)}}}
\end{minipage}
\caption{\footnotesize{
\textbf{(a)} The negative part of the value function $v_-$.
\textbf{(b)}
The distortion function $g_1$.
Both are constructed using the analytic center $({\bm v}_-^{\rm AC},{\bm \sigma}_1^{\rm AC})$.}}
\label{fig:analyticminus}
\end{figure}

\underline{Distortion function $g_1$}. 
Consider a set of breakpoints of $g_1$ denoted by ${\cal S}_1={\cal S}_2$.
Let $N_1=|{\cal S}_1|=25$ and let ${\bm \sigma}_1$ denote the vector of increments of $g_1$.
Then we define $\widetilde{\bm \sigma}_1\in \R^{23}$.\\
\noindent
\underline{Value function $v_-$}.
Let ${\cal X}_-^0=\{-0.5,-0.4,-0.3,-0.2,-0.1,0\}$ be a set of breakpoints of the value function $v_-$ and $J_-^0=|{\cal X}_-^0|=6$.
Let ${\bm v}_-$ denote the vector of increments of $v_-$.
Then $\widetilde{\bm v}_-\in \R^4$.
Let $({\bm v}_-^{\rm AC},{\bm \sigma}_1^{\rm AC})$ denote the analytic center of the polyhedron $\Delta_-^{M_1}$. 
The corresponding value and distortion functions are given by
$v_-^{\rm AC}(x)=({\bm R}{\bm v}_-^{\rm AC})^\top{\bm h}^{{\cal X}_-^{M_1}}(x)-1$
and $g_1^{\rm AC}(t)=({\bm R}{\bm \sigma}_1^{\rm AC})^\top{\bm h}^{{\cal S}_1}(t)$.\\
\noindent Figure~\ref{fig:analyticminus} depicts the change of the value function over $[-0.5,0]$ induced by the analytic center,
and the behavior of the distortion function $g_1$ over $[0,1]$.

\subsection{Estimating the Value Function \texorpdfstring{$v$}{v}}
\label{sec:estimation_v}
\begin{figure}
\centering
\begin{minipage}[t]{0.49\textwidth}
\centering
\includegraphics[width=3.0in]{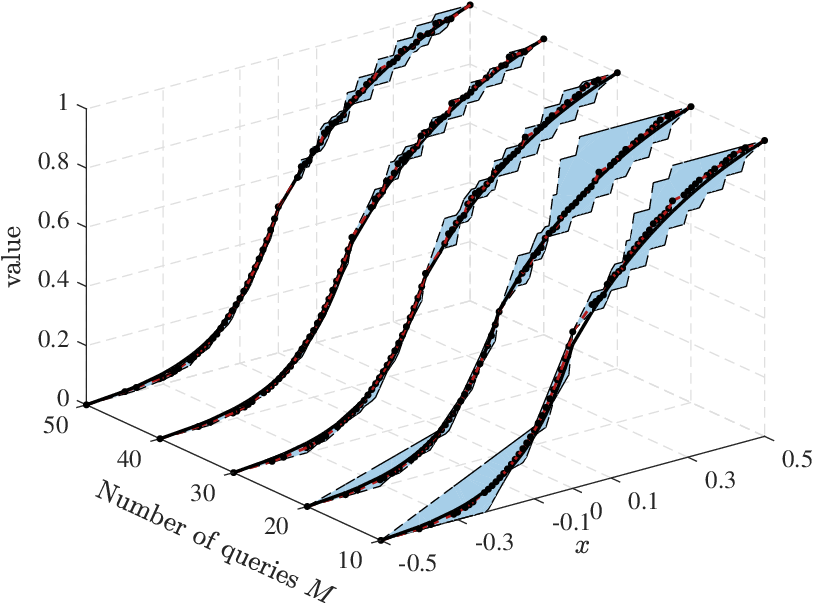}
\caption{\footnotesize Change of the Range of the Estimated Value Functions over
$[-0.5,0.5]$ as $M_4$ Increases from 10 to~50}
\label{fig:value}
\end{minipage}
\hfill
\begin{minipage}[t]{0.49\textwidth}
\centering
\includegraphics[width=3.0in]{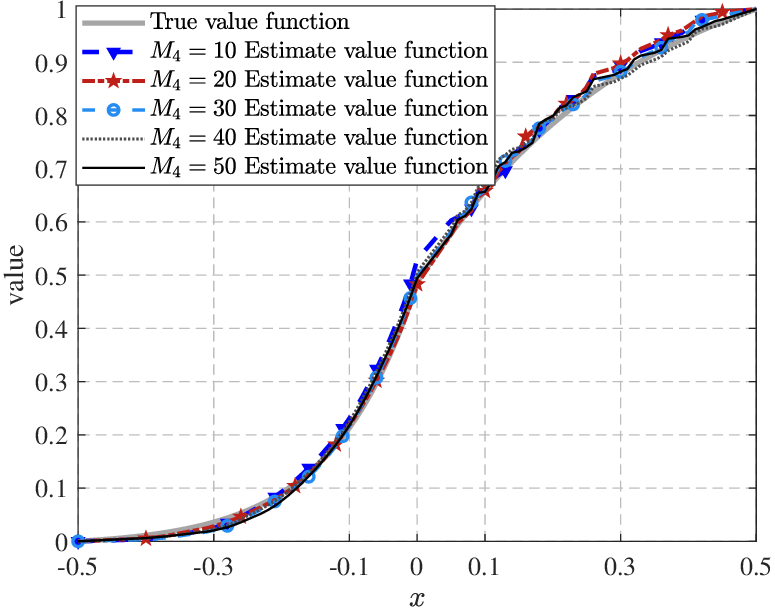}
\caption{\footnotesize Evolution of the Estimated Value Function Constructed via
the Analytic Center of
$\widehat{\Delta}_{M_1+M_2+M_3+M_4}$}
\label{fig:CPT_PLA-4}
\end{minipage}
\end{figure}

Let ${\cal X}_0 = {\cal X}_-^{M_1} \cup {\cal X}_+^{M_2}$.
Figure~\ref{fig:value} depicts the range of the estimated ambiguity set of piecewise linear value functions corresponding to $\widehat{\Delta}_{M_1+M_2+M_3+M_4}$, 
where the blue area is the pointwise range of these functions.
Figure~\ref{fig:CPT_PLA-4} depicts the change of the value function identified by the analytic center of $\widehat{\Delta}_{M_1+M_2+M_3+M_4}$ defined in \eqref{eq:Delta1234} and the estimate of the ratio $\widehat{r}$.
For Question Type 4, we generate $M_4$ queries and set $M_4 = M_2 = M_1$ as $M_1$ increases from $10$ to $50$.
The estimated distortion functions $\widehat{g}_1$ and $\widehat{g}_2$,
defined in~\eqref{eq:hat_g1g2}, are constructed after $M_1$ and $M_2$ queries of Types 1 and 2, respectively. 

\subsection{In-Sample and Out-of-Sample Performance and SAA Convergence of PRGSR-CPT Measure}
\label{sec:sensitivity}
\begin{figure}[!ht]
\centering
\begin{minipage}[t]{0.32\linewidth}
\centering
\includegraphics[width=2.1in]{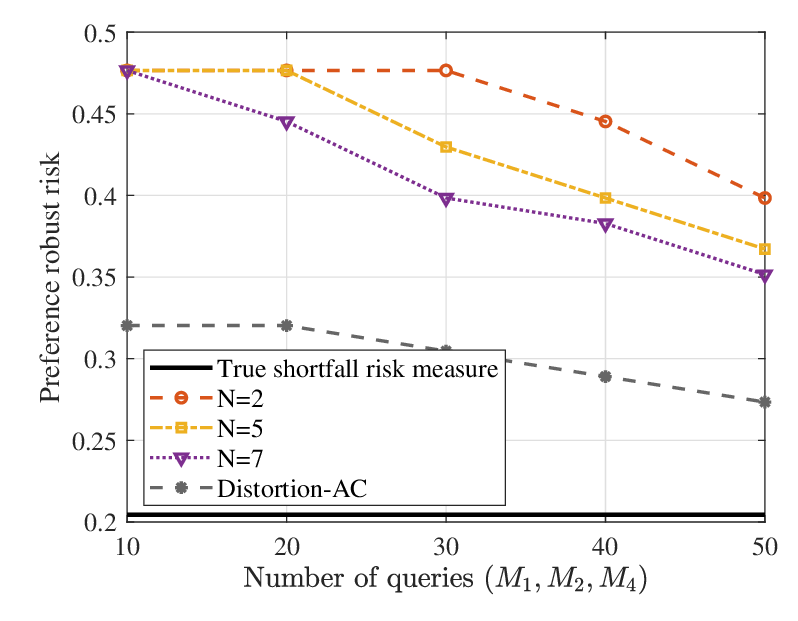} 
\par{\hspace{0.7cm} \footnotesize\textbf{(A)}} 
\end{minipage}
\hspace{0.01\linewidth}
\begin{minipage}[t]{0.64\linewidth}
\centering
\begin{minipage}[t]{0.48\linewidth}
\centering
\includegraphics[width=2.1in]{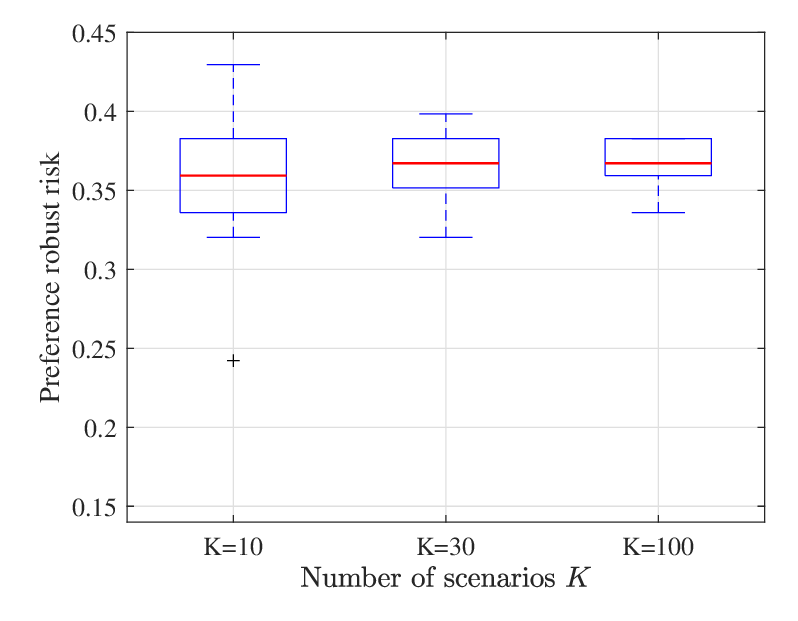}
\par{\hspace{0.7cm}  \footnotesize\textbf{(B1)}}
\end{minipage}
 \hfill
\begin{minipage}[t]{0.48\linewidth}
\centering
\includegraphics[width=2.1in]{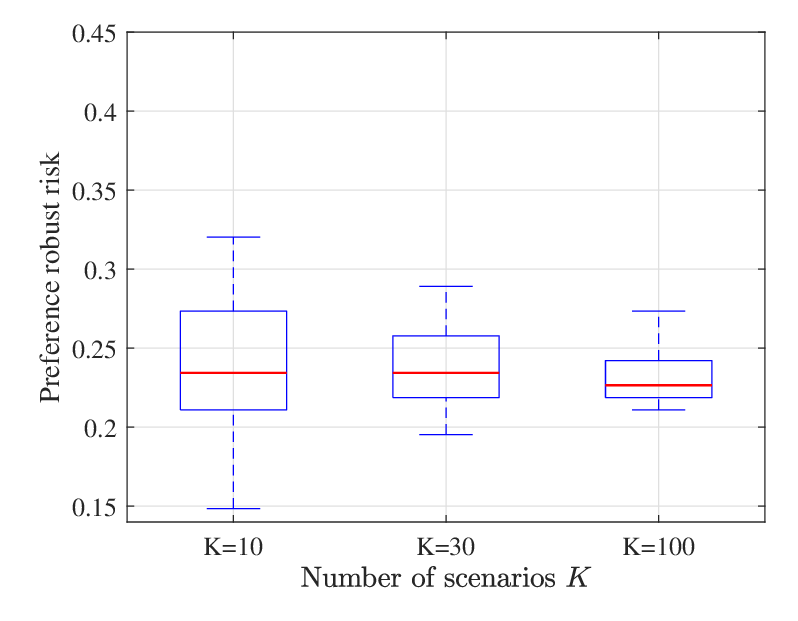}
 \par{\hspace{0.7cm} \footnotesize\textbf{(B2)}}
\end{minipage}
\end{minipage}
\caption{\footnotesize \textbf{(A)} Approximate PRGSR-CPT values obtained by solving the inner problem
\eqref{eq:trilinear_relax_2} with $K=10$, where the number of McCormick segments varies
from $N=2$ to $N=7$.
\textbf{(B)} Boxplots of approximate PRGSR-CPT values using SAA for different inner problem
formulations and sampling sizes $K=10,30,100$; (B1) inner problem \eqref{eq:bilinear} with fixed
${\bm\sigma}_i={\bm\sigma}_i^{\rm AC}$; (B2) inner problem \eqref{eq:trilinear_relax_2} with the number of segments of piecewise McCormick relaxation $N=2$.}
 \label{fig:SAA}
 \end{figure}

\noindent{\bf In-Sample Performance and SAA Convergence under the PRGSR-CPT Measure} In all of the preceding numerical studies, we used a sample size of $K=10$. 
Figure~\ref{fig:SAA} (A) depicts the convergence of the PRGSR-CPT measure computed by Algorithm~\ref{algorithm:thm4-1}, where the inner problem 
$\sup_{(v,g_1,g_2)\in (\mathbb{V}\mathbb{G}_1\mathbb{G}_2)_M} \bbe_{g_1g_2}[v(\xi-x)]$ is approximated via the MILP formulation~\eqref{eq:trilinear_relax_2}. 
Specifically, we illustrate the convergence behavior of the PRGSR-CPT measure as the number of queries $M_4$ (with $M_1 = M_2 = M_4$) increases from $10$ to $50$, and the discretization level $N$ in PMR varies from $2$ to $7$. 
The gray dotted line represents the PRGSR-CPT measure obtained by solving the inner problem using the bilinear objective function in~\eqref{eq:bilinear}, 
with ${\bm \sigma}_1$ and ${\bm \sigma}_2$ fixed at ${\bm \sigma}_1^{\rm AC}$ and ${\bm \sigma}_2^{\rm AC}$, respectively.
To illustrate the SAA convergence, we investigate the performance of the robust shortfall with inner problem \eqref{eq:trilinear_relax_2} by increasing $K$ to $30$ and $100$.  
We continue to assume that $\xi$ is uniformly distributed over $[0,0.5]$.
For each fixed sample size, we generate $20$ sets of iid samples  and then solve 
$\rho_{(\mathbb{V}\mathbb{G}_1\mathbb{G}_2)_{M}}(\xi)$ for each set of samples. 
Figure~\ref{fig:SAA} (B1)-(B2) depicts the change of boxplot of the $20$ PRGSR-CPT values as $K$ increases from $10$ to $30$ and $100$.
We observe that the PRGSR-CPT measure appears to converge. 
To ensure feasibility of the MILP~\eqref{eq:trilinear_relax_2},
we impose shape constraints on the value functions over fixed breakpoint sets, 
namely concavity of $v_+$ over ${\cal X}_+^0$ and convexity of $v_-$ over negative breakpoints in ${\cal X}_{M_4}$, instead of adopting the updated breakpoints ${\cal X}_-^{M_1}$ obtained by adding $\{\xi_k-x\}$ as in Corollary~\ref{coro:inner_trac}.

\noindent{\bf Out-of-Sample Performance.}
We investigate the quality of the worst-case value and distortion functions $v_J^{\rm worst}$, $g_{1N_1}^{\rm worst}$, 
and $g_{2N_2}^{\rm worst}$ of problem \eqref{eq:trilinear_relax_2} by examining the out-of-sample performance of the resulting preference functional $\bbe_{{g}_{1N_1}^{\rm worst}{g}_{2N_2}^{\rm worst}}[v_J^{\rm worst}(\cdot)]$.
We do so by comparing the consistency of inequalities:
$
E^L(m):=  \bbe_{{g}_{1N_1}^{\rm worst} {g}_{2N_2}^{\rm worst}}[v_J^{\rm worst}(Z_{1,m})]
 -\bbe_{{g}_{1N_1}^{\rm worst} {g}_{2N_2}^{\rm worst}}[v_J^{\rm worst}(Z_{2,m})] \leq (\geq) \;0$,
 and
$E^*(m):= \bbe_{{g}_1^*{g}^*_2}[v^*(Z_{1,m})]
 -\bbe_{{g}^*_1 {g}^*_2}[v^*(Z_{2,m})] \leq (\geq)\; 0$, 
where $Z_{1,m}$ takes the values $r_1^m$ and $r_3^m$
with probabilities $1-p^m$ and $p^m$, respectively,
and $Z_{2,m}=r_2^m$,
with $r_1^m,r_2^m,r_3^m$ being uniformly generated from $[-0.5,0.5]$
and $p^m$ being uniformly generated from $(0,1)$ for $m=1,\cdots,100$.
When the directions of the inequalities ${E}^L(m)$ and ${E}^*(m)$
are consistent, we count it as a ``hit''.
By convention (see \citep{McC02}), 
we define the {\em hit rate} as the ratio of the total number of hits to the sample size. 
We select $100$ out-of-sample queries, 
and denote the number of hits by $n$.
Then the hit rate is ${n}/{100}$.
    
\begin{figure}[!ht]
\centering
\begin{minipage}[t]{0.32\linewidth}
\centering
 \includegraphics[width=2.2in]{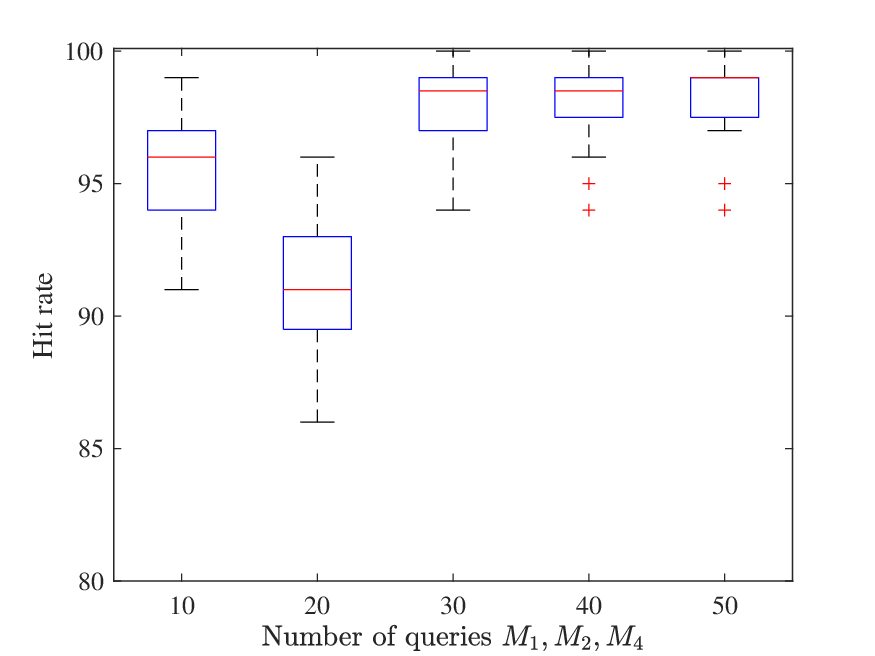}
 \text{\footnotesize \hspace{0.7cm} {\textbf{(a)}}}
\end{minipage}
\begin{minipage}[t]{0.32\linewidth}
\centering
\includegraphics[width=2.2in]{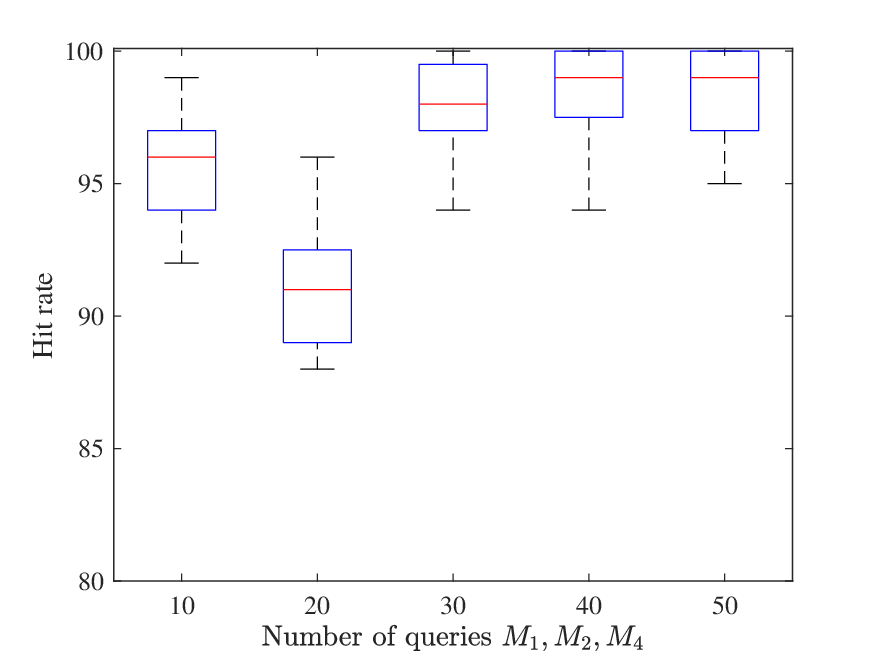}
\text{\footnotesize \hspace{0.7cm} {\textbf{(b)}}}
\end{minipage}
\begin{minipage}[t]{0.32\linewidth}
\centering
\includegraphics[width=2.2in]{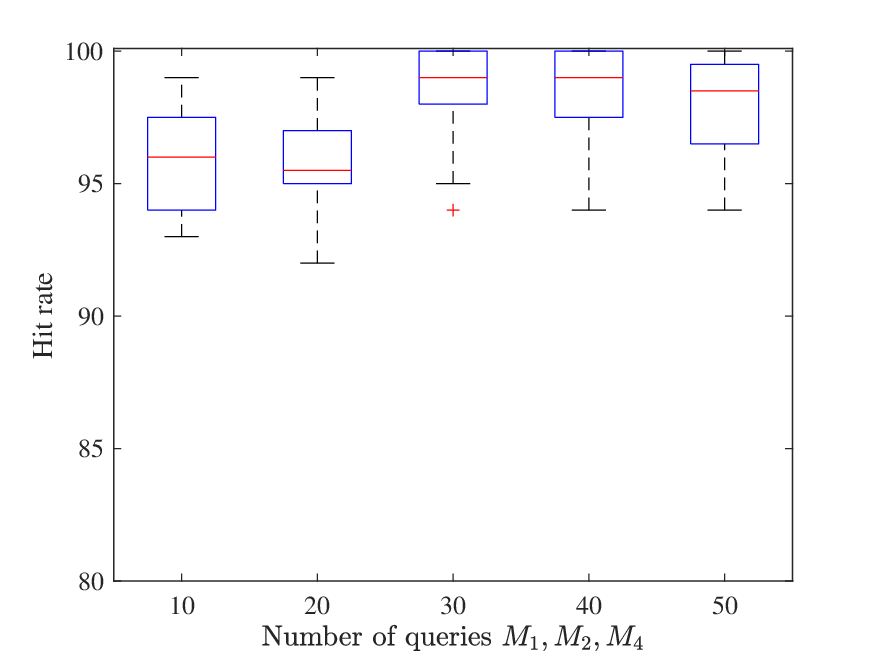}
 \text{\footnotesize{ \hspace{0.7cm} \textbf{(c)}}}
\end{minipage}
\caption{\footnotesize
\textbf{(a)} Boxplots of out-of-sample hit rates when ${\bm \sigma}_i = {\bm \sigma}_i^{\rm AC}$, $i=1,2$. 
\textbf{(b)} Boxplots of out-of-sample hit rates based on the piecewise McCormick relaxation (PMR) approach with $N=2$.
\textbf{(c)} Boxplots of out-of-sample hit rates based on the PMR approach with $N=5$.
The hit rates ($n/100$) are reported as $M_4$ varies, with $M_1 = M_2 = M_4$, illustrating the impact of the number of segments in the PMR.}
\label{fig:hit_rate}
\end{figure}

We have run $20$ simulations for each fixed number of queries $M_4$. 
Figure~\ref{fig:hit_rate} (a) presents the boxplots of the resulting hit rates 
when ${\bm \sigma}_1$ and ${\bm \sigma}_2$ are fixed at ${\bm \sigma}_1^{\rm AC}$ and ${\bm \sigma}_2^{\rm AC}$, respectively, 
in the bilinear objective of problem~\eqref{eq:bilinear}.
We can see that as $M_4$ increases to $50$, the hit rate reaches $98\%$,
which means the out-of-sample performance of the estimated preference function is very good.
Figure~\ref{fig:hit_rate} (b) depicts the hit rate when the number of pieces is $N=2$ in the piecewise McCormick relaxation (PMR) with inner problem \eqref{eq:trilinear_relax_2}, 
whereas Figure~\ref{fig:hit_rate} (c) depicts the case when $N=5$. 
We can see that in both cases, the hit rate is not affected much.

\section{Concluding Remarks}
\label{sec:conclusion}
In this paper, we consider a robust CPT-based shortfall to simultaneously address ambiguity in both the distortion function and the value function in CPT
and extend the well-known polyhedral method for eliciting linear utility functions to piecewise linear value and distortion functions. 
The work effectively advances the current research on preference robust optimization (PRO) from expected utility theory-based decision-making models (both  the von Neumann and Morgenstern EUT and Yaari's dual expected utility) to CPT-based models, 
and paves the way for further exploration of adaptive elicitation methods for nonlinear and nonconvex preference functions.
There are a few points that we would like to highlight including 
(i) the limitation of the developed models and preference elicitation based on CPT;
(ii) application to one-stage CPT-based PRO models.
We refer readers to Appendix~\ref{SC:conclu-rem} for details.

\newpage 
\appendix 
\section*{Appendix}

\section{Properties of the Preference Functional \texorpdfstring{$\bbe_{{g_1}{g_2}}[v(\cdot)]$}{}}
Before proving Proposition~\ref{prop:PGSR-robust-equal} in Section~\ref{sec:proof_Prop2}, 
we establish two preliminary lemmas analyzing the properties of the preference functional $\bbe_{g_1g_2}[v(\cdot)]$ defined as in \eqref{eq:defi-bbe-wwv}.

\begin{lemma}
\label{lemma:monotone-of-bbe-xieta}
Let $\bbe_{{g_1}{g_2}}[v(\cdot)]$ be defined as in \eqref{eq:defi-bbe-wwv}. 
Then
\begin{equation*}
\bbe_{{g_1}{g_2}}[v(\eta)]\leq \bbe_{{g_1}{g_2}}[v(\xi)],\; \forall \eta,\xi \in {\cal L}^{\infty}(\R) \inmat{ with } \eta(\omega)\leq \xi(\omega),\; \forall \omega \in \Omega.
\end{equation*}
\end{lemma}
\proof{Proof.}
By definition (see \eqref{eq:defi-bbe-wwv}),
\begin{equation*}
\bbe_{{g_1}{g_2}}[v(\xi)] =\int_{0}^{\infty}{g_2}(\mathbb P(v(\xi)>t))dt +\int_{0}^{\infty} -{g_1}(\mathbb P(-v(\xi){\geq } t))dt.
\end{equation*}
Since $\eta(\omega)\leq \xi(\omega)$ for all $\omega\in \Omega$ and $v(\cdot)$ is monotonically increasing, 
$v(\eta(\omega))\leq v(\xi(\omega))$ for all $\omega\in \Omega$, and consequently 
\begin{equation*}
\mathbb P(v(\eta)>t)\leq \mathbb P(v(\xi)>t) \quad \inmat{and} \quad \mathbb P(-v(\eta) \geq t)\geq \mathbb P(-v(\xi) \geq t).
\end{equation*}
Moreover, since ${g_1},{g_2}$ are increasing, 
we have 
\begin{equation*}
{g_2}(\mathbb P(v(\eta)>t))\leq {g_2}(\mathbb P(v(\xi)>t))\quad \text{and}\quad -{g_1}(\mathbb P(-v(\eta)\geq t))\leq -{g_1}(\mathbb P(-v(\xi) \geq t)).
\end{equation*}
A combination of the inequalities gives rise to $\bbe_{{g_1}{g_2}}[v(\eta)]\leq \bbe_{{g_1}{g_2}}[v(\xi)]$.
The proof is complete.

\hfill
\Box

\begin{lemma}
\label{lemma:monotone-of-bbewv}
Let $\bbe_{{g_1}{g_2}}[v(\cdot)]$ be defined as in \eqref{eq:defi-bbe-wwv} and $\xi\in {\cal L}^{\infty}(\R)$.
Then the following assertions hold.
\begin{itemize}
\item[(i)] 
For any $a, b\in \R$ with $a\leq b$,
\begin{equation*}
\bbe_{{g_1}{g_2}}[v(\xi-a)]\geq \bbe_{{g_1}{g_2}}[v(\xi-b)].
\end{equation*}
\item[(ii)]
For any $a, b\in \R$ with $a< b$,
\begin{equation*}
\bbe_{{g_1}{g_2}}[v(\xi-a)]>
\bbe_{{g_1}{g_2}}[v(\xi-b)].
\end{equation*}
\end{itemize}
\end{lemma}

\proof{Proof.}
Part (i) follows directly from Lemma~\ref{lemma:monotone-of-bbe-xieta}.

Part~(ii). Let $\Omega_1:=\{\omega\in \Omega: \xi(\omega)< b\}$.
For a real number $a\in \R$, let $(a)_+:=\max\{a,0\}$ and $(a)_-:=\min\{a,0\}$.
For any $a< b$, we have 
\begin{equation*}
\xi(\omega)-b< [\xi(\omega)-a]_-, \; \forall \omega\in \Omega_1 
\quad \inmat{and}\quad
\xi(\omega)-b< [\xi(\omega)-a]_+, \; \forall \omega\in \Omega \backslash{\Omega_1}.
\end{equation*}
Since the value function $v\in \mathscr V$ is strictly increasing, 
\begin{equation*}
v^-(\xi(\omega)-b)< v^-([\xi(\omega)-a]_-), \;  \forall \omega\in \Omega_1 \quad \inmat{and}\quad
v^+(\xi(\omega)-b)< v^+([\xi(\omega)-a]_+), \; \forall \omega\in \Omega\backslash{\Omega_1}.
\end{equation*}
By the definition of $\bbe_{{g_1}{g_2}}[v(\cdot)]$, we have 
\begin{equation*}
\bbe_{g_1 g_2}[v(\xi-b)]
= \int_{-\infty}^b v^-(u-b)\, d\big(g_1(F_\xi(u))\big) 
+ \int_b^{\infty} v^+(u-b)\, d\big(1-g_2(1-F_\xi(u))\big).
\end{equation*}
Since we have either $\mathbb{P}(\Omega_1)>0$ or $\mathbb{P}(\Omega\backslash \Omega_1)>0$ or both, 
\begin{equation*}
\begin{aligned}
\bbe_{{g_1}{g_2}}[v(\xi-b)] < & 
\int^{b}_{-\infty} v^-([u-a]_-) d{(g_1(F_{\xi}(u)))} 
+\int^{\infty}_{b} v^+([u-a]_+) d{(1-g_2(1-F_{\xi}(u)))}  \\
\leq &\int^{a}_{-\infty}
v^-(u-a) d{(g_1(F_{\xi}(u)))} 
+\int^{\infty}_{a} v^+(u-a) d{(1-g_2(1-F_{\xi}(u)))}\\
= &\bbe_{{g_1}{g_2}}[v(\xi-a)].
\end{aligned}
\end{equation*}
The proof is complete.
\hfill \Box

\section{Proof of Proposition~\ref{prop:PGSR-robust-equal}}
\label{sec:proof_Prop2}
\proof{Proof of Proposition~\ref{prop:PGSR-robust-equal}.}
To simplify the exposition, we call $x$ admissible for \eqref{eq:rho-CPT} if $\bbe_{g_1g_2}[v(\xi-x)]\leq 0$, and admissible for \eqref{eq:PRGRS-CPT} if $\sup_{(v,g_1,g_2)\in\mathbb{V}\mathbb{G}_1\mathbb{G}_2} \bbe_{g_1g_2}[v(\xi-x)]\le0$.
We begin by  showing that the equality \eqref{eq:PRGSR-equivalent} 
holds when one of the quantities is $+\infty$.
First, consider the case that $\rho_{\mathbb{V}\mathbb{G}_1\mathbb{G}_2}
(\xi) =+\infty$.
Then the admissible set corresponding to \eqref{eq:PRGRS-CPT} is empty.
Assume for the sake of a contradiction that $\hat{x}:=\Pi_{\mathbb{V}\mathbb{G}_1\mathbb{G}_2}(\xi)<+\infty$.
Then $\rho_{v,g_1,{g_2}}(\xi)\leq \hat x<+\infty$ for all $(v,g_1,{g_2})\in {\mathbb{V}\mathbb{G}_1\mathbb{G}_2}$.
Thus, for any positive number $\epsilon$, it follows by Lemma~\ref{lemma:monotone-of-bbewv}~(i) that 
\begin{equation*}
\bbe_{g_1{g_2}}[v(\xi-\hat x-\epsilon)] \leq \bbe_{g_1{g_2}}[v(\xi-\rho_{v,g_1,{g_2}}(\xi)-\epsilon)] \leq 0, \; \forall (v,g_1,{g_2})\in {\mathbb{V}\mathbb{G}_1\mathbb{G}_2}.
\end{equation*}
The latter implies $\hat x+\epsilon$ is admissible for \eqref{eq:PRGRS-CPT}, a contradiction.

Conversely, consider the case that $\Pi_{\mathbb{V}\mathbb{G}_1\mathbb{G}_2}(\xi)=+\infty$.
Assume for a contradiction that $x^*:=\rho_{\mathbb{V}\mathbb{G}_1\mathbb{G}_2}(\xi)<+\infty$.
Then for any positive number $\epsilon$, $x^*+\epsilon$ is admissible for \eqref{eq:PRGRS-CPT}, that is,
\begin{equation*}
\sup_{(v,g_1,{g_2})\in \mathbb{V}\mathbb{G}_1\mathbb{G}_2} \bbe_{g_1{g_2}}[v(\xi-x^*-\epsilon)]\leq 0,
\end{equation*}
which implies
\begin{equation*}
\bbe_{g_1{g_2}}[v(\xi-x^*-\epsilon)]\leq 0,\;\forall (v,g_1,{g_2})\in \mathbb{V}\mathbb{G}_1\mathbb{G}_2.
\end{equation*}
The latter guarantees that $\rho_{v,g_1,{g_2}}(\xi)\leq x^*+\epsilon, \forall (v,g_1,{g_2})\in \mathbb{V}\mathbb{G}_1\mathbb{G}_2$ and hence $\Pi_{\mathbb{V}\mathbb{G}_1\mathbb{G}_2}(\xi)\leq x^*+\epsilon<+\infty$, a contradiction.

We now turn to discuss the case that the risk measures for $\xi$ in \eqref{eq:PRGRS-CPT} and \eqref{eq:PRGRS-CPT-equiv} are both finite.
Let $x$ be admissible for \eqref{eq:PRGRS-CPT}.
Then $x$ is also admissible for \eqref{eq:rho-CPT} for any $(v,g_1, {g_2})\in \mathbb{V}\mathbb{G}_1\mathbb{G}_2$, which implies $x^*\geq \hat{x}$.

Conversely, for any fixed $(v,g_1, {g_2})\in \mathbb{V}\mathbb{G}_1\mathbb{G}_2$.
Then by definition $\hat{x} \geq  \rho_{v,g_1, {g_2}}(\xi)$.
It follows by Lemma~\ref{lemma:monotone-of-bbewv}~(i) that
\begin{equation*}
\mathbb{E}_{g_1{g_2}}[v(\xi-\hat{x})] \leq \mathbb{E}_{g_1{g_2}}[v(\xi-\rho_{v,g_1, {g_2}}(\xi))] \leq 0.
\end{equation*}
This shows $\hat{x}$ is admissible for \eqref{eq:PRGRS-CPT}.
This shows $\hat{x}\geq x^*$. 
\hfill \Box

\section{An Example of PRGSR-CPT with a Finite Ambiguity Set and Its Advantages}
\label{exa:finite}
\begin{example}[Worst CPT-Based Shortfall vs. PRGSR-CPT]
\label{exa:Calculate_worst}
Assume that the DM's preference can be represented by the CPT-based shortfall $\rho_{v,g_1,{g_2}}(\cdot)$.
For two prospects $A=(p_1,x_1;\ldots;p_n,x_n)$ and $B=(q_1,y_1;\ldots;q_n,y_n)$, we have $\rho_{v,g_1,{g_2}}(A)< \rho_{v,g_1,{g_2}}(B)$.

Consider a simple case that the DM's preference may be represented by three triples of value and distortion functions $\{(v^i,{g}_1^i,{g}_2^i): i=1,2,3\}$.
Suppose that $(v^1,{g}_1^1,{g}_2^1)$ and $(v^2,{g}_1^2,{g}_2^2)$ satisfy $\rho_{v^j,g_1^j,{g_2^j}}(A)< \rho_{v^j,g_1^j,{g_2^j}}(B)$, for $j=1,2$, but $(v^3,{g}_1^3,{g}_2^3)$ does not.
Then we define the ambiguity set as $\mathbb{ V}\mathbb{G}_1\mathbb{G}_2
=\{ (v^i,{g}_1^i,{g}_2^i): i=1,2\}$.
Let $\xi$ be a random variable representing financial loss.
Then the worst (largest) CPT-based shortfall of $\xi$ is defined as
\begin{equation}
\label{eq:PRGRS-CPT-equiv-ex}
\Pi_{\mathbb{ V}\mathbb{G}_1\mathbb{G}_2}(\xi)
:= \max_{i=1,2} \rho_{v^i,{g}_1^i,{g}_2^i}(\xi)
\end{equation}
and  the PRGSR-CPT is defined as $\rho_{\mathbb{ V}\mathbb{G}_1\mathbb{G}_2}(\xi):= \inf\limits_{x\in \R}  \left\{x: \mathop{\max}\limits_{i=1,2} \mathbb E_{{g}_1^i{g}_2^i}[v^i(\xi-x)]\leq 0 \right\}.$ 
\end{example}

\begin{remark}[Advantages of PRGSR-CPT]
\label{rem:advantage}
\begin{itemize}
\item[(i)]
In practice, the ambiguity set $\mathbb{ V}\mathbb{G}_1\mathbb{G}_2$ is not a finite set, and consequently, $\Pi_{\mathbb{ V}\mathbb{G}_1\mathbb{G}_2}(\xi)$ is computationally challenging to evaluate.
To alleviate this difficulty, 
we use the equivalence between $\Pi_{\mathbb{ V}\mathbb{G}_1\mathbb{G}_2}(\xi)$ and $\rho_{\mathbb{ V}\mathbb{G}_1\mathbb{G}_2}(\xi)$ in Proposition~\ref{prop:PGSR-robust-equal} to reformulate the constraint function of PRGSR-CPT into a bilinear program when the ambiguity sets admit a linear structure (see Section~\ref{sec:tractability}).

\item[(ii)]
We confine our study to the case that $[\underline{x},\bar{x}]$ satisfies $\underline{x}<0<\bar{x}$, 
since we choose $v(0)$ as the maximal tolerable distorted expectation of $\xi-x$ in \eqref{eq:rho-CPT}.
Suppose that the DM has a reference point or benchmark in wealth denoted by $B\neq 0$.
In this case, the distorted expectation of the value function is   $\bbe_{g_1g_2}^B[v(\xi)]:=\int_{B}^{+\infty} v(x-B) d(-g_2(1-F_{\xi}(x)))
+\int_{-\infty}^B v(x-B) d(g_1(F_{\xi}(x)))$, 
where $v\in \mathscr{V}$, $g_i\in \mathscr{G}$, for $i=1,2$.
If we are able to estimate a bound  $c>0$ for $|\xi-x|$, then we can set  $[\underline{x},\bar{x}]=[-B-c,-B+c]$. 
Subsequently, we may replace $v(0)$ in the right-hand side of \eqref{eq:rho-CPT} 
with $v(-B)\neq 0$.
\end{itemize} 
\end{remark}

\section{Proof of Proposition~\ref{Prop:Convex}}
\proof{Proof of Proposition~\ref{Prop:Convex}.}
By Proposition~\ref{prop:PGSR-robust-equal},
\begin{equation*}
\rho_{\mathbb{V}\mathbb{G}_1\mathbb{G}_2^{\rm con}}(\xi) = \sup_{(v,g_1,{g_2})\in \mathbb{V}\mathbb{G}_1\mathbb{G}_2^{\rm con}} \rho_{v,g_1,{g_2}}(\xi).
\end{equation*}
By \citet[Corollary 2.8]{MaC18}, $\rho_{v,g_1,{g_2}}(\cdot)$ is a convex risk measure if and only if $g_1(p)$ and $1-{g_2}(1-p)$ are convex in $p\in [0,1]$, and $v$ is convex on both $\R_-$ and $\R_+$, and $(v,g_1,{g_2})$ satisfies the inequality in \eqref{eq:constraint-vww}.
Hence, ${\cal R}_{\rm con}=\{\rho_{v,g_1,{g_2}} :(v,g_1,g_2)\in\mathbb{V}\mathbb{G}_1\mathbb{G}_2^{\rm con}\}$.
Taking the supremum over identical index sets yields
\begin{equation*}
\rho_{\mathbb{V}\mathbb{G}_1\mathbb{G}_2^{\rm con}}(\xi) = \rho_{{\cal R}_{\rm con}}(\xi).
\end{equation*}
The proof is complete.
\hfill \Box

\section{Proof of Proposition~\ref{prop:robust-RM-Pi}}
\proof{Proof of Proposition~\ref{prop:robust-RM-Pi}.}
By \citet[Theorem~2.10]{MaC18}, $\rho_{v,g_1,{g_2}}(\cdot)$ is coherent if and only if $v\in {\cal V}_{\rm coh}$ and $g_1(p)$ and $1-{g_2}(1-p)$ are convex,
and $\Upsilon(p):=a_1(1-g_2(1-p))-a_2g_1(p)$
is nondecreasing on $[0,1]$.
Hence, ${\cal R}_{\rm coh}=\{\rho_{v,g_1,{g_2}}:(v,g_1,g_2)\in{\mathbb{V}\mathbb{G}_1\mathbb{G}_2^{\rm coh}}\}$.
Moreover, it follows by Proposition~\ref{prop:PGSR-robust-equal} that $\rho_{{\mathbb{V}\mathbb{G}_1\mathbb{G}_2^{\rm coh}}}(\xi)$ can be rewritten as 
$\sup_{(v,g_1,{g_2})\in {\mathbb{V}\mathbb{G}_1\mathbb{G}_2^{\rm coh}}}
\rho_{v,g_1,{g_2}}(\xi)$.
Taking the supremum over identical index sets yields
\begin{equation*}
\Pi_{\mathbb{V}\mathbb{G}_1\mathbb{G}_2^{\rm coh}}
(\xi)=\rho_{{\cal R}_{\rm coh}}(\xi).
\end{equation*}
The proof is complete. \hfill \Box

\section{Explanation of Strong Risk Taking in Remark~\ref{rem:risk-taking}}
\label{SC:risk-seeking}
We will show that (i) the relation in \eqref{eq:ineq5} is equivalent to $\frac{v_+'(0)}{v_-'(0)}\geq \sup_{s\in (0,1)} \frac{g'_{1,-}(s)}{g'_{2,-}(1-s)}$;
(ii) the positive homogeneity of CPT implies that the value function $v$ is positively homogeneous.
\begin{itemize}
\item[(i)] 
Note that the CPT value for the prospect $(p_1,x_1;\ldots;p_n,x_n)$ is evaluated by the functional ${\rm CPT}(p_1,x_1;\cdots;p_n,x_n)=\sum_{j=1}^n \pi_j v(x_j)$,
where $v$ is the value function, and $\pi_j$, $j=1,\ldots,n$, are defined as
\begin{equation*}
\pi_j=\left\{\begin{aligned}
g_{2}(p_1+\cdots+p_j)-g_{2}(p_1+\cdots+p_{j-1}) & \; \inmat{if}\;j \leq k,\\
g_{1}(p_j+\cdots+p_n)-g_{1}(p_{j+1}+\cdots+p_n) & \; \inmat{if}\;j > k,
\end{aligned}
\right.
\end{equation*}
where $k\in \{0,\cdots,n\}$ satisfies $x_1\geq \cdots \geq x_k \geq 0 > x_{k+1}\geq \cdots \geq x_n$.
The relation in \eqref{eq:ineq5} is equivalent to the inequality expressed in terms of CPT-values:
\begin{equation*}
\begin{aligned}
&  \sum_{j=1}^k{(g_2(p_1+\cdots+p_j)-g_2(p_{1}+\cdots+p_{j-1}))}v(x_j) \\
& \quad +\sum_{j={k+1}}^n{(g_1(p_j+\cdots+p_n)-g_1(p_{j+1}+\cdots+p_{n}))}v(x_j)\\
&\quad \leq   \sum_{l=1,l\neq i}^k(g_2(p_1+\cdots+p_l)-g_2(p_{1}+\cdots+p_{l-1}))v(x_l)  \\
&\quad \quad  +\sum_{l={k+1},l\neq j}^n(g_1(p_l+\cdots+p_n)-g_1(p_{l+1}+\cdots+p_{n}))v(x_l)\\
&\quad  \quad  +(g_2(p_1+\cdots+p_i)-g_2(p_{1}+\cdots+p_{i-1}))v\left(x_i+\frac{\delta}{p_i}\right) \\
&\quad \quad  +(g_1(p_j+\cdots+p_n)-g_1(p_{j+1}+\cdots+p_{n}))v\left(x_j-\frac{\delta}{p_j}\right),
\end{aligned}
\end{equation*}
which reduces to
\begin{equation*}
\begin{aligned}
0\leq& (g_2(p_1+\cdots+p_i)-g_2(p_{1}+\cdots+p_{i-1}))\left(v\left(x_i+\frac{\delta}{p_i}\right) -v(x_i)\right)\\
&+(g_1(p_j+\cdots+p_n)-g_1(p_{j+1}+\cdots+p_{n}))\left(v\left(x_j-\frac{\delta}{p_j}\right)-v(x_j)\right).
\end{aligned}
\end{equation*}
Let $p=p_i$, $q=p_j$, $r=p_1+\ldots+p_{i-1}$ and $s=p_{j+1}+\ldots+p_n$. 
Then the inequality becomes 
\begin{equation*}
(g_1(s+q)-g_1(s))\left(v(x_j)-v\left(x_j-\frac{\delta}{q}\right)\right)
\leq (g_2(p+r)-g_2(r))\left(v\left(x_i+\frac{\delta}{p}\right) -v(x_i)\right),
\end{equation*}
which can be rewritten as
\begin{equation}
\label{eq:relation_vg}
\frac{v(x_i+\frac{\delta}{p}) -v(x_i)}{v(x_j)-v(x_j-\frac{\delta}{q})}
\geq \frac{g_1(s+q)-g_1(s)}{g_2(p+r)-g_2(r)}.
\end{equation}
Since $v$ is convex over $\R_+$ and $\R_-$,
$\frac{v(x_i+\frac{\delta}{p}) -v(x_i)}{{\delta}/{p}}$ decreases as $x_i\downarrow 0$ and $\delta\downarrow0$, attaining its minimum $v'_+(0)$ at $0$, 
and $\frac{v(x_j)-v(x_j-\frac{\delta}{q})}{{\delta}/{q}}$ increases as $x_j\uparrow 0$ and $\delta\downarrow0$, attaining the maximum $v'_-(0)$ at $0$. 
Thus, the left-hand side of \eqref{eq:relation_vg} attains its minimum of $\frac{v_+'(0)}{v_-'(0)} \frac{q}{p}$.
Since $g_2$ is concave, $\frac{g_2(p+r)-g_2(r)}{p}\geq g'_{2,-}(p+r)$. 
Likewise, since $g_1$ is convex, we have $\frac{g_1(s+q)-g_1(s)}{q}\leq g'_{1,-}(s+q)$, 
where $g'_{2,-}(p+r)$ and $g'_{1,-}(s+q)$ denote the left derivatives of $g_2$ at $p+r$ and of $g_1$ at $s+q$, respectively. 
Thus, the right-hand side of \eqref{eq:relation_vg} attains its maximum equal to
$\frac{g'_{1,-}(s+q)}{g'_{2,-}(p+r)} \frac{q}{p}$.
Summarizing the above discussion, 
we conclude from \eqref{eq:relation_vg} that
\begin{equation}
\label{eq:ineq6}
\frac{v_+'(0)}{v_-'(0)}\geq \frac{g'_{1,-}(s+q)}{g'_{2,-}(p+r)}, \; \forall s+q\in (0,1),
\end{equation}
which is equivalent to 
\begin{equation}
\label{eq:SRT-condition}
\frac{v_+'(0)}{v_-'(0)}\geq \sup_{\eta\in (0,1)} \frac{g'_{1,-}(\eta)}{g'_{2,-}(1-\eta)}.
\end{equation}
To establish sufficiency, observe that $p+r\leq1-(s+q)$.
By the concavity of $g_2$,
$g'_{2,-}(p+r)\geq g'_{2,-}(1-(s+q))$.
Hence, \eqref{eq:SRT-condition} implies \eqref{eq:ineq6}.
To establish necessity, consider adjacent outcomes with $j=i+1$,
for which $p+r=1-(s+q)$, and let $p,q\downarrow0$.
The resulting derivative ratios yield \eqref{eq:SRT-condition}.
Finally,
$g'_{2,-}(1-\eta)=(1-g_2(1-\cdot))'_+(\eta)$.

\item[(ii)] 
Observe that $r{\rm CPT}(p_1,x_1;\cdots;p_n,x_n)={\rm CPT}(p_1,rx_1;\cdots;p_n,rx_n)$ can be rewritten in terms of CPT values:
\begin{equation*}
\begin{aligned}
& \sum_{j=1}^k{(g_2(p_1+\cdots+p_j)-g_2(p_{1}+\cdots+p_{j-1}))}rv(x_j)  \\
& \quad   +\sum_{j={k+1}}^n{(g_1(p_j+\cdots+p_n)-g_1(p_{j+1}+\cdots+p_{n}))}rv(x_j)\\
& \quad = \sum_{j=1}^k(g_2(p_1+\cdots+p_j)-g_2(p_{1}+\cdots+p_{j-1}))v(rx_j) \\
&\quad \quad  +\sum_{j={k+1}}^n(g_1(p_j+\cdots+p_n)-g_1(p_{j+1}+\cdots+p_{n}))v(rx_j).
\end{aligned}
\end{equation*}
Considering the certain prospect $(1,x_j)$ for each $j=1,\ldots,n$,
the positive homogeneity of CPT gives
\begin{equation*}
\begin{aligned}
 rv(x_j)=v(rx_j), \; \forall j=1,\cdots,n.
\end{aligned}
\end{equation*}
Thus, $v$ is positively homogeneous.
\end{itemize}

\section{An Example of Calculating \texorpdfstring{$\bbe_{g_1g_2}[v(\xi-x)]$}{}}
\label{SC:exa-bbe}
We now illustrate the computation of $\bbe_{g_1g_2}[v(\xi-x)]$ presented in Section~\ref{Sec:computation_gsr}.
\begin{example}
\label{exa:compute_bbe}
Consider the value function $v$ and the distortion functions $g_1,g_2$ as in Figure~\ref{fig:exa_value_distortion}.

\begin{figure}[!ht]
\centering
\begin{minipage}[t]{0.48\linewidth}
\centering
{
\includegraphics[width=\linewidth]{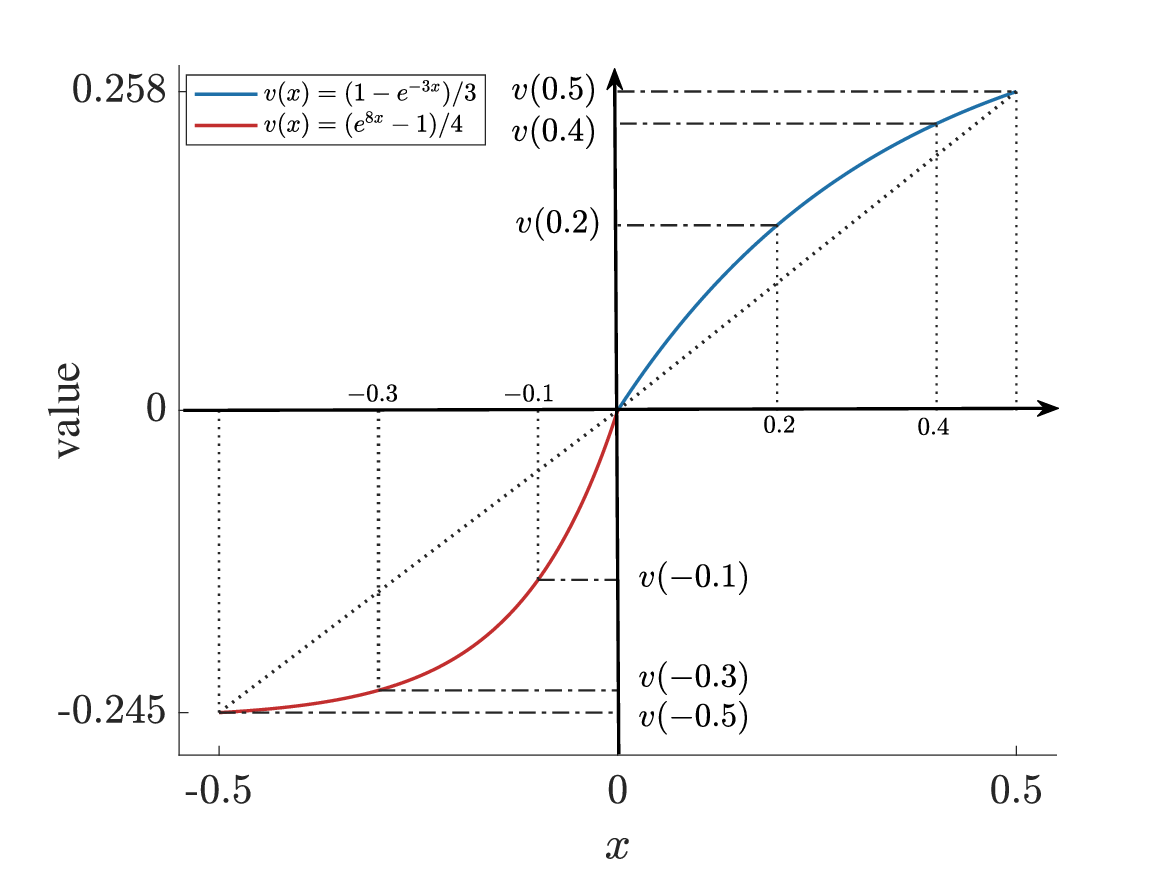}
 \text{\hspace{0.5cm} \footnotesize\bf{(a)}} 
}
\end{minipage}
\hfill
\begin{minipage}[t]{0.48\linewidth}
\centering
{
\includegraphics[width=\linewidth]{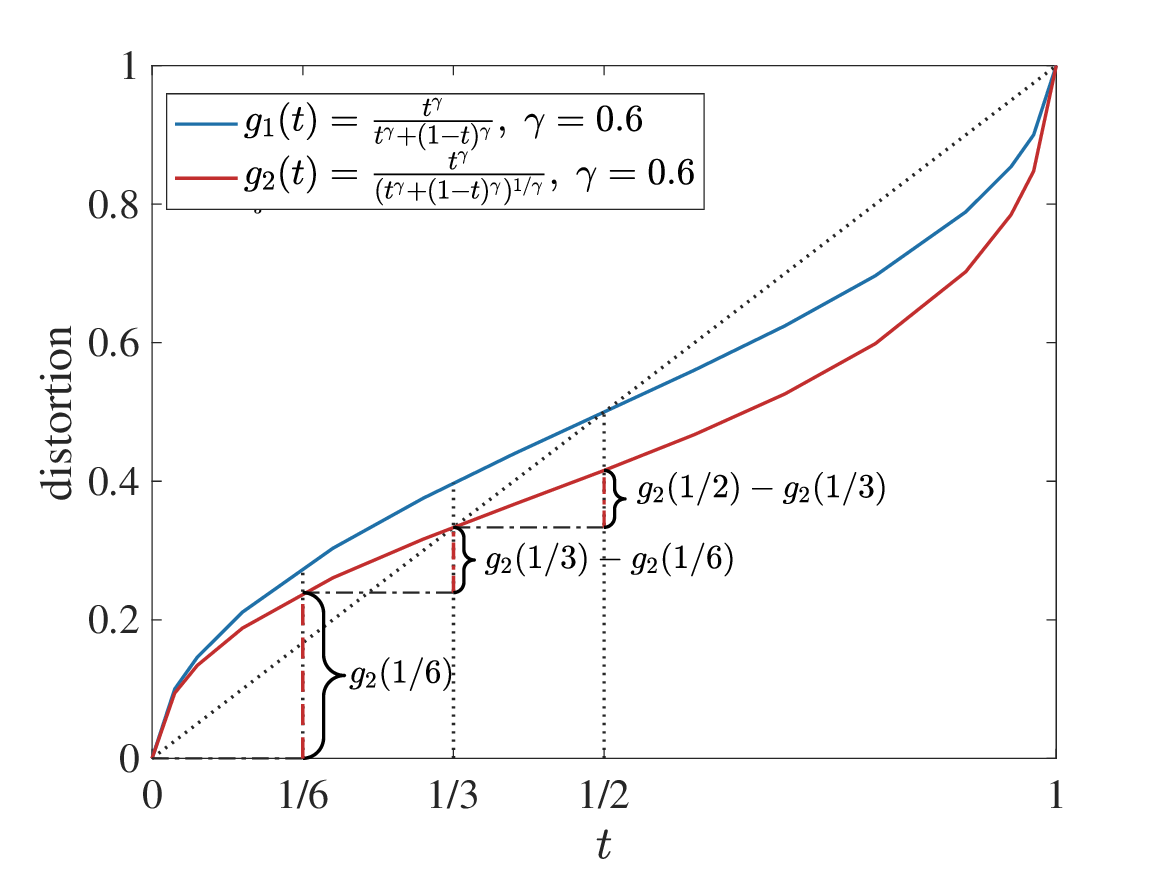}
\text{\hspace{0.5cm} \footnotesize\bf{(b)}}
}
\end{minipage}
\caption{\footnotesize\textbf{(a)} Example of the value function (nonsmooth at $0$).
\textbf{(b)} Example of the distortion functions.}
\label{fig:exa_value_distortion}
\end{figure}

To see how the distorted expectation $\bbe_{g_1g_2}[v(\xi-x)]$ is computed,
let $\xi$ be a random variable with $P(\xi=-0.1)=1/6$, $P(\xi=0.2)=1/6$, $P(\xi=-0.3)=1/6$, $P(\xi=0.4)=1/6$, $P(\xi=-0.5)=1/6$, $P(\xi=0.5)=1/6$ and $x=0$.
Then from Tables~\ref{tal:h1}--\ref{tal:h2},
\begin{equation*}
\begin{aligned}
h_{g_1g_2}^1[v(\xi-x)] =& v(-0.5)g_1({1}/{6})+v(-0.3)\Big(g_1({1}/{3})-g_1({1}/{6})\Big)+v(-0.1)\Big(g_1({1}/{2})-g_1({1}/{3})\Big)\\
=& -0.0676-0.027-0.0141=-0.1087,\\
h_{g_1g_2}^2[v(\xi-x)]=&v(0.2)\Big(g_2({1}/{2})-g_2({1}/{3})\Big)+v(0.4)\Big(g_2({1}/{3})-g_2({1}/{6})\Big)+v(0.5)g_2({1}/{6})\\
=& 0.0123+0.0219+0.0619=0.0961.
\end{aligned}
\end{equation*}
Then $\bbe_{g_1g_2}[v(\xi-x)]=h_{g_1g_2}^1[v(\xi-x)]+h_{g_1g_2}^2[v(\xi-x)]=-0.0126$.

\begin{table}[!ht]
\centering
\caption{Distorted probabilities for Losses}
{%
\begin{tabular}{cccc}
\hline
 Order No. $k$ & 1 & 2 & 3 \\
\hline
 Loss $\xi_k<0$ & $-0.5$ & $-0.3$ & $-0.1$ \\
$q_k$ & $q_1=1/6$ & $q_2=1/3$ & $q_3=1/2$ \\
\makecell[c]{Distorted \\ probability}
& \makecell[c]{$g_1(q_1)-g_1(q_0)$\\
$=g_1(1/6)$}
& \makecell[c]{$g_1(q_2)-g_1(q_1)$\\
$=g_1(1/3)-g_1(1/6)$}
& \makecell[c]{$g_1(q_3)-g_1(q_2)$\\
$=g_1(1/2)-g_1(1/3)$}\\
\hline
\end{tabular}
}
\label{tal:h1}
\end{table}

\begin{table}[!ht]
\centering
\caption{Distorted Probabilities for Gains}
{%
\begin{tabular}{cccc}
\hline
 Order No. $k$ & 4 & 5 & 6 \\
\hline
Gain $\xi_k>0$ & $0.2$ & $0.4$ & $0.5$ \\
$q_k$ & $q_4=2/3$ & $q_5=5/6$ & $q_6=1$ \\
\makecell[c]{Distorted \\ probability}
& \makecell[c]{$g_2(1-q_3)-g_2(1-q_4)$\\
$=g_2(1/2)-g_2(1/3)$}
& \makecell[c]{$g_2(1-q_4)-g_2(1-q_5)$\\
$=g_2(1/3)-g_2(1/6)$}
& \makecell[c]{$g_2(1-q_5)-g_2(1-q_6)$\\
$=g_2(1/6)$}\\
\hline
\end{tabular}
}
\label{tal:h2}
\end{table}

\end{example}

\section{Bisection Method for Computing the CPT-based Shortfall \texorpdfstring{\eqref{eq:rho-CPT}}{}}
\label{sec:Bisection-norobust}
We present a general form of the bisection method.
Suppose that based on the minimal and maximal values of $\xi$ and the monotonicity of the mapping $t\mapsto \bbe_{g_1g_2}[v(\xi-t)]$, 
we can obtain a lower bound $t^0_l$ and an upper bound $t^0_u$ for the CPT-based shortfall risk measure $t^*$: 
$\rho_{v,g_1,g_2}(\xi)=\inf\{t\in \R: \bbe_{g_1g_2}[v(\xi-t)]\leq 0\}$, 
i.e., $\rho_{v,g_1,g_2}(\xi)=t^*$, where $t^0_l\leq t^* \leq t^0_u$.
Then we can carry out the bisection search to find $t^*$.
\begin{algorithm}[!ht]
\caption{Bisection method for computing $\rho_{v,g_1,g_2}(\xi)$}
\label{algorithm:GS}
\begin{itemize}[labelwidth=3.5em,labelsep=0.5em,leftmargin=4em]
\item[Step 0.]
Set $k:=0$. Set $t_l :=t^0_l$ and $t_u :=t^0_u$.
\item[Step 1.]
Set $t_k :=\frac{t_l+t_u}{2}$ and compute $s :=\bbe_{{g_1}{g_2}}[v(\xi-t_k)]$.
If $s <0$, then update $t_u := t_k$ {\color{black}and $t_l:=t_l$; 
if $s>0$, then update $t_l:= t_k$ and $t_u:=t_u$;
stop if $s=0$.}

\item[Step 2.] {\color{black} Set $k:=k+1$ and go to Step 1.}
\end{itemize}
\end{algorithm}

\section{Proof of Proposition~\ref{prop:ambiguity_N_opt}}
\label{EC:Prop4}

\proof{Proof of Proposition~\ref{prop:ambiguity_N_opt}.}
Observe that $r_1^m,r_2^m,r_3^m\in {\cal J}$ and $1-p^m\in {\cal S}_1$,
and $p^m\in {\cal S}_2$  for $m\in [M]$.
By construction in \eqref{eq:PLA_v_g1g2}, 
$v_J$ coincides with $v$ on ${\cal J}$ and $g_{iN_i}$ coincides with $g_i$ on ${\cal S}_i$.
It follows by \eqref{eq:distorted_exp-2} and \eqref{eq:PLA_v_g1g2} that
\begin{equation}
\label{eq:bbe-equal}
\begin{aligned}
\bbe_{g_{1{N}_1} g_{2{N}_2}}[v_{{J}}(Z_1^m)]
=& v_{{J}}(r_1^m)\big(g_{1{N}_1}(1-p^m)\mathds{1}_{r_1^m<0}+(1-g_{2{N}_2}(p^m))\mathds{1}_{r_1^m\geq  0}\big) \\
&  +v_{{J}}(r_3^m)\big((1-g_{1{N}_1}(1-p^m))\mathds{1}_{r_3^m<0}+g_{2{N}_2}(p^m)\mathds{1}_{r_3^m\geq 0}\big)  \\
=& v(r_1^m)\big(g_{1}(1-p^m)\mathds{1}_{r_1^m<0}+(1-g_{2}(p^m))\mathds{1}_{r_1^m\geq 0}\big) \\
&  +v(r_3^m)\big((1-g_{1}(1-p^m))\mathds{1}_{r_3^m<0}+g_{2}(p^m)\mathds{1}_{r_3^m\geq 0}\big)
=\bbe_{g_{1} g_{2}}[v(Z_1^m)], 
\end{aligned}
\end{equation}
and $\bbe_{g_{1{N}_1}g_{2{N}_2}}[v_{{J}}(Z_2^m)]=v_{J}(r_2^m)=v(r_2^m)=\bbe_{g_{1}g_{2}}[v(Z_2^m)]$, $m\in [M]$.
The relation means that under the specific structures of $v_J$ and $g_{iN_i}$,
the distorted expected value of $v_J(Z_i^m)$ based on $g_{iN_i}$ coincides with the distorted expected value of $v(Z_i^m)$ based on $g_1, g_2$.
This shows $(v_{J},g_{1N_1},g_{2N_2})\in\mathbb{V}_J\mathbb{G}_{N_1}\mathbb{G}_{N_2}$ as desired.
\hfill \Box

\section{Three Types of Inconsistency}
\label{Sec:inconsistency}

Let ${\rm CPT}(X)=\bbe_{g_1g_2}[v(X)]$ for a random variable $X\in {\cal L}^\infty(\mathbb{R})$.
We distinguish three sources of inconsistency in preference elicitation.
\begin{itemize}
\item[(1)] 
The evaluation of the CPT values of the random variables (lotteries) contains noise. Instead of having ${\rm CPT}(X_m) -{\rm CPT}(Y_m)\leq 0$,
the errors in the evaluation of ${\rm CPT}(X_m)$ and/or ${\rm CPT}(Y_m)$ may lead to ${\rm CPT}(X_m) -{\rm CPT}(Y_m)> 0$.
To address the issue, we may introduce a positive number $\gamma_m>0$ such that if ${\rm CPT}(X_m)- {\rm CPT}(Y_m)\leq \gamma_m$, 
then we would still regard $X_m$ as preferable to $Y_m$.
By setting a bound on the total errors, 
we may consider the ambiguity set 
\begin{equation*}
\left\{(v,g_1,g_2)\in \mathscr{V}\times \mathscr{G}\times \mathscr{G}:\exists\, \gamma_m\ge 0, \;
{\rm CPT}(X_m)- {\rm CPT}(Y_m)\leq \gamma_m,\;m\in[M],\;
\sum_{m=1}^M \gamma_m\leq \Gamma\right\}.
\end{equation*}
This kind of approach is widely used in the literature of preference robust optimization (PRO, see \citep{ArD15,GXZ23}). 

\item[(2)] 
The second type of error occurs in the perception of the outcomes of lotteries \citep{ArD15}. 
Specifically, the random variable $X$ is perceived as $X+\delta$.
If $\delta_{\min}\leq \delta \leq \delta_{\max}$ almost surely, 
then we could replace the condition ${\rm CPT}(X_m) -{\rm CPT}(Y_m)\leq 0$ with ${\rm CPT}(X_m+\delta_{\min}) -{\rm CPT}(Y_m+\delta_{\max})\leq 0$.
Let $\mathscr U:=\mathscr{V}_{\rm S}\times \mathscr{G}_{\rm invS}^1\times \mathscr{G}_{\rm invS}^2$.
In the case when $\delta_{\min}<0$ and $\delta_{\max}>0$, 
the latter is a relaxation of the former due to the monotonicity of the CPT functional.
In our formulations $\bbe_{g_{1}g_{2}}[v(X_m)] \leq \bbe_{g_{1}g_{2}}[v(Y_m)]$, 
for $m\in[M]$, we would need to replace $X_m$ with $X_m+\delta_{\min}$ and $Y_m$ with $Y_m+\delta_{\max}$, 
and the ambiguity set is replaced by
\begin{equation*}
(\mathbb{V}\mathbb{G}_1\mathbb{G}_2)_M
=\left\{(v,g_1,g_2)\in  \mathscr U: \bbe_{g_{1}g_{2}}[v(X_m+\delta_{\min})] 
\leq \bbe_{g_{1}g_{2}}[v(Y_m+\delta_{\max})], \;
m\in[M] \right\}.
\end{equation*}

\item[(3)]
The third type of error potentially arises from the DM's wrong responses to the elicitation query.  
In the spirit of \cite{BeO13}, we may allow at most $\epsilon M$ of the total $M$ pairwise comparisons to be mistaken (receiving wrong answers from the DM), 
in which case at most $100(1-\epsilon)\%$ of the $M$  pairwise comparisons are correct.
In that case, we can introduce binary variables $\delta_{m}$ with $\delta_{m}=1$ if the DM is mistaken about lottery $m$ and $0$ otherwise.
Consequently, we can replace the condition ${\rm CPT}(X_m)\leq {\rm CPT}(Y_m)$ with the two inequalities
\begin{equation*}
{\rm CPT}(X_m)\leq {\rm CPT}(Y_m)+\delta_m {\cal M} 
\quad 
\inmat{and}
\quad 
{\rm CPT}(Y_m)\leq {\rm CPT}(X_m)+(1-\delta_m) {\cal M},
\end{equation*}
where ${\cal M}$ is a large constant (``big $M$''), and we require $\sum_{m=1}^M\delta_m\leq \epsilon M$.
\end{itemize}

Next, we present an example which shows that the true preference function is excluded from the constructed ambiguity set when we use the wrong model (EUT, DRM or RDEU) to construct the set of preference functions (value function and/or distortion function), whereas CPT describes the DM's preference.
\begin{example}[Preference inconsistency with EUT, DRM and RDEU]
\label{exa:PRE_INC}
We consider a risky lottery $Z_1$ and a certain lottery $Z_2$ as defined in \eqref{eq:questions} with parameters $r_1\in [-0.5,0)$, $r_3\in (0,0.5]$,
$r_2\in [-0.5,0.5]$, and $p\in (0,1)$.
The true utility function $v^*$ and distortion functions $g_1^*, g_2^*$
are the same as those in Example~\ref{exa:comparison_three}.
The preference relation between $Z_1$ and $Z_2$ under the preference functional $\rho(\cdot)$ can be quantified by the inequality $h(\rho(Z_1)-\rho(Z_2))\leq 0$,
where $h=1$ when $Z_1 \preceq Z_2$, and $h=-1$ otherwise.
Here $\rho(\cdot)$ may represent $\rho_{\rm EUT}(\cdot):=\bbe[v^*(\cdot)]$, $\rho_{\rm DRM}(\cdot):=\bbe_{g_1^*}[\cdot]$, $\rho_{\rm RDEU}(\cdot):=\bbe_{g_1^*}[v^*(\cdot)]$, and $\rho_{\rm CPT}(\cdot):=\bbe_{g_1^*g_2^*}[v^*(\cdot)]$, under EUT, DRM, RDEU, and CPT, respectively.
We assume that the preference of the DM follows CPT.
Preference inconsistency may occur when we use EUT, DRM, or RDEU to describe the DM's preference in the following three cases.

\begin{enumerate}
\item[1.]
Take $(r_1,r_2,r_3,p)=(-0.15,0.01,0.25,0.4)$.
Under CPT, $\bbe_{g_1^*g_2^*}[v^*(Z_1)]-\bbe_{g_1^*g_2^*}[v^*(Z_2)]=0.0046\geq 0$, and hence the DM's response to this query is $Z_2\preceq Z_1$, $h=-1$, and 
\begin{subequations}
\begin{align}
& \underline{\rho_{\rm EUT}(Z_1)-\rho_{\rm EUT}(Z_2)= -0.0052 \leq 0}
\quad \inmat{and} \quad  \rho_{\rm DRM}(Z_1)-\rho_{\rm DRM}(Z_2)= 0.0043 \geq 0, \label{eq:HE2-a} \\
& \rho_{\rm RDEU}(Z_1)-\rho_{\rm RDEU}(Z_2)= 0.0716 \geq 0 \quad \inmat{and} \quad \rho_{\rm CPT}(Z_1)-\rho_{\rm CPT}(Z_2)=0.0046 \geq 0.
\end{align}
\end{subequations}
Based on the responses of the DM and  under the framework of the EUT, 
we can use the inequality $\rho_{\rm EUT}(Z_1)-\rho_{\rm EUT}(Z_2) \geq 0$ to construct the set $S_1$ of utility functions, 
i.e., $S_1=\{v: v(r_1)(1-p)+v(r_3)p-v(r_2)\geq 0\}$.
However, the true utility satisfies the inequality \eqref{eq:HE2-a}, 
which has the opposite sign to the inequality defining $S_1$. 
Thus, the constructed set $S_1$ does not contain the true utility.

\item[2.]
Take $(r_1,r_2,r_3,p)=(-0.5,0.01,0.1,0.8)$.
Under CPT, $\bbe_{g_1^*g_2^*}[v^*(Z_1)]-\bbe_{g_1^*g_2^*}[v^*(Z_2)]=0.0182\geq 0$, and hence the DM's response to this query is $Z_2\preceq Z_1$, $h=-1$, and 
\begin{subequations}
\begin{align}
& \rho_{\rm EUT}(Z_1)-\rho_{\rm EUT}(Z_2)= 0.0391\geq 0
\quad \inmat{and} \quad 
\underline{\rho_{\rm DRM}(Z_1)-\rho_{\rm DRM}(Z_2)=-0.0453\leq 0}, \label{eq:HE4-b} \\
& \rho_{\rm RDEU}(Z_1)-\rho_{\rm RDEU}(Z_2)= 0.0217\geq 0
\quad \inmat{and} \quad 
\rho_{\rm CPT}(Z_1)-\rho_{\rm CPT}(Z_2)=0.0182\geq 0.
\end{align}
\end{subequations}
Based on the responses of the DM and under the framework of DRM, 
we can use the inequality $\rho_{\rm DRM}(Z_1)-\rho_{\rm DRM}(Z_2)\geq 0$ to construct the set $S_2$ of distortion functions, 
i.e., $S_2=\{g_1:r_1 g_1(1-p)+r_3 (1-g_1(1-p))-r_2\geq 0\}$.
However, the true distortion satisfies the inequality \eqref{eq:HE4-b}, 
which has the opposite sign to the inequality defining $S_2$.
Thus, the constructed set $S_2$  does not contain the true distortion.

\item[3.] Take $(r_1,r_2,r_3,p)=(-0.25,0,0.25,0.4)$.
Under CPT, $\bbe_{g_1^*g_2^*}[v^*(Z_1)]-\bbe_{g_1^*g_2^*}[v^*(Z_2)]=-0.0086\leq 0$, then the DM's response to this query is $Z_1\preceq Z_2$, $h=1$, and 
\begin{subequations}
\begin{align}
& \rho_{\rm EUT}(Z_1)-\rho_{\rm EUT}(Z_2)= -0.0190\leq 0
\quad \inmat{and} \quad 
\rho_{\rm DRM}(Z_1)-\rho_{\rm DRM}(Z_2)=-0.0446\leq 0, \\
& \underline{\rho_{\rm RDEU}(Z_1)-\rho_{\rm RDEU}(Z_2)= 0.0493\geq 0}
\quad \inmat{and} \quad 
\rho_{\rm CPT}(Z_1)-\rho_{\rm CPT}(Z_2)=-0.0086\leq 0. \label{eq:HE1-c}
\end{align}
\end{subequations}
Based on the responses of the DM and under the framework of RDEU, 
we can use the inequality $\rho_{\rm RDEU}(Z_1)-\rho_{\rm RDEU}(Z_2)\leq 0$ to construct the set $S_3$ of utility functions and distortion functions, 
i.e., $S_3=\{(v,g_1): v(r_1) g_1(1-p)+v(r_3)(1-g_1(1-p))-v(r_2) \leq 0\}$.
However, the true utility and distortion function satisfy the inequality \eqref{eq:HE1-c}, which has the opposite sign to the inequality defining $S_3$. 
Thus, the constructed set $S_3$ does not contain the true utility and true distortion functions.
\end{enumerate}
\end{example}

\section{Proof of Theorem~\ref{thm:reformulation-constraint-function-VW}}
\label{EC:Theorem1}

\proof{Proof of Theorem~\ref{thm:reformulation-constraint-function-VW}.}
\underline{Part (i)}.
Since $\mathbb{V}_J\mathbb{G}_{N_1}\mathbb{G}_{N_2}\subset\mathbb{V}\mathbb{G}_{1}\mathbb{G}_{2}$,  
by the definition of $h_{\mathbb{V}\mathbb{G}_{1}\mathbb{G}_{2}}$ in \eqref{eq:constraint_fun},
for any fixed $x\in \R$,
\begin{equation}
\label{eq:hleqh-N}
h_{\mathbb{V}\mathbb{G}_{1}\mathbb{G}_{2}}(\xi-x)
\geq h_{\mathbb{V}_J\mathbb{G}_{N_1}\mathbb{G}_{N_2}}(\xi-x).
\end{equation}
It suffices to show that $h_{\mathbb{V}\mathbb{G}_{1}\mathbb{G}_{2}}(\xi-x)\leq h_{\mathbb{V}_J\mathbb{G}_{N_1}\mathbb{G}_{N_2}}(\xi-x)$.
Let $\varepsilon$ be a small positive number, 
and $(v_{\varepsilon}^*,g_{1,\varepsilon}^{*},g_{2,\varepsilon}^{*})\in {\mathbb{V}\mathbb{G}_{1}\mathbb{G}_{2}}$ be such that 
\begin{equation}
\label{eq:epsilon}
\bbe_{g_{1,\varepsilon}^*g_{2,\varepsilon}^*}[v_{\varepsilon}^*(\xi-x)]
\geq 
\sup_{(v,g_1,g_2)\in \mathbb{V}\mathbb{G}_{1}\mathbb{G}_{2}} \bbe_{g_1g_2}[v(\xi-x)]
-\varepsilon.
\end{equation}
By Proposition~\ref{prop:ambiguity_N_opt}, 
we can find a tuple of PLFs $(v_{J,\varepsilon}^*,g_{1N_1,\varepsilon}^{*},g_{2N_2,\varepsilon}^{*})\in \mathbb{V}_J\mathbb{G}_{N_1}\mathbb{G}_{N_2}$ such that $v^*_{J,\varepsilon} (\xi_k-x)= v^{*}_{\varepsilon}(\xi_k-x)$,  $k=1,\cdots,K$, $g_{1N_1,\varepsilon}^*(q_k)=g_{1,\varepsilon}^*(q_k)$ for $k=1,\cdots,K_0(x)-1$,
and $g_{2N_2,\varepsilon}^*(1-q_k)=g_{2,\varepsilon}^*(1-q_k)$ for $k=K_0(x)-1,\cdots,K$.
Consequently, we have
\begin{equation}
\label{eq:equal_bbe}
\begin{aligned}
&\bbe_{g_{1N_1,\varepsilon}^*g_{2N_2,\varepsilon}^*}[v_{J,\varepsilon}^*(\xi-x)]  \\
&\qquad =  \displaystyle \sum_{k=1}^{K_0(x)-1} v_{J,\varepsilon}^*(\xi_k-x)({g_{1N_1,\varepsilon}^*}(q_k)-{g_{1N_1,\varepsilon}^*}(q_{k-1}))  \\
&\qquad \quad   +\sum_{k=K_0(x)}^{K}
v_{J,\varepsilon}^*(\xi_k-x)({g_{2N_2,\varepsilon}^*}(1-q_{k-1})-{g_{2N_2,\varepsilon}^*}(1-q_k)) \\
&\qquad =  \displaystyle \sum_{k=1}^{K_0(x)-1} v_{\varepsilon}^*(\xi_k-x)(g_{1,\varepsilon}^*(q_k)-g_{1,\varepsilon}^*(q_{k-1}))
+\sum_{k=K_0(x)}^{K}
v_{\varepsilon}^*(\xi_k-x)(g_{2,\varepsilon}^*(1-q_{k-1})-g_{2,\varepsilon}^*(1-q_k)) \\
&\qquad = \bbe_{g_{1,\varepsilon}^*g_{2,\varepsilon}^*}[v_{\varepsilon}^*(\xi-x)].
\end{aligned}
\end{equation}
It follows by \eqref{eq:epsilon} and \eqref{eq:equal_bbe} that
\begin{equation*}
\begin{aligned}
\sup_{(v,g_1,g_2)\in \mathbb{V}_J\mathbb{G}_{N_1}\mathbb{G}_{N_2}} \bbe_{g_1g_2}[v(\xi-x)]
\geq& \displaystyle \bbe_{g_{1N_1,\varepsilon}^{*} g_{2N_2,\varepsilon}^{*}}[v_{J,\varepsilon}^*(\xi-x)]
= \bbe_{g_{1,\varepsilon}^*g_{2,\varepsilon}^*}[v_{\varepsilon}^*(\xi-x)] \\
\geq& \sup_{(v,g_1,g_2)\in \mathbb{V}\mathbb{G}_{1}\mathbb{G}_{2}}\bbe_{g_1g_2}[v(\xi-x)] -\varepsilon,
\end{aligned}
\end{equation*}
and hence
\begin{equation}
\label{eq:converse} h_{\mathbb{V}_J\mathbb{G}_{N_1}\mathbb{G}_{N_2}}(\xi-x) \geq h_{\mathbb{V}\mathbb{G}_{1}\mathbb{G}_{2}}(\xi-x),
\end{equation}
because $\varepsilon$ can be arbitrarily small.
A combination of \eqref{eq:hleqh-N} and \eqref{eq:converse} yields
\begin{equation*}
h_{\mathbb{V}\mathbb{G}_{1}\mathbb{G}_{2}}(\xi-x)
=h_{\mathbb{V}_J\mathbb{G}_{N_1}\mathbb{G}_{N_2}}(\xi-x).
\end{equation*}

\underline{Part (ii).} 
We begin by giving a finite-dimensional  program representation of $h_{\mathbb{V}_J\mathbb{G}_{N_1}\mathbb{G}_{N_2}}(\xi-x)$.
For a fixed $x\in \R$, we can write $h_{\mathbb{V}_J\mathbb{G}_{N_1}\mathbb{G}_{N_2}}(\xi-x)$ as
\begin{equation}
\label{eq:h_VG1G2-N}
\begin{aligned}
h_{\mathbb{V}_J\mathbb{G}_{N_1}\mathbb{G}_{N_2}}(\xi-x)
= \sup_{v_J,g_{1N_1},{g_{2N_2}}} & \bbe_{g_{1N_1}g_{2N_2}}[v_{J}(\xi-x)]  \\
{\rm s.t.}\;\;\;\;\;\; & h_m(\bbe_{g_{1N_1}g_{2N_2}}[v_J(Z_1^m)]-v_J(r_2^m))\leq 0, \; m\in [M], \\
& v_J \in \mathscr{V}_{\rm S}^J,\; g_{iN_i}\in \mathscr{G}_{\rm invS}^{iN_i}, \; i=1,2. 
\end{aligned}
\end{equation}
Recall that $v_{J}(x)=(x{\bm a}+{\bm b})^\top{\bm \phi}(x)$, 
and $g_{iN_i}(t)={\bm \sigma}_{i}^\top{\bm h}^{{\cal S}_i}(t)$, $i=1,2$ (see \eqref{eq:PLA_utility} and \eqref{eq:PLA_distortion}).
Based on Part (i), 
it suffices to reformulate problem \eqref{eq:h_VG1G2-N} in terms of ${\bm a}, {\bm b}, {\bm \sigma}_1$ and ${\bm \sigma}_2$.
By \eqref{eq:distorted_exp-2}, we obtain
\begin{equation*}
\begin{aligned}
\bbe_{g_{1N_1} g_{2N_2}}[v_J(Z_1^m)]
=& (r_1^m{\bm a}+{\bm b})^\top{\bm \phi}(r_1^m)\Big({\bm \sigma}_1^\top{\bm h}^{{\cal S}_1}(1-p^m)\mathds{1}_{r_1^m<0}+(1-{\bm \sigma}_2^\top{\bm h}^{{\cal S}_2}(p^m))\mathds{1}_{r_1^m\geq 0}\Big)\\
&+(r_3^m{\bm a}+{\bm b})^\top{\bm \phi}(r_3^m)\Big((1-{\bm \sigma}_1^\top{\bm h}^{{\cal S}_1}(1-p^m))\mathds{1}_{r_3^m<0}+{\bm \sigma}_2^\top{\bm h}^{{\cal S}_2}(p^m)\mathds{1}_{r_3^m\geq 0}\Big)\\
=& R^m({\bm a},{\bm b},{\bm \sigma}_1,{\bm \sigma}_2),
\end{aligned}
\end{equation*}
and $\bbe_{g_{1N_1}g_{2N_2}}[v_J(Z_2^m)]=v_J(r_2^m)= (r_2^m {\bm a}+{\bm b})^\top {\bm \phi}(r_2^m)$.
Thus, the inequalities $h_m(\bbe_{g_{1N_1} g_{2N_2}}[v_J(Z_1^m)]- \bbe_{g_{1N_1}g_{2N_2}}[v_J(Z_2^m)])\leq 0$ for $m\in [M]$ 
are equivalent to $h_m(R^m({\bm a},{\bm b},{\bm \sigma}_1,{\bm \sigma}_2)-(r_2^m{\bm a}+{\bm b})^\top{\bm \phi}(r_2^m))\leq 0$ for $m\in [M]$, which is constraint \eqref{robust_Wrj}.
On the other hand, by \eqref{eq:saa_ww}, the objective function of $h_{\mathbb{V}_J\mathbb{G}_{N_1}\mathbb{G}_{N_2}}(\xi-x)$ takes the form
\begin{equation*}
\begin{aligned} 
\bbe_{g_{1N_1}g_{2N_2}}[v_J(\xi-x)] 
=& \sum_{k=1}^{K_0(x)-1} v_J(\xi_k-x)({g_{1N_1}}(q_k)-{g_{1N_1}}(q_{k-1})) \\
& \displaystyle +\sum_{k=K_0(x)}^{K} v_J(\xi_k-x)({g_{2N_2}}(1-q_{k-1})-{g_{2N_2}}(1-q_k))\\
=& \sum_{k=1}^{K_0(x)-1}((\xi_k-x) {\bm a}+{\bm b})^\top{\bm \phi}(\xi_k-x){\bm \sigma}_1^\top ({\bm h}^{{\cal S}_1}(q_k)-{\bm h}^{{\cal S}_1}(q_{k-1})) \\
& +\sum_{k=K_0(x)}^{K}((\xi_k-x) {\bm a}+{\bm b})^\top{\bm \phi}(\xi_k-x)
{\bm \sigma}_2^\top ({\bm h}^{{\cal S}_2}(1-q_{k-1})-{\bm h}^{{\cal S}_2}(1-q_{k})).
\end{aligned}
\end{equation*}
Thus, problem \eqref{robust_v_o} provides an equivalent finite-dimensional reformulation of $h_{\mathbb{V}_J\mathbb{G}_{N_1}\mathbb{G}_{N_2}}(\xi-x)$ and the proof is complete.
\hfill \Box

\section{Bisection Method for Computing  PRGSR-CPT \texorpdfstring{\eqref{eq:PRGRS-CPT}}{}}
\label{sec:bisection-robust}

\begin{algorithm}[!ht]
\caption{Bisection method for computing PRGSR-CPT \eqref{eq:PRGRS-CPT}}
\label{algorithm:thm4-1}
\begin{itemize}[labelwidth=3.5em,labelsep=0.5em,leftmargin=4em]
\item[Step 0.]
Set $k:=0$. Set $t_l :=t^0_l$ and $t_u :=t^0_u$.
\item[Step 1.]
Set $t_k:=\frac{t_l+t_u}{2}$ and compute $s :=\sup_{(v,g_1,g_2)\in \mathbb{V}\mathbb{G}_1\mathbb{G}_2}\bbe_{{g_1}{g_2}}[v(\xi-t_k)]$.
If $s<0$, then update $t_u := t_k$; otherwise, update $t_l:= t_k$.
Set $k:=k+1$. 
\end{itemize}
\end{algorithm}

\section{Proof of Proposition~\ref{prop:pS_2}}

\proof{Proof of Proposition~\ref{prop:pS_2}.}
Since ${\bm c}_2^m>0$ and ${\bm e}^\top{\bm c}_2^m<1$, $({\bm c}_2^m,1-{\bm e}^\top{\bm c}_2^m)>0$, 
which implies that the function $p \mapsto ({\bm R}{\bm c}_2^{m})^\top {\bm h}^{{\cal S}_2}(p)$ is strictly increasing on $[0,1]$.
Moreover, by the assumption that $(\widetilde{\bm \sigma}_2^{1m},1-{\bm e}^\top\widetilde{\bm \sigma}_2^{1m})$ has no zero component, 
it follows that $\{p^{m,i}\}$ is determined by $({\bm R}{\bm c}_2^{m})^\top {\bm h}^{{\cal S}_2}(p^{m,i})=D_i$ in \eqref{eq:gene_m_Q-1}, for $i\in [N_2-2]$.
Note that ${\bm h}^{{\cal S}_2}(p^{m,i})=\Big(\underbrace{1,\cdots,1}_{j_i-1},\frac{p^{m,i}-t_{j_i}}{t_{j_i+1}-t_{j_i}},\underbrace{0,\cdots,0}_{N_2-1-j_i}\Big) \in \R^{N_2-1}$, for $p^{m,i}\in (t_{j_i},t_{j_i+1}]$ with some $j_i\in [N_2-1]$.
Observe that
\begin{equation}
\label{eq:Rg-a}
({\bm R}{\bm c}_2^m)^\top {\bm h}^{{\cal S}_2}(p^{m,i})
=\sum_{l=1}^{j_i-1}c_{2,l}^m+c_{2,j_i}^m\frac{p^{m,i}-t_{j_i}}{t_{j_i+1}-t_{j_i}},\;\inmat{for}\; p^{m,i} \in (t_{j_i},t_{j_i+1}] \; \inmat{with some}\; j_i\leq N_2-2.
\end{equation}
Let $s\in [N_2-2]$.
If $D_i=\sum_{j=1}^s c_{2,j}^m$ and $p^{m,i}$ satisfies $({\bm R}{{\bm c}_2^m})^\top {\bm h}^{{\cal S}_2}(p^{m,i})=D_i$, then \eqref{eq:Rg-a} is equivalent to
\begin{equation*}
\sum_{l=1}^{j_i-1}c_{2,l}^m+c_{2,j_i}^m\frac{p^{m,i}-t_{j_i}}{t_{j_i+1}-t_{j_i}}
=\sum_{j=1}^s c_{2,j}^m,
\end{equation*}
which means $p^{m,i}=t_{s+1}$.
The proof is complete.
\hfill \Box

\section{RRUS versus the Modified Polyhedral Method}
\label{sec:RRUS_MPH}

\begin{example}[RRUS vs. the Modified Polyhedral Method] 
\label{exa:comparison}
Consider the true utility function:
$u^*(t)=2(e^{8t}-1)/8, \;\; \inmat{for}\;\; t\in [-0.5,0]$
and $u^*(t)=(1-e^{-3t})/3 \;\; \inmat{for}\;\; t\in [0,0.5]$.
\begin{itemize}
\item[1.] {\bf Random Relative Utility Split Scheme}.
We randomly generate $10$ queries with a risky lottery $Z_1^m$ and a deterministic lottery $Z_2^m$:
\begin{equation*}
Z_1^m=\left\{\begin{aligned}
r_1^m & \inmat{w.p.}\; 1-p^m \\
r_3^m & \inmat{w.p.}\; p^m
\end{aligned}
\right.
\quad 
\inmat{and}
\quad 
Z_2^m=r_2^m,
\end{equation*}
where the parameters $(r_1^m,r_2^m,r_3^m)$ are given by
\begin{equation*}
\begin{aligned}
m=1: &(-0.39389,-0.12723,0.13943);&& m=2: (-0.34508,0.00218,0.34945) \\
m=3: &(-0.02082,0.08383,0.18848); &&
m=4:  (-0.47790,-0.12964,0.21861) \\
m=5: &(-0.09662,0.08570,0.26803); &&
m=6: (-0.37031,-0.34015,-0.30999) \\
m=7: &(-0.31495,-0.27677,-0.23860); &&
m=8: (-0.36184,-0.09444,0.17296) \\
m=9: & ( -0.35396,-0.05848,0.23700); &&
m=10: (-0.46137,-0.38929,-0.31721).
\end{aligned}
\end{equation*}
Next, we use RRUS to determine the parameter $p^m$.
Specifically, let ${\cal X}=\{-0.5,0, 0.5,r_1^m,r_2^m,r_3^m: m=1,\ldots,10\}$
be the set of breakpoints with $N=|{\cal X}|=33$ and $\mathscr{U}_N$ be the set of all increasing piecewise linear utility functions defined over $[-0.5,0.5]$ with breakpoints ${\cal X}$. 
Let ${\cal U}_N^0=\mathscr{U}_N$.
For $(r_1^1,r_2^1,r_3^1)$, let $u_N(r_1^1)=0$ and $u_N(r_3^1)=1$.
Then 
\begin{equation*}
\bbe[u_N(Z_1^1)]-\bbe[u_N(Z_2^1)]=
(1-p^1)u_N(r_1^1)+p^1 u_N(r_3^1)-u_N(r_2^1)
=p^1 -u_N(r_2^1).
\end{equation*}
Let $I_1^1=\min_{u_N\in {\cal U}_N^0}u_N(r_2^1)$ and $I_2^1=\max_{u_N\in {\cal U}_N^0}u_N(r_2^1)$.
Then $p^1=0.5(I_1^1+I_2^1)$, and 
\begin{equation*}
{\cal U}_N^1=\{u_N\in {\cal U}_N^0: h_1((1-p^1)u_N(r_1^1)+p^1 u_N(r_3^1)-u_N(r_2^1))\leq 0\},
\end{equation*}
where $h_1=1$ if $Z_1^1 \preceq Z_2^1$ and $h_1=-1$ otherwise.
We consider three cases.
\begin{itemize}
\item[(i)] {No information about the shape and the Lipschitz modulus}.
We obtain $p^1=\ldots=p^{10}=0.5$.

\item[(ii)] {No information about the shape but with complete information about  the Lipschitz modulus $L=2$ of $u^*$}. 
$p^1=p^2=p^3=0.5$, $p^4=0.35374$, $p^5=0.62686$, $p^6=p^7=0.5$, $p^8=0.31244$, $p^9=0.34520$ and $p^{10}=0.5$.
 
\item[(iii)] {No information about Lipschitz modulus but with complete information about the $S$-shape with an inflection point $0$}.
We obtain $p^1=0.33849$, $p^2=0.50312$, $p^3=0.73959$, $p^4=0.19973$, $p^5=0.70754$, $p^6=0.25000$, $p^7=0.25000$, $p^8=0.30879$, $p^9=0.31161$ and $p^{10}=0.25000$.
\end{itemize}
Then ${\cal U}_N^{10}=\left\{u_N\in \mathscr{U}_N: h_m((1-p^m)u_N(r_1^m)+p^mu_N(r_3^m)-u_N(r_2^m))\leq 0, \;m=1,\ldots,10 \right\}$. 
Let
\begin{equation*}
\bar{u}_i=\max_{u_N\in {\cal U}_N^{10},u_N(-0.5)=0,u_N(0.5)=1}u_N(x_i)
\quad \inmat{and} \quad 
\underline{u}_i=\min_{u_N\in {\cal U}_N^{10},u_N(-0.5)=0,u_N(0.5)=1}u_N(x_i),\; x_i\in {\cal X}.
\end{equation*}
Figure~\ref{fig:RRUS_M10}~(A1) depicts the curve passing through $\{(x_i,u_N({x}_i))\}_{i=1}^{N}$ and the ``\textbf{upper curve}'' passing through $\{(x_i,\bar{u}_i)\}_{i=1}^{N}$, the ``\textbf{lower curve}'' through $\{(x_i,\underline{u}_i)\}_{i=1}^{N}$ in case (i).
Figure~\ref{fig:RRUS_M10}~(A2) shows the lower and upper curves in case (ii).
Figure~\ref{fig:RRUS_M10}~(A3) depicts the lower and upper curves in case (iii).

\item[2.] {\bf Modified Polyhedral Method with No Information about the Shape and the Lipschitz Modulus}.
Let us develop the modified polyhedral method for the utility function under the von Neumann and Morgenstern EUT.
We begin with an initial polyhedron for the vector of increments of the utility function: ${\Delta}_0^{\rm u}:=\{ {\bm v}\in \R^{J_0^{\rm u}-2}: {\bm v}\geq 0, {\bm e}^\top{\bm v}\leq 1\}$, where ${\cal X}_0^{\rm u}=\{-0.5,0,0.5,r_1^k,r_2^k,r_3^k:k=1,2,3\}$ is the set of breakpoints of utility functions, where $\{(r_1^{k},r_2^{k},r_3^{k})\}_{k=1}^3$ are obtained from the first three RRUS queries. 
Let $J_0^{\rm u}:=|{\cal X}_0^{\rm u}|=12$.
We generate questionnaires based on the preference function $\bbe[u(\cdot)]$.
In {Steps 1--2}, let ${\bm c}^m$ be the analytic center of $\Delta_{m-1}^{\rm u}$ and ${\bm v}^m_1$ and ${\bm v}^m_2$ be the two intersection points between the longest axis of Sonnevend's inner ellipsoid of the polyhedron $\Delta_{m-1}^{\rm u}$ and the boundary of the polyhedron.
In {Step 3}, we generate the $m$-th pair of queries by solving the following two programs:
\begin{subequations}
\label{eq:ploy-EUT}
\begin{align}
r_2^{m,i} & \in \arg\max_{x} \Big\{ \rho_1(x;{\bm v}_1^m):
x \in [\underline{x},\bar{x}],\; \rho_1(x;{\bm c}^m)\leq D_i \Big\},\\
(r_1^{m,i},r_3^{m,i},p^{m,i}) & \in \arg\max_{(x,y,p)} \left\{ \rho_2((x,y,p);{\bm v}_2^m): 
\begin{aligned}
& x\in [\underline{x},0), \; y \in (0,\bar{x}], \; p\in (0,1)\\
&\rho_2((x,y,p);{\bm c}^m)\leq D_i
\end{aligned}
\right\},
\end{align}
\end{subequations}
where $D_i\in \{\sum_{j=1}^i c_j^m: i=1,\ldots,J_{m-1}^{\rm u}-2\}$, $\rho_1(x;{\bm v})=({\bm R}{\bm v})^\top{\bm h}^{{\cal X}_{m-1}^{\rm u}}(x)$ and $\rho_2((x,y,p);{\bm v})=({\bm R}{\bm v})^\top((1-p){\bm h}^{{\cal X}_{m-1}^{\rm u}}(x)+p{\bm h}^{{\cal X}_{m-1}^{\rm u}}(y))$.
Note that in the formulations $Z_2^m$ is uniquely determined by $r_2^{m}$ and $Z_1^m$ is uniquely determined by $(r_1^m,r_3^m,p^m)$.

We select $(r_{1}^{m},r_2^m,r_3^m):=(r_{1}^{m,s_0},r_2^{m,s_0},r_3^{m,s_0})$ from $\{(r_{1}^{m,i},r_2^{m,i},r_3^{m,i})\}_{i=1}^{J_{m-1}^{\rm u}-2}$ and set $p^m:=p^{m,s_0}$, where
$s_0 \in \arg\min_{s \in \{1,\ldots,J_{m-1}^{\rm u}-2\}} \big|({\bm R}\bm c^m)^\top((1-p^{m,s}){\bm h}^{{\cal X}_{m-1}^{\rm u}}(r_1^{m,s}) + p^{m,s} {\bm h}^{{\cal X}_{m-1}^{\rm u}}(r_3^{m,s}) - {\bm h}^{{\cal X}_{m-1}^{\rm u}}(r_2^{m,s}))\big|$.
This selects the cut hyperplane whose deviation from the analytic center is minimal.
We then update the set of breakpoints: ${\cal X}_m^{\rm u}:={\cal X}_{m-1}^{\rm u}\cup \{r_{1}^{m},r_3^m\}$ and $J_m^{\rm u}:=|{\cal X}_m^{\rm u}|$.
The ambiguity set of the utility functions is then constructed by
\begin{equation*}
{\cal U}_m=\Big\{u_{J_m^{\rm u}}:[\underline{x},\bar{x}]\to [0,1]: h_l(\bbe[u_{J_m^{\rm u}}(Z_1^l)]-\bbe[u_{J_m^{\rm u}}(Z_2^l)])\leq 0,\; \forall l\in [m]\Big\},
\end{equation*}
where $h_l=1$ if $\bbe[u(Z_1^l)]\leq \bbe[u(Z_2^l)]$, and $h_l=-1$ otherwise.

Figure~\ref{fig:RRUS_M10}~(B) depicts the results of the modified polyhedral method for eliciting the true utility function using $10$ queries.
We emphasize that the polyhedral method is more efficient for utility functions with no prior information on their shape or Lipschitz modulus.

\end{itemize}
Moreover, we can generate $M\in \{30,50,100\}$ queries with parameters $(r_1^m,r_2^m,r_3^m)$ with $m=1,\ldots,M$. We then determine parameters $p^m$ by random relative utility split scheme.
Figures~\ref{fig:RRUS_M30}--\ref{fig:RRUS_M100} display the performances of utility function elicitation and the upper and lower curves generated by the resulting feasible set ${\cal U}_{M}$ with $J_M^{\rm u}=3M+3$.
We can see that the modified polyhedral method performs much better than 
RRUS in all of the cases.

\begin{figure}[!ht]
    \centering
\begin{minipage}[t]{0.24\linewidth}
\centering
\includegraphics[width=1.7in]{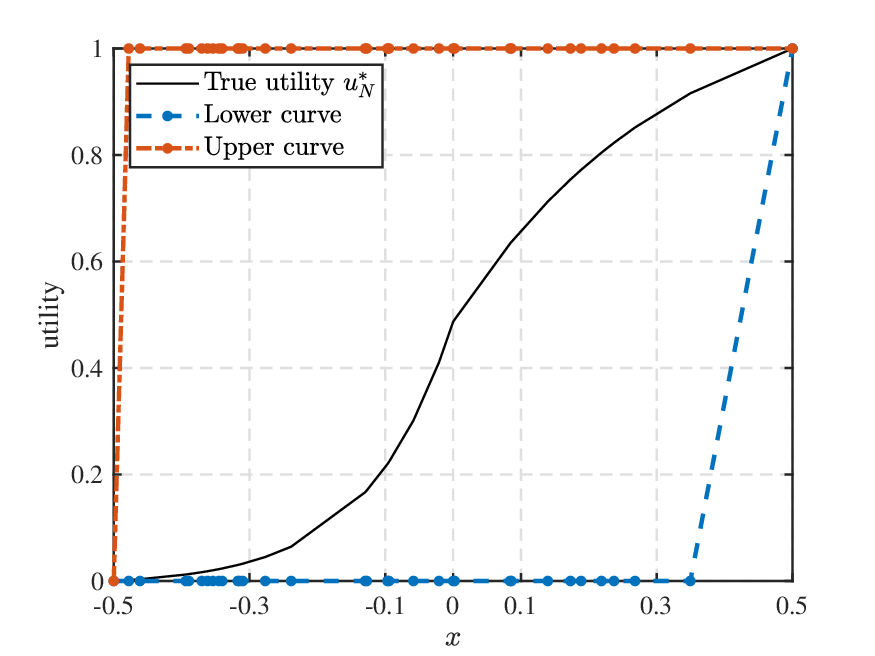}
\text{\hspace{0.5cm}\footnotesize{\bf{(A1)}}}
\end{minipage}
\hfill
\begin{minipage}[t]{0.24\linewidth}
\centering
\includegraphics[width=1.7in]{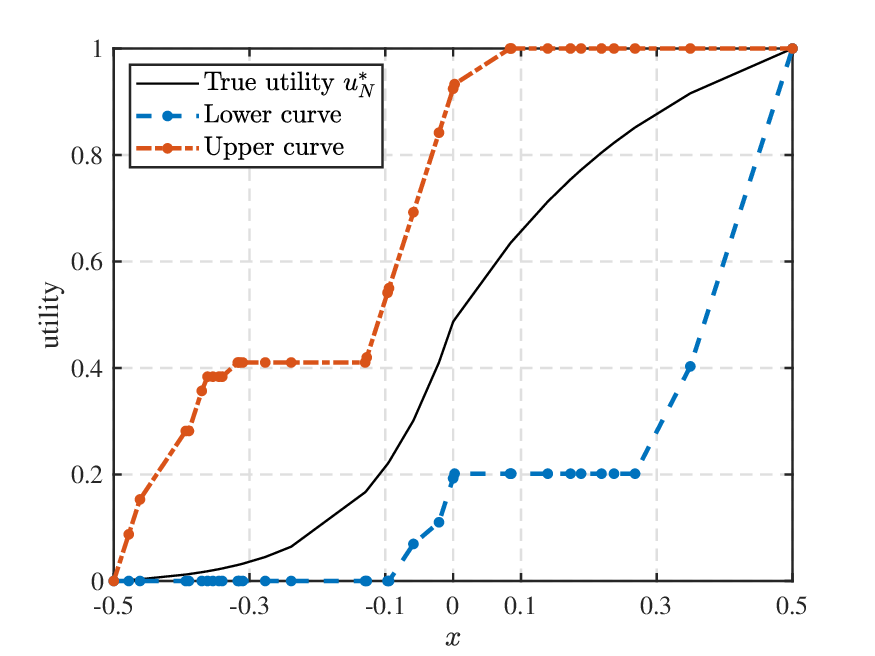}
\text{\hspace{0.5cm}\footnotesize \bf{(A2)}}
\end{minipage}
\hfill
\begin{minipage}[t]{0.24\linewidth}
\centering
\includegraphics[width=1.7in]{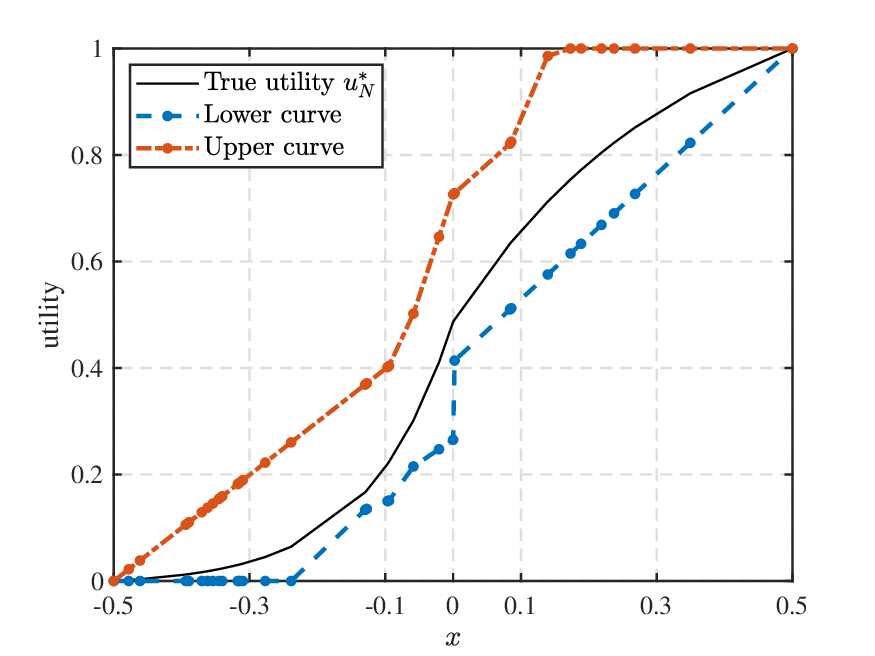}
\text{\hspace{0.5cm}\footnotesize\bf{(A3)}}
\end{minipage}
\hfill
\begin{minipage}[t]{0.24\linewidth}
\centering
\includegraphics[width=1.7in]{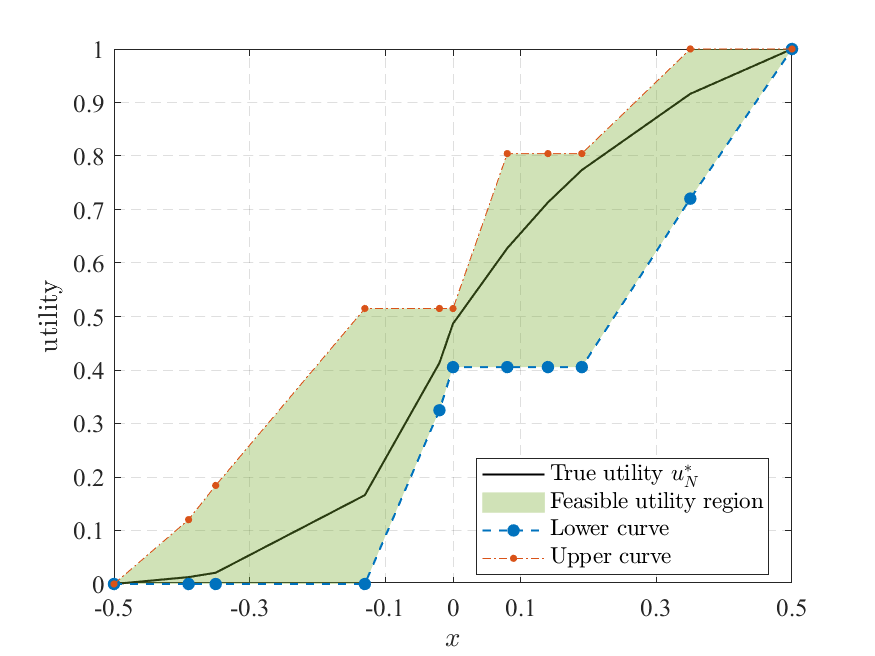}
\text{\hspace{0.5cm}\footnotesize \textbf{(B)}}
\end{minipage}
\caption{\footnotesize{\textbf{(A)} 
RRUS scheme for preference elicitation in the S-shaped utility setting with the range of utility functions 
$\mathcal{U}_{N}^{M}$ for $M=10$.
(A1) 
Case (i): no shape information and no Lipschitz constant $L$. 
The domain of the 
utility function is $[-0.5,0.5]$ ($x$-axis).
The vertical axis marks the range of the utility functions  ($[0,1]$). 
The blue curve corresponds to $\{(x_i,\underline{u}_i)\}_{i=1}^N$. 
The orange curve corresponds to $\{(x_i,\bar{u}_i\}_{i=1}^N)$. 
(A2) Case (ii): with Lipschitz constant $L$ but without the S-shape condition. 
(A3) Case (iii): with the S-shape condition but without Lipschitz constant $L$. 
\textbf{(B)} 
Polyhedral method for preference elicitation in the S-shaped utility setting with the range of utility functions $\mathcal{U}_{M}$ for $M=10$ without prior information.}}
\label{fig:RRUS_M10}
\end{figure}

\begin{figure}[!ht]
    \centering
\begin{minipage}[t]{0.24\linewidth}
\centering
 \includegraphics[width=1.7in]{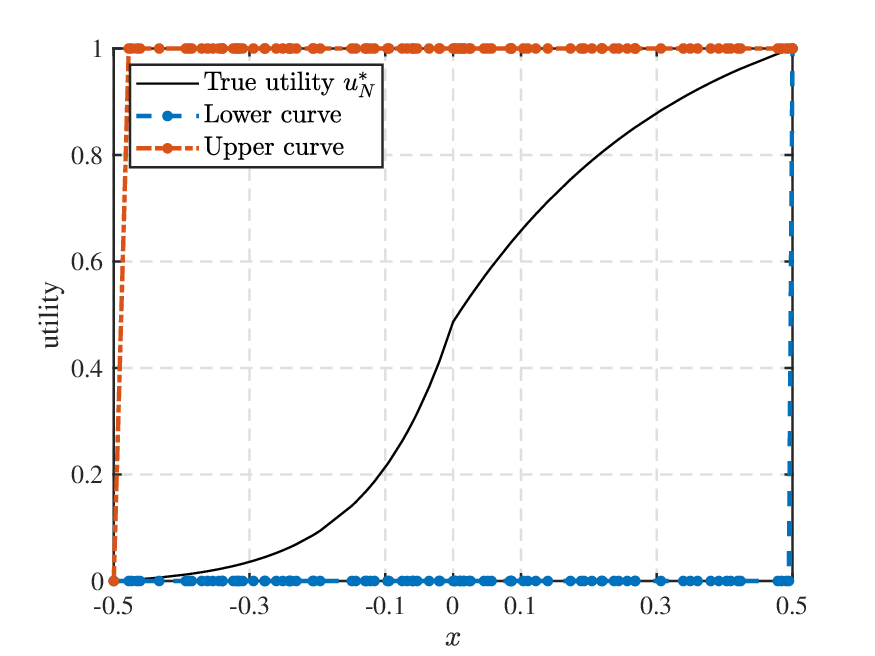}
\text{\hspace{0.5cm}\footnotesize \bf{(A1)}}
\end{minipage}
\hfill
\begin{minipage}[t]{0.24\linewidth}
\centering
\includegraphics[width=1.7in]{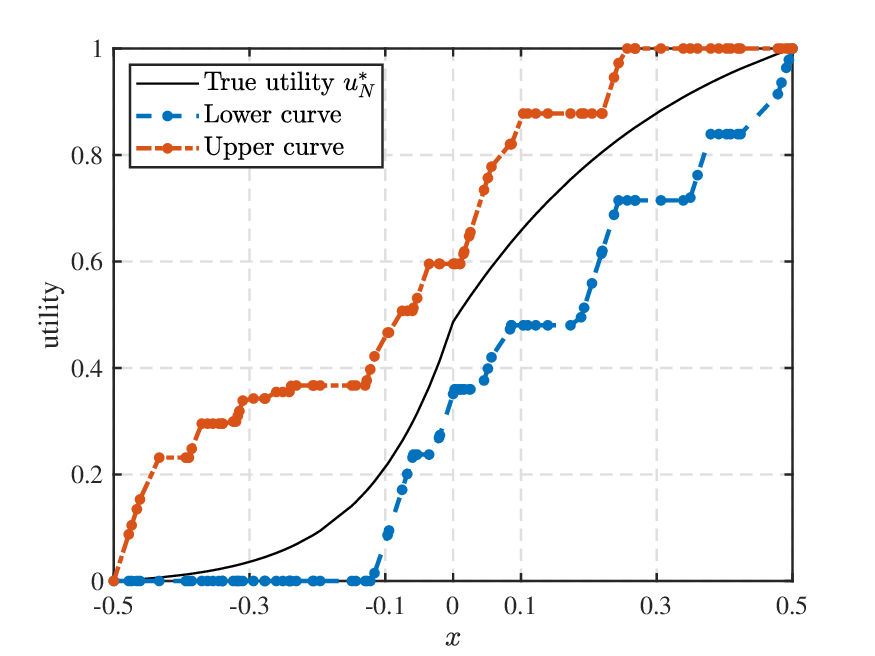}
\text{\hspace{0.5cm} \footnotesize \bf{(A2)}}
\end{minipage}
\hfill
\begin{minipage}[t]{0.24\linewidth}
\centering
\includegraphics[width=1.7in]{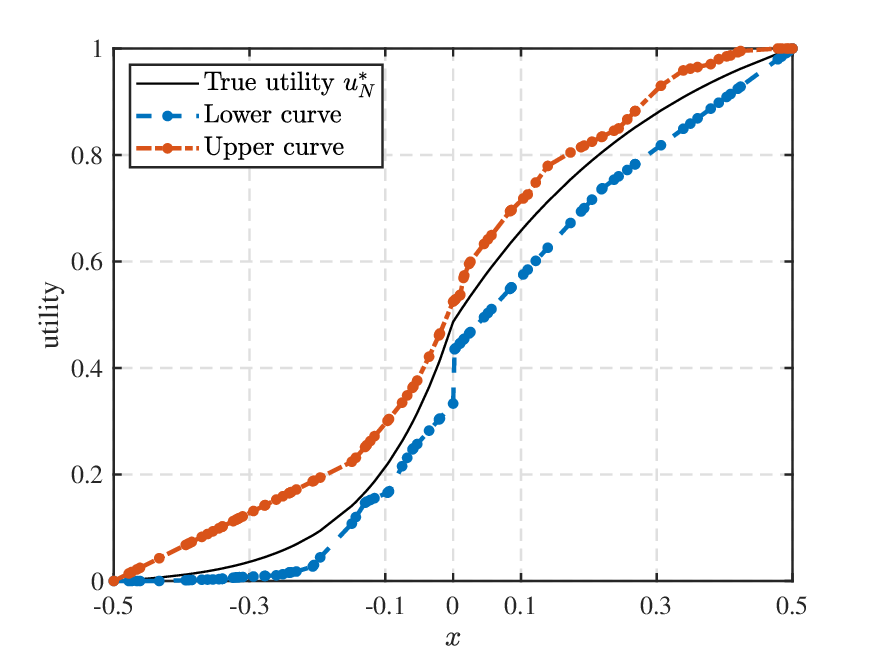}
\text{\hspace{0.5cm}\footnotesize \bf{(A3)}}
\end{minipage}
\hfill
\begin{minipage}[t]{0.24\linewidth}
\centering
\includegraphics[width=1.7in]{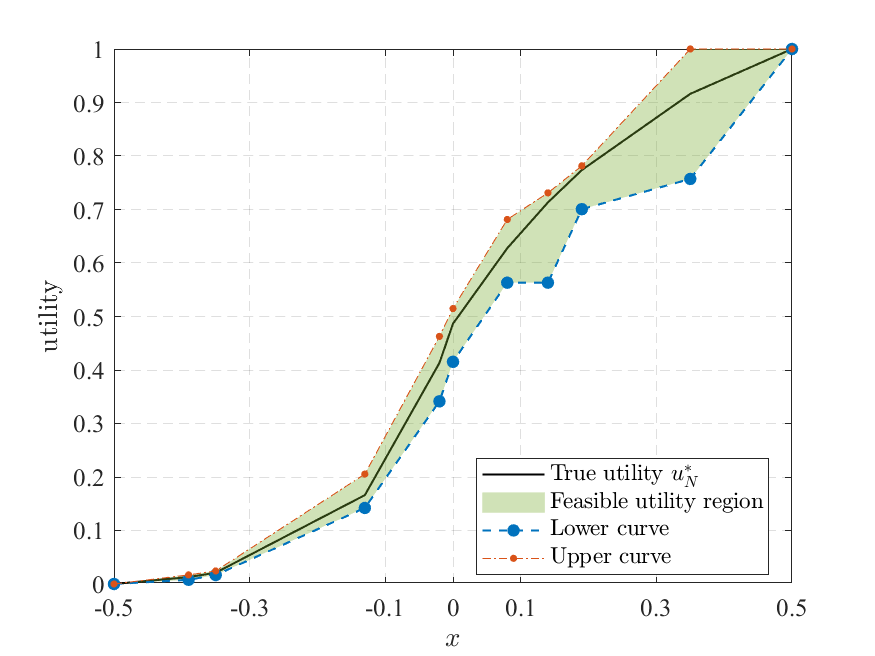}
\text{\hspace{0.5cm}\footnotesize \bf{(B)}}
\end{minipage}
\caption{\footnotesize{\textbf{(A)}
RRUS scheme for preference elicitation in the S-shaped utility setting with the range of utility functions 
$\mathcal{U}_{N}^{M}$ for $M=30$.  
\textbf{(B)}  Polyhedral method for preference elicitation in the S-shaped utility setting with the range of utility functions 
$\mathcal{U}_{M}$ for $M=30$ without prior information.
The labels, graphs and lines have the same interpretation as in Figure~\ref{fig:RRUS_M10}.
}}
\label{fig:RRUS_M30}
\end{figure}

\begin{figure}[!ht]
    \centering
\begin{minipage}[t]{0.24\linewidth}
\centering
 \includegraphics[width=1.7in]{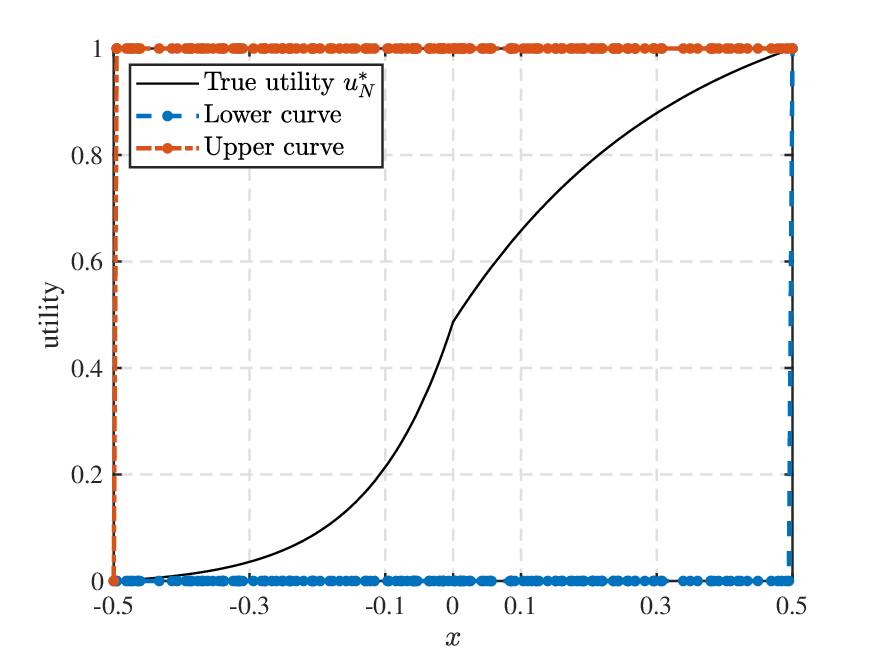}
\text{\hspace{0.5cm} \footnotesize \bf{(A1)}}
\end{minipage}
\hfill
\begin{minipage}[t]{0.24\linewidth}
\centering
\includegraphics[width=1.7in]{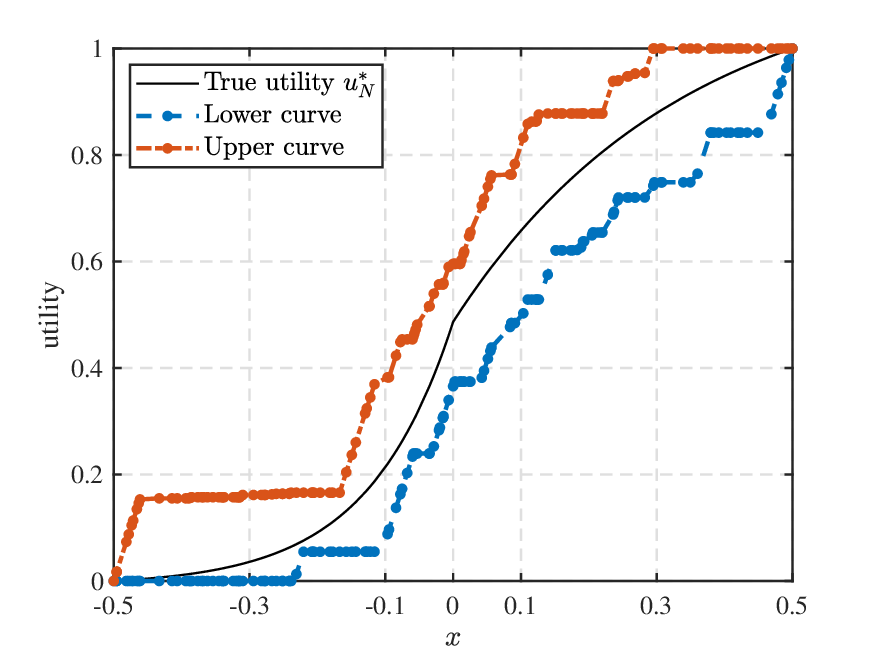}
\text{\hspace{0.5cm} \footnotesize \bf{(A2)}}
\end{minipage}
\hfill
\begin{minipage}[t]{0.24\linewidth}
\centering
 \includegraphics[width=1.7in]{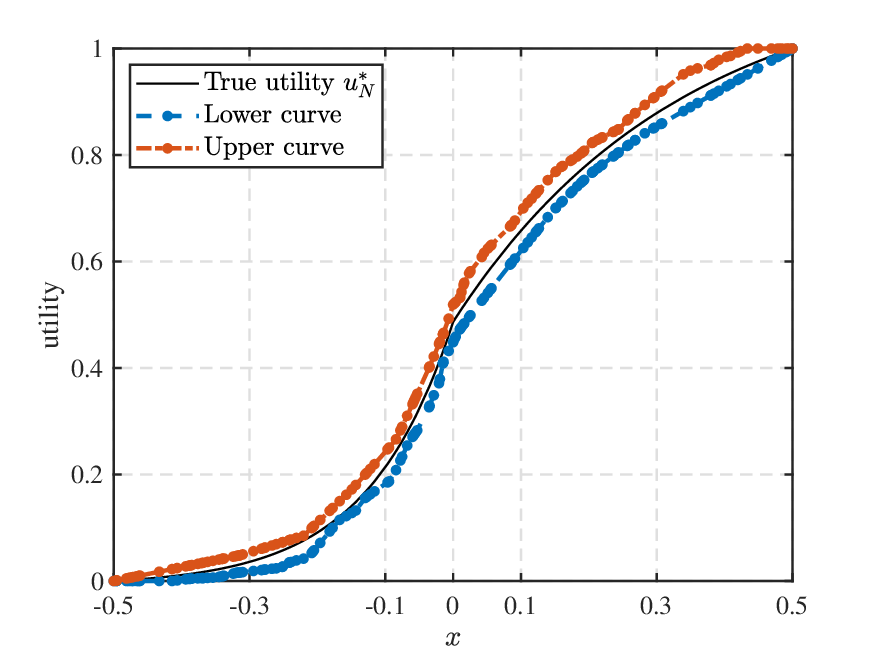}
\text{\hspace{0.5cm} \footnotesize \bf{(A3)}}
\end{minipage}
\hfill
\begin{minipage}[t]{0.24\linewidth}
\centering
\includegraphics[width=1.7in]{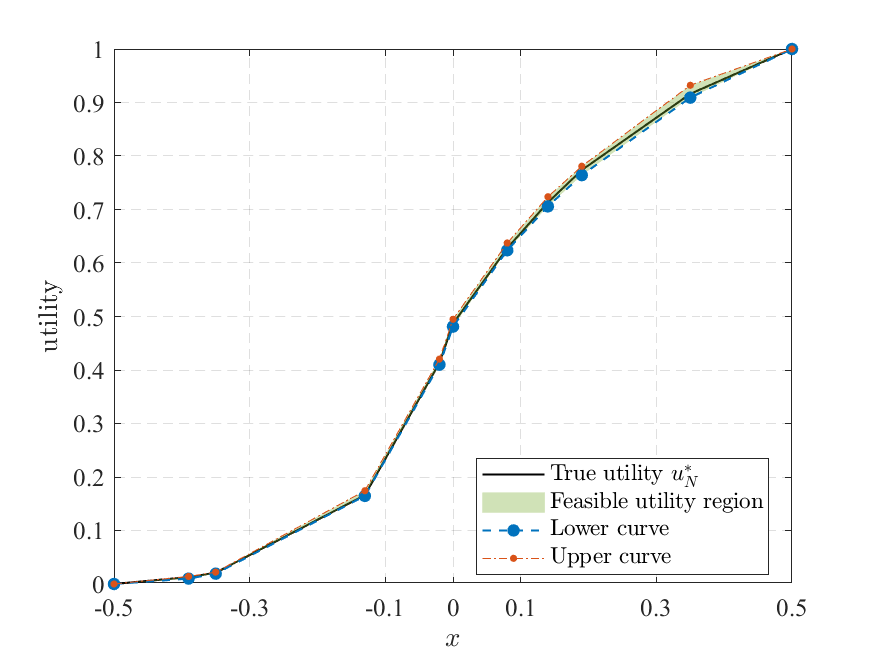}
\text{\hspace{0.5cm} \footnotesize \bf{(B)}}
\end{minipage}
\caption{\footnotesize{\textbf{(A)}
RRUS scheme for preference elicitation in the S-shaped utility setting with the range of utility functions 
$\mathcal{U}_{N}^{M}$ for $M=50$.   
\textbf{(B)}  Polyhedral method for preference elicitation in the S-shaped utility setting with the range of utility functions 
$\mathcal{U}_{M}$ for $M=50$ without prior information.
The labels, graphs and lines have the same interpretation as in Figure~\ref{fig:RRUS_M10}.
}}
\label{fig:RRUS_M50}
\end{figure}

\begin{figure}[!ht]
    \centering
\begin{minipage}[t]{0.24\linewidth}
\centering
 \includegraphics[width=1.7in]{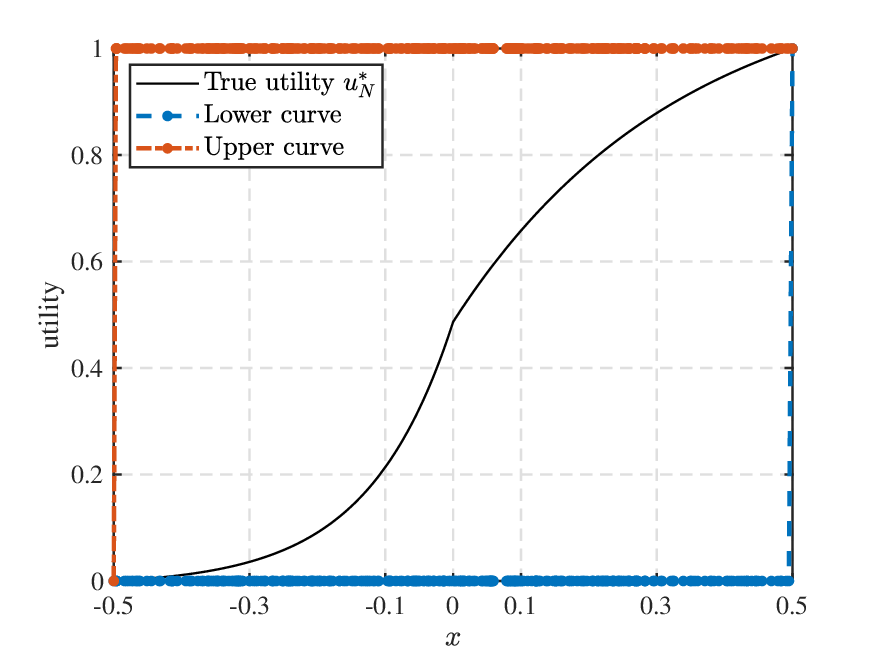}
\text{\hspace{0.5cm} \footnotesize \bf{(A1)}}
\end{minipage}
\hfill
\begin{minipage}[t]{0.24\linewidth}
\centering
 \includegraphics[width=1.7in]{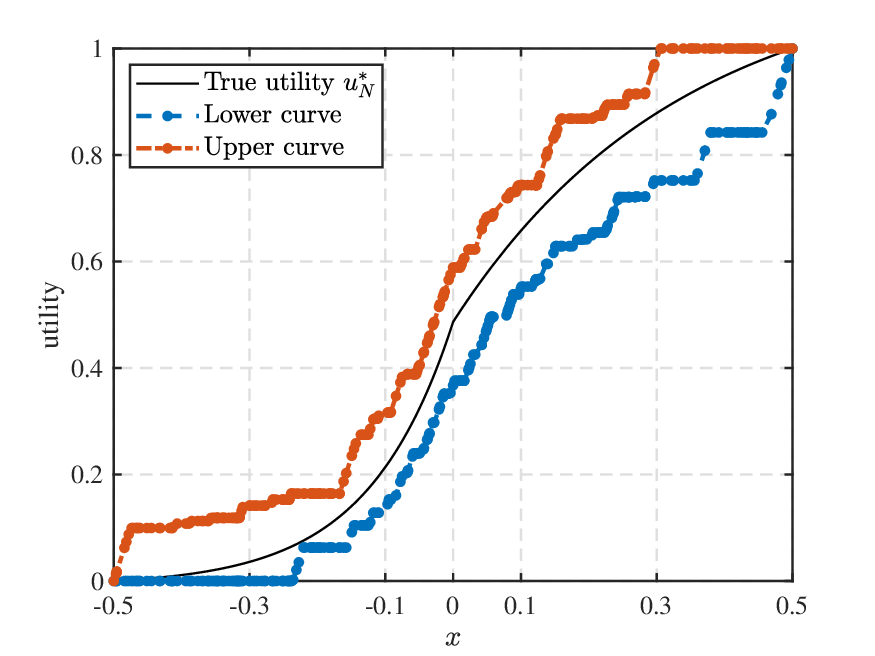}
\text{\hspace{0.5cm} \footnotesize \bf{(A2)}}
\end{minipage}
\hfill
\begin{minipage}[t]{0.24\linewidth}
\centering
\includegraphics[width=1.7in]{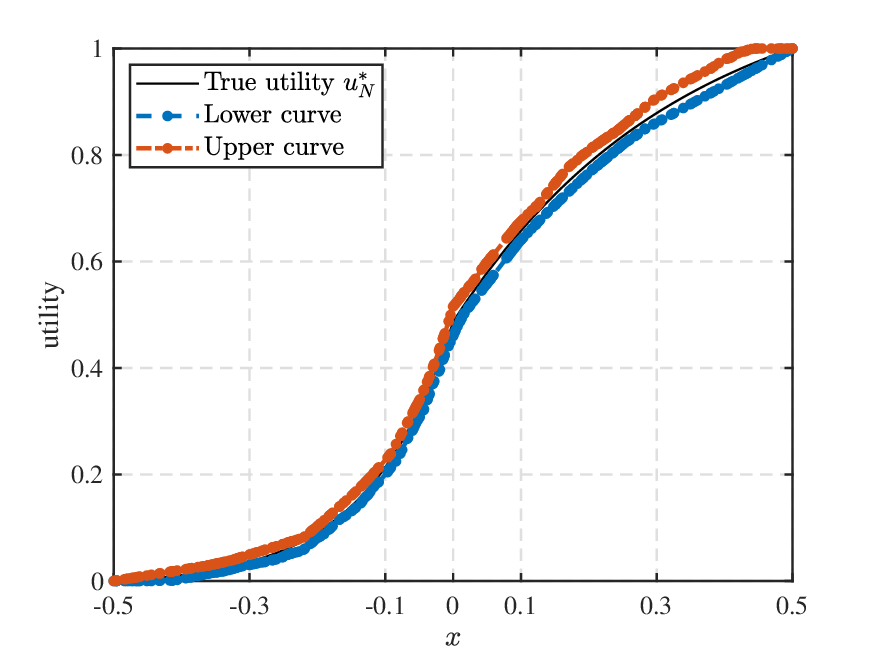}
\text{\hspace{0.5cm} \footnotesize \bf{(A3)}}
\end{minipage}
\hfill
\begin{minipage}[t]{0.24\linewidth}
\centering
\includegraphics[width=1.7in]{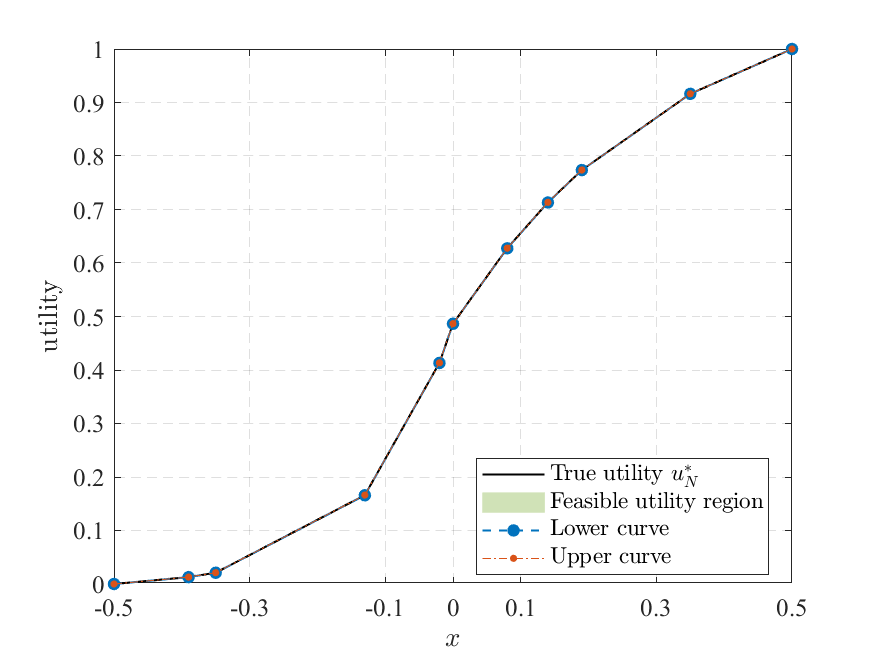}
\text{\hspace{0.5cm} \footnotesize \bf{(B)}}
\end{minipage}
\caption{\footnotesize{\textbf{(A)}
RRUS scheme for preference elicitation in the S-shaped utility setting with the range of utility functions 
$\mathcal{U}_{N}^{M}$ for $M=100$.  
\textbf{(B)}  Polyhedral method for preference elicitation in the S-shaped utility setting with the range of utility functions 
$\mathcal{U}_{M}$ for $M=100$ without prior information.
The labels, graphs and lines have the same interpretation as in Figure~\ref{fig:RRUS_M10}.
}}
\label{fig:RRUS_M100}
\end{figure}

\end{example}

\section{Proof of Proposition~\ref{Prop:CB}}

\proof{Proof of Proposition~\ref{Prop:CB}.}
Analogous to the arguments in the proof of Proposition~\ref{prop:pS_2},
we have ${\bm c}^m=({\bm c}_1^m,{\bm c}_2^m)$ with ${\bm R}{\bm c}_1^m>0$
and  ${\bm R}{\bm  c}_2^m>0$.
By the assumption that $\widetilde{\bm \sigma}_2^{1m}$ and $\widetilde{\bm v}_+^{2m}$ have no zero components, it follows that $p^{m,i}$ and $r_{2+}^{m,i}$ satisfy  
\begin{equation*}
\left[\begin{matrix}
{\bm R}{\bm c}_1^{m}\\
{\bm R}{\bm c}_2^{m}
\end{matrix}\right]^\top
\left[\begin{matrix}
0 \\
{\bm h}^{{\cal S}_2}(p^{m,i})
\end{matrix}\right]= D_i
\quad \inmat{and}\quad
\left[\begin{matrix}
{\bm R}{\bm c}_1^{m}\\
{\bm R}{\bm c}_2^{m}
\end{matrix}\right]^\top
\left[\begin{matrix}
{\bm h}^{{\cal X}_+^{m-1}}(r_{2+}^{m,i})\\
0
\end{matrix}\right] =D_i
\end{equation*}
for $i\in [N_2-2]$.
Then $l_{\rm cut,+}^{m,m-1}({\bm c}^m)=h_m(({\bm R}{\bm c}_2^{m})^\top{\bm h}^{{\cal S}_2}(p^{m,i}) -({\bm R}{\bm c}_1^{m})^\top {\bm h}^{{\cal X}_+^{m-1}}(r_{2+}^{m,i}))
=h_m(D_i-D_i)=0$ for $i\in [N_2-2]$.
Hence the cut hyperplane $l_{\rm cut,+}^{m,m-1}(\cdot)$ passes through the analytic center ${\bm c}^m$.
The proof is complete.
\hfill
\Box

\section{Proof of Proposition~\ref{prop:decrease_AC}}
\proof{Proof of Proposition~\ref{prop:decrease_AC}.}
Recall that the analytic center of the polyhedron that is defined by a set of linear inequalities $\tilde{\bm a}_i^\top {\bm x} \leq \tilde{b}_i,~~i\in [I]$
is defined as an optimal solution of the problem:
\begin{equation*}
\min_{{\bm x}}\; -\sum_{i=1}^I \log(\tilde{b}_i-\tilde{\bm a}_i^\top {\bm x}).
\end{equation*}
If we assume that $\|\tilde{\bm a}_i\| = 1$, for $i\in [I]$, then the slack 
\begin{equation*}
\tilde{b}_i-\tilde{\bm a}_i^\top {\bm x} =\inmat{the distance of}~~{\bm x} \inmat{~to~the~ hyperplane} ~~{\cal H}_i = \{{\bm x}: \tilde{\bm a}_i^\top {\bm x} = \tilde{b}_i\}.
\end{equation*}
The analytic center is the point that maximizes the product of distances to the defining hyperplanes $\tilde{\bm a}_i^\top {\bm x} = \tilde{b}_i$, for $i\in [I]$ (see \citep[p.~420]{BoV04}).
For the case of ${\Delta}_+^m$, to ease the exposition, we set the updated initial set $\Delta_+^0=\{(\widetilde{\bm v}_+,\widetilde{\bm \sigma}_2)\in \R^{J_++N_2-4}: {\bm A}(\widetilde{\bm v}_+,\widetilde{\bm \sigma}_2) \leq {\bm b}\}$, where ${\bm A}\in \R^{(2J_++N_2-4)\times (J_++N_2-4)}$ and ${\bm b}\in \R^{2J_++N_2-4}$.
The analytic center ${\bm c}_+^{m}$ of polyhedron ${\Delta}_+^{m-1}$ is the optimal solution of the following problem w.r.t.~$(\widetilde{\bm v}_+,\widetilde{\bm \sigma}_2)$:
\begin{equation*}
\begin{aligned}
\max\limits_{\widetilde{\bm v}_+,\widetilde{\bm \sigma}_2,{\bm s}}~~& \sum_{i=1}^{2J_+^{m-1}+N_2-4+m-1} \log(s_i) \\
{\rm s.t.}\quad & {\bm A}(\widetilde{\bm v}_+,\widetilde{\bm \sigma}_2)+{\bm s}_{[2J_+^{m-1}+N_2-4]}={\bm b}, \; {\bm s}\geq 0,\\
& h_l\big(-({\bm R}\widetilde{\bm \sigma}_2)^\top{\bm h}^{{\cal S}_2}(p^l)- ({\bm R}\widetilde{\bm v}_+)^\top({\bm h}^{{\cal X}_+^{m-1}}(r_{2+}^l)-{\bm e})\big)+s_{2J_+^{m-1}+N_2-4+l}=0,  \; l\in[m-1], 
\end{aligned}
\end{equation*}
where ${\bm s}_{[J_+^{m-1}+N_2-4]}=(s_1,\cdots,s_{J_+^{m-1}+N_2-4})$ and ${\bm s}:=(s_1,\cdots,s_{J_+^{m-1}+N_2-4+m-1})\in \R^{J_+^{m-1}+N_2-4+m-1}$.
For the ambiguity set ${\Delta}_+^{m-2}$, the analytic center, denoted by ${\bm c}_+^{m-1}$, maximizes 
\begin{equation*}
\begin{aligned}
\vt_{m-1}(\widetilde{\bm v}_+,\widetilde{\bm \sigma}_2):= & \sum_{j=1}^{2J_+^{m-2}+N_2-4}{\rm log}({b}_j-{\bm a}_j^\top(\widetilde{\bm v}_+,\widetilde{\bm \sigma}_2))\\
& +\sum_{l=1}^{m-2}{\rm log}\Big(-h_l\big(({\bm R}\widetilde{\bm \sigma}_2)^\top{\bm h}^{{\cal S}_2}(p^l)- ({\bm R}\widetilde{\bm v}_+)^\top{\bm h}^{{\cal X}_+^{m-2}}(r_{2+}^l)\big)\Big).
\end{aligned}
\end{equation*}
For the ambiguity set ${\Delta}_{+}^{m-1}$ with ${\cal X}_+^{m-1}={\cal X}_+^{m-2}$, the analytic center ${\bm c}_+^{m}$ maximizes 
\begin{equation*}
\begin{aligned}
\vt_{m}(\widetilde{\bm v}_+,\widetilde{\bm \sigma}_2)
=&\sum_{j=1}^{2J_+^{m-1}+N_2-4}{\rm log}({b}_j-{\bm a}_j^\top(\widetilde{\bm v}_+,\widetilde{\bm \sigma}_2))\\
& +\sum_{l=1}^{m-1}{\rm log}\Big(-h_l\big(({\bm R}\widetilde{\bm \sigma}_2)^\top{\bm h}^{{\cal S}_2}(p^l)- ({\bm R}\widetilde{\bm v}_+)^\top{\bm h}^{{\cal X}_+^{m-1}}(r_{2+}^l)\big)\Big)\\
:=&\vt_{m-1}(\widetilde{\bm v}_+,\widetilde{\bm \sigma}_2)+{\rm log}\Big(-h_{m-1}\big(({\bm R}\widetilde{\bm \sigma}_2)^\top{\bm h}^{{\cal S}_2}(p^{m-1})- ({\bm R}\widetilde{\bm v}_+)^\top{\bm h}^{{\cal X}_+^{m-1}}(r_{2+}^{m-1})\big)\Big).
\end{aligned}
\end{equation*}
Note that ${\bm c}_+^m\in {\Delta}_+^{m-1}$, ${\bm c}_+^{m-1}\in \Delta_+^{m-2}$ make the objective functions $\vt_m(\widetilde{\bm v}_+,\widetilde{\bm \sigma}_2)$ and $\vt_{m-1}(\widetilde{\bm v}_+,\widetilde{\bm \sigma}_2)$
attain the optimal values, respectively.

By the fact that $\delta_{m-1}:=-h_{m-1}\big(({\bm R}\widetilde{\bm \sigma}_2)^\top{\bm h}^{{\cal S}_2}(p^{m-1})- ({\bm R}\widetilde{\bm v}_+)^\top{\bm h}^{{\cal X}_+^{m-1}}(r_{2+}^{m-1})\big)\in (0,1)$, 
we have ${\rm log}(\delta_{m-1})<0$ and subsequently,
\begin{equation*}
\vt_{m}({\bm c}^{m}_+)< \vt_{m-1}({\bm c}^{m}_+)\leq \vt_{m-1}({\bm c}^{m-1}_+).
\end{equation*}
Note that the above inequality still holds when the vector ${\bm a}_i$ for $i\in [2J_+^{m-1}+N_2-4]$ of matrix ${\bm A}$ and the coefficient vector of $-h_l\big(({\bm R}\widetilde{\bm \sigma}_2)^\top{\bm h}^{{\cal S}_2}(p^l)- ({\bm R}\widetilde{\bm v}_+)^\top{\bm h}^{{\cal X}_+^{m-2}}(r_{2+}^l)\big)$, for $l\in [m-1]$, are normalized.
By construction, $\Delta_+^{m-1}\subseteq\Delta_+^{m-2}$,
so the feasible region shrinks as the number of queries increases.
The proof is complete.
\hfill
\Box

\section{An Example of the Modified Polyhedral Method}
\label{SC:exa-polyhedral}
In this section, we present an example of the modified polyhedral method introduced in Section~\ref{sec:modified_PM}.
\begin{example}[Illustration of the Polyhedral Method for $\widehat{g}_2$]
\label{exa:polyhdedral}
We use a simple example to explain the steps in the polyhedral method, where the true utility function is $v_+^*(x)=(1-e^{-3x})/3$ for $x>0$ defined over $[0,0.5]$ and the true distortion function is $g_2^*(t)=t^{\gamma}/(t^{\gamma}+(1-t)^{\gamma})$ for $\gamma=0.6$ defined over $[0,1]$.
The piecewise linear value function has three breakpoints, 
i.e., $J_+=3$ including the two ends of the interval and the one within it is $0.25$, that is, $x_1=0$, $x_2=0.25$, $x_3=0.5$.
The piecewise linear distortion function has four breakpoints,
i.e., $N_2=4$ including the two ends of $[0,1]$ and the two interior breakpoints are $0.3$ and $0.6$, that is, $t_1=0$, $t_2=0.3$, $t_3=0.6$, $t_4=1$.
Consider the normalized value function on $[0,0.5]$.
The true increment vectors of the value function and the distortion function are given by ${\bm v}^*_+=(0.679,0.321)$, ${\bm \sigma}^*_2=(0.375,0.185,0.439)$,
and $(\widetilde{\bm v}_+^*,\widetilde{\bm \sigma}_2^*)=(v_{+,1},\sigma_{2,1},\sigma_{2,2})=(0.679,0.375,0.185)$ (these are all column vectors: we write them in row format in running text).
We impose the concavity of the value function $v_+$. To ease the exposition, we do not update the set of breakpoints of the value function.
The parameter $D_i$ in \eqref{eq:gene_m_Q-1} is set to a randomly generated constant $D$.
The $l$-th cutting hyperplane is constructed as $\tilde{\bm \sigma}_2^\top{\bm h}^{{\cal S}_2}_{[N_2-1]}(p^l) -\widetilde{\bm v}_+^\top{\bm h}^{{\cal X}_+}_{[J_+-1]}(r_{2+}^l)$, where ${\cal S}_2=\{t_1,t_2,t_3,t_4\}$, ${\cal X}_+=\{x_1,x_2,x_3\}$, ${\bm h}^{{\cal S}_2}_{[N_2-1]}(p^l)$ and ${\bm h}^{{\cal X}_+}_{[J_+-1]}(r_{2+}^l)$ denote the truncated vectors consisting of the first $N_2-1$ and $J_+-1$ entries.
We focus on the iteration of the polyhedral method.
\begin{equation*}
\Delta_+^0=\{(v_{+,1},\sigma_{2,1},\sigma_{2,2})\in \R^{1+2}:
v_{+,1}\leq 1,\; v_{+,1}\geq 0,\; \sigma_{2,1}+\sigma_{2,2}\leq 1, \; \sigma_{2,i}\geq0,i=1,2\}.
\end{equation*}

\underline{Iteration 1}. {\underline{\bf Step 1}:
The analytic center of {$\Delta_+^0$} is ${\bm c}^1=(0.789,0.333,0.333)$.
\underline{\bf Step 2}:
$(v^{11}_{+,1},\sigma^{11}_{2,1},\sigma^{11}_{2,2})=(0.789,0.667,0)$, $(v^{21}_{+,1},\sigma^{21}_{2,1},\sigma^{21}_{2,2})=(0.789,0,0.667)$.}
\underline{\bf Step 3}: $p^1=0.9496$, $r_{2+}^1=0.4503$ ($D=0.958$).
Then the cutting hyperplane is $-v_{+,1}+\sigma_{2,1}+\sigma_{2,2}=0$, $h_1=1$. The $1$st polyhedron is
\begin{equation*}
{\Delta_+^1=\Delta_+^0}\cap \{(v_{+,1},\sigma_{2,1},\sigma_{2,2})\in \R^{1+2}: -v_{+,1}+\sigma_{2,1}+\sigma_{2,2}\leq 0\}.
\end{equation*}

\underline{Iteration 2}.
\underline{\bf Step 1}: The analytic center of {$\Delta_+^1$} is ${\bm c}^2=(0.848,0.228,0.228)$.
\underline{\bf Step 2}: $(v^{12}_{+,1},\sigma^{12}_{2,1},\sigma^{12}_{2,2})=(0.848, 0.456,0)$, $(v^{22}_{+,1},\sigma^{22}_{2,1},\sigma^{22}_{2,2})=(0.848,0,0.456)$.
\underline{\bf Step 3}:
$p^2=0.6030$, $r_{2+}^2=0.1355$ ($D=0.46$).
Then the cutting hyperplane is $-0.542v_{+,1}+\sigma_{2,1}+\sigma_{2,2}=0$, $h_2=-1$. The $2$nd polyhedron is
\begin{equation*}
{\Delta_+^2=\Delta_+^1}\cap \{(v_{+,1},\sigma_{2,1},\sigma_{2,2})\in \R^{1+2}: -0.542v_{+,1}+\sigma_{2,1}+\sigma_{2,2}\geq 0\}.
\end{equation*}

\underline{Iteration 3}.
\underline{\bf Step 1}:
The analytic center of {$\Delta_+^2$} is {${\bm c}^3=(0.847,0.328,0.328)$.}
\underline{\bf Step 2}:  $(v^{13}_{+,1},
\sigma^{13}_{2,1},
\sigma^{13}_{2,2})
=(0.847,0.656,0)$, 
$(v^{23}_{+,1},
\sigma^{23}_{2,1},
\sigma^{23}_{2,2})
=(0.847,
0,
0.656)$.
\underline{\bf Step 3}: $p^3=0.6049$, $r_{2+}^3= 0.1948$ ($D=0.66$).
Then the cutting hyperplane is $-0.779v_{+,1}+\sigma_{2,1}+\sigma_{2,2}=0$, $h_3=1$. The $3$rd polyhedron is
\begin{equation*}
{\Delta_+^3=\Delta_+^2}\cap \{(v_{+,1},\sigma_{2,1},\sigma_{2,2})\in \R^{1+2}: -0.779v_{+,1}+\sigma_{2,1}+\sigma_{2,2}\leq 0\}.
\end{equation*}

\underline{Iteration 4}.
\underline{\bf Step 1}:
The analytic center of {$\Delta_+^3$} is ${\bm c}^4=(0.875,0.284,0.284)$. 
\underline{\bf Step 2}: $(v^{14}_{+,1},\sigma^{14}_{2,1},\sigma^{14}_{2,2})=(0.875,0.567,0)$, $(v^{24}_{+,1},\sigma^{24}_{2,1},\sigma^{24}_{2,2})=(0.875,0,0.567)$.
\underline{\bf Step 3}:
$p^4=0.1903$, $r_{2+}^4=0.0514$ ($D=0.18$).
Then the cutting hyperplane is $-0.205v_{+,1}+0.634\sigma_{2,1}=0$, $h_4=-1$. The $4$th polyhedron is
\begin{equation*}
{\Delta_+^4=\Delta_+^3}\cap \{(v_{+,1},\sigma_{2,1},\sigma_{2,2}) \in \R^{1+2}: -0.206v_{+,1}+0.634\sigma_{2,1}\geq 0\}.
\end{equation*}
    
The cutting hyperplanes and the polyhedra $\Delta_+^i$ for $i=1,2,3,4$ are presented in Figure~\ref{fig:polyhedral}.

\begin{figure}[!ht]
\centering
\begin{minipage}[t]{0.45\linewidth}
\centering
\includegraphics[width=2.8in]{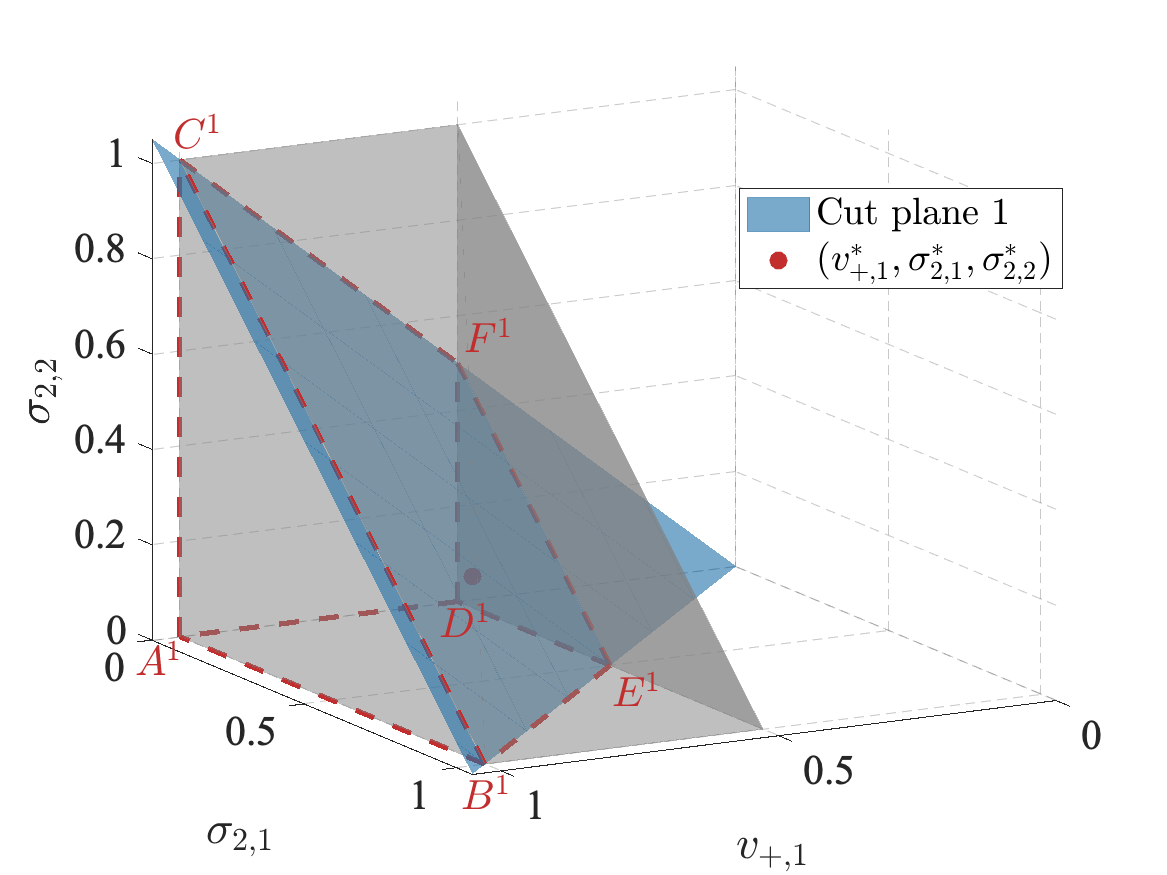}
\text{\hspace{0.5cm}\footnotesize \textbf{(a)}}
\end{minipage}
\begin{minipage}[t]{0.45\linewidth}
\centering
\includegraphics[width=2.8in]{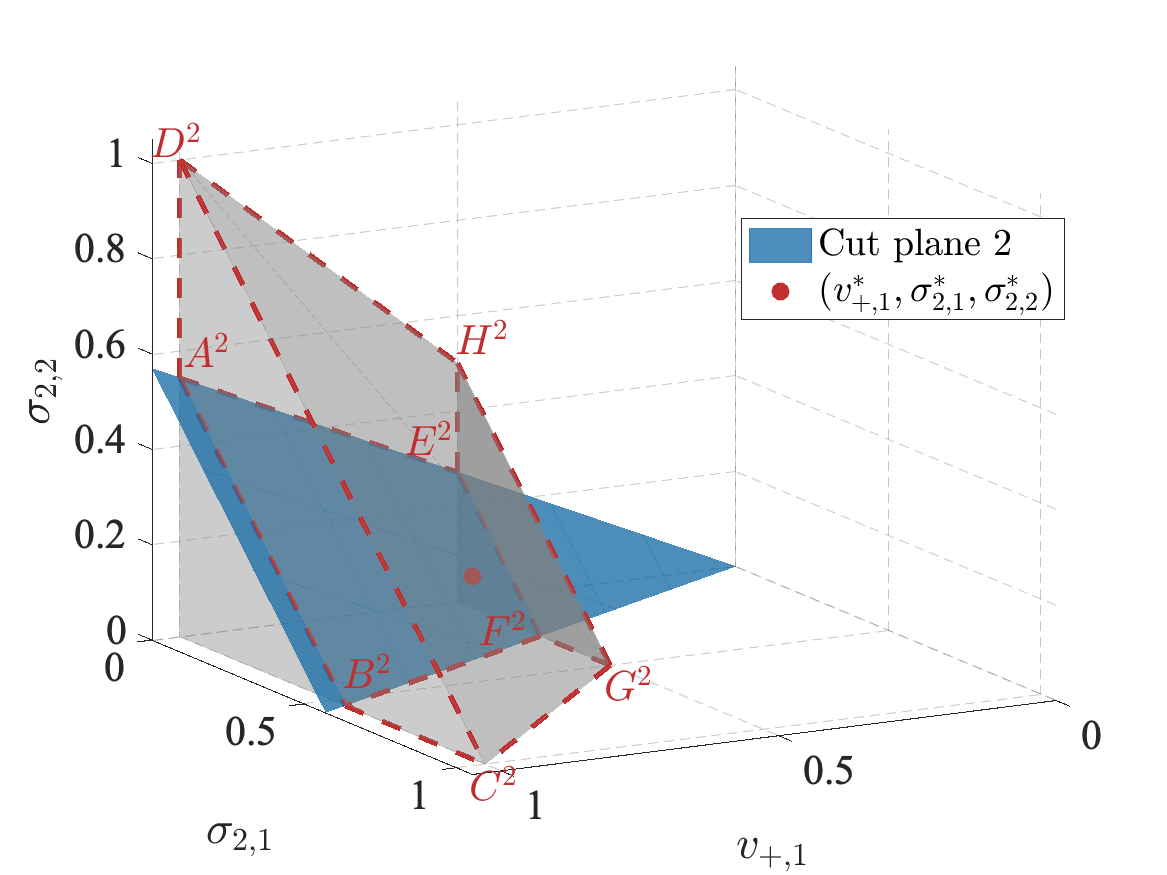}
\text{\hspace{0.5cm}\footnotesize \textbf{(b)}}
\end{minipage}
\begin{minipage}[t]{0.45\linewidth}
\centering
\includegraphics[width=2.8in]{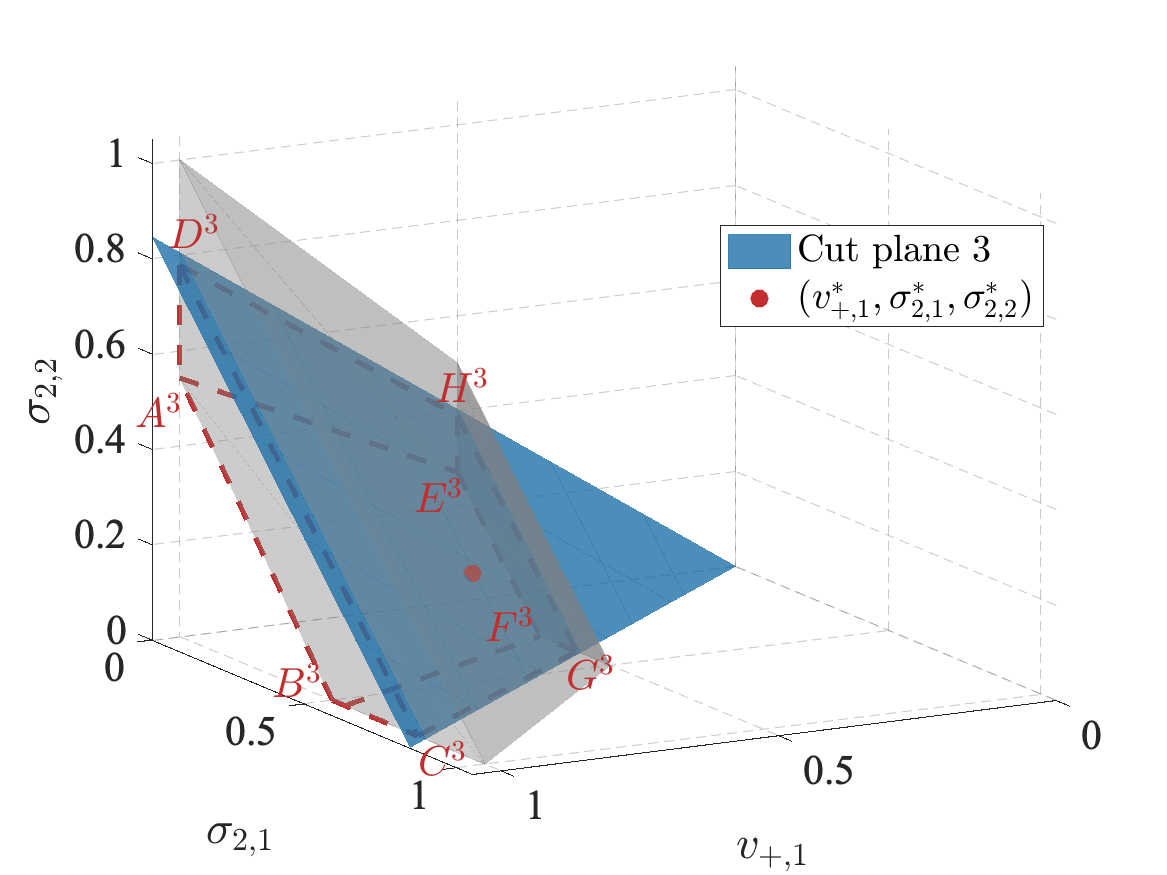}
\text{\hspace{0.5cm}\footnotesize \textbf{(c)}}
\end{minipage}
\begin{minipage}[t]{0.45\linewidth}
\centering
\includegraphics[width=2.8in]{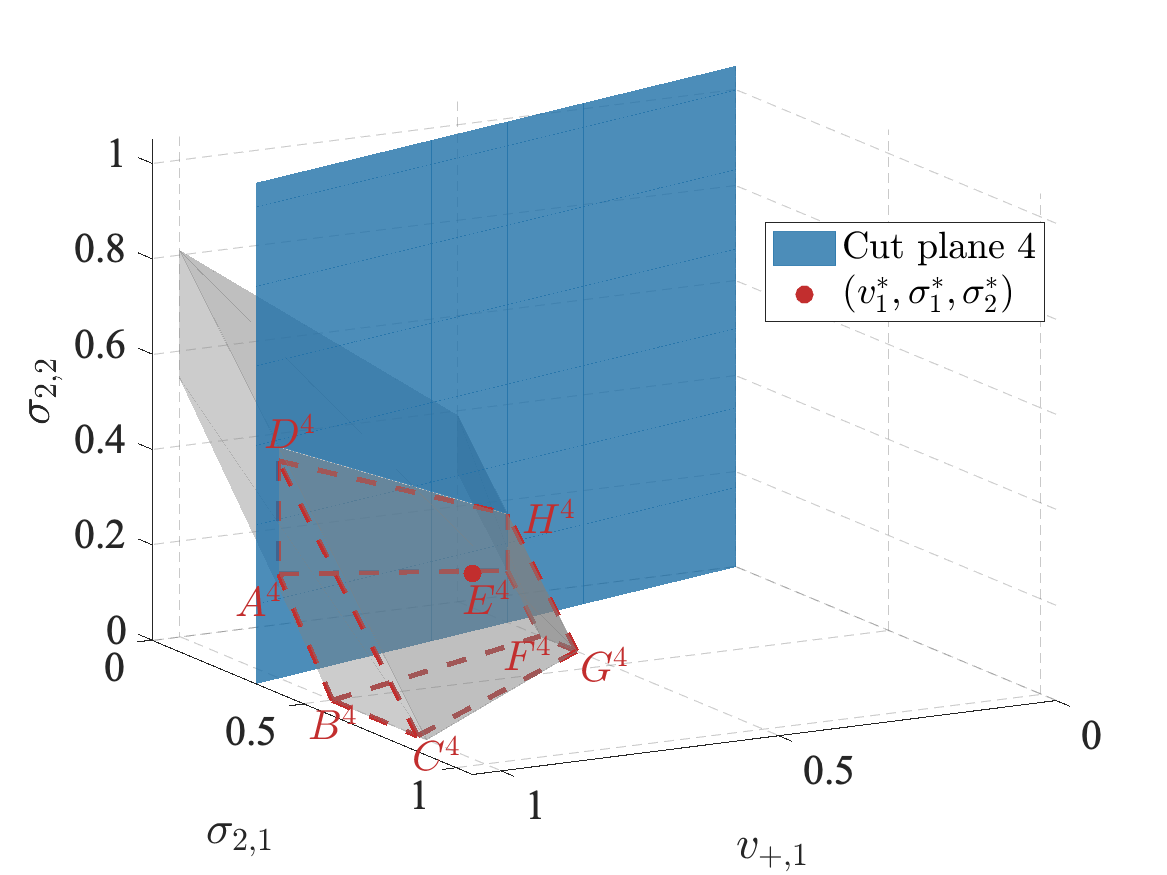}
\text{\hspace{0.5cm}\footnotesize \textbf{(d)}}
\end{minipage}
 \caption{\footnotesize{Visualization of the cut planes and the resulting polyhedron $\Delta_+^{m}$
for $m=1,2,3,4$.
The red point is $(v_{+,1}^*,\sigma_{2,1}^*,\sigma_{2,2}^*)=(0.679,0.375, 0.185)$.
\textbf{(a)} The grey polyhedron represents $\Delta_+^0$ and the one framed by the red dashed lines is $\Delta_+^1$ ($A^1$-$B^1$-$C^1$-$D^1$-$E^1$-$F^1$) after the first cut. 
The cut plane is $-v_{+,1}+\sigma_{2,1}+\sigma_{2,2}=0$. $h_1=1$.
\textbf{(b)} The resulting polyhedron $\Delta_+^2$: $A^2$-$B^2$-$C^2$-$D^2$-$E^2$-$F^2$-$G^2$-$H^2$ corresponding to $h_2=-1$.
\textbf{(c)} The resulting polyhedron $\Delta_+^3$: $A^3$-$B^3$-$C^3$-$D^3$-$E^3$-$F^3$-$G^3$-$H^3$ corresponding to $h_3=1$.
\textbf{(d)} The resulting polyhedron $\Delta_+^4$: $A^4$-$B^4$-$C^4$-$D^4$-$E^4$-$F^4$-$G^4$-$H^4$ corresponding to $h_4=-1$.}}
\label{fig:polyhedral}
\end{figure}
\end{example}

\section{Piecewise McCormick Relaxation: Implementation Details}
\label{sec:Piecewise_Mc}
In this section, we present a piecewise linearization of the bilinear terms
$v_{+,i}\sigma_{2,j}$ and $v_{-,i}\sigma_{1,j}$ appearing in the objective function~\eqref{eq:objective-2}.

\noindent{{\bf Piecewise Linearization of Bilinear Term $v_{+,i}\sigma_{2,j}$}}.
For each $j\in [N_2-1]$, 
we divide the domain $[\sigma_{2,j}^L,\sigma_{2,j}^U]$ of $\sigma_{2,j}$ into $N$ equally sized, mutually exclusive sub-regions. 
Next, we calculate the lower and upper bounds of $\sigma_{2,j}$ in each sub-region:
\begin{equation}
\label{eq:sigma_2}
\sigma_{2,j,n}^L:=\sigma_{2,j}^L+\frac{(\sigma_{2,j}^U-\sigma_{2,j}^L)(n-1)}{N}
\quad 
\inmat{and}
\quad
\sigma_{2,j,n}^U:=\sigma_{2,j}^L+\frac{(\sigma_{2,j}^U-\sigma_{2,j}^L)n}{N}
,\;n\in [N].
\end{equation}
Let ${\bm \sigma}_{2,j}^L:=({\sigma}_{2,j,1}^L,\cdots,{\sigma}_{2,j,N}^L)$
and ${\bm \sigma}_{2,j}^U:=({\sigma}_{2,j,1}^U,\cdots,{\sigma}_{2,j,N}^U)$, for $j\in [N_2-1]$.
In the piecewise McCormick relaxation, we define $w_{ij}^+:=v_{+,i}\sigma_{2,j}$,
for $i\in [J_{+}^{M_2}-1]$ and $j\in [N_2-1]$.
By introducing binary variables ${\bm y}_j^{g_2}=(y_{j,1}^{g_2},\cdots,y_{j,N}^{g_2})$ to represent that $\sigma_{2,j}$ lies within sub-region $[\sigma_{2,j,n}^L,\sigma_{2,j,n}^U]$ or not, we have
\begin{subequations}
\label{eq:wMcC}
\begin{align}
& w_{ij}^+ \geq v_{+,i}^L\sigma_{2,j}+(({\bm \sigma}_{2,j}^L)^\top{\bm y}_{j}^{g_2}) v_{+,i}-
v_{+,i}^L(({\bm \sigma}_{2,j}^L)^\top{\bm y}_{j}^{g_2}), \label{eq:wMcC-a}\\
& w_{ij}^+ \geq v_{+,i}^U\sigma_{2,j}+(({\bm \sigma}_{2,j}^U)^\top{\bm y}_{j}^{g_2}) v_{+,i}-
v_{+,i}^U (({\bm \sigma}_{2,j}^U)^\top{\bm y}_{j}^{g_2}), \\
& w_{ij}^+\leq v_{+,i}^U\sigma_{2,j}+(({\bm \sigma}_{2,j}^L)^\top{\bm y}_{j}^{g_2}) v_{+,i}-
v_{+,i}^U
(({\bm \sigma}_{2,j}^L)^\top{\bm y}_{j}^{g_2}), \\
& w_{ij}^+\leq v_{+,i}^L \sigma_{2,j}+(({\bm \sigma}_{2,j}^U)^\top{\bm y}_{j}^{g_2}) v_{+,i}-
v_{+,i}^L
(({\bm \sigma}_{2,j}^U)^\top{\bm y}_{j}^{g_2}) \label{eq:wMcC-b},\\
& {\bm e}^\top{\bm y}_{j}^{g_2}=1, \label{eq:wMcC-c}
\end{align}
\end{subequations}
where ${\bm y}_{j}^{g_2},{\bm \sigma}_{2,j}^L,{\bm \sigma}_{2,j}^U\in \R^{N}$ for $j\in [N_2-1]$, and constraint \eqref{eq:wMcC-c} ensures that $\sigma_{2,j}$ lies within a single sub-region.
There are bilinear terms $y_{j,n}^{g_2}v_{+,i}$ in constraints \eqref{eq:wMcC-a}--\eqref{eq:wMcC-b}.
Furthermore, we apply the McCormick relaxation to each bilinear term $y_{j,n}^{g_2}v_{+,i}$. 
Specifically, let $\widehat{V}_{j,n,i}^+:=y_{j,n}^{g_2} v_{+,i}$,
where $y_{j,n}^{g_2}\in \{0, 1\}$, $\sigma_{2,j,n}^L\leq \sigma_{2,j} \leq \sigma_{2,j,n}^U$, for $i\in [J_{+}^{M_2}-1]$, $j\in [N_2-1]$, $n\in [N]$.
Then 
\begin{equation}
\label{eq:GVY}
 \widehat{V}_{j,n,i}^+ \geq y_{j,n}^{g_2} v_{+,i}^L, \;\;\;
\widehat{V}_{j,n,i}^+\geq v_{+,i}+y_{j,n}^{g_2}v_{+,i}^U-v_{+,i}^U, \;\;\;
\widehat{V}_{j,n,i}^+ \leq v_{+,i}+y_{j,n}^{g_2}v_{+,i}^L-v_{+,i}^L \;\;\;\inmat{and}\;\;\;
\widehat{V}_{j,n,i}^+ \leq y_{j,n}^{g_2} v_{+,i}^U,
\end{equation}
for $i\in [J_{+}^{M_2}-1]$, $j\in [N_2-1]$, $n\in [N]$.
Consequently, by introducing $\sigma_{2,j,n}^L$ and $\sigma_{2,j,n}^U$ from \eqref{eq:sigma_2} along with the definition $\widehat{V}_{j,n,i}^+ =y_{j,n}^{g_2} v_{+,i}$,
we can reformulate the inequality system \eqref{eq:wMcC} as 
\begin{subequations}
\label{eq:wMcC-2}
\begin{align}
& w_{ij}^+ \geq v_{+,i}^L \sigma_{2,j} +\sum_{n=1}^N\left(\sigma_{2,j,n}^L
\widehat{V}_{j,n,i}^+\right)- v_{+,i}^L(({\bm \sigma}_{2,j}^L)^\top{\bm y}_{j}^{g_2}), \\
& w_{ij}^+ \geq v_{+,i}^U \sigma_{2,j}
+\sum_{n=1}^N\left(\sigma_{2,j,n}^U\widehat{V}_{j,n,i}^+\right)
- v_{+,i}^U (({\bm \sigma}_{2,j}^U)^\top{\bm y}_{j}^{g_2}), \\
& w_{ij}^+\leq  v_{+,i}^U \sigma_{2,j}
+\sum_{n=1}^N\left(\sigma_{2,j,n}^L\widehat{V}_{j,n,i}^+\right)
- v_{+,i}^U (({\bm \sigma}_{2,j}^L)^\top{\bm y}_{j}^{g_2}), \\
& w_{ij}^+\leq v_{+,i}^L \sigma_{2,j} +\sum_{n=1}^N\left(\sigma_{2,j,n}^U
\widehat{V}_{j,n,i}^+\right) -v_{+,i}^L (({\bm \sigma}_{2,j}^U)^\top{\bm y}_{j}^{g_2}),\\
& {\bm e}^\top{\bm y}_{j}^{g_2}=1,
\;{\bm y}_{j}^{g_2}\in \{0,1\}^{N},\; \inmat{constraints}\;\eqref{eq:GVY}.
\end{align}
\end{subequations}
In this way, we linearize the bilinear term $v_{+,i}\sigma_{2,j}$ in \eqref{eq:wMcC-2} by using piecewise McCormick relaxation.

\noindent{\bf Piecewise Linearization of the Bilinear Term $v_{-,i}\sigma_{1,j}$}.
Analogous to the discussions about $v_{+,i}\sigma_{2,j}$, 
we divide the range of each $\sigma_{1,j}$ into $N$ sub-regions and then apply the McCormick relaxation to each sub-region.
Specifically, for $\sigma_{1,j}$, $j\in {[N_1-1]}$, define $\sigma_{1,j,n}^L$ and $\sigma_{1,j,n}^U$ accordingly.
Let $w_{ij}^-:= v_{-,i}\sigma_{1,j}$ for $i\in [J_-^{M_1}-1]$ and $j\in [N_1-1]$ and introduce binary variables ${\bm y}_j^{g_1}:=(y_{j,1}^{g_1},\ldots,y_{j,N}^{g_1})$, for $j\in [N_1-1]$ to represent $\sigma_{1,j}$ which lies within sub-region $[\sigma_{1,j,n}^L,\sigma_{1,j,n}^U]$ or not.
Next, we apply the McCormick relaxation to the bilinear term $y_{j,n}^{g_1} v_{-,i}$ over each sub-region.
Denote $\widehat{V}_{j,n,i}^-:=y_{j,n}^{g_1} v_{-,i}$ for $i\in [J_-^{M_1}-1]$ and $j\in [N_1-1]$.
Analogous to system \eqref{eq:wMcC-2}, 
we can develop constraints for $({\bm w}^-,{\bm \sigma}_{1},\widehat{\bm V}^-, {\bm y}^{g_1})$:
\begin{subequations}
\label{eq:wMcC-2-minus}
\begin{align}
& w_{ij}^- \geq v_{-,i}^L \sigma_{1,j} +\sum_{n=1}^N\left(\sigma_{1,j,n}^L
\widehat{V}_{j,n,i}^-\right) - v_{-,i}^L(({\bm \sigma}_{1,j}^L)^\top{\bm y}_{j}^{g_1}), \\
& w_{ij}^- \geq v_{-,i}^U \sigma_{1,j} +\sum_{n=1}^N\left(\sigma_{1,j,n}^U
\widehat{V}_{j,n,i}^-\right) - v_{-,i}^U (({\bm \sigma}_{1,j}^U)^\top{\bm y}_{j}^{g_1}), \\
& w_{ij}^-\leq  v_{-,i}^U \sigma_{1,j} +\sum_{n=1}^N\left(\sigma_{1,j,n}^L
\widehat{V}_{j,n,i}^-\right) - v_{-,i}^U (({\bm \sigma}_{1,j}^L)^\top{\bm y}_{j}^{g_1}),\;\;  \\
& w_{ij}^-\leq v_{-,i}^L \sigma_{1,j}+\sum_{n=1}^N\left(\sigma_{1,j,n}^U
\widehat{V}_{j,n,i}^-\right) -v_{-,i}^L (({\bm \sigma}_{1,j}^U)^\top{\bm y}_{j}^{g_1}),\\
& {\bm e}^\top{\bm y}_{j}^{g_1}=1, \; {\bm y}_{j}^{g_1}\in \{0,1\}^{N},\\
& \widehat{V}_{j,n,i}^- \geq y_{j,n}^{g_1} v_{-,i}^L, \;\;
\widehat{V}_{j,n,i}^-\geq v_{-,i}+y_{j,n}^{g_1}v_{-,i}^U-v_{-,i}^U,\;\;
\widehat{V}_{j,n,i}^- \leq v_{-,i}+y_{j,n}^{g_1}v_{-,i}^L-v_{-,i}^L,\;\;
\widehat{V}_{j,n,i}^- \leq y_{j,n}^{g_1} v_{-,i}^U,
\end{align}
\end{subequations}
where $\sigma_{1,j,n}^L=\sigma_{1,j}^L+\frac{(\sigma_{1,j}^U-\sigma_{1,j}^L)(n-1)}{N}$ and  $\sigma_{1,j,n}^U=\sigma_{1,j}^L+\frac{(\sigma_{1,j}^U-\sigma_{1,j}^L)n}{N}$ for $n\in [N]$.
In this way, we have linearized the bilinear term $v_{-,i}\sigma_{1,j}$ in \eqref{eq:wMcC-2-minus} by using piecewise McCormick relaxation.

\section{Remarks on the PRO Model \texorpdfstring{$\rho_{\mathbb{V}\mathbb{G}_1\mathbb{G}_2}(f({\bm z},{\bm \zeta}))$}{} in Section~\ref{sec:computation}}
\label{SC:extention_models}

We apply the proposed robust shortfall risk measure to a multivariate decision-making problem with a loss function
$$
\rho_{{\mathbb{V}\mathbb{G}_1\mathbb{G}_2}}(f({\bm z},{\bm \zeta}))
:= \inf\limits_{x\in \R} \left\{ x:  \sup_{(v,g_1,g_2)\in \mathbb{V}\mathbb{G}_1\mathbb{G}_2} \mathbb
E_{g_1{g_2}}[v(f({\bm z},{\bm \zeta})-x)]\leq v(0)\right\},
$$
where ${\bm z} \in \R^n$ is a decision vector in the feasible set $Z\subset \R^n$, ${\bm \zeta}:\Omega\to \R^n$ is a random vector, 
$f({\bm z},{\bm \zeta})={\bm z}^\top{\bm \zeta}$ represents a financial loss.
In this case, we cannot directly the bisection method to solve this problem because we have to tackle ${\bm z}$ and $x$ simultaneously and the constraint function is multivariate which is not necessarily monotonic w.r.t. $(x,{\bm z})$.
Let $\mathbb{V}_J\mathbb{G}_{1N_1}\mathbb{G}_{2N_2}$ be defined as in \eqref{eq:Question_VJ}.
We can use $\mathbb{V}_J\mathbb{G}_{1N_1}\mathbb{G}_{2N_2}$ to approximate the ambiguity set $\mathbb{V}\mathbb{G}_{1}\mathbb{G}_{2}$, 
and then consider the approximate preference robust problem: 
\begin{align}
\label{eq:approximate_port}
 \inf\limits_{ {\bm z}\in Z}\; \rho_{\mathbb{V}_J\mathbb{G}_{1N_1}\mathbb{G}_{2N_2}}(f({\bm z},{\bm \zeta}))
:= \inf\limits_{x\in \R, {\bm z}\in Z} \left\{x :\sup_{(v,g_1,g_2)\in \mathbb{V}_J\mathbb{G}_{1N_1}\mathbb{G}_{2N_2}} \mathbb E_{g_1{g_2}}[v(f({\bm z},{\bm \zeta})-x)]\leq v(0)\right\}. \qquad 
\end{align}
For any fixed $(x,{\bm z})$, we order the scenarios so that
$f({\bm z},{\bm\zeta}_1)\leq\cdots\leq
f({\bm z},{\bm\zeta}_K)$,
with the probabilities $p_k$ reordered accordingly,
and reformulate the inner maximization problem as follows:
\begin{subequations}
\label{eq:PRGSR_V_port}
\begin{align}
\sup_{{\bm v}_-,{\bm v}_+,{\bm \sigma}_1,{\bm \sigma}_2}  & \left\{
\sum_{k=1}^{K_0({\bm z},x)-1} \big(({\bm  v}_-)^\top{\bm h}^{{\cal X}_-^{M_1}}(f({\bm z},{\bm \zeta}_k)-x)-1 \big) {\bm \sigma}_1^\top ({\bm h}^{{\cal S}_1}(q_k)-{\bm h}^{{\cal S}_1}(q_{k-1})) \right. \nonumber \\
& \;\;\; \left. +
\widehat{r} \sum_{k=K_0({\bm z},x)}^{K} 
({\bm  v}_+)^\top{\bm h}^{{\cal X}_+^{M_2}}(f({\bm z},{\bm \zeta}_k)-x) {\bm \sigma}_2^\top({\bm h}^{{\cal S}_2}(1-q_{k-1})-{\bm h}^{{\cal S}_2}(1-q_{k})) \right\}  \\
 {\rm s.t.}\quad & \eqref{eq:bilinear-b}-\eqref{eq:bilinear-h},
\end{align}
\end{subequations}
where $q_{0}=0$, $q_k=\sum_{j=1}^{k}p_j$ for $k\in [K]$, we label ${\bm \zeta}_k$ satisfying $f({\bm z},{\bm \zeta_k})-x<0$ with index $k\in \{1,\cdots,K_0({\bm z},x)-1\}$, 
and for those satisfying $f({\bm z},{\bm \zeta_k})-x\geq 0$ by $k\in \{K_0({\bm z},x),\cdots,K\}$.
For any fixed $(x,{\bm z})$, problem \eqref{eq:PRGSR_V_port} is a bilinear program with linear constraints.
The piecewise McCormick relaxation described in
Section~\ref{sec:Piecewise_Mc} can be applied to approximate this problem.

\subsection{The Worst-Case Value Function and the Evaluation of PRGSR-CPT}
\label{sec:worst-case-v}
\begin{figure}[!ht]
\begin{minipage}[t]{0.3\linewidth}
\centering
\includegraphics[width=2.1in]{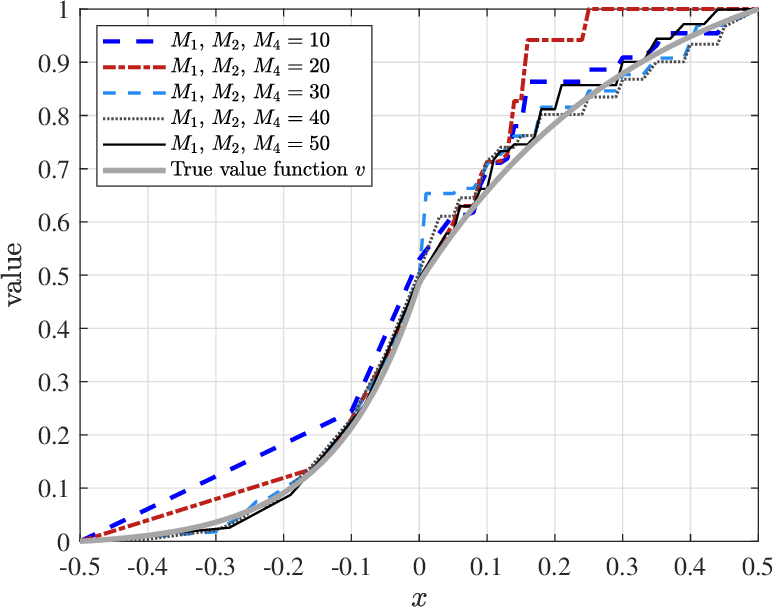}
\text{\hspace{0.9cm}\footnotesize \textbf{(a)}}
\end{minipage}
\hfill
\begin{minipage}[t]{0.3\linewidth}
\centering
\includegraphics[width=2.1in]{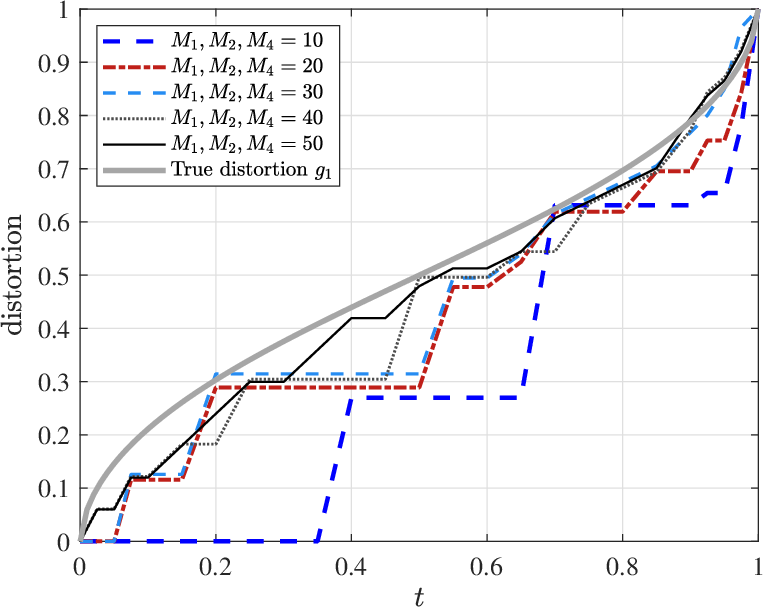}
\text{\hspace{0.9cm}\footnotesize \textbf{(b)}}
\end{minipage}
\hfill
\begin{minipage}[t]{0.3\linewidth}
\centering
\includegraphics[width=2.1in]{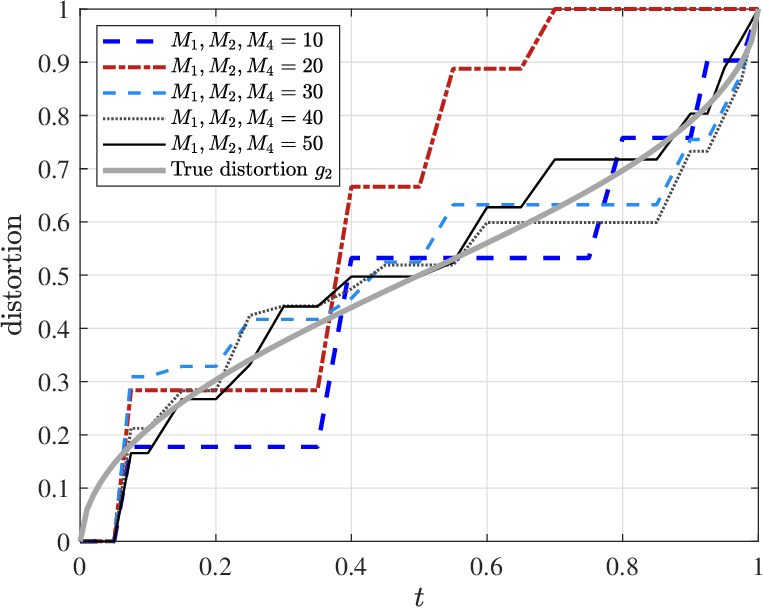}
\text{\hspace{0.9cm}\footnotesize   \textbf{(c)}}
\end{minipage}
\hfill
\caption{\footnotesize{\textbf{(a)} The worst-case value function $v$. 
\textbf{(b)} The worst-case distortion function $g_1$. 
\textbf{(c)} The worst-case distortion function $g_2$. 
In (a)--(c), the numbers of queries satisfy $M_1=M_2=M_4\in\{10,20,30,40,50\}$. 
The worst-case functions are obtained from the optimal solution of problem~\eqref{eq:bilinear}, where ${\bm \sigma}_1$ and ${\bm \sigma}_2$ in the objective function are fixed at their analytic centers ${\bm \sigma}_1^{\rm AC}$ and ${\bm \sigma}_2^{\rm AC}$, respectively.
}}
\label{fig:uti_1-2}
\end{figure}

\begin{figure}[!ht]
\begin{minipage}[t]{0.3\linewidth}
\centering
\includegraphics[width=2.1in]{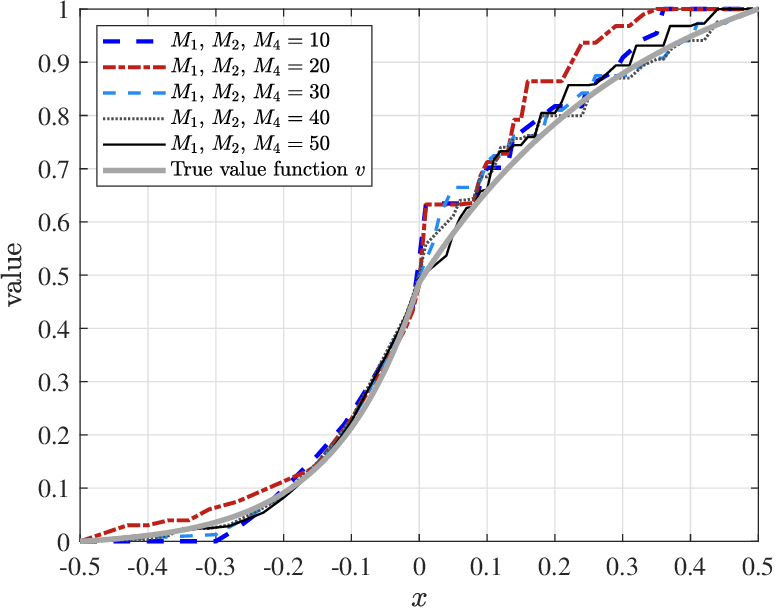}
\text{\hspace{0.9cm}\footnotesize   \textbf{(a)}}
\end{minipage}
\hfill
\begin{minipage}[t]{0.3\linewidth}
\centering
\includegraphics[width=2.1in]{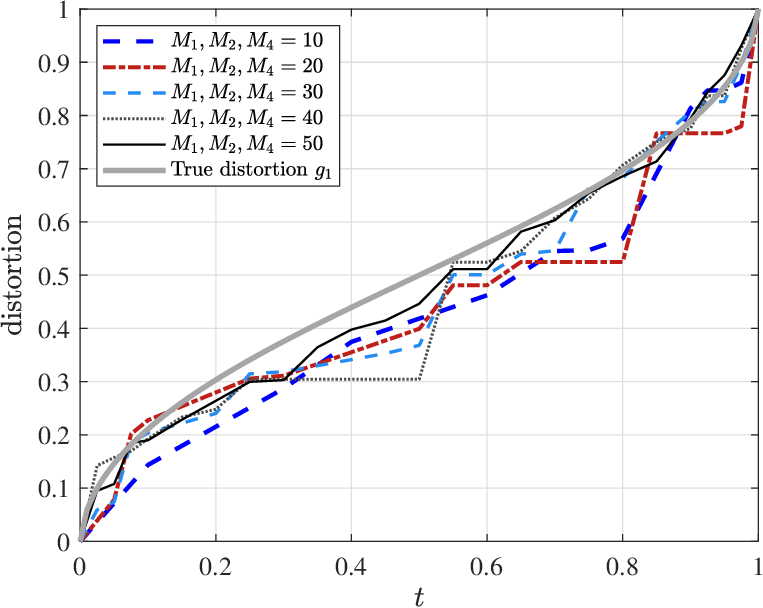}
\text{\hspace{0.9cm}\footnotesize \textbf{(b)}}
\end{minipage}
\hfill
\begin{minipage}[t]{0.3\linewidth}
\centering
\includegraphics[width=2.1in]{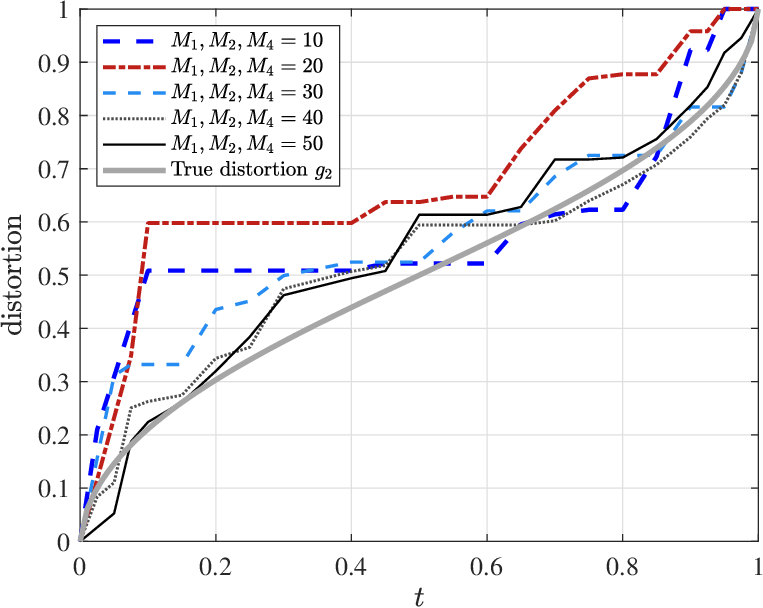}
\text{\hspace{0.9cm}\footnotesize  \textbf{(c)}}
\end{minipage}
\hfill
\caption{\footnotesize{\textbf{(a)} The worst-case value function $v$. 
\textbf{(b)} The worst-case distortion function $g_1$. 
\textbf{(c)} The worst-case distortion function $g_2$. 
In (a)--(c), the numbers of queries satisfy $M_1=M_2=M_4\in\{10,20,30,40,50\}$. 
The worst-case functions are obtained from the optimal solution of problem~\eqref{eq:trilinear_relax_2}, where the number of segments in the PMR is $N=5$.}}
\label{fig:N20}
\end{figure}

We can evaluate $h_{(\mathbb{V}\mathbb{G}_1\mathbb{G}_2)_{M}}(\xi-x)$ with different numbers of queries, and obtain the PRGSR-CPT measure $\rho_{(\mathbb{V}\mathbb{G}_1\mathbb{G}_2)_{M}}(\xi)$ and the worst-case value and distortion functions, 
where the worst-case functions are obtained by solving subproblem~\eqref{eq:trilinear_relax_2} with the optimal solution $x$ of the corresponding PRGSR-CPT model~\eqref{eq:PRGRS-CPT}.
Figure~\ref{fig:value}~(a) depicts the holistic value functions when $M_1=M_2=M_4=10,20,30,40,50$ ($M=M_1+M_2+M_3+M_4$).
This shows that the worst-case value function converges to the true value function.
For the change of optimal value functions of $\rho_{(\mathbb{V}\mathbb{G}_1\mathbb{G}_2)_{M}}(\xi)$, 
as the number of queries increases, 
we show the optimal value function
and the optimal distortion functions  for the case where ${\bm \sigma}_i={\bm \sigma}_i^{\rm AC}$, $i=1,2$, 
in the objective function of problem \eqref{eq:bilinear} in Figure~\ref{fig:uti_1-2}.
The worst-case value function and distortion functions of problem \eqref{eq:trilinear_relax_2} with $N=5$ are shown in Figure~\ref{fig:N20}.
Note that the CPT-based shortfall risk of $\xi$ with true functions $v^*$, $g_1^*$ and $g_2^*$ is $0.2044$.

Figure~\ref{fig:uti_1-2} shows the change of the optimal value function and 
the optimal distortion functions as the number of queries ($M_1,M_2,M_4$) increases from $10$ to $50$, 
when ${\bm \sigma}_i={\bm \sigma}_i^{\rm AC}$,
$i=1,2$, of problem \eqref{eq:bilinear}.
Figure~\ref{fig:N20} depicts the worst-case value function and distortion functions of problem \eqref{eq:trilinear_relax_2} with the number of pieces of McCormick relaxation being $N=5$.

\section{Implications of Using an EUT-Based Polyhedral Method
to Elicit CPT Preferences}
\label{sec:illustration_Inconsis}
Consider the case that the DM's preference is characterized by the CPT with true value function $v^*$ and the distortion functions $g_1^*$ and $g_2^*$ defined as in Section~\ref{sec:setup}.
Suppose that EUT is mistakenly used to model the DM’s preference.
We develop a modified polyhedral method, analogous to that in Section~\ref{sec:ratior}, 
to elicit the unknown utility function under the assumption that the distortion function is the identity.
The ambiguity set of utility functions is defined as
\begin{equation*}
{\cal U}_M=\left\{u_{J_M^{\rm u}}:[\underline{x},\bar{x}]\to [0,1]: h_l(\bbe[u_{J_M^{\rm u}}(Z_1^l)]-\bbe[u_{J_M^{\rm u}}(Z_2^l)])\leq 0,\; l\in [M]\right\},
\end{equation*}
where $h_l=1$ if $\rho_{\rm CPT}(Z_1^l)\leq \rho_{\rm CPT}(Z_2^l)$, and $h_l=-1$ otherwise.
We now use this ambiguity set to compute the corresponding robust utility-based shortfall of $\xi$:
\begin{equation}
\label{eq:Shortfall_EUT}
\inf \left\{x\in \R: \sup_{u\in {\cal U}_M} (\bbe[u(\xi-x)]-u(0))\leq 0\right\}.
\end{equation}
Let $\widetilde{\cal X}_M^{\rm u} = {\cal X}_M^{\rm u} \cup \{\xi_k-x\}_{k=1}^K$.
Analogous to the computation of PRGSR-CPT in Section~\ref{sec:bisection-robust}, we apply a bisection method. 
For each fixed $x$, the inner problem can be reformulated as
\begin{equation*}
\begin{aligned}
\sup_{{\bm v}}\; &\sum_{k=1}^K({\bm R}{\bm v})^\top {\bm h}^{\widetilde{\cal X}_M^{\rm u}}(\xi_k-x) p_k-({\bm R}{\bm v})^\top {\bm h}^{\widetilde{\cal X}_M^{\rm u}}(0) \\
{\rm s.t.}\; & {\bm v}\geq 0, \; {\bm e}^\top {\bm v}\leq 1, \\
& h_m(({\bm R}{\bm v})^\top ((1-p^m){\bm h}^{\widetilde{\cal X}_M^{\rm u}}(r_1^m)+p^m {\bm h}^{\widetilde{\cal X}_M^{\rm u}}(r_3^m)-{\bm h}^{\widetilde{\cal X}_M^{\rm u}}(r_2^m)))\leq 0,\;\forall m\in [M].
\end{aligned}
\end{equation*}
The utility-based shortfall can be calculated by the bisection method as described in Section~\ref{sec:bisection-robust}.

\begin{figure}[!ht]
\begin{minipage}[t]{0.49\linewidth}
\centering
\includegraphics[width=3.1in]{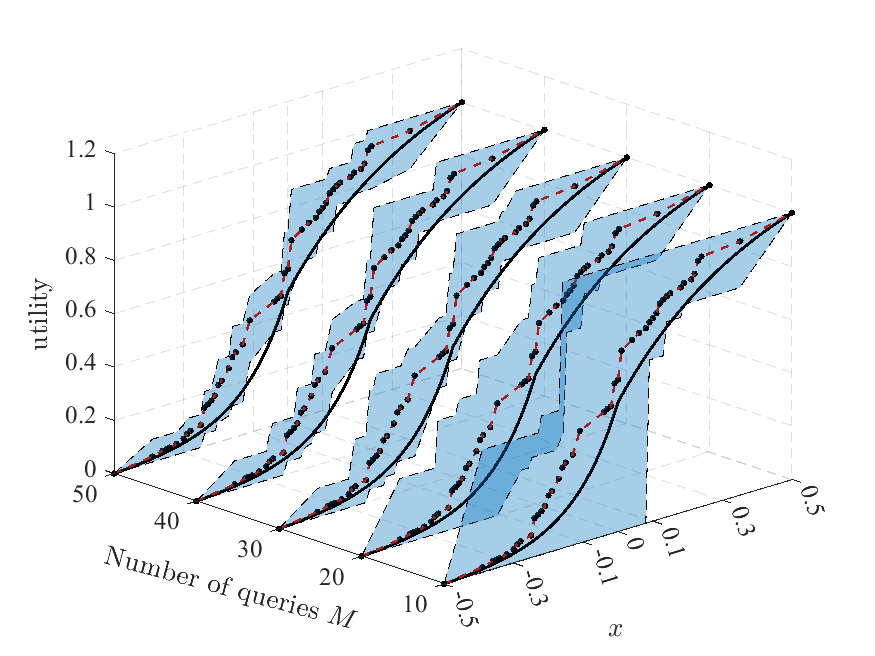}
\text{\footnotesize \textbf{(a)}}
\end{minipage}
\hfill
\begin{minipage}[t]{0.49\linewidth}
\centering
\includegraphics[width=2.9in]{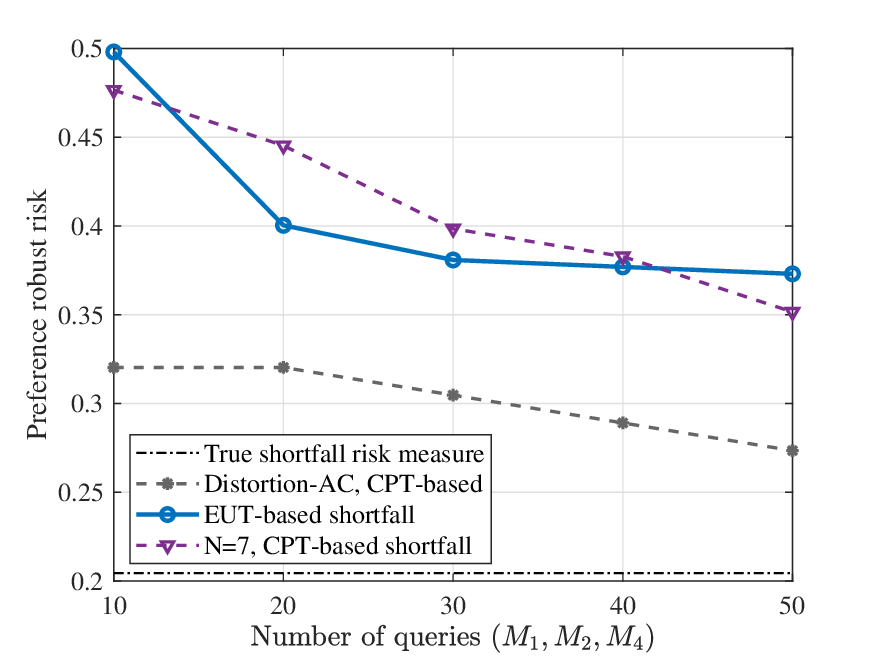}
\text{\footnotesize \textbf{(b)}}
\end{minipage}
\caption{\footnotesize{\textbf{(a)} Range of utility functions.
The solid line represents the true value function $v^*$. 
The shaded region shows the envelope of admissible utility functions over $[-0.5,0.5]$ induced by ${\cal U}_M$ for $M=10,20,30,40,50$. 
\textbf{(b)} Robust shortfall under EUT and CPT. 
The curves show the shortfall risk measure as the number of queries $M_4$ (with $M_1 = M_2 = M_4$) increases with the number of segments in the PMR being $N=7$. 
The horizontal dotted line indicates the true shortfall value, 
illustrating the higher shortfall risk value using EUT.}}
\label{fig:EUT}
\end{figure}

Figure \ref{fig:EUT} (a) shows that when the distortion function is incorrectly fixed as the identity (i.e., under EUT),
the true value function lies outside the admissible set of utility functions inferred from the queries when $M\geq 40$.
This indicates a misspecification of the preference model.
Figure \ref{fig:EUT} (b) further illustrates the impact of this misspecification on the resulting robust shortfall risk measure.
When the distortion function is fixed as the identity, 
the resulting EUT-based robust shortfall is higher than the CPT-based robust shortfall for larger values when $M_4=50$ ($M=50$ for the EUT-based shortfall).
In addition, it is also higher than the gray dotted line, 
which represents the PRGSR-CPT measure obtained by solving the inner problem using the bilinear objective function in~\eqref{eq:bilinear}, 
with ${\bm \sigma}_1$ and ${\bm \sigma}_2$ fixed at ${\bm \sigma}_1^{\rm AC}$ and ${\bm \sigma}_2^{\rm AC}$, respectively.
These observations suggest that ignoring probability distortion leads to a systematic overestimation of the robust shortfall risk. 
In contrast, the CPT-based approach, 
which allows for flexible distortion functions, 
provides estimates that are closer to the true risk level.
Overall, this explains why rejecting the true value function under a misspecified model (EUT) does not imply better estimation accuracy; 
instead, it reflects model misspecification that leads to biased (conservative) risk evaluation.

\section{Qualitative Summary of Convergence of the Estimated Value Function and Distortion Functions}
\label{sec:quantatative}
In this section, we provide a qualitative summary and supporting numerical evidence for the convergence behavior of the proposed elicitation procedure.
In Section~\ref{sec:quantitative1}, 
we quantify the distance from the analytic center to the boundaries of the corresponding polyhedra $\Delta_-^{M_1}$, 
$\Delta_+^{M_2}$, and $\Delta_{M_4}$ defined in \eqref{eq:Delta_minus}, \eqref{eq:poly_Delta+}, 
and \eqref{eq:whole_uti} with $m=M_4$.
In Section~\ref{sec:quantitative2}, 
we evaluate the mean squared distance between the estimated value and distortion functions and their true counterparts.

\subsection{Distance from the Analytic Center to the Defining Hyperplanes}
\label{sec:quantitative1}
\begin{figure}[!ht]
\begin{minipage}[t]{0.32\linewidth}
\centering
\includegraphics[width=2.0in]{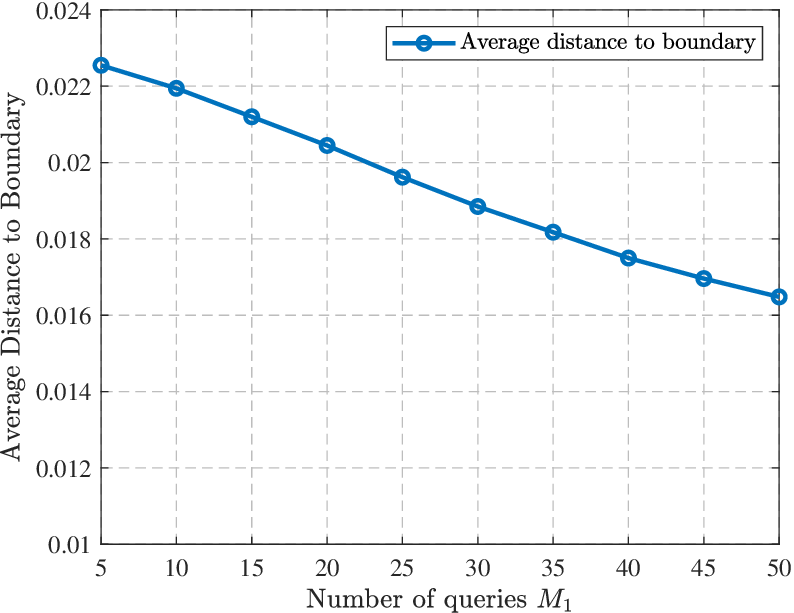}
\text{\footnotesize \hspace{0.5cm} \textbf{(a)}}
\end{minipage}
\begin{minipage}[t]{0.32\linewidth}
\centering
\includegraphics[width=2.0in]{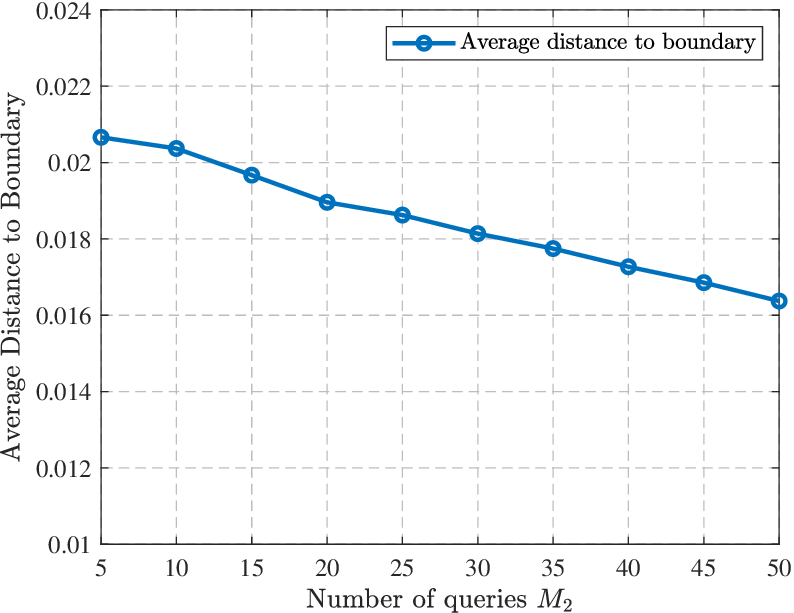}
 \text{\footnotesize \hspace{0.5cm}  \textbf{(b)}} 
\end{minipage}
\begin{minipage}[t]{0.32\linewidth}
\centering
\includegraphics[width=2.0in]{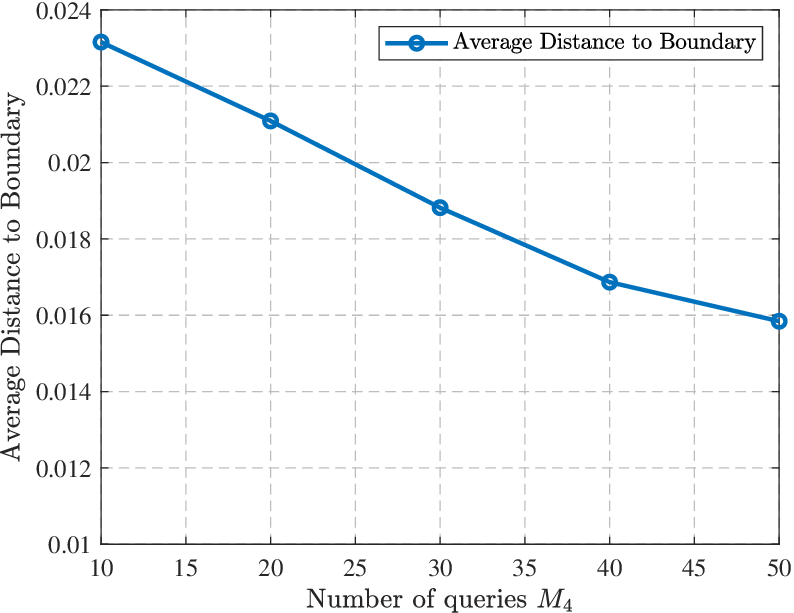}
 \text{\footnotesize \hspace{0.5cm}  \textbf{(c)}}
\end{minipage}
\caption{\footnotesize{\textbf{(a)} Distance from the elicited analytic center to the defining hyperplanes of $\widehat{\Delta}_{-}^{M_1}$ as the number of queries $M_1$ increases. 
\textbf{(b)} Distance from the elicited analytic center to the defining hyperplanes of 
$\widehat{\Delta}_{+}^{M_2}$ as the number of queries $M_2$ increases. 
\textbf{(c)} Distance from the elicited analytic center to the defining hyperplanes of
$\widetilde{\Delta}_{M_1+M_2+M_3+M_4}^v$ as the number of queries $M_1=M_2=M_4$ increases.
}}
\label{fig:distance}
\end{figure}

Let $\widehat{\Delta}_-^{M_1}$ denote the normalized version of the polyhedron $\Delta_-^{M_1}$ obtained by scaling each constraint so that the corresponding normal vector has unit norm. 
For simplicity, we write $\widehat{\Delta}_-^{M_1} = \{x : A_i^- x \le b_i^-,\ i \in {\cal I}_-\}$.
The analytic center of $\widehat{\Delta}_-^{M_1}$ is obtained by solving
$\max_x \sum_{i\in {\cal I}_-} \log(b_i^- - A_i^- x)$.
Let $x_-^{\rm AC}$ denote the resulting analytic center. 
The average distance from $x_-^{\rm AC}$ to the defining hyperplanes of $\widehat{\Delta}_-^{M_1}$ is defined as
\begin{equation*}
d_{\rm Ave}^{M_1} = \frac{1}{|{\cal I}_-|}\sum_{i\in{\cal I}_-} (b_i^- - A_i^- x_-^{\rm AC}).
\end{equation*}
For the normalized polyhedron $\widehat{\Delta}_+^{M_2}=\{x:A_i^+x\leq b_i^+, i\in {\cal I}_+\}$ with analytic center $x_+^{\rm AC}$, 
let
\begin{equation*}
d_{\rm Ave}^{M_2} = \frac{1}{|{\cal I}_+|}\sum_{i\in{\cal I}_+}(b_i^+ - A_i^+ x_+^{\rm AC}).
\end{equation*}
For the normalized polyhedron $\widetilde{\Delta}_{M}^v=\{x:A_ix\leq b_i, i\in {\cal I}\}$ with $M=M_1+M_2+M_3+M_4$ and analytic center $x^{\rm AC}$, 
let 
\begin{equation*}
d_{\rm Ave}^{M} =\frac{1}{|{\cal I}|}\sum_{i\in{\cal I}} (b_i - A_i x^{\rm AC}).
\end{equation*}
Figure~\ref{fig:distance} depicts the convergence of these distances as the number of queries increases. 
Specifically, Figure~\ref{fig:distance}(a) shows the convergence of the distances from the analytic center to the defining hyperplanes  of $\widehat{\Delta}_-^{M_1}$, 
while Figure~\ref{fig:distance}(b) and (c) show the corresponding distances for $\widehat{\Delta}_+^{M_2}$ and $\widetilde{\Delta}_M^v$.

\subsection{Mean Squared Distance to the True Functions}
\label{sec:quantitative2}
To evaluate the estimation accuracy, we measure the mean squared distance 
between the estimated functions and their true counterparts using randomly 
generated evaluation points.
For a function $f$ and its estimate $\hat f$, the error is defined by 
\begin{equation*}
\dd(f,\hat f)=\frac{1}{N}\sum_{i=1}^{N}
\left|f(x_i^m)-\hat f(x_i^m)\right|^2,
\end{equation*}
where $\{x_i^m\}_{i=1}^N$ are uniformly generated evaluation points.
Throughout the experiments, we set $N=1000$.
Let $v_-^*$, $v_+^*$, and $v^*$ denote the true normalized value functions over $[-0.5,0]$, $[0,0.5]$, and $[-0.5,0.5]$, respectively.

\begin{figure}[!ht]
\begin{minipage}[t]{0.31\linewidth}
\centering
\includegraphics[width=2.1in]{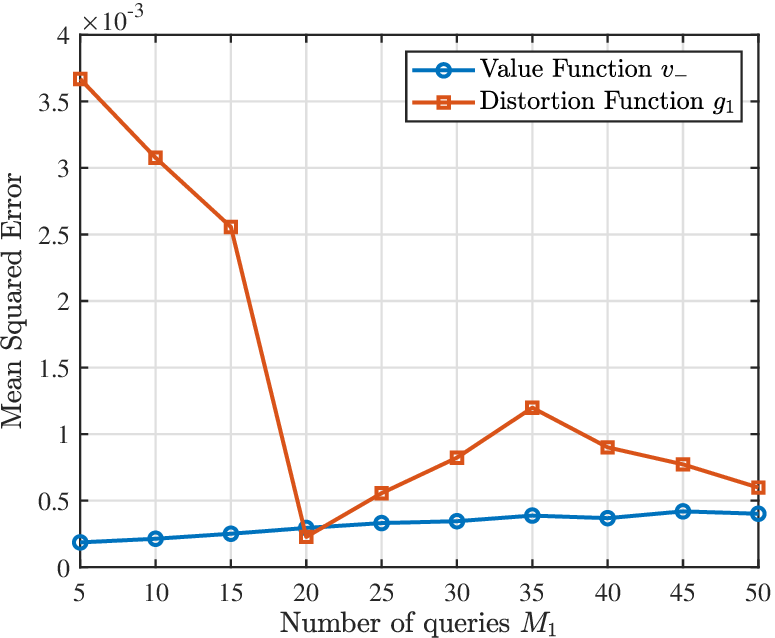}
\text{\footnotesize \hspace{0.5cm} \textbf{(a)}} 
\end{minipage}
\hfill
\begin{minipage}[t]{0.31\linewidth}
\centering
\includegraphics[width=2.1in]{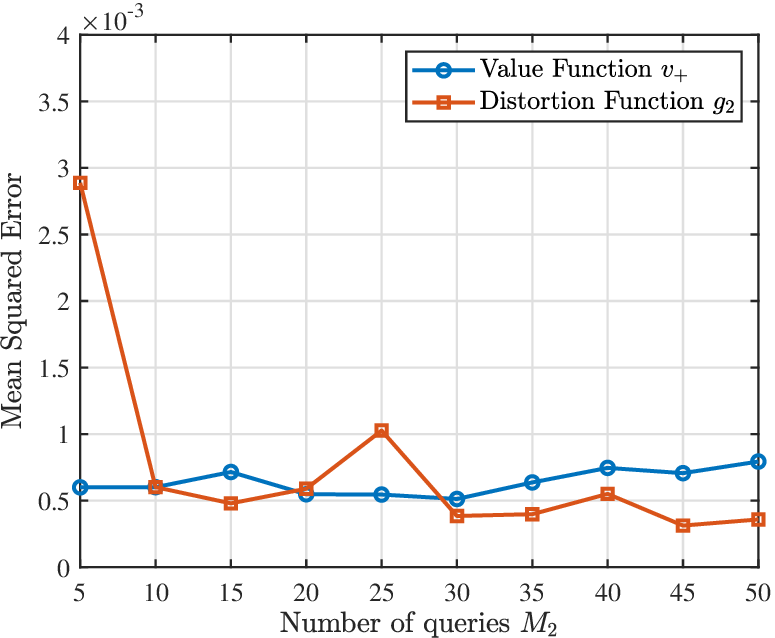}
\text{\footnotesize \hspace{0.5cm} \textbf{(b)}}
\end{minipage}
\hfill
\begin{minipage}[t]{0.31\textwidth}
\centering
\includegraphics[width=2.1in]{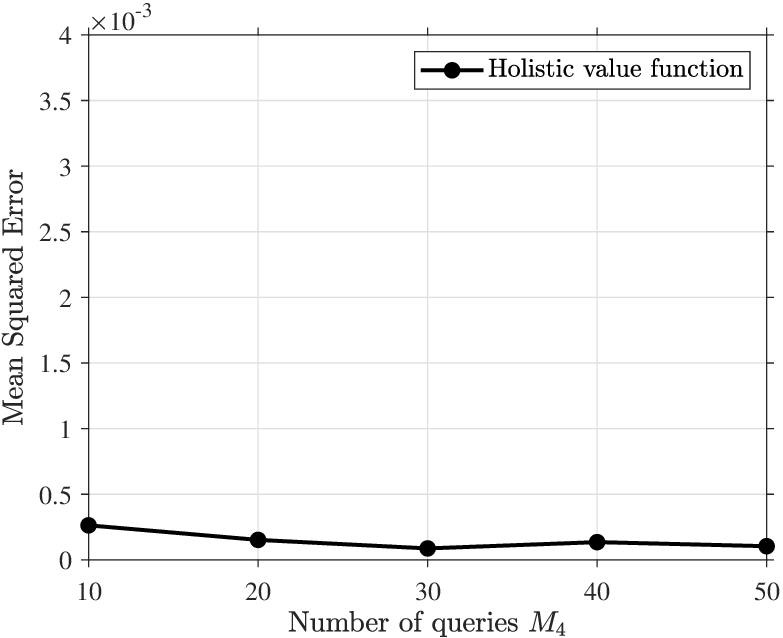}
\text{\footnotesize \hspace{0.5cm} \textbf{(c)}}
\end{minipage}
\caption{\footnotesize{\textbf{(a)} Mean squared errors of the estimated distortion function $g_1$ obtained from the analytic center of polyhedron 
$\widehat{\Delta}_{-}^{M_1}$. 
\textbf{(b)} Mean squared errors of the estimated distortion function $g_2$ obtained from the analytic center of polyhedron $\widehat{\Delta}_{+}^{M_2}$. 
\textbf{(c)} Mean squared error of the estimated holistic value function $v$ obtained from the analytic center of polyhedron $\widetilde{\Delta}_{M_1+M_2+M_3+M_4}^{v}$. 
Here, $M_1$, $M_2$, $M_3$, and $M_4$ denote the numbers of queries of Types 1--4, respectively. 
In particular, for $M_4 \in \{10,20,30,40,50\}$, we set $M_1=M_2=M_4$. 
The construction of $\widetilde{\Delta}_{M_1+M_2+M_3+M_4}^{v}$ relies on 
$\widehat{\Delta}_{M_1+M_2+M_3+M_4}^{v-v+}$,
where the estimators $\widehat{g}_1$ and $\widehat{g}_2$ are defined as in \eqref{eq:hat_g1g2}.}}
\label{fig:uti_1}
\end{figure}

\underline{Negative value and distortion functions.}
For the negative part of the value function, evaluation points are generated over $[-0.5,0]$, while for the distortion function they are generated over $(0,1)$. 
Let $({\bm c}_{1,-}^{M_1},{\bm c}_{2,-}^{M_1})$ denote the analytic center of $\Delta_-^{M_1}$. 
The corresponding estimated value and distortion functions are given by
$v_-^{\rm AC}(x) = ({\bm R}{\bm c}_{1,-}^{M_1})^\top {\bm h}^{{\cal X}_-^{M_1}}(x)-1$ and $g_1^{\rm AC}(p) = ({\bm R}{\bm c}_{2,-}^{M_1})^\top {\bm h}^{{\cal S}_1}(p)$.
Figure~\ref{fig:uti_1} (a) shows the evolution of the corresponding errors 
$\dd(v_-^{\rm AC},v_-^*)$ and $\dd(g_1^{\rm AC},g_1^*)$ as the number of 
queries increases from $5$ to $50$.
Overall, the estimation error of the distortion function decreases consistently, whereas the improvement for the negative part of the value function is relatively moderate.

\underline{Positive value and distortion functions.}
Let $({\bm c}_{1,+}^{M_2},{\bm c}_{2,+}^{M_2})$ denote the analytic center of $\Delta_+^{M_2}$. 
The corresponding estimated value and distortion functions are given by
$v_+^{\rm AC}(x) = ({\bm R}{\bm c}_{1,+}^{M_2})^\top {\bm h}^{{\cal X}_{+}^{M_2}}(x)$ and $g_2^{\rm AC}(p) = ({\bm R}{\bm c}_{2,+}^{M_2})^\top {\bm h}^{{\cal S}_2}(p)$.
Figure~\ref{fig:uti_1} (b) illustrates the convergence of the errors 
$\dd(v_+^{\rm AC},v_+^*)$ and $\dd(g_2^{\rm AC},g_2^*)$. 

\underline{Holistic value function.}
Let $({\bm c}_1^{M},{\bm c}_2^M)$ denote the analytic center of the polyhedron $\widehat{\Delta}_{M_1+M_2+M_3+M_4}^v$. 
To evaluate the error of the overall value function $v^{\rm AC}(x)=\left(\frac{1}{1+\widehat{r}}{\bm R}{\bm c}_1^{M},\frac{\widehat{r}}{1+\widehat{r}}{\bm R}{\bm c}_2^{M}\right)^\top{\bm h}^{{\cal X}_{M_4}}(x)$, 
we generate  evaluation points uniformly over $[-0.5,0.5]$. 
The resulting error $\dd(v^{\rm AC}, v^*)$ is shown in Figure~\ref{fig:uti_1} (c).

The numerical results indicate that the proposed polyhedral elicitation 
method successfully reduces estimation error but also appears to stabilize after some number of queries.
For $\Delta_-^{M_1}$, the estimation error appears to become stable after approximately 35 questions of Type 1, 
while it is around 20  questions for the estimation of the holistic value function.
These results suggest that while our elicitation may eventually produce high-accuracy functions, further reductions in estimation error may require additional queries.

\section{Comparative Analysis of Trade-offs with Standard Parametric Approach}
\label{sec:comparition_para}
Suggested by a reviewer, we have also undertaken some comparative analysis between our method and a standard parametric approach where the true value/distortion functions are not necessarily parametric.  
We consider a parametric family commonly used in CPT (see Table~\ref{tab:Distortion}). Since the true preference functions may not belong to this parametric family (a situation we refer to as the model being misspecified), 
the fitted functions may deviate from the true ones even when the parameters are optimally estimated.

\noindent\underline{Specification of a parametric family}.
We consider the following family of parametric value functions:
\begin{equation}
\label{eq:value_hat}
\widetilde{v}(x;\alpha,\beta,\lambda)= \left\{\begin{aligned}
& x^{\alpha}, && x \ge 0, \\
&-\lambda (-x)^{\beta}, && x < 0,
\end{aligned}
\right.
\end{equation}
and distortion functions:
\begin{equation}
\label{eq:g_hat}
\widetilde{g}_1(t;\delta_1,\gamma_1)=\exp\!\left(-\delta_1(-\ln t)^{\gamma_1}\right),
\qquad
\widetilde{g}_2(t;\delta_2,\gamma_2)=\exp\!\left(-\delta_2(-\ln t)^{\gamma_2}\right),
\end{equation}
where $\alpha, \beta \in [0,1]$, $\delta_k \in [0,2]$, and $\gamma_k \in [0,2]$ for $k=1,2$.

\noindent\underline{Estimation of the Parameter}. 
To estimate the parameters of the value function, we first normalize the 
elicited value function over the interval $[-0.5,0.5]$. 
Preference queries provide observations $v_i^M = \left(\frac{1}{1+\widehat{r}}{\bm R}{\bm c}_1^{M},\frac{\widehat{r}}{1+\widehat{r}}{\bm R}{\bm c}_2^{M}\right)^\top{\bm h}^{{\cal X}_{M_4}}(x_i)$ at selected points $x_i \in {\cal X}_{M}$ obtained by the analytic center $({\bm c}_1^{M},{\bm c}_2^M)$ of the polyhedron $\widehat{\Delta}_{M_1+M_2+M_3+M_4}^v$.
The parameters $(\alpha,\beta,\lambda)$ are estimated by solving the following mean squared error problem:
\begin{equation}
\label{eq:least_SQ_value}
\begin{aligned}
\min_{\alpha,\beta \in [0,1],
\tau\in [0,2], \lambda\tau\in [0,2]} 
& \sum_{i=1}^{J}\left(\tau \widetilde{v}(x_i;\alpha,\beta,\lambda) - v_i^M
\right)^2 \\
\text{s.t.}\qquad \quad \;\; & \tau\big(\widetilde{v}(\bar{x})-\widetilde{v}(\underline{x})\big)=1,
\end{aligned}
\end{equation}
where the scaling parameter $\tau$ is introduced to align the parametric value function with the normalized elicited value function.
Similarly, the parameters $(\delta_k,\gamma_k)$ of the distortion functions are estimated via nonlinear least squares.
The elicited values are obtained from the analytic center $({\bm c}_{-,1}^{M_1},{\bm c}_{-,2}^{M_1})$ of the polyhedron $\Delta_-^{M_1}$ as $g_{1,j}^{M_1}=({\bm R}{\bm c}_{-,2}^{M_1})^\top{\bm h}^{{\cal S}_1}(p_j)$.
Likewise, the elicited values of $g_{2,j}^{M_2}=({\bm R}{\bm c}^{M_2}_{+,2})^\top{\bm h}^{{\cal S}_2}(p_j)$ are obtained from the analytic center $({\bm c}_{+,1}^{M_2},{\bm c}_{+,2}^{M_2})$ of polyhedron $\Delta_+^{M_2}$. We solve the following regression problems:
\begin{equation}
\label{eq:para_estimate}
\begin{aligned}
\min_{\delta_1,\gamma_1\in [0,2]} 
& \sum_{j=1}^{S}
\left(
\exp\!\left(-\delta_1(-\ln p_j)^{\gamma_1}\right)-g_{1,j}^{M_1}
\right)^2, \\
\min_{\delta_2,\gamma_2\in [0,2]} 
& \sum_{j=1}^{S}
\left(
\exp\!\left(-\delta_2(-\ln p_j)^{\gamma_2}\right)-g_{2,j}^{M_2}
\right)^2.
\end{aligned}
\end{equation}
Let $(\alpha^*,\beta^*,\lambda^*)$ and $(\delta_k^*,\gamma_k^*)$, $k=1,2$, 
denote optimal solutions to problems \eqref{eq:least_SQ_value}  and \eqref{eq:para_estimate}, respectively.

\noindent\underline{Evaluation metric}.
To quantify the discrepancy between the estimated parametric functions and the true preference functions, we use the mean squared error evaluated over
randomly generated points.
For the distortion functions, the error is defined as
\begin{equation*}
\dd(\widetilde{g}_{k}(\cdot;\delta_k^*,\gamma_k^*),g_k)
= \frac{1}{N}\sum_{i=1}^{N}
|\widetilde{g}_{k}(t_i^m;\delta_k^*,\gamma_k^*)-g_k(t_i^m)|^2,
\; k=1,2,
\end{equation*}
where $\{t_i^m\}_{i=1}^N$ are uniformly drawn from $(0,1)$.
For the value function, let
\begin{equation*}
\dd(\widetilde{v}(\cdot;\alpha^*,\beta^*,\lambda^*),v^*)
= \frac{1}{J}\sum_{i=1}^{J}
|\widetilde{v}(x_i^m;\alpha^*,\beta^*,\lambda^*)-v^*(x_i^m)|^2,
\end{equation*}
where $\{x_i^m\}_{i=1}^{J}$ are uniformly drawn from $[-0.5,0.5]$.
Throughout the experiments, we set $N=1000$ and $J=1000$.
For comparison, we also implement a parametric benchmark, 
where the functional forms of the value and distortion functions are assumed to be correctly specified (see Section~\ref{sec:setup}).
The parameters are estimated by solving \eqref{eq:least_SQ_value} and \eqref{eq:para_estimate}, 
which yields the minimum achievable approximation error under the assumed parametric specification.

\begin{figure}[!ht]
\centering
\begin{minipage}{0.32\textwidth}
\centering
\includegraphics[width=\linewidth]{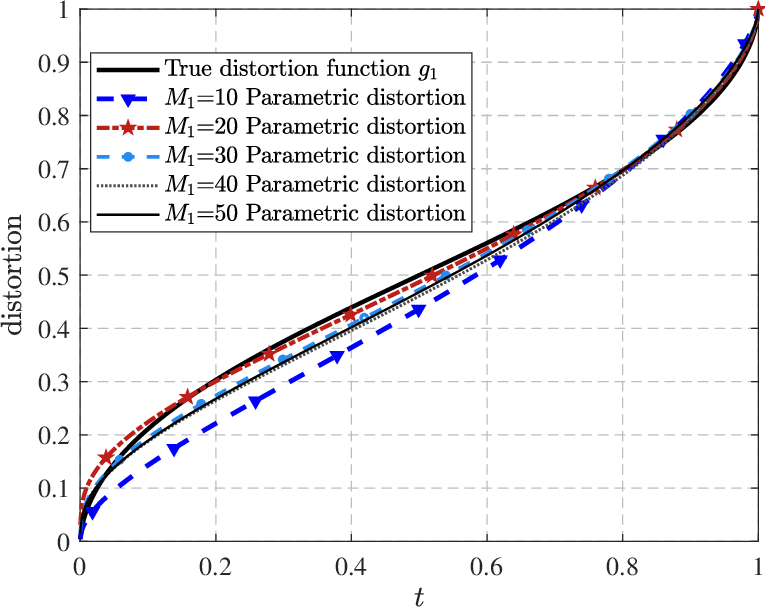}
\par\smallskip
{\footnotesize \hspace{0.2cm} \textbf{(a)}}
\end{minipage}
\hfill
\begin{minipage}{0.32\textwidth}
\centering
\includegraphics[width=\linewidth]{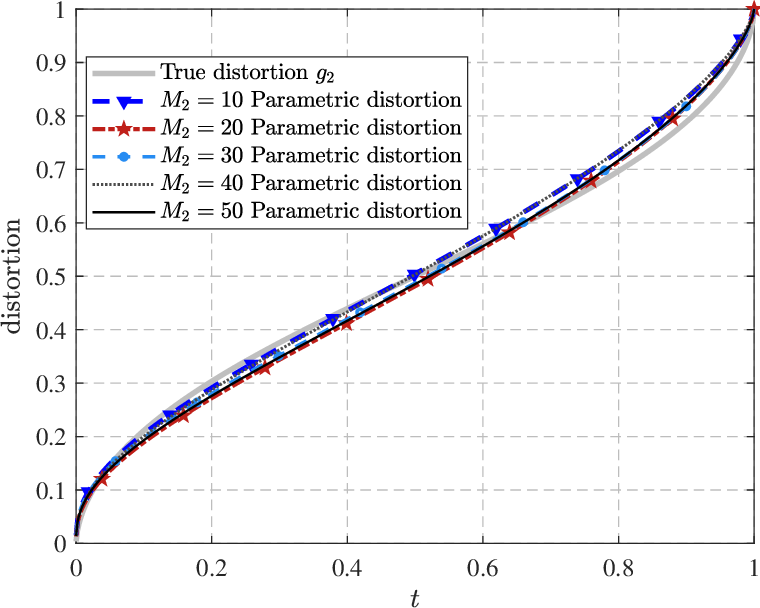}
\par\smallskip
{\footnotesize \hspace{0.2cm} \textbf{(b)}}
\end{minipage}
\hfill
\begin{minipage}{0.32\textwidth}
\centering
\includegraphics[width=\linewidth]{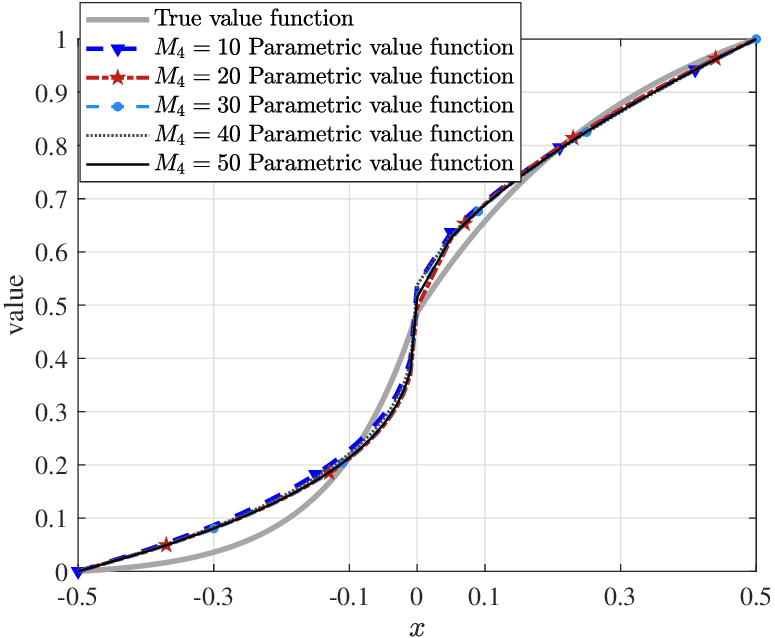}
\par\smallskip
{\footnotesize \hspace{0.2cm} \textbf{(c)}}
\end{minipage}
 \caption{\footnotesize \textbf{(a)} Parametric Estimation of $\widetilde{g}_1(\cdot;\delta_1,\gamma_1)$ using function values derived from the
analytic center of the polyhedron $\Delta_-^{M_1}$.
 \textbf{(b)} Parametric Estimation of $\widetilde{g}_2(\cdot;\delta_2,\gamma_2)$ using function values derived from the
analytic center of the polyhedron $\Delta_+^{M_2}$.
 \textbf{(c)} Parametric Estimation of $\widetilde{v}(\cdot;\alpha,\beta,\lambda)$ using the function values derived from the analytic center of the polyhedron $\widehat{\Delta}_{M}^{v}$ (the projection of $\widehat{\Delta}_{M}$ onto the space of
$({\bm R}\widetilde{\bm v}_-, {\bm R}\widetilde{\bm v}_+)$).
In (a)--(c),
the black solid curve represents the true functions,
while the dashed curves represent the estimated parametric functions.}
\label{fig:dist_minus}
\end{figure}

\noindent\underline{Comparative analysis}.
Figure~\ref{fig:dist_minus} depicts the estimated distortion and 
the value functions based on $\Delta_-^{M_1}$, $\Delta_+^{M_2}$, 
and $\Delta_{M_1+M_2+M_3+M_4}^v$, respectively. 
The parametric CPT specification produces smooth functional forms but may deviate from the true functions due to model misspecification. 
In contrast, the nonparametric polyhedral approach more closely tracks 
the true preference functions, 
although small local fluctuations may 
appear as a consequence of the polyhedral elicitation procedure.

\begin{figure}[!ht]
\centering
\begin{minipage}{0.32\textwidth}
\centering
\includegraphics[width=\linewidth]{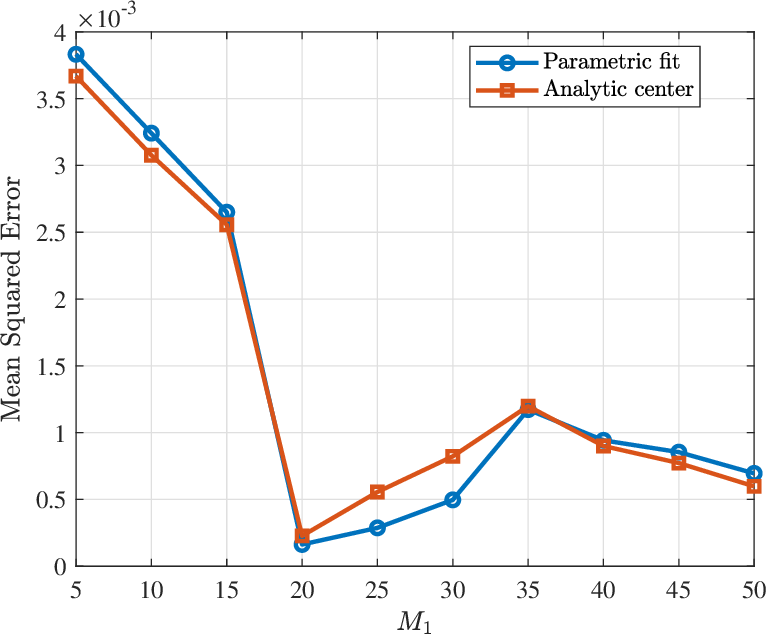}
\par\smallskip
{\footnotesize \hspace{0.2cm} \textbf{(a)}}
\end{minipage}
\hfill
\begin{minipage}{0.32\textwidth}
\centering
\includegraphics[width=\linewidth]{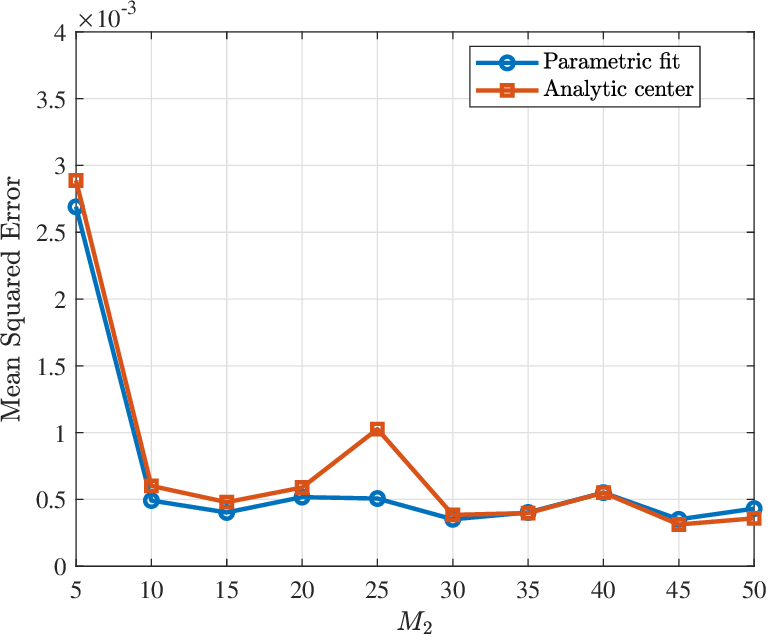}
\par\smallskip
{\footnotesize \hspace{0.2cm} \textbf{(b)}}
\end{minipage}
\hfill
\begin{minipage}{0.32\textwidth}
\centering
\includegraphics[width=\linewidth]{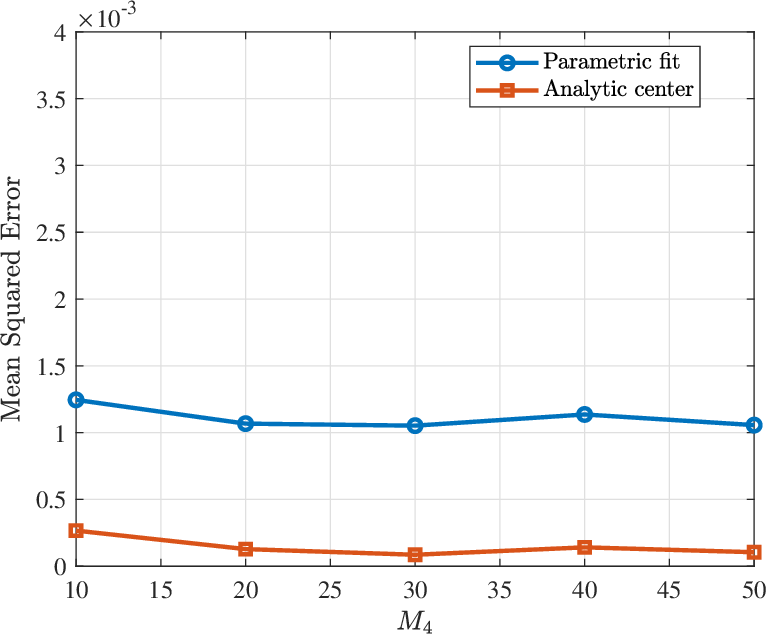}
\par\smallskip
{\footnotesize \hspace{0.2cm} \textbf{(c)}}
\end{minipage}
\caption{\footnotesize \textbf{(a)} Evolution of the mean squared error for $g_1$ as the number of 
queries $M_1$ increases.
\textbf{(b)} Evolution of the mean squared error for $g_2$ as the number of 
queries $M_2$ increases.
\textbf{(c)} Evolution of the mean squared error for $v$ as the number of 
queries $M_4$ increases.
In (a)-(c), the blue curve corresponds to the 
parametric estimation and the orange curve corresponds to the nonparametric 
(polyhedral) estimation.
Both estimators improve as the number of queries increases. 
The nonparametric estimator generally achieves lower error in the large-sample regime, as illustrated when $M_1 \geq 35$, $M_2 \geq 40$, and across all tested values of $M_4$.
}
\label{fig:dist_distortion2}
\end{figure}

Figure~\ref{fig:dist_distortion2} further reports the evolution of the mean squared error as the number of queries increases. 
The results show that the parametric method achieves similar performance in terms of distortion function quality,
while demonstrating a clear improvement in terms of the quality of the holistic value functions $v$. 
We attribute the latter loss in performance to the misspecification of the parametric family for $\tilde{v}$.
Based on the test results, 
we conclude that the performance of the parametric approach is subject to model misspecification and estimation error, 
highlighting the advantages of the nonparametric polyhedral method.

\section{Extended Concluding Remarks}
\label{SC:conclu-rem}
In this section, we provide additional details on the implementation and limitations of the proposed models and CPT‑based preference elicitation approaches, 
as well as their application to optimal decision‑making problems.

\begin{itemize}
\item[(i)] In practice, the only directly observable data are historical realizations of $\xi$ and ${\bm \zeta}$, 
whereas the support of the random variable $\xi$ and random vector ${\bm \zeta}$ are not directly observable. 
However, the support set of $\xi$ or ${\bm \zeta}$ can often be inferred from domain-specific knowledge (e.g., the random loss $\xi$ is nonnegative and is upper bounded, and the random return rates ${\bm \zeta}$ can be supposed bounded); 
it may also be constructed from historical data (e.g., the convex hull of all historical realizations of ${\xi}$ and ${\bm \zeta}$ respectively).

\item[(ii)] There is a limitation in the developed models and preference elicitation approaches,
that is, we are unable to use the CPT-based shortfall risk measure to elicit the DM’s risk preference directly; 
rather we have to use the distorted expected values of the random prospects, 
that is, $\bbe_{g_1g_2}[v(\cdot)]$ to elicit the DM’s preference via pairwise comparison questionnaires.

\item[(iii)] Weak preferences are unavoidable under the current questionnaire design.
When the expected values of two lotteries are very close, the decision maker may not show a strong preference.
The weak preferences may affect the correctness of the added inequality.
Due to rounding of the generated cut points $r_1,r_2,r_3$, 
the sign of the newly added constraint in the ambiguity set may be reversed. 
To tackle the issue, one may incorporate the robust approaches \cite{BeO13}.
We leave this for future exploration.

\item[(iv)] 
The robust model, preference elicitation approaches and the formulations
discussed in this paper can be easily applied to one-stage decision-making problems with a CPT-based objective:
\begin{equation}
\label{eq:PRO-SXZ-CPT}
\inmat{(PRO-CPT)}  \quad \;\;
\max_{{\bm z}\in Z} \min_{(v,g_1,g_2)\in \mathbb{VG}_1\mathbb{G}_2}
 \bbe_{{g}_1{g}_2}[v(f({\bm z},{\bm \zeta}))],
\end{equation}
where $\mathbb{VG}_1\mathbb{G}_2$ is the ambiguity set of $(v,g_1,g_2)$.
Specifically, we may begin by obtaining an ambiguity set $(\mathbb{V}_-\mathbb{G}_1)_M$ and $(\mathbb{V}_+\mathbb{G}_2)_M$ of $({v}_-,{g}_1)$ and $({v}_+,{g}_2)$ of the true value functions respectively using the modified
polyhedral method presented in Section~\ref{sec:modified_PM}.
Then we can use the by-product, estimate $\widehat{g}_1$ and $\widehat{g}_2$ to construct the ambiguity set of the normalized value function $v$ over $[\underline{x},\bar{x}]$ with $M$ pairs of queries to elicit the range of the potential ratio $r$ (see Section~\ref{sec:ratior}).
Consequently, we use $(\mathbb{V}_-\mathbb{G}_1)_M$, $(\mathbb{V}_+\mathbb{G}_2)_M$ and the range of $r$ to construct $(\mathbb{V}\mathbb{G}_1\mathbb{G}_2)_M$.
Then we can consider the following simplified preference robust optimization model:
\begin{equation}
\label{eq:PRO-SXZ-CPT-2}
\inmat{(PRO-CPT-S)}  \quad \;\;
\vt:=\max_{{\bm z}\in Z} \min_{(v,g_1,g_2)\in (\mathbb{V}\mathbb{G}_1\mathbb{G}_2)_M}   \bbe_{{g}_1{g}_2}[v(f({\bm z},{\bm \zeta}))].
\end{equation}
Analogous to discussions in Section~\ref{SC:extention_models}, 
we can use $(\mathbb{V}_J\mathbb{G}_{N_1}\mathbb{G}_{N_2})_M$ to approximate the ambiguity set $(\mathbb{V}\mathbb{G}_1\mathbb{G}_2)_M$, and then consider the approximate preference robust optimization model: 
\begin{equation}
\label{eq:approximate_CPT}
\inmat{(PRO-CPT-S-J)}  \quad \;\;\vt_{JN}
:= \max_{{\bm z}\in Z}\min_{(v,g_1,g_2)\in (\mathbb{V}_J\mathbb{G}_{N_1}\mathbb{G}_{N_2})_M} \mathbb
E_{g_1{g_2}}[v(f({\bm z},{\bm \zeta}))].\qquad 
\end{equation}
For each fixed decision vector ${\bm z}$, the inner minimization problem as a bilinear program with linear constraints and can be approximated using the piecewise McCormick relaxation described in Section~\ref{sec:Piecewise_Mc}.
The error arising from the approximation can be quantified as in \citep{GXZ23}. 
In comparison with the expected utility maximin problem in \citep{HZXZ22,ArD15,GXZ23} and the preference robust distortion risk measure in \cite{WaX23}, 
the preference robust optimization model \eqref{eq:PRO-SXZ-CPT} covers a larger class of problems where the DM’s preference is represented by CPT.
\end{itemize}


\begin{thebibliography}{34}
\providecommand{\natexlab}[1]{#1}
\providecommand{\url}[1]{\texttt{#1}}
\providecommand{\urlprefix}{URL }

\bibitem[{Abdellaoui(2000)}]{Abd00}
Abdellaoui M (2000) Parameter-free elicitation of utility and probability  weighting functions. \emph{Management Science} 46(11):1497--1512.

\bibitem[{Abdellaoui(2002)}]{Abd02}
Abdellaoui M (2002) A genuine rank‐dependent generalization of the {von} Neumann‐Morgenstern expected utility theorem. \emph{Econometrica} 70(2):717--736.

\bibitem[{Abdellaoui et~al.(2007)Abdellaoui, Bleichrodt, \protect\BIBand{} Paraschiv}]{ABP07}
Abdellaoui M, Bleichrodt H, Paraschiv C (2007) Loss aversion under prospect theory: A parameter-free measurement. \emph{Management Science} 53(10):1659--1674.

\bibitem[{Allais(1953)}]{Al53}
Allais M (1953) Le comportement de l'homme rationnel devant le risque: critique
des postulats et axiomes de l'{\'e}cole am{\'e}ricaine. \emph{Econometrica: Journal of the Econometric Society} 503--546.

\bibitem[{Armbruster \protect\BIBand{} Delage(2015)}]{ArD15}
Armbruster B, Delage E (2015) Decision making under uncertainty when preference
information is incomplete. \emph{Management Science} 61(1):111--128.

\bibitem[{Baillon et~al.(2022)Baillon, Bleichrodt, Emirmahmutoglu, Jaspersen, \protect\BIBand{} Peter}]{BBEJ22}
Baillon A, Bleichrodt H, Emirmahmutoglu A, Jaspersen J, Peter R (2022)
When risk perception gets in the way: Probability weighting and underprevention. \emph{Operations Research} 70(3): 1371--1392.

\bibitem[{Benartzi \protect\BIBand{} Thaler(1995)}]{BeT95}
Benartzi S, Thaler RH (1995) Myopic loss aversion and the equity premium puzzle. \emph{The Quarterly Journal of Economics} 110(1):73--92.

\bibitem[{Bertsimas\protect\BIBand{} O'Hair(2013)}]{BeO13}
Bertsimas D,  O'Hair A (2013)
Learning preferences under noise and loss aversion: An optimization approach. 
\emph{Operations Research} 61(5): 1190--1199.

\bibitem[{Bleichrodt \protect\BIBand{} Pinto(2000)}]{BlP00}
Bleichrodt H, Pinto JL (2000)
A parameter-free elicitation of the probability weighting function in medical decision analysis. \emph{Management Science} 46(11): 1485--1496.

\bibitem[{Boyd \protect\BIBand{} Vandenberghe(2004)}]{BoV04}
Boyd S, Vandenberghe L (2004) \emph{Convex Optimization}
(Cambridge University Press).

\bibitem[{Burden et~al.(2015)Burden, Faires, \protect\BIBand{} Burden}]{BFB15}
Burden RL, Faires JD, Burden AM (2015) \emph{Numerical Analysis} (Cengage learning).

\bibitem[{Camerer \protect\BIBand{} Ho(1994)}]{CaH94}
Camerer CF, Ho TH (1994) Violations of the betweenness axiom and nonlinearity in probability. \emph{Journal of Risk and Uncertainty} 8(2): 167--196.

\bibitem[{Chen \protect\BIBand{} Liu(2026)}]{ChL26}
Chen B, Liu J (2026) Eliciting Von Neumann-Morgenstern utility from discrete choices with response error. \emph{arXiv preprint arXiv:2603.26065}.

\bibitem[{Delage et~al.(2022)Delage, Guo, \protect\BIBand{} Xu}]{DGX22}
Delage E, Guo S, Xu H (2022) Shortfall risk models when information on loss function is incomplete. \emph{Operations Research} 70(6):3511--3518.

\bibitem[{Ellsberg(1961)}]{El1961}
Ellsberg D (1961) Risk, ambiguity, and the savage axioms. \emph{The Quarterly
Journal of Economics} 75(4):643--669.

\bibitem[{F{\"o}llmer \protect\BIBand{} Schied(2002)}]{FoS02}
F{\"o}llmer H, Schied A (2002) Convex measures of risk and trading constraints.
\emph{Finance and Stochastics} 6(4):429--447.

\bibitem[{F{\"o}llmer \protect\BIBand{} Schied(2004)}]{FoS04}
F{\"o}llmer H, Schied A (2004)
\emph{Stochastic Finance: An Introduction in Discrete Time} (Walter de Gruyter, Berlin).

\bibitem[{Giesecke et~al.(2008)Giesecke, Schmidt, \protect\BIBand{} Weber}]{GSW08}
Giesecke K, Schmidt T, Weber S (2008) Measuring the risk of large losses. \emph{Journal of Investment Management} 6(4): 1--15.

\bibitem[{Gonzalez \protect\BIBand{} Wu(1999)}]{GoW99}
Gonzalez R, Wu G (1999) On the shape of the probability weighting function. \emph{Cognitive Psychology} 38(1): 129--166.

\bibitem[{Guo \protect\BIBand{} Xu(2019)}]{GuX18}
Guo S, Xu H (2019) Distributionally robust shortfall risk optimization model
and its approximation. \emph{Mathematical Programming} 174(1-2):473--498.

\bibitem[{Guo et~al.(2024)Guo, Xu, \protect\BIBand{} Zhang}]{GXZ23}
Guo S, Xu H, Zhang S (2024) Utility preference robust optimization with moment-type information structure. \emph{Operations Research} 72(5): 2241–-2261.

\bibitem[{Haskell et~al.(2016)Haskell, Fu, \protect\BIBand{} Dessouky}]{HFD16}
Haskell WB, Fu L, Dessouky M (2016) Ambiguity in risk preferences in robust stochastic optimization. \emph{European Journal of Operational Research} 254(1):214--225.

\bibitem[{Hu \protect\BIBand{} Stepanyan(2017)}]{HuS17}
Hu J, Stepanyan G (2017) Optimization with reference-based robust preference
constraints. \emph{SIAM Journal on Optimization} 27(4):2230--2257.

\bibitem[{Kahneman \protect\BIBand{} Tversky(1979)}]{KaT79}
Kahneman D, and Tversky A (1979) Prospect theory: an analysis of decision under risk. \emph{Econometrica} 47(2): 263--291.

\bibitem[{Konno(1976)}]{Kon76}
Konno H (1976) A cutting plane algorithm for solving bilinear programs. \emph{Mathematical Programming} 11(1): 14--27.

\bibitem[{Hu et~al.(2025)Hu, Zhang, Xu, \protect\BIBand{} Zhang}]{HZXZ22}
Hu J, Zhang D, Xu H, Zhang S (2025)
Distributional utility preference robust optimization models in multi-attribute decision making. \emph{Mathematical Programming} 212: 519–-565.

\bibitem[{Kothiyal et~al.(2011)Kothiyal, Spinu, \protect\BIBand{} Wakker}]{KSW11}
Kothiyal A, Spinu V, Wakker PP (2011) Prospect theory for continuous distributions: A preference foundation. \emph{Journal of Risk and Uncertainty} 42:195--210.

\bibitem[{Lattimore et~al.(1992)Lattimore, Baker, \protect\BIBand{} Witte}]{LBW92}
Lattimore PK, Baker JR, Witte AD (1992) The influence of probability on risky choice: A parametric examination. \emph{Journal of Economic Behavior and Organization} 17(3): 377--400.

\bibitem[{Mao \protect\BIBand{} Cai(2018)}]{MaC18}
Mao T, Cai J (2018) Risk measures based on behavioural economics theory. \emph{Finance and Stochastics} 22:367--393.

\bibitem[{Markle et~al.(2018)Markle, Wu, White, \protect\BIBand{} Sackett}]{markle2018goals}
Markle A, Wu G, White R, Sackett A (2018) Goals as reference points in marathon running: A novel test of reference dependence. \emph{Journal of Risk and Uncertainty} 56:19--50.

\bibitem[{McCormick(1976)}]{McC76}
McCormick GP (1976) Computability of global solutions to factorable nonconvex programs: Part I—Convex underestimating problems. \emph{Mathematical Programming} 10(1): 147--175.

\bibitem[{McCullough(2002)}]{McC02}
McCullough D (2002) A user's guide to conjoint analysis. \emph{Marketing Research} 14(2):18--23.

\bibitem[{Pesenti et~al.(2025)Pesenti, Wang, \protect\BIBand{} Wang}]{PWW20}
Pesenti S, Wang Q, Wang R (2025) Optimizing distortion riskmetrics with distributional uncertainty.  \emph{Mathematical Programming} 213: 51–-106.
  
\bibitem[{Prelec(1998)}]{Pre98}
Prelec D (1998) The probability weighting function. \emph{Econometrica} 497--527.

\bibitem[{Quiggin(1982)}]{Qui82}
Quiggin J (1982) A theory of anticipated utility. \emph{Journal of Economic Behavior $\&$ Organization} 3(4):323--343.

\bibitem[{Quiggin(2012)}]{Qui12}
Quiggin J (2012) \emph{Generalized Expected Utility Theory: The Rank-Dependent Model} (Springer Science \& Business Media).

\bibitem[{Rieger \protect\BIBand{} Wang(2006)}]{RiW06}
Rieger MO, Wang M (2006) Cumulative prospect theory and the St. Petersburg paradox. \emph{Economic Theory} 28: 665--679.

\bibitem[{Rothschild \protect\BIBand{} Stiglitz(1970)}]{RoS70}
Rothschild M, Stiglitz J E (1970) Increasing risk: I. A definition. \emph{Journal of Economic Theory} 2(3): 225--243.

\bibitem[{Saur{\'e} \protect\BIBand{} Vielma(2019)}]{SaV19}
Saur{\'e} D, Vielma JP (2019) Ellipsoidal methods for adaptive choice-based
conjoint analysis. \emph{Operations Research} 67(2):315--338.

\bibitem[{Schmeidler(1989)}]{Sch89}
Schmeidler D (1989). Subjective probability and expected utility without additivity. \emph{Econometrica: Journal of the Econometric Society}
57(3):571--587.

\bibitem[{Schmidt \protect\BIBand{} Zank(2008)}]{ScZ08}
Schmidt U, Zank H (2008) Risk aversion in cumulative prospect theory. \emph{Management Science} 54(1): 208--216.

\bibitem[{Shapiro et~al.(2009)Shapiro, Dentcheva, \protect\BIBand{} Ruszczynski}]{SDR14}
Shapiro A, Dentcheva D, Ruszczynski A (2009) \emph{Lectures on Stochastic Programming} (SIAM).

\bibitem[{Sherali \protect\BIBand{} Alameddine(1992)}]{ShA92}
Sherali HD, Alameddine A (1992) A new reformulation-linearization technique for bilinear programming problems. \emph{Journal of Global Optimization} 2: 379--410.

\bibitem[{Sherali \protect\BIBand{} Shetty(1980)}]{ShS80}
Sherali HD, Shetty CM (1980) A finitely convergent algorithm for bilinear programming problems using polar cuts and disjunctive face cuts. \emph{Mathematical Programming} 19(1): 14--31.

\bibitem[{Stott(2006)}]{Sto06}
Stott HP (2006) Cumulative prospect theory's functional menagerie. \emph{Journal of Risk and Uncertainty} 32:101--130.

\bibitem[{Thieu(1988)}]{Thi88}
Thieu TV (1988) A note on the solution of bilinear programming problems by reduction to concave minimization. \emph{Mathematical Programming} 41: 249--260.

\bibitem[{Toubia et~al.(2004)Toubia, Hauser, \protect\BIBand{} Simester}]{THS04}
Toubia O, Hauser JR, Simester DI (2004) Polyhedral methods for adaptive choice-based conjoint analysis. \emph{Journal of Marketing Research} 41(1):116--131.  

\bibitem[{Toubia et~al.(2003)Toubia, Simester, Hauser, \protect\BIBand{} Dahan}]{TSHD03}
Toubia O, Simester DI, Hauser JR, Dahan E (2003) Fast polyhedral adaptive conjoint estimation. \emph{Marketing Science} 22(3):273--303.

\bibitem[{Toubia et~al.(2007)Toubia, Hauser,\protect\BIBand{} Garcia}]{THG07}
Toubia O, Hauser J, Garcia R (2007) Probabilistic polyhedral methods for adaptive choice-based conjoint analysis: Theory and application. \emph{Marketing Science} 26(5):596--610.

\bibitem[{Tuthill \protect\BIBand{} Frechette(2002)}]{TuF02}
Tuthill JW, Frechette DL (2002) Non-expected utility theories: Weighted
  expected, rank dependent, and cumulative prospect theory utility. \emph{2002
  Conference, St. Louis, Missouri, April 22--23} $\;$(1266-2016-101906):24, \urlprefix{\footnotesize \url{https://doi.org/10.22004/ag.econ.19073}}.

\bibitem[{Tversky \protect\BIBand{} Fox(1995)}]{TvF95}
Tversky A, Fox CR (1995) Weighing risk and uncertainty. \emph{Psychological
 Review} 102(2):269.

\bibitem[{Tversky \protect\BIBand{} Kahneman(1992)}]{TvK92}
Tversky A, Kahneman D (1992) Advances in prospect theory: Cumulative representation of uncertainty. \emph{Journal of Risk and Uncertainty} 5:297--323.

\bibitem[{Tversky \protect\BIBand{} Wakker(1995)}]{TvW95}
Tversky A, Wakker P (1995) Risk attitudes and decision weights. \emph{Econometrica: Journal of the Econometric Society} 63(6):1255--1280.

\bibitem[{Vayanos et~al.(2020)Vayanos, Ye, McElfresh, Dickerson,
  \protect\BIBand{} Rice}]{VYMDR20}
Vayanos P, Ye Y, McElfresh D, Dickerson J, Rice E (2020) Robust active preference elicitation. \emph{arXiv preprint arXiv:2003.01899}.
Available at {\footnotesize \url{https://doi.org/10.48550/arXiv.2003.01899}} (Accessed 2020). 

\bibitem[{von~Neumann \protect\BIBand{} Morgenstern(1947)}]{VNM47}
von~Neumann J, Morgenstern O (1947) \emph{Theory of games and economic behavior, 2nd rev} (Princeton university press).

\bibitem[{Wakker \protect\BIBand{} Tversky(1993)}]{WaT93}
Wakker P, Tversky A (1993) An axiomatization of cumulative prospect theory. \emph{Journal of Risk and Uncertainty} 7:147--175.

\bibitem[{Wang \protect\BIBand{} Xu(2023)}]{WaX23}
Wang W, Xu H (2023) Preference robust distortion risk measure and its application. \emph{Mathematical Finance} 33(2):389--434.  

\bibitem[{Weber(2006)}]{We06}
Weber S (2006) Distribution-invariant risk measures, information, and dynamic consistency. \emph{Mathematical Finance} 16(2):419--441.

\bibitem[{White(1992)}]{Whi92}
White DJ (1992) A linear programming approach to solving bilinear programmes. \emph{Mathematical Programming} 56: 45--50.

\bibitem[{Wu \protect\BIBand{} Gonzalez(1996)}]{WuG96}
Wu G, Gonzalez R (1996) Curvature of the probability weighting function. \emph{Management Science} 42(12):1676--1690.

\bibitem[{Wu \protect\BIBand{} Gonzalez(1998)}]{WuG98}
Wu G, Gonzalez R (1998) Common consequence conditions in decision making under risk. \emph{Journal of Risk and Uncertainty} 16: 115--139.

\bibitem[Zhang \protect\BIBand{} Xu(2022)]{ZhX21}
Zhang S and Xu H (2022) Insurance premium-based shortfall risk measure induced by cumulative prospect theory. \emph{Computational Management Science} 19(4): 703--738.

\bibitem[{Zhang et~al.(2022)Zhang, Xu, \protect\BIBand{} Wang}]{ZXW22}
Zhang Y, Xu H, Wang W (2022) Preference robust models in multivariate utility-based shortfall risk minimization. \emph{Optimization Methods and Software} 37(2): 712--752.

\bibitem[{Zhang et~al.(2025)Zhang, Guo, Sim,   \protect\BIBand{} Xu}]{ZGSX25}
Zhang, S., Guo, S., Sim, M. and Xu, H. (2025) Modified Polyhedral Method for Elicitation of Shape-Free Utility and Conservatism Reduction in Robust Optimization. arXiv preprint arXiv:2503.23269. Available at {\footnotesize \url{https://doi.org/10.48550/arXiv.2003.01899}} (Accessed 2025).
\end{thebibliography}
\end{document}